\documentclass[10pt]{amsart}

\usepackage{a4wide}
\usepackage[utf8]{inputenc}
\usepackage{amsmath, amssymb, bbm}
\usepackage{amsthm}
\usepackage{bbm}
\usepackage{blkarray}
\usepackage[matrix,arrow,curve]{xy}
\usepackage{tikz}
\usepackage{graphicx}
\usepackage{color}
\usepackage[normalem]{ulem}
\usepackage{enumerate}
\usepackage{enumitem}
\definecolor{darkblue}{rgb}{0,0,0.6}
\usepackage[ocgcolorlinks,colorlinks=true, citecolor=darkblue, filecolor=darkblue, linkcolor=darkblue, urlcolor=darkblue]{hyperref}
\usepackage[capitalize,noabbrev]{cleveref}
\usepackage[mathscr]{euscript}
\usepackage{stmaryrd}

\usepackage{sansmath}

\newcommand{\cocolon}{\nobreak \mskip6mu plus1mu \mathpunct{}\nonscript\mkern-\thinmuskip {:}\mskip2mu \relax}

\newcommand*{\coloneqq}{\mathrel{\vcenter{\baselineskip0.5ex \lineskiplimit0pt
                     \hbox{\scriptsize.}\hbox{\scriptsize.}}}%
                     =}
                     
\usepackage{pict2e,picture}
\makeatletter
\DeclareRobustCommand{\bbDelta}{{\mathpalette\bb@Delta\relax}}
\newcommand{\bb@Delta}[2]{%
  \begingroup
  \sbox\z@{$\m@th#1\Delta$}%
  \dimendef\Dht=6 \dimendef\Dwd=8
  \setlength{\Dwd}{\wd\z@}%
  \setlength{\Dht}{\ht\z@}%
  \begin{picture}(\Dwd,\Dht)
  \put(0,0){$\m@th#1\Delta$}
  \put(.52\Dwd,.77\Dht){\line(-13,-26){.35\Dht}}
  \end{picture}%
  \endgroup
}
\makeatother

\DeclareFontFamily{T1}{cbgreek}{}
\DeclareFontShape{T1}{cbgreek}{m}{n}{<-6>  grmn0500 <6-7> grmn0600 <7-8> grmn0700 <8-9> grmn0800 <9-10> grmn0900 <10-12> grmn1000 <12-17> grmn1200 <17-> grmn1728}{}
\DeclareSymbolFont{quadratics}{T1}{cbgreek}{m}{n}
\DeclareMathSymbol{\qoppa}{\mathord}{quadratics}{19}
\DeclareMathSymbol{\Qoppa}{\mathord}{quadratics}{21}

\newtheoremstyle{introthms}
	{}{}{\itshape}{}{\bfseries }{}{ }
	{\thmname{#1} \thmnumber{#2}. \thmnote{\bfseries{(#3)}}}

\newtheoremstyle{thms}
	{}{}{\itshape}{}{\bfseries }{}{ }
	{\thmname{#1}\thmnumber{#2}. \thmnote{\bfseries{[#3]}}}
\newtheoremstyle{thms2}
	{}{}{\itshape}{}{\bfseries }{}{ }
	{\thmname{#1}\thmnumber{#2}. \thmnote{\bfseries{[#3]}}}

\newtheoremstyle{name}
	{}{}{\itshape}{}{\bfseries }{}{ }
	{\thmname{#1}\thmnumber{#2}\thmnote{\bfseries{[#3]}}}
\newtheoremstyle{defs}
	{}{12pt}{\normalfont}{}{\bfseries }{}{ }
	{\thmname{#1} \thmnumber{#2}. \thmnote{\bfseries{(#3)}}}
\newtheoremstyle{defs2}
	{}{12pt}{\normalfont}{}{\bfseries }{}{ }
	{\thmname{#1}\thmnumber{#2}. \thmnote{\bfseries{(#3)}}}
\newtheoremstyle{rmk}
	{}{}{\normalfont}{}{\itshape }{}{ }{\thmname{#1}. \thmnote{#3}}
\newtheoremstyle{claim}
	{}{}{\normalfont}{}{\itshape}{}{ }{\thmname{#1} \thmnumber{#2}. \thmnote{#3}}

\theoremstyle{introthms}
\newtheorem{introthm}{Theorem}

\swapnumbers
	
\theoremstyle{thms}
\newtheorem{proposition}{Proposition}[subsection]
\newtheorem{theorem}[proposition]{Theorem}
\newtheorem{add}[proposition]{Addendum}
\newtheorem{lemma}[proposition]{Lemma}
\newtheorem{corollary}[proposition]{Corollary}

\newtheorem*{ucorollary}{Corollary}

\newtheorem{aproposition}{Proposition}[subsection]
\newtheorem{alemma}[aproposition]{Lemma}
\newtheorem{atheorem}[aproposition]{Theorem}
\newtheorem{acorollary}[aproposition]{Corollary}

\theoremstyle{defs}
\newtheorem{definition}[proposition]{Definition}
\newtheorem{construction}[proposition]{Construction}
\newtheorem{example}[proposition]{Example}
\newtheorem{remark}[proposition]{Remark}

\newtheorem{observation}[proposition]{Observation}
\newtheorem{convention}[proposition]{Convention}

\newtheorem{adefinition}[aproposition]{Definition}
\newtheorem{aconstruction}[aproposition]{Construction}
\newtheorem{aremark}[aproposition]{Remark}

\theoremstyle{defs2}

\theoremstyle{rmk}

\theoremstyle{claim}

\renewcommand{\L}{\mathrm L}               
\renewcommand{\l}{\mathcal {L}}
\newcommand{\K}{\mathrm K}                 
\renewcommand{\k}{\mathcal K}                 
\newcommand{\GW}{\mathrm{GW}}              
\newcommand{\KR}{\mathrm{KR}}
\newcommand{\kr}{\mathrm{kr}}
\newcommand{\gw}{\mathcal{GW}}              
\newcommand{\Q}{\mathcal Q}                
\newcommand{\Poinc}{\mathrm{Pn}}            

\newcommand{\spcforms}{\mathrm{Fm}} 			
\newcommand{\catforms}{\mathrm{He}}            
\newcommand{\op}{^\mathrm{op}}              
\newcommand{\rev}{\mathrm{rev}}
\newcommand{\dec}{\mathrm{dec}}
\newcommand{\Bil}{\mathrm B}               
\newcommand{\Lin}{\Lambda}               
\newcommand{\Dual}{\mathrm D}              
\newcommand{\fib}{\mathrm{fib}}            
\newcommand{\cof}{\mathrm{cof}}            
\newcommand{\nerv}{\mathrm N}              
\newcommand{\asscat}{\mathrm{asscat}}

\newcommand{\pcone}{\mathrm{lpc}}

\newcommand{\Hom}{\mathrm{Hom}}            
\renewcommand{\hom}{\mathrm{hom}}          
\newcommand{\Met}{\operatorname{Met}}

\newcommand{\sd}{\operatorname{TwAr}}
\newcommand{\Twar}{\operatorname{TwAr}}
\newcommand{\colim}{\operatorname{colim}}
\newcommand{\Hyp}{\operatorname{Hyp}}
\newcommand{\hyp}{\operatorname{hyp}}
\newcommand{\fin}{\mathrm{fin}}
\newcommand{\grp}{\mathrm{grp}}
\newcommand{\fgt}{\mathrm{fgt}}
\newcommand{\srg}{\mathrm{srg}}

\newcommand{\heart}{\mathrm{ht}}
\newcommand{\Sq}{\operatorname{Sq}}
\newcommand{\refl}{\mathrm{refl}}
\newcommand{\ev}{\mathrm{ev}}
\newcommand{\Qclass}{\Q_{\mathrm{cl}}}  
\newcommand{\Qhclass}{\Q^\mathrm{h}_{\mathrm{cl}}} 
\newcommand{\Ho}{\mathrm{Ho}}
\newcommand{\Ffl}{\mathrm{F}^\fpm}     
\newcommand{\Horn}{\mathrm{H}}

\newcommand{\Ab}{\mathrm{Ab}}
\newcommand{\Stab}{\mathrm{Stab}}          
\newcommand{\Cath}{\mathrm{Cat^h}}         
\newcommand{\Catp}{\mathrm{Cat}^\mathrm{p}}         
\newcommand{\Catwic}{\mathrm{Cat}^{\mathrm{add}}_{\flat}}
\newcommand{\Catpad}{\mathrm{Cat}^{\mathrm{ap}}_{\flat}}
\newcommand{\Spa}{\mathrm{Sp}}              
\newcommand{\Sps}{\mathrm{An}}              
\newcommand{\sSps}{{\mathrm{s}\Sps}}

\newcommand{\Cat}{\mathrm{Cat}_\infty}            
\newcommand{\Catx}{\mathrm{Cat}^\mathrm{st}}      
\newcommand{\Catw}{\Catx_\mathrm{wt}}    
\newcommand{\Catpw}{\Catp_\mathrm{wt}}    
\newcommand{\Cata}{\mathrm{Cat}^\mathrm{add}}    
\newcommand{\Cataq}{\mathrm{Cat}^\mathrm{aq}}              
\newcommand{\Fun}{\mathrm{Fun}}            

\newcommand{\qshift}[1]{^{[#1]}} 
\newcommand{\wshift}[1]{^{\langle #1\rangle}}

\newcommand{\inj}{\mathrm{inj}}
\newcommand{\const}{\mathrm{const}}
\newcommand{\Dperf}{\mathcal D^{\mathrm{p}}}

\newcommand{\PProj}{\mathcal{P}}
\newcommand{\Modp}{\mathrm{Mod}^\omega}
\newcommand{\Mod}{{\mathrm{Mod}}}

\newcommand{\C}{\mathcal C}                
\newcommand{\D}{\mathcal D}                
\newcommand{\QF}{\Qoppa}               

\newcommand{\pr}{\mathrm{pr}}              
\newcommand{\id}{\mathrm{id}}              

\newcommand{\Eone}{\mathbb E_1}
\newcommand{\Ct}{{\mathrm{C}_2}}
\newcommand{\hC}{{\mathrm{hC}_2}}           
\newcommand{\tC}{{\mathrm{tC}_2}}           
\newcommand{\sym}{\mathrm{s}}
\newcommand{\fw}{\mathrm{fw}}
\newcommand{\bw}{\mathrm{bw}}
\newcommand{\cobmor}{\rightsquigarrow}
\newcommand{\Cob}{\mathsf{Cob}}
\newcommand{\Qcl}{\mathcal Q_{\mathrm{cl}}}
\newcommand{\Qhcl}{\mathcal Q_{\mathrm{cl}}^{\mathrm h}}
\newcommand{\Cobsd}{\Cob^\sgd}
\newcommand{\sgd}{\mathrm{sd}}

\newcommand{\cc}{\mathrm{sc}}
\newcommand{\Cobfw}{\Cob^{\mathrm{fw}}}
\newcommand{\Cobbw}{\Cob^{\mathrm{bw}}}

\newcommand{\lcone}[1]{^\triangleleft {#1}}
\newcommand{\rcone}[1]{{#1}^\triangleright}
\newcommand{\Ar}{\operatorname{Ar}}
\newcommand{\E}{\mathbb E}

\newcommand{\mc}{m} 
\newcommand{\oc}{p} 

\newcommand{\Qsq}{\Q^{\square}}
\newcommand{\SC}{\operatorname{SC}}
\newcommand{\SSC}{\operatorname{SSC}}

\newcommand{\Span}{\mathrm{Span}}
\newcommand{\coll}{\mathfrak C}

\newcommand{\Psd}{\Poinc_\sgd}
\newcommand{\Pmsd}[1]{\Poinc_{{#1}\text{-}\sgd}}

\renewcommand{\dec}{\operatorname{dec}}

\newcommand{\DCob}{\Cob^{(2)}}
\newcommand{\SurgCob}{\mathrm{SurgCob}}
\newcommand{\SurgCobfw}{\mathrm{SurgCob}_\fw}
\newcommand{\SurgCobbw}{\mathrm{SurgCob}_\bw}

\renewcommand{\Cobsd}{\Cob_\sgd}
\renewcommand{\Cobfw}{\Cob_\fw}
\renewcommand{\Cobbw}{\Cob_\bw}
\newcommand{\Cobmsd}[1]{\Cob_{{#1}\text{-}\sgd}}
\newcommand{\Cobmfw}[1]{\Cob_{{#1}\text{-}\fw}}
\newcommand{\Cobmbw}[1]{\Cob_{{#1}\text{-}\bw}}

\newcommand{\chifw}{\chi_\fw}
\newcommand{\chibw}{\chi_\bw}

\newcommand{\lcube}[2]{[#1]^{#2}}
\newcommand{\powerar}{\lcube a {r+1}}
\newcommand{\Pccsd}[2]{\Poinc_{[{#1}]\text{-}\mathrm{sd}^{[{#2}]}}}
\newcommand{\Cobccsd}[2]{\Cob_{[{#1}]\text{-}\mathrm{sd}^{[{#2}]}}}

\newcommand{\sk}{\operatorname{sk}}
\newcommand{\Dinj}{\bbDelta_\inj}
\newcommand{\eDinj}{{\widehat\bbDelta}_\inj}
\newcommand{\CCbdd}{\mathrm{Ch}^{\mathrm{bd}}}

\makeatletter
\providecommand{\leftsquigarrow}{%
  \mathrel{\mathpalette\reflect@squig\relax}%
}
\newcommand{\reflect@squig}[2]{%
  \reflectbox{$\m@th#1\rightsquigarrow$}%
}
\makeatother

\newcommand{\grpcr}{\mathrm{cr}}
\newcommand{\Surg}{\mathrm{Surg}}

\newcommand{\Deltainj}{\bbDelta_\mathrm{inj}}
\newcommand{\short}{\mathrm{sht}}
\newcommand{\axesar}{{[a]^{r+1}_{\mathrm{ax}}}}
\newcommand{\msrg}{\mathrm{msrg}}

\newcommand{\Db}{\mathrm{db}}
\newcommand{\core}{\grpcr}

\newcommand{\rsrg}{\mathrm{rsrg}}
\newcommand{\sks}{{-\s}}
\newcommand{\skq}{{-\q}}
\newcommand{\g}{\mathrm{g}}
\newcommand{\Proj}{\mathrm{Proj}}

\newcommand{\fpm}{\lambda}
\newcommand{\Unimod}{\mathrm{Unimod}}
\newcommand{\q}{\mathrm{q}}
\newcommand{\s}{\mathrm{s}}

\newcommand{\vis}{\mathrm{v}}
\renewcommand{\H}{\mathrm{H}}

\newcommand{\gs}{\mathrm{gs}}

\newcommand{\gr}{{\mathrm{g}\fpm}}
\newcommand{\Einf}{\mathbb E_\infty}

\newcommand{\pl}{\oplus}
\newcommand{\A}{\mathcal{A}}
\newcommand{\B}{\mathcal{B}}
\newcommand{\Grp}{\operatorname{Grp}}  
\newcommand{\x}{x} 
\newcommand{\y}{y}
\newcommand{\z}{z}
\newcommand{\cob}{w}  
\newcommand{\hrar}{\hookrightarrow}
\newcommand{\hlar}{\hookleftarrow}

\newcommand{\st}{\stackrel}
\newcommand{\F}{\mathcal{F}}
\newcommand{\G}{\mathcal{G}}
\newcommand{\Hlgy}{\operatorname{Hlgy}} 
\newcommand{\Seq}{\operatorname{Seq}}

\newcommand{\lto}{\longrightarrow}
\newcommand{\vphi}{\varphi}
\newcommand{\bS}{\mathrm{S}}
\newcommand{\Sig}{\Sigma}
\newcommand{\sig}{\sigma}
\newcommand{\Del}{\bbDelta}
\newcommand{\Om}{\Omega}
\newcommand{\Lam}{\Lambda}

\newcommand{\QFD}{\Phi}
\newcommand{\NN}{\mathbb{N}}
\newcommand{\Poset}{\mathrm{Poset}}       
\newcommand{\Beta}{\mathfrak{P}}
\newcommand{\R}{\mathrm{R}}
\newcommand{\cosk}{\mathrm{cosk}}
\newcommand{\I}{\mathcal{I}} 
\newcommand{\scc}{\mathrm{scc}}

\newcommand{\Lag}{\mathcal{E}}
\newcommand{\wic}{\flat}
\newcommand{\Kspace}{\k}
\newcommand{\GWspace}{\gw}
\newcommand{\Lspace}{\l}

\newcommand{\Sur}{\Surg}

\title[Stable Moduli spaces of hermitian forms]{Stable Moduli spaces of hermitian forms}

\author[Fabian Hebestreit]{Fabian Hebestreit}
\address{Universität Bielefeld; Fakultät für Mathematik, Bielefeld, Germany}
\email{hebestreit@math.uni-bielefeld.de}
\author[Wolfgang Steimle]{Wolfgang Steimle}
\address{Universit\"at Augsburg; Institut f\"ur Mathematik, Augsburg, Germany}
\email{wolfgang.steimle@math.uni-augsburg.de}

\makeatletter
\let\@wraptoccontribs\wraptoccontribs
\makeatother

\contrib[with an Appendix by]{Yonatan Harpaz}
\address{Université Paris Cité, Sorbonne Université, Paris, France}
\email{harpaz@imj-prg.fr}

\date{}
\dedicatory{In memory of Bruce Williams.}
\begin{document}
\setcounter{tocdepth}{1}

\begin{abstract}
We prove that Grothendieck-Witt spaces of Poincar\'e categories are, in many cases, group completions of certain moduli spaces of hermitian forms. This, in particular, identifies Karoubi's classical hermitian and quadratic $\K$-groups with the genuine Grothendieck-Witt groups from our joint work with Calm\`es, Dotto, Harpaz, Land, Moi, Nardin and Nikolaus, and thereby completes our solution of several conjectures in hermitian K-theory. 

The method of proof is abstracted from work of Galatius and Randal-Williams on cobordism categories of manifolds using the identification of the Grothendieck-Witt space of a Poincar\'e category as the homotopy type of the associated cobordism category. 
\end{abstract}

\maketitle
\tableofcontents


\section{Introduction}

Grothendieck-Witt theory, or hermitian $K$-theory, is the study of moduli spaces of unimodular forms after group completion. More specifically, following Quillen's monumental development of higher algebraic K-theory in \cite{QuillenHigherK,QuillenHigherK2}, Karoubi introduced the following definition in \cite{karoubi-quillen}: For simplicity let us at first restrict attention to a commutative ring $R$ 
and consider the groupoids
\[\Unimod^\fpm(R),\qquad  \fpm\in \{\s, \sks, \q, \skq\},\]
whose objects are pairs $(P, q)$ consisting of a finitely generated projective $R$-module $P$ and (according to the respective value of $\fpm$) a unimodular symmetric, skew-symmetric, quadratic or skew-quadratic form $q$, and whose morphisms are the isometries. While these groupoids are of great import in many areas of mathematics, ranging from number theory to manifold topology, they are not understood even for $R=\mathbb Z$. Per construction the components of $\Unimod^\fpm(R)$ form the set of isometry classes of unimodular forms of type $\fpm$ over $R$, and the homology of (the nerve of) this groupoid describes the group homology of the corresponding isometry groups:
\[\H_*(\Unimod^\fpm(R);\mathbb Z) \cong \bigoplus_{[P,q]} \H_*(\mathrm O(P,q);\mathbb Z),\]
where the direct sum runs over all such isometry classes.

The operation of orthogonal sum induces a symmetric monoidal structure on $\Unimod^\fpm(R)$ and thus a coherently associative and commutative addition, i.e.\ an $\Einf$-structure, on its nerve. As a tractable simplification of $\Unimod^\fpm(R)$ Karoubi defined the Grothendieck-Witt space of $R$ as the group completion of this $\Einf$-monoid,
\[\gw^\fpm_\mathrm{cl}(R) = \Unimod^\fpm(R)^\grp,\]
mimicking Quillen's definition of the algebraic K-space of $R$,
\[\k(R) = \Proj(R)^\grp,\]
where $\Proj(R)$ is the groupoid of finitely generated projective $R$-modules (see \cite{segal} for the translation of the original definition into the present language). Via Grothendieck's homotopy hypothesis, i.e. Joyal's theorem \cite{joyal}, the passage to the underlying $\Einf$-monoid in either case simply views a symmetric monoidal groupoid as a symmetric monoidal $\infty$-groupoid, making $\Einf$-groups as above the homotopical version of Picard groupoids in algebraic geometry.

The abelian group of components $\pi_0(\gw_{\mathrm{cl}}^\fpm(R))$ is the  classical Grothendieck-Witt group, i.e. the ordinary group completion of the monoid of isomorphism classes of unimodular forms of type $\fpm$, variously denoted by
\begin{align*}
\mathrm{KO}_0(R), \mathrm{KH}_0(R) \ \text{or}\ \mathrm{GW}_0(R) & \quad \quad \text{for } \fpm = \s \\
\mathrm{KSp}_0(R) \ \text{or} \ \mathrm{GW}_0(R,-) & \quad \quad\text{for } \fpm = -\s \\
\mathrm{KO}^\q_0(R), \mathrm{KQ}_0(R) \ \text{or} \ {}_1\mathcal{L}(R) & \quad \quad\text{for } \fpm = \q \\
\mathrm{KSp}^\q_0(R) \ \text{or} \ {}_{-1}\mathcal{L}(R) & \quad \quad\text{for } \fpm = -\q
\end{align*}
in the literature; we emphatically warn the reader that we shall follow conventions in geometric topology and use the letter $\L$ for Witt-, rather than Grothendieck-Witt groups. These group completions are much easier to access than the original monoids: For example it is known that
\[\pi_0(\gw_\mathrm{cl}^\q(\mathbb Z))
= 8 \mathbb Z \oplus \mathbb Z \subseteq \mathbb Z \oplus \mathbb Z = 
\pi_0(\gw_\mathrm{cl}^\s(\mathbb Z))
\]
via signature and rank of the positive definite part, despite the classification of unimodular forms over the integers being a wide open field. 

The effect of group completion on homology is also understood by the group completion theorem of McDuff and Segal \cite{mcduff-segal}: For example, whenever $R$ has no $2$-torsion, we have
\[\H_*(\gw^\q_\mathrm{cl}(R)_0; \mathbb Z) = \colim_{g \in \mathbb N} \H_*(\mathrm{O}_{g,g}(R);\mathbb Z),\]
where $\mathrm{O}_{g,g}(R)$ is the automorphism group of the standard (symmetric) hyperbolic form on $R^{2g}$ (and the subscript $0$ on the left denotes the unit path component); for general $R$ one has to either take the subgroup $\mathrm{O}_{g,g}^\q(R)$ on the right that further preserves the canonical quadratic refinement or replace the Grothendieck-Witt space on the left by the version associated the form parameter $\lambda=\mathrm{ev}$ which we shall discuss later. For nice enough $R$ such as $R=\mathbb Z$ one similarly finds
\[\H_*(\gw^\s_\mathrm{cl}(R)_0; \mathbb Z) = \colim_{g \in \mathbb N} \H_*(\mathrm{O}_{\langle g,g\rangle}(R);\mathbb Z),\]
where $\mathrm{O}_{\langle g,g\rangle}(R)$ is the automorphism group of the sum of $g$-copies of the standard positive and negative definite forms on $R^g$ each. By definition the homotopy groups of $\gw^\fpm_{\mathrm{cl}}(R)$ are the higher Grothendieck-Witt-groups. While there is no direct description of these groups in terms of unimodular forms and their automorphisms, they contain much subtle arithmetic information about $R$, that is highly non-trivial to extract, just as in the case of higher K-groups, see \emph{e.g.}\ \cite{weibel, BKSO-fixed-point, rondigs}.
 
Primarily by the work of Karoubi and Schlichting \cite{karoubi-quillen, Karoubi-Le-theoreme-fondamental, hornbostel-schlichting, Schlichting2010, schlichting-derived}, Grothendieck-Witt spaces are fairly well-understood relative to Quillen's algebraic K-space if one assumes that $2$ is a unit in $R$. By the tremendous progress in the understanding of the latter object, through the advent of topological cyclic homology, motivic homotopy theory, and the systematic study of assembly maps, this has led to a good understanding of Grothendieck-Witt theory of rings in which $2$ is invertible. However, a lot of interesting subtlety is lost in this regime: For example, the forgetful maps
\[\Unimod^\q(R) \longrightarrow \Unimod^\s(R) \quad \mathrm{and} \quad  \Unimod^{\skq}(R)\longrightarrow \Unimod^{\sks}(R)\]
taking a quadratic form to its polarisation are then equivalences.

\subsection{Genuine Grothendieck-Witt spaces}\label{sec:genuine_GW}

In recent joint work with Calm\`es, Dotto, Harpaz, Land, Moi, Nardin and Nikolaus \cite{9authI, 9authII, 9authIII, 9authIV, 9authV}, to which this paper is a companion, we introduced and investigated a new framework for Grothendieck-Witt theory, located in the modern realm of stable $\infty$-categories. We showed that the Grothendieck-Witt spaces arising from this new definition can always be decomposed into an algebraic K- and an L-(or Witt-)theoretic part, regardless of the invertibility of $2$. By addressing these pieces separately, such a decomposition allows for both structural and computational control not previously available without assuming $2$ a unit.  The principal goal of the present paper is to show that the classical Grothendieck-Witt spaces can be expressed as instances of this new theory, thereby rendering them accessible to computation.

The input for our new Grothendieck-Witt functor is Lurie's notion of a Poincar\'e category consisting of a stable $\infty$-category $\C$ together with a certain type of functor $\QF \colon \C\op \rightarrow \Spa$ to the $\infty$-category of spectra, which we call a Poincar\'e structure. It can be thought of as a derived version of the functor 
$\Ffl \colon \PProj(R)\op\to \Ab$ that sends $P$ to the abelian group of forms of type $\fpm$ on $P$; here $\PProj(R)$ denotes the full subcategory of $\mathrm{Mod}(R)$ spanned by the finitely generated projective $R$-modules, so that $\Proj(R)$ is its groupoid core. The output is a Grothendieck-Witt-space, defined by a version of the hermitian $\Q$-construction, which we denote $\gw(\C, \QF)$.

Now, by the results of \cite{9authI}, there are unique Poincar\'e structures $\QF^{\g\fpm}\colon \Dperf(R)\op\to \Spa$ on the perfect derived category of $R$ such that 
\[\pi_i\QF^{\g\fpm}(P) \cong \begin{cases} \Ffl(P) & i= 0 \\
 0 & \text{else} \end{cases}\]	
naturally in $P\in \PProj(R)\subseteq \Dperf(R)$. This identification provides a canonical map
\[\Unimod^\fpm(R) \longrightarrow \gw(\Dperf(R),\QF^{\gr})\]
of $\Einf$-monoids, which factors through the classical Grothendieck-Witt space, since the target is an $\Einf$-group.

\begin{introthm}
\label{thm:main_special}
Let $R$ be a commutative ring and $\fpm\in \{\s, \sks, \q, \skq\}$. Then the canonical map
\[\gw^\fpm_\mathrm{cl}(R) \longrightarrow \gw(\Dperf(R),\QF^{\g\fpm})\]
is an equivalence.
\end{introthm}

The conjecture that this should be so is originally due to Nikolaus, and served as one of the main motivations for our series. Combining Theorem \ref{thm:main_special} with the main results of \cite{9authII, 9authIII}, we for example find that:

\begin{ucorollary}
For any Dedekind ring $R$ there is a fibre sequence
\[\k(R)_\hC \xrightarrow{\hyp} \gw^{\s}_\mathrm{cl}(R) \xrightarrow{\mathrm{bord}} \l^{\s}(R),\]
that canonically splits after inverting $2$; here, the action on $\k(R)$ is induced by dualising a projective module and $\l^s(R)$ denotes Ranicki's symmetric L- or Witt-space. Furthermore, if the field of fractions of $R$ is a number field, the map
\[\fgt \colon \gw^\s_\mathrm{cl}(R) \longrightarrow \k(R)^\hC\]
is a $2$-adic equivalence.
\end{ucorollary} 

Along with its skew-symmetric and quadratic companions this solves several long standing conjectures of Berrick, Hesselholt, Karoubi, Madsen, Thomason, Williams and others regarding the Grothendieck-Witt groups of such rings, and allowed us to perform essentially complete calculations of all four flavours of Grothendieck-Witt groups of the integers in \cite{9authIII}. We refer the reader to the introductions of \cite{9authII, 9authIII} for a more detailed account of the history of such results.

The fibre sequence of the corollary also gives access to the homotopy type of Grothendieck-Witt spaces, and thus the stable cohomology of orthogonal groups. For example one can use this fibre sequence to determine the homotopy type of $\gw^\s_{\mathrm{cl}}(\mathbb Z)$, as explained in \cite{nik-lec} $2$-adically; this method will be used for systematic cohomology computations in forthcoming joint work with Land and Nikolaus \cite{homologyo}. As a sample calculation we here include:

\begin{ucorollary}
Certain classes $w_i,v_i \in \H^i(\mathrm{O}_{\langle g,g\rangle}(\mathbb Z);\mathbb F_2)$ and $a_i \in \H^{2i-1}(\mathrm{O}_{\langle g,g\rangle}(\mathbb Z);\mathbb F_2)$ induce a ring isomorphism
\[\mathbb F_2[w_i,v_i, a_{i}\mid i \geq 1] \longrightarrow \lim_{g \in \mathbb N\op}\H^*(\mathrm{O}_{\langle g,g\rangle}(\mathbb Z);\mathbb F_2).\]
where $\mathrm{O}_{\langle g,g\rangle}(\mathbb Z)$ is the isometry group of the symmetric form
\[\left(\mathbb Z^{2g},\begin{pmatrix} \mathrm{I}_g & 0 \\ 0&  -\mathrm I_g\end{pmatrix}\right).\]
\end{ucorollary}

Given the results of \cite{9authIII} and the present paper, this particular computation can also be deduced quickly from the main results of \cite{berrick-karoubi} by base-change to $\mathbb Z[\frac12]$, and we give that proof here. 

It is curious that to the best of our knowledge no result on the homological stability for the system formed by the groups $\mathrm{O}_{\langle g,g\rangle}(\mathbb Z)$ is contained in the literature as of yet (outside the case of rational coefficients); see \ref{homstabsymmrmk} below for a more detailed discussion of this point.

\subsection{Grothendieck-Witt spaces and weight structures}

The techniques we employ to prove Theorem \ref{thm:main_special} are not limited to the set of examples considered above: They extend to prove versions of Theorem \ref{thm:main_special} for arbitrary rings, invertible modules with involution as coefficients (which includes the classical set-up of rings with involution) and in fact arbitrary form parameters as in \cite{Bak} or \cite{Schlichting2019}, which for example contains the case of even forms. Moreover, they are not limited to such classical settings at all. For instance, they apply to the homotopy quadratic Poincar\'e structure $\QF^\q$ on $\Dperf(R)$, which by fundamental work of Ranicki captures the classical surgery obstruction groups $\L_*^\q(R)$ of Wall \cite{Ranickiblue}. In addition, our results also seamlessly extend to (connective) ring spectra and various flavours of parametrised spectra.
All of these examples fit into Bondarko's concept of a weight structure from \cite{Bondarko2010weights}, in which we frame our general results.

A weight structure on a stable $\infty$-category $\C$ can be defined by two subcategories $\C_{[0,\infty]}$ and $\C_{[-\infty,0]}$, closed under finite colimits and limits, respectively, and retracts, such that the mapping spectra $\hom_\C(X,Y)$ are connective, whenever $X \in \C_{[-\infty,0]}$ and $Y \in \C_{[0,\infty]}$ and 
such that every object $X \in \C$ lies in a fibre sequence 
\[Y \rightarrow X \rightarrow Z\]
where $Y \in \C_{[-\infty,0]}$ and $Z\qshift{-1} \in \C_{[0,\infty]}$; here and in the following, $(-)\qshift i$ denotes the $i$-fold shift, i.e.\ suspension, in $\C$. The heart of such a weight structure is the subcategory 
\[\C^\heart = \C_{[-\infty,0]} \cap \C_{[0,\infty]}.\]
If the stable hull of $\C^\heart$ is $\C$ itself, then we call the weight structure exhaustive, in which case it is entirely determined by its heart. As an example $\mathrm{Proj}(R) \subseteq \Dperf(R)$ is the heart of an exhaustive weight structure as is more generally the category spanned by retracts of $R^n, n \in \mathbb N$ inside $\Modp_{R}$ for any connective $\Eone$-ring $R$, where $\Modp_{R} \subseteq \Mod{R}$ is the category of compact $R$-modules.

To explain our general version of Theorem \ref{thm:main_special}, let us mention that any Poincar\'e structure $\QF$ on $\C$ in particular gives rise to a duality equivalence
\[\Dual_\QF\colon \C\op\to \C\]
on $\C$, and to an abstract notion of unimodular form, which following Ranicki is called a Poincar\'e object. It consists of a pair $(X,q)$ with $X \in \C$ and $q \in \Omega^\infty \QF(X)$, subject to a nondegeneracy condition, which in particular implies $X \simeq \Dual_\QF X$. In the four examples of $(\Dperf(R),\QF^\gr)$ the duality is simply given by $\mathbb R\mathrm{Hom}_R(-,R) \colon \Dperf(R)\op \rightarrow \Dperf(R)$. Just as for unimodular forms, there is an orthogonal sum of Poincar\'e objects, which assembles them into an $\Einf$-monoid $\Poinc(\C,\QF)$ and as such it comes equipped with a tautological map
\[\Poinc(\C,\QF) \longrightarrow \gw(\C,\QF).\]
If we denote by $\Poinc^\heart(\C,\QF)$ the full subspace of $\Poinc(\C,\QF)$ spanned by those Poincar\'e objects $(X,q)$ with $X \in \C^\heart$, we show the following weight theorem for Grothendieck-Witt spaces:

\begin{introthm}\label{thm:resolution}
Let $(\C,\QF)$ a Poincar\'e category equipped with an exhaustive weight structure, such that $\QF(X)$ is connective for each $X \in \C^\heart$ and $\Dual_\QF$ preserves the heart of $\C$. Then the canonical map
\[\Poinc^\heart(\C,\QF)^\grp \longrightarrow \gw(\C,\QF)\]
is an equivalence.
\end{introthm}

This result implies Theorem \ref{thm:main_special}, since
\[\Poinc^\heart(\Dperf(R),\QF^{\g\fpm}) \simeq \Unimod^\fpm(R)\]
for the weight structure mentioned above, and we also find immediately the evident variant of Theorem \ref{thm:main_special} for arbitrary form parameters.
In the case of Ranicki's quadratic Poincar\'e structure 
\[\QF^\q(C) = \hom_R(C \otimes_R C,R)_\hC\] 
on $\Dperf(R)$ the forgetful map
\[\Poinc^\heart(\Dperf(R),\QF^\q) \rightarrow \Unimod^\q(R)\]
is an isomorphism at the level of path components, but the automorphism groups on the left hand side can have many higher homotopy groups, so that the Theorem  \ref{thm:resolution} expresses $\gw(\Dperf(R), \QF^\q)$ as the group completion of a less classical $\Einf$-monoid. Similarly, let us mention the visible (or universal) Poincar\'e structure $\QF^\mathrm v$ on $\Spa^\fin$, the category of finite spectra, constructed in \cite{9authI}. It satisfies 
\[\GW(\Spa^\fin,\QF^\mathrm v) \simeq \mathrm{LA}^\mathrm v(*)\]
for the $\mathrm{LA}$-spectrum introduced by Weiss and Williams in their pursuit of a direct description of homeomorphism groups of closed manifolds. In this case Theorem \ref{thm:resolution} applies to describe the right hand side as the group completion of $\Poinc^\heart(\Spa^\fin,\QF^\mathrm v)$, which consists of orthogonal groups over the sphere spectrum, in analogy with Waldhausen's description of $\k(\Spa^\fin) = \Omega^\infty \mathrm A(*)$ in terms of spherical general linear groups. To the best of our knowledge this result is also new.

In fact, Waldhausen's result, and many similar ones, are also contained in Theorem \ref{thm:resolution} via the Poincar\'e categories $\Hyp(\C)=(\C\times \C\op, \hom_\C)$, whenever $\C$ is equipped with an exhaustive weight structure; in this case $\C\op$ and thus $\Hyp(\C)$ inherit weight structures and the duality 
is given by $(X,Y) \mapsto (Y,X)$ so preserves the heart. We thus find 
\[\grpcr(\C^\heart)^\grp \simeq \Poinc^\heart\Hyp(\C)^\grp \simeq \gw(\Hyp(\C)) \simeq \k(\C),\]
where $\grpcr$ denotes the groupoid core and the outer two equivalences are immediate from the definitions. This equivalence recovers a combination of Quillen's `+-equals-$\Q$'-theorem with the group completion theorem of Segal and McDuff and the work of Gillet and Waldhausen, mediating between the classical set-up of exact categories and their derived variants. The full statement has been a folklore result for some time, we originally learned about it from unpublished work of Clausen. We will give a brief discussion of the literature surrounding this result in the body of the text.

In distinction with existing proofs in the non-hermitian setting, we do not have to resort to any notion of exact or Waldhausen categories as an intermediary step; given the state of development of hermitian K-theory for such categories, it is in fact not clear whether such a proof can be lifted from algebraic to hermitian K-spaces.\footnote{During the revision process of the present paper Marlowe and Schlichting did give a proof of Theorem \ref{thm:main_special} along these more classical lines in \cite{SchlichtingII,SchlichtingIII}. The generality of their result is, however, incomparable to Theorem \ref{thm:resolution} in that it concerns only $\mathbb Z$-linear categories (thus ruling out appliciations to honest ring spectra), but allows for non-split exact structures.} Our proof is more direct and instead takes inspiration from cobordism theory.

\subsection{Parametrised algebraic surgery}
Theorem \ref{thm:resolution} is the zero-dimensional case of a more general result which we formulate next. Both statement and method of proof are greatly inspired by the work of Galatius, Madsen, Randal-Williams, Tillmann and Weiss  \cite{MW, GMTW, GRW} on the cobordism category of manifolds, so let us pause to recall some of the relevant results.

The geometric cobordism category $\Cob_\theta$, where $\theta$ is some $d$-dimensional vector bundle, has as objects closed $d-1$-dimensional $\theta$-oriented manifolds and as morphisms $\theta$-oriented cobordisms between such; for technical reasons --- that luckily do not persist to the algebraic setting --- one furthermore has to require all objects to be equipped with a chosen embedded disk $\mathrm{D}^{d-1}$, and allow only those cobordisms that are cylindrical over these embedded discs. The higher structure of this category is arranged so that $\Hom_{\Cob_\theta}(M,N)$ is the moduli space of all cobordisms between $M$ and $N$, so that (non-canonically) 
\[\Hom_{\Cob_\theta}(M,N) \simeq \bigsqcup_{W} \mathrm{BDiff}^\theta_\partial(W,\mathrm{D}^d)\]
where the notation refers to diffeomorphisms that fix the embedded copy of $\mathrm{D}^d$ in the cylindrical part, and $W$ ranges over all such diffeomorphism classes of cobordisms between $M$ and $N$.

In the simplest case Galatius and Randal-Williams then consider in \cite{GRW} (employing our numbering conventions) the subcategories 
\[\Cob_\theta^{\mc,\oc} \subseteq \Cob_\theta^{\mc} \subseteq \Cob_\theta\]
with the middle entry spanned by all morphisms $\xymatrix{W \colon M \ar@{~>}[r] & N}$ such that the inclusion $M \subseteq W$ is a $\mc$-connective, and the right by all objects $M$ with are, furthermore, $\oc$-connective. By performing what they termed parametrised surgery, they show that these inclusions realise to equivalences provided, $\mc, \oc\lesssim d/2$ and the base space of $\theta$ is suitably connected. Their method ultimately shows that for $d=2n$ the map
\[\Hom_{\Cob^{n-1}_{\theta}}(S^{2n-1},S^{2n-1})^\grp \simeq \Omega|\Cob_\theta^{n-1,n}| \longrightarrow \Omega|\Cob_{\theta}|\]
is an equivalence, whenever $\theta$ is a (once-stable) bundle over  an $(n+1)$-connective base space, as the sphere is the only $n$-connective closed $(2n-1)$-dimensional $\theta$-manifold in the base point component of $\vert \Cob_\theta\vert$.

This result can be combined with the theorem of Galatius, Madsen, Tillmann and Weiss \cite{GMTW} that 
\[|\Cob_\theta| \simeq \Omega^{\infty-1}\mathrm{MT}\theta,\]
where $\mathrm{MT}\theta$ is the Thom spectrum of the stable vector bundle $-\theta$. Applying the group completion theorem of McDuff and Segal one concludes 
\[\colim_{g \in \mathbb N} \mathrm{H}_*(\mathrm{Diff}((\mathrm{S}^{n} \times \mathrm{S}^n)^{\#g},\mathrm{D}^{2n}); \mathbb Z) \simeq \mathrm{H}_*(\Omega^{\infty}_0
\mathrm{MT}\theta_n; \mathbb Z),\]
in close correspondence with the description of the homology of $\gw^{\fpm}_{\mathrm{cl}}(R)$ in terms of orthogonal groups given above; here $\theta_n$ is the vector bundle classified by the map $\tau_{\geq n+1}\mathrm{BO}(2n) \rightarrow \mathrm{BO}(2n)$, and the homology on the left is continuous group homology.

In the present paper we show that their method can be transported to the algebraic situation at hand. Indeed, as part of \cite{9authII} we have constructed for every Poincar\'e category $(\C,\QF)$ an algebraic cobordism ($\infty$-)category $\Cob(\C,\QF)$, which satisfies
\[|\Cob(\C,\QF)| \simeq \Omega^{\infty-1}\GW(\C,\QF)\]
for the Grothendieck-Witt spectrum associated to $(\C,\QF)$ and thus in particular
\[\Omega|\Cob(\C,\QF)| \simeq \gw(\C,\QF).\] 
The objects of $\Cob(\C, \QF)$ are the Poincar\'e objects $X$ in $(\C,\QF\qshift{1})$, where we denote by $-\qshift{i}$ the shift in $\Spa$ or indeed any stable $\infty$-category. The morphisms $\xymatrix{W \colon X \ar@{~>}[r] & Y}$ are given by Poincar\'e cobordisms $X \leftarrow W \rightarrow Y$, Ranicki style. This relates to the geometric situation by applying cochains. 

For a Poincar\'e category equipped with a weight structure, we can thus consider a similar sequence of subcategories
\[\Cob^{\mc,\oc}(\C,\QF) \subseteq \Cob^{\mc}(\C,\QF) \subseteq \Cob(\C,\QF),\]
where all morphisms $X\leftarrow W \to Y$ in the middle category are required to have $W \rightarrow X$ $\mc$-connective in the weight structure, and all objects on the left $\oc$-connective. For example, if $\Dual_\QF$ preserves the heart of $\C$ as in Theorem \ref{thm:resolution} above, we find
\[\Poinc^\heart(\C,\QF) \simeq \Hom_{\Cob^{0}(\C,\QF)}(0,0)\] 
and thus
\[\Poinc^\heart(\C,\QF)^\grp \simeq \Omega|\Cob^{0,1}(\C,\QF)|,\]
since $0 \in \Cob^{0,1}(\C,\QF)$ is the only object. 

In order to state the full result we prove, let us say that a Poincar\'e category with an exhaustive weight structure has dimension $d$ if $\Dual_\QF(X)\qshift{-d} \in \C^\heart$ for every $X \in \C^\heart$. Furthermore, we need to recall that to any Poincar\'e structure is associated its linear part $\Lin_\QF \colon \C\op \rightarrow \Spa$, defined as the Goodwillie derivative of $\QF$. In a vague sense it plays the role of the bundle $\theta$. Let us for now simply mention that $\Lin_{\QF^{\g\fpm}}(P)$ is connective for all $P \in \Dperf(R)^\heart=\Proj(R)$. Replacing the geometric surgery used by Galatius and Randal-Williams with Ranicki's algebraic surgery we show:

\begin{introthm}\label{thm:main}
Let $(\C,\QF)$ be a Poincar\'e category of dimension $d$. Then the inclusions
\[\Cob^{\mc,\oc}(\C,\QF) \subseteq \Cob^{\mc}(\C,\QF) \subseteq \Cob(\C,\QF)\]
become equivalences upon realisation provided
\begin{enumerate}
\item $2\mc\leq d+1$;
\item $2\oc \leq d+1$;
\item $\oc \leq \mc+1$, and
\item $\Lin_\QF(X)$ is $(\oc-1)$-connective for every $X \in \C^\heart$.
\end{enumerate}
In particular, if $d=2\oc$,
\[|\Cob^{\oc,\oc}(\C,\QF)| \rightarrow |\Cob(\C,\QF)|\]
is an equivalence if $\Lin_\QF(X)$ is $(\oc-1)$-connective for every $X \in \C^\heart$. In this case, and if $\Lin_\QF(X)$ is even $\oc$-connective for all $X \in \C^\heart$, then furthermore the inclusion
\[\Cob^{\oc,\oc+1}(\C,\QF) \subseteq \Cob^{\oc,\oc}(\C,\QF)_0\]
becomes an equivalence upon realisation, where the subscript denotes the component of $0$.
\end{introthm}

Let us immediately remark that the index shift in the statement as compared to the result of Galatius and Randal-Williams stems from our use of cochains (rather than chains) when translating from geometry to algebra.

Via the equivalence
\[\Poinc^\heart(\C,\QF)^\grp \simeq \Omega|\Cob^{0,1}(\C,\QF)|,\]
Theorem \ref{thm:main} in particular contains Theorem \ref{thm:resolution} and thus also Theorem \ref{thm:main_special} by choosing $d=0$. The statement for $d=2\oc$ can of course be reduced to this case by considering $\C\wshift{\oc}$, the $\oc$-fold shifted weight structure on $\C$, resulting in 
\[\Poinc^\heart(\C\wshift{\oc},\QF)^\grp \simeq \gw(\C,\QF)\]
whenver $(\C,\QF)$ has dimension $2\oc$, and $\Lin_\QF$ takes $\oc$-connective values on $\C^\heart$. 

Theorem \ref{thm:main} is also interesting in odd dimensions, however: For instance, if $(\C, \QF)$ is a Poincar\'e category of dimension 0 and $\QF(X)$ is discrete for all $X \in \C^\heart$, then direct inspection shows that $\Cob^{0,0}(\C, \QF\qshift{-1})$ is an ordinary category. For $(\C,\QF) = (\Dperf(R),\QF^{\epsilon \gs})$ it is equivalent to the hermitian $\Q$-construction $\Qhcl(\PProj(R),\Hom_R(-,R),\epsilon)$ in the context of exact categories with duality from \cite{Schlichting2010}, where $\PProj(R)$ is endowed with the split exact structure; similarly in the case of a general form parameter. 

Another example is provided by hyperbolic Poincar\'e categories, where one additionally has an equivalence $\Hyp(\C)\simeq \Hyp(\C)\qshift{-1}$. We thus obtain a new, albeit rather inefficient, proof of Quillen's equivalence
\[\k(R) = \Proj(R)^\grp\simeq \Omega \vert \Qcl(\PProj(R))\vert\]
to his $\Q$-construction from \cite{QuillenHigherK}. Combining these two statements with (a rotation of) the Bott-Genauer sequence
\[\GW(\C,\QF) \longrightarrow \K(\C) \longrightarrow \GW(\C,\QF\qshift{1})\]
from \cite{9authII}, we conclude the existence of a fibre sequence
\[\Unimod^{\epsilon \s}(R)^\grp \longrightarrow \vert \Qhcl(\PProj(R), \Hom_R(-,R),\epsilon)\vert \longrightarrow \vert \Qcl(\PProj(R))\vert\]
and similarly for general form parameters. The fibre of the right hand map is classically used to define the Grothendieck-Witt space of the exact category with duality $\PProj(R)$, and for this reason the fibre sequence was long sought after: Assuming $2$ invertible in $R$, it was first established by Schlichting in \cite{schlichting-giffen}, correcting an error in an earlier attempt in \cite{CL}. In the symmetric case this extra assumption was then removed by Hesselholt and Madsen in \cite{HM}, and the case of general form parameters was finally obtained in \cite{Schlichting2019} by a rather different route.

Theorem \ref{thm:main} already has content at the level of path components. For example, assuming again that $(\C,\QF)$ is a Poincar\'e category of dimension $0$ with $\QF$ taking connective values on $\C^\heart$, we find
\[|\Cob^{-n-1,-n}(\C,\QF\qshift{-n-1})| \simeq |\Cob(\C,\QF\qshift{-n-1})|\]
for $n \geq 0$. Passing to components on both sides gives an isomorphism
\[\L_n^\mathrm{sht}(\C,\QF) \longrightarrow \L_n(\C,\QF)\]
where the left hand sides denotes the group of Poincar\'e objects in $(\C,\QF\qshift{-n})$ that lie in $\C_{[-n,0]}$ modulo bordisms that lie in $\C_{[-n-1,0]}$. Using the stronger
\[|\Cob^{0,0}(\C,\QF\qshift{-1})| \simeq |\Cob(\C,\QF\qshift{-1})|\] 
one even finds that
\[\L_0^\heart(\C,\QF) \longrightarrow \L_0(\C,\QF)\]
is an isomorphism, where we denote by $\L_0^\heart(\C,\QF)$ the abelian group obtained by taking Poincar\'e objects in the heart of $\C$ as cycles, and dividing them by the congruence relation generated by declaring all strictly metabolic forms equivalent to $0$, where a strictly metabolic form $M$ is one that admits a Lagrangian (i.e. a null bordism) which also lies in the heart of $\C$. For $(\Dperf(R),\QF^{\g\fpm})$, the left hand side is nothing but the usual Witt group of $R$ for forms of type $\fpm$. Both of these computations were obtained as part of \cite{9authIII} with more direct proofs.

As a final application we combine this analysis with Theorem \ref{thm:resolution} to compute the geometric fixed points of the direct sum real $\K$-spectrum $\kr(\mathcal A,\Dual)$ of an additive $\infty$-category with duality $(\mathcal A,\Dual)$ introduced by Heine, Lopez-Avila and Spitzweck in \cite{Spitzweckgroupcompletion}: The geometric fixed points $\kr(\mathcal A,\Dual)^{\varphi\Ct}$ are connective and for $k \geq 0$ we produce a canonical isomorphism
\[\pi_k(\kr(\mathcal A,\Dual)^{\varphi\Ct}) \longrightarrow \L^{\mathrm{sht}}_k(\Stab(\mathcal A),\QF^\sym_\Dual),\]
where $\Stab(\mathcal A)$ is the stabilisation of $\mathcal A$, which is modelled by bounded chain complexes in $\mathcal A$, whenever $\mathcal A$ is an ordinary category. \\

The paper is accompanied by an appendix. In it Yonatan Harpaz explains how an $\L$-theoretic version of Theorem \ref{thm:resolution} can be obtained by a substantial simplification of our methods. With this result in place, one can attempt to deduce Theorem \ref{thm:resolution}, independently from the body of the text, by splitting the statement into its $\K$- and $\L$-theoretic parts. He shows that this can be achieved using the existing literature on algebraic and hermitian K-theory (most notably the Gillet-Waldhausen theorem and results of Schlichting from \cite{Schlichting2019}), whenever $\C^\heart$ is an ordinary category and $\QF\colon (\C^\heart)\op \rightarrow \Spa$ takes discrete values. In particular, this provides an alternate proof of Theorem \ref{thm:main_special}.

In contrast, these inputs are rederived from Theorem \ref{thm:resolution} in the body of the text (even without the restriction to ordinary categories).

\subsection{Organisation of the paper} In Section \ref{sec:recoll} we give a brief overview of those parts of the Grothendieck-Witt theory of Poincar\'e categories  that are relevant for the present paper, in particular, expounding the genuine Poincar\'e structures associated to general form parameters. In Section \ref{sec:weight} we recount the definition of a weight structure, record its interaction with Poincar\'e structures and place all examples mentioned so far in that framework. In Section \ref{sec:para} we explain the analogues of the geometric surgery moves of Galatius and Randal-Williams in our algebraic setting, that we then use in Section \ref{sec:cplx} to give an outline of the proof of Theorem \ref{thm:main} and treat its recurring parts. The main steps of the proof then occupy Sections \ref{sec:objects} and \ref{sec:morph}, where we show that the inclusions
\[\Cob^{\oc+1}(\C,\QF)\longrightarrow \Cob^{\oc}(\C,\QF) \quad \text{and}\quad\Cob^{\mc+1,\oc}(\C,\QF) \longrightarrow \Cob^{\mc,\oc}(\C,\QF),\]
respectively, become equivalences after realisation under appropriate assumptions. The final Section \ref{sec:appl} deduces Theorems \ref{thm:main_special}, \ref{thm:resolution} and \ref{thm:main} as \ref{main_special2}, \ref{cor:res2} and \ref{thm:main2}, respectively, the two corollaries of Theorem \ref{thm:main_special} as \ref{hlgyo} and \ref{cor:2corintro}, and also explains the additional applications sketched above.

\subsection{Acknowledgements}
We first of all wish to thank Baptiste Calmès, Emanuele Dotto, Markus Land, Kristian Moi, Denis Nardin and Thomas Nikolaus for the exciting collaboration, developing in particular the original vision for algebraic cobordism categories beyond our wildest dreams, and especially Yonatan Harpaz for (in addition) inspecting our arguments very carefully, catching several oversights and of course producing the appendix. 

We also thank Mauricio Bustamante, Dustin Clausen, Bert Heuts, Manuel Krannich, Jacob Lurie, Vikram Nadig, Nathan Perlmutter, Georgios Raptis, Ulrike Tillmann, Michael Weiss, Christoph Winges and Lior Yanovski for a number of illuminating discussions.

Finally, we acknowledge a tremendous intellectual debt to S\o{}ren Galatius and Oscar Randal-Williams, and want to express our gratitude for their continued support and interest in this project. \\

During the preparation of this article FH was a member of the Hausdorff Center for Mathematics at the University of Bonn funded by the German Research Foundation (DFG), grant no.\ EXC 2047 390685813, and received support from the European Research Council (ERC) through the grant ``Moduli spaces, Manifolds and Arithmetic'', grant no.\ 682922. Both authors were supported by the Engineering and Physical Sciences Research Council through the program ``Homotopy harnessing higher structures'' at the Isaac Newton Institute for Mathematical Sciences, grants no.\ EP/K032208/1 and EP/R014604/1. WS was further supported by the DFG through the priority program ``Geometry at Infinity'' at the University of Augsburg, grant no.\ SPP 2026. YH was supported by the ERC through the grant ``Foundations of Motivic Real K-Theory'', grant no.\ 949583, and by the French National Research Agency  through the grant “Chromatic Homotopy and K-theory”, grant no. 16-CE40-0003.

During the revision process FH was further supported by the DFG through the centre ``Integral structures in Geometry and Representation Theory", grant no.\ TRR 358-4913924, at the University of Bielefeld.

\subsection{Notations and conventions}
We shall frequently refer to $\infty$-categories simply as categories when there is no risk of confusion; we will explicitly say that $\C$ is an ordinary category whenever we need to make the distinction. We will denote the shift functor, i.e. suspension, in stable $\infty$-categories by $(-)\qshift{1}$. We use $\Hom_\C$ to refer to the \emph{space} of homomorphisms, and $\hom_\C$ for its spectral enhancement, in case $\C$ is stable, and $\grpcr(\C)$ denotes the groupoid core of a category $\C$. The $\infty$-category of spaces/$\infty$-groupoids will be denoted $\Sps$ (for animae, following a suggestion of Clausen and Scholze), and we will use all three terms interchangeably.

\section{Recollections on Poincar\'e categories}\label{sec:recoll}

In the present section we provide a summary of the parts of \cite{9authI,9authII} that we will need throughout the paper, and then give an account of the examples that we will consider. Knowledge of the examples is not needed for the proof of Theorems \ref{thm:resolution} or \ref{thm:main}; outside illustrations they will only be required in Section \ref{sec:appl}. We include them here with the aim of providing a coherent account and to give the reader a flavour of the theory.
Let us also mention that we have separated out the treatment of algebraic surgery, which we recall in greater detail in Section \ref{sec:para}.

\subsection{Poincaré categories and Poincaré objects}
A hermitian structure on a small stable $\infty$-category $\C$ is a reduced, quadratic functor $\QF \colon \C\op \rightarrow \Spa$. A pair $(\C,\QF)$ consisting of this data is called a hermitian category. These organise into an $\infty$-category $\Cath$ whose morphisms consist of hermitian functors, that is pairs $(f,\eta)$ where $f \colon \C \rightarrow \C'$ is an exact functor and $\eta \colon \QF \Rightarrow \QF' \circ f\op$ is a natural transformation, see \cite[Section 1.2]{9authI}.

To such a hermitian category is associated its category of hermitian forms $\catforms(\C,\QF)$: Objects are pairs $(X,q)$ where $X \in \C$ and $q$ is a $\QF$-hermitian form on $X$, i.e.\ a point in $\Omega^\infty \QF(X)$. Morphisms are maps in $\C$ preserving the hermitian forms. We shall denote the groupoid core of $\catforms(\C,\QF)$ by $\spcforms(\C,\QF)$. In order to impose a non-degeneracy condition on the forms in $\catforms(\C,\QF)$, one needs a non-degeneracy condition on the hermitian \(\infty\)-category $(\C,\QF)$ itself. To this end recall the classification of quadratic functors from Goodwillie calculus, see \cite[Section 1.3]{9authI}: Any reduced quadratic functor uniquely extends to a cartesian diagram
\[\xymatrix{\QF(X) \ar[r] \ar[d] & \Lin_\QF(X) \ar[d]\\
\Bil_\QF(X,X)^\hC \ar[r] & \Bil_\QF(X,X)^\tC
}\]
where $\Lin_\QF \colon \C\op \rightarrow \Spa$ is linear (i.e.\ exact) and $\Bil_\QF \colon \C\op \times \C\op \rightarrow \Spa$ is bilinear (i.e.\ exact in each variable) and symmetric (i.e.\ comes equipped with a refinement to an element of $\Fun(\C\op \times \C\op,\Spa)^\hC$, with $\Ct$ acting by flipping the input variables).

A hermitian structure $\QF$ is called Poincaré if there exists an equivalence $\Dual \colon \C\op \rightarrow \C$ such that
\[
\Bil_\QF(X,Y) \simeq \hom_\C(X,\Dual Y)
\]
naturally in $X,Y \in \C\op$. By Yoneda's lemma, such a functor $\Dual$ is uniquely determined if it exists, so we refer to it as $\Dual_\QF$. By the symmetry of $\Bil_\QF$ the functor $\Dual_\QF$ then automatically satisfies $\Dual_\QF \circ \Dual_\QF\op \simeq \id_\C$. Any hermitian functor $(F,\eta) \colon (\C,\QF) \rightarrow (\C',\QF')$ between Poincaré-categories (i.e.\ hermitian categories whose hermitian structure is Poincaré) induces a tautological map
\[
F \circ \Dual_\QF \Longrightarrow \Dual_{\QF'} \circ F\op,
\]
and we say that $(F,\eta)$ is a Poincaré functor if this transformation is an equivalence. Poincaré categories together with Poincaré functors form a (non-full) subcategory $\Catp$ of $\Cath$. \\

Now, if $(\C,\QF)$ is Poincaré, then to any hermitian form $(X,q) \in \catforms(\C,\QF)$ there is tautologically associated a map
\[
q^\sharp \colon X \longrightarrow \Dual_\QF X
\]
as the image of $q$ under
\[
\Omega^\infty \QF(X) \longrightarrow \Omega^\infty\Bil_\QF(X,X) \simeq \Hom_\C(X,\Dual_\QF X)
\]
and we say that $(X,q)$ is Poincaré if $q^\sharp$ is an equivalence, see \cite[Section 2.1]{9authI}. The full subspace of the core of $\catforms(\C,\QF)$ spanned by the Poincaré forms is denoted by $\Poinc(\C,\QF)$ and provides a functor
\[
\Poinc \colon \Catp \rightarrow \Sps.
\]

\subsection{The hermitian $\Q$-construction}

To a Poincar\'e category are associated its Grothendieck-Witt space or spectrum and an algebraic cobordism category. Both are defined using the hermitian $\Q$-construction which we now recall. For $K \in \Cat$ let $\Twar(K)$ denote the twisted arrow category of $K$. We define $\Q_K(\C) \subseteq \Fun(\Twar(K),\C)$ for stable $\C$ to be the full subcategory spanned by all diagrams $F$ for which the squares
\[\xymatrix{F(i \rightarrow l) \ar[d]\ar[r] & F(j \rightarrow l) \ar[d]\\
            F(i \rightarrow k) \ar[r] & F(j \rightarrow k)}\]
are (co)cartesian for every functor $[3] \rightarrow K$ with value $i \rightarrow j \rightarrow  k \rightarrow l$; for example this pertains to the middle square in $\Twar[2]$:
\[\xymatrix@-1pc{ & & F(0 \leq 2) \ar[ld]\ar[rd]& & \\
 & F(0 \leq 1) \ar[ld]\ar[rd]& & F(1 \leq 2) \ar[ld]\ar[rd]& \\
F(0 \leq 0) && F(1 \leq 1) && F(2 \leq 2)}\]
and by pasting pullbacks it always suffices to consider chains of the form $i\to k \xrightarrow\id k\to l$. If $\C$ is equipped with a hermitian structure $\QF$ then
\[\QF^{\Twar(K)} \colon \Fun(\Twar(K),\C)\op \longrightarrow \Spa, \quad F \longmapsto \lim_{\Twar(K)\op} \QF \circ F\op\]
defines a hermitian structure on $\Fun(\Twar(K),\C)$. If $\QF$ is Poincar\'e, then so is the restriction of $\QF^{\Twar(K)}$ to $\Q_K(\C)$ (though $\QF^{\Twar(K)}$ is not even Poincar\'e for $K = [2]$). The construction is functorial in $K$ yielding the desired hermitian $\Q$-construction
\[\Q \colon \Cat\op \times \Catp \longrightarrow \Catp,\]
see \cite[Section 2.2]{9authII}. It preserves limits in both variables, so fixing the input Poincar\'e category to be $(\C,\QF)$, its restriction to $\bbDelta\op \subseteq \Cat\op$ defines a complete Segal object of $\Catp$, and thus $\Poinc\Q(\C,\QF) \colon \bbDelta\op \rightarrow \Sps$ is a complete Segal space, since $\Poinc\colon \Catp \rightarrow \Sps$ also preserves limits. 

Denote now by
\[\asscat \colon \sSps \longleftrightarrow \Cat \cocolon \nerv\]
Rezk's nerve adjunction (i.e $\asscat$ is the colimit preserving extension of the inclusion $\bbDelta \subseteq \Cat$ and $\nerv$, the nerve, its right adjoint), whose right adjoint is fully faithful with essential image the complete Segal spaces see e.g.\ \cite{JoyalTierney, HSteine}. We note that for a simplicial space $X$, its realisation $|X|\coloneqq \colim_{\bbDelta\op} X$ is canonically equivalent to the realisation $|\asscat X|$ of the associated category since adjoints compose and spaces have constant nerve (the realisation functor of categories being the left adjoint to the inclusion $\Sps\to \Cat$).
We put
\[\Cob(\C,\QF) = \asscat \Poinc\Q(\C,\QF\qshift{1}).\]
Unwinding the definitions objects of $\Cob(\C,\QF)$ are exactly Poincar\'e objects in $(\C,\QF\qshift{1})$, and a morphism $\xymatrix{(X,q) \ar@{~>}[r] & (X',q')}$ is an algebraic cobordism \`a la Ranicki, i.e.\ given by a span
\[\xymatrix{X & \ar[l]_f W \ar[r]^g & X'}\]
together with an equivalence $\eta \colon f^*q \sim g^*q'$ in $\Omega^\infty\QF\qshift{1}(W)$, which satisfies Lefschetz duality in the sense that the induced map
\[\eta_\# \colon \fib(W \rightarrow X) \longrightarrow \Dual_\QF(\fib(W \rightarrow X'))\]
is an equivalence. Here indexing conventions adhere to the geometric definitions where $\Cob_d$ has $d$-dimensional cobordisms as morphisms. In particular, we find $\pi_0|\Cob(\C,\QF)|$ given by the cobordism group of Poincar\'e objects in $(\C,\QF\qshift{1})$, which is isomorphic to $\L_0(\C,\QF\qshift{1}) = \L_{-1}(\C,\QF)$. For the purposes of this paper we shall treat it as the definition of the $\L$-groups on the left (and more generally of $\L_n(\C,\QF\qshift{n+1})$), and refer the reader to \cite[Section 2.3]{9authII} for a proof that this agrees with the Witt-style definition, giving by taking the quotient of $\pi_0(\Poinc(\C,\QF\qshift{1}))$ by metabolic forms (which we will recall below). Finally, we recall that there is a tautological equivalence
\[\Poinc\Q_K(\C,\QF\qshift{1}) \simeq \Hom_{\Cat}(K,\Cob(\C,\QF))\]
that we shall exploit throughout.\\

The Grothendieck-Witt space of $(\C,\QF)$ can be defined in a number of ways. Closest in spirit to classical constructions for exact categories is
\[\gw(\C,\QF) = \fib(|\Poinc\Q(\C,\QF)| \xrightarrow{\fgt} |\core\Q(\C)|).\] 
but as part of the Genauer-Bott sequence described at the end of this subsection we also established an equivalence
\[\gw(\C,\QF) = \Omega|\Cob(\C,\QF)|\]
in \cite[Section 4.1]{9authII} and showed that $\gw \colon \Catp \rightarrow \Sps$ is the initial functor equipped with a transformation from $\Poinc \colon \Catp \rightarrow \Sps$ that is additive and group-like; this is arguably the cleanest definition, but as we won't need it we refrain from expounding the meaning of these terms here. 

Any of the definitions of the Grothendieck-Witt space can be extended to that of a Grothendieck-Witt spectrum $\GW(\C,\QF)$ with
\[\Omega^\infty \GW(\C,\QF) \simeq \gw(\C,\QF) \quad \text{and} \quad \Omega^{\infty-i}\GW(\C,\QF) \simeq |\Poinc \Q^{(i)}(\C,\QF\qshift{i})|\]
for all $i > 0$, where $\Q^{(i)}\colon \Catp\to \Fun((\bbDelta\op)^i,\Catp)$ denotes the $i$-fold iteration of $\Q$, see \cite[Section 4.2]{9authII}. In particular, we find
\[|\Cob(\C,\QF)| \simeq \Omega^{\infty-1}\GW(\C,\QF)\]
in analogy with the theorem of Galatius, Madsen, Tillmann and Weiss. While this spectrum, which in contrast to the situation for $\K(\C)$ is generally not bounded below, will play essentially no role in the present paper, the iterated $\Q$-construction on the right hand side does. We note that 
\[\Poinc\Q^{(i)}(\C,\QF\qshift{i}) \colon (\bbDelta\op)^i \longrightarrow \Sps\]
is an $i$-fold complete Segal space, and thus models an $i$-fold (or $i$-tuple) category. We call it $\Cob^{(i)}(\C,\QF)$ in analogy with the higher cobordism categories in manifolds topology; of course $i$-fold complete Segal spaces are a perfectly good \emph{definition} of $i$-fold categories in which case the above definition merely amounts to an abbreviation. In the double category $\Cob^{(2)}(\C,\QF)$, a morphism square is given by a diagram
\[\xymatrix@-1pc{X & W \ar[r]\ar[l]& X' \\
            V \ar[u]\ar[d]& \ar[r]\ar[l]\ar[u]\ar[d]U & V' \ar[u]\ar[d]\\
            Y & W' \ar[r]\ar[l]& Y'}\]
in which the corners are Poincar\'e objects of $(\C,\QF\qshift{2})$, the outer spans are cobordisms between these, which can be glued to a single Poincar\'e object in $(\C,\QF\qshift{1})$, of which $U$ is then a nullbordism. \\

In a different direction, restricting the Poincar\'e structure on $\Q_1(\C,\QF)$ to the fibre of the functor
\[d_0 \colon \Q_1(\C) \longrightarrow \C\]
gives a Poincar\'e structure on the arrow category $\mathrm{Ar}(\C)$, which we call the metabolic Poincar\'e category $\Met(\C,\QF)$ associated to $(\C,\QF)$. Its Poincar\'e objects can be thought of either as the equivalent of unimodular forms equipped with a Lagrangian, or as Poincar\'e objects with boundary in $(\C,\QF\qshift{-1})$, see \cite[Section 2.3]{9authI}. The latter point of view nicely fits with the null-composite sequence 
\[(\C,\QF\qshift{-1}) \longrightarrow \Met(\C,\QF) \xrightarrow{\mathrm{met}} (\C,\QF)\]
of Poincar\'e categories, where $\mathrm{met}$ is the restriction of $d_1 \colon \Q_1(\C,\QF) \rightarrow (\C,\QF)$ and the left hand arrow takes $X$ to $X \rightarrow 0$. Applying the Grothendieck-Witt functor results in the fibre sequence
\[\GW(\C,\QF\qshift{-1}) \xrightarrow{\fgt} \K(\C) \xrightarrow{\hyp} \GW(\C,\QF)\]
of spectra, which we termed the Bott-Genauer sequence, as it is analogous to Genauer's fibre sequence form geometric topology, see \cite[Section 4.3]{9authII}.

\subsection{Examples of Poincaré categories}
As mentioned above the material of this somewhat lengthy subsection is not strictly speaking necessary to read most of the present paper, but we shall freely reference the examples explained here throughout the text.

We start out by discussing hermitian structures given as the animations of form parameters in the sense of Schlichting \cite{Schlichting2019}, which includes form parameters in the sense of Bak \cite{Bak} as a special case. Given a ring $R$, a form parameter $\fpm$ consists of an $R \otimes_\mathbb Z R$-module $M$, equipped with an involution $\sigma$ that is linear over the flip of $R \otimes_\mathbb Z R$, together with a factorisation 
\[M_\Ct \xrightarrow{\tau} Q \xrightarrow{\rho} M^\Ct\]
of the norm map $[m] \mapsto m + \sigma(m)$ that satisfies condition ($*$) below and is quadratic over $R$ in the following sense: $Q$ is an abelian group equipped with an action of the multiplicative monoid of $R$, such that multiplication by any $r \in R$ is an additive map $Q \rightarrow Q$, but such that in lieu of distributivity one only has that 
\[\Bil(r,s)(x) = (r+s)x - rx - sx\]
defines a bilinear map $R \times R \rightarrow \Hom_\Ab(Q,Q)$. $M_\Ct$ and $M^\Ct$ also carry such actions by restricting the $R \otimes_\mathbb Z R$-module structure on $M$ along the diagonal
\[R \longrightarrow R \otimes R, \quad r \longmapsto r \otimes r.\]
In addition to being additive and $R$-equivariant for these actions we require that $\tau$ and $\rho$ satisfy
\begin{equation*}\tag{$\ast$}
\Bil(r,s)(q) = \tau[(r \otimes s)\rho(q)],
\end{equation*}
which in particular prescribes the composition $\tau \circ \rho \colon Q \rightarrow Q$. This requirement is geared to make the notion of a $\fpm$-form work: It consists of a finitely generated projective $R$-module $P$ together with a $\Ct$-equivariant  $R \otimes R$-linear map $b \colon P \otimes P \rightarrow M$ and an $R$-equivariant function $q\colon P \rightarrow Q$ that satisfy
\[q(p+p') -q(p)-q(p') = \tau[b(p,p')] \quad \text{and}\quad \rho q(p) = b(p,p).\]
For commutative $R$ the examples of $M$-valued (skew-)symmetric and (skew-)quadratic forms for some $R$-module $M$ correspond to the choices $\rho = \id_{M^\Ct}$ and $\tau = \id_{M_\Ct}$, respectively, with $\sigma = \pm \id_M$ and the $R \otimes_\mathbb Z R$-module structure on $M$ induced by the multiplication of $R$. A ring with skew-involution similarly fits into this framework by taking $M=R$ 
and letting the skew-involution provide both the $R \otimes_\mathbb Z R$-module structure as well as $\sigma$.

Sending a projective module $P$ to the abelian group of $\fpm$-forms begets the functor
\[\mathrm F^\fpm \colon \PProj(R)\op \longrightarrow \Ab\]
that already appeared in the introduction. It uniquely extends to a hermitian structure $\QF^{\g\fpm}_M$ on $\Dperf(R)$ that makes the diagram
\[\xymatrix{
\PProj(R)\op \ar[r]^-{\mathrm F_\fpm}\ar[d] & \Ab \ar[d]^{\mathrm H} 
\\
\Dperf(R)\op \ar[r]^-{\QF^{\g\fpm}_M} & \Spa
}
\]
commute, see \cite[Section 4.2]{9authI}. This hermitian structure is Poincar\'e if and only if $M$ is perfect as an $R$-module, when pulled back along the left, say, inclusion $R \rightarrow R \otimes_\mathbb Z R$, and in addition the canonical map
\[R \longrightarrow \mathbb R\Hom_R(M,M), \quad r \longmapsto (1 \otimes r) \cdot -\]
is an equivalence. In this case the duality associated to $\QF^{\g\fpm}_M$ is then given by \[\Dual_{\QF^{\g\fpm}_M} = \mathbb R\Hom_R(-,M).\] For most of our results, in particular for Theorem \ref{thm:resolution} to be applicable to $(\Dperf(R),\QF_M^{\g\fpm})$, we will furthermore have to assume that $M$ is not just perfect but projective as an $R$-module, and for notational ease will therefore refer to $\lambda$ as invertible only if it satisfies this stronger condition (we will discuss this point further in \ref{examplesofpoincweight} below). In the simplest example of a commutative ring $R$ and the two $R$-module structures on $M$ agreeing, this exactly amounts to $M$ being $\otimes_R$-invertible. At any rate, to put the description of the linear part $\Lambda_{\QF^{\g\fpm}_M}$ into context, let us briefly recall the classification of hermitian structures on module categories from \cite[Section 3.2]{9authI}. \\

Fix an $\Eone$-algebra $A$ over a base $\Einf$-ring spectrum $k$ and a subgroup $c \subseteq \K_0(A)$. Consider then the category of compact $A$-module spectra $\Modp_{A}$ or more generally its full subcategory $\Mod^c_A$ spanned by all those $X \in \Modp_{A}$ with $[X] \in c \subseteq \K_0(A)$. For example $\Mod^{\langle A \rangle}_A$ is the stable subcategory of $\Mod(A)$ spanned by $A$ itself. For the reader mostly interested in the applications to discrete rings as discussed above, we recall that any discrete ring $R$ gives rise to such data, via the Eilenberg-Mac Lane functor $\mathrm{H} \colon \Ab \longrightarrow \Spa$, which is lax symmetric monoidal and therefore induces a functor
\[
\mathrm{Ring} \longrightarrow \mathrm{Alg}_{\Eone}(\Mod_{\mathrm{H}\mathbb Z})
\] 
In this way, any discrete ring may be regarded as an $\Eone$-algebra over $\mathrm{H}\mathbb Z$. There are, furthermore, equivalences
\[
\Modp_{\mathrm{H} R} \simeq \Dperf(R) \quad \text{and} \quad \Mod^{\langle \mathrm H R \rangle} _{\mathrm{H} R} \simeq \D^\mathrm f(R),
\]
where $\Dperf(R)$ denotes the full subcategory of the derived $\infty$-category $\D(R)$ of $R$ spanned by the perfect complexes, i.e.\ finite chain complexes of finitely generated projective $R$-modules and $\D^\mathrm f(R)$ is the full subcategory spanned by the finite chain complexes of finite free $R$-modules. In this regime the reader should keep in mind, that terms such as $\otimes_{\mathrm{H} \mathbb Z}$ or $\Hom_{\mathrm{H} R}$ will evaluate to the functors $\otimes^\mathbb L_\mathbb Z$ and $\mathbb R\Hom_R$. \\

Hermitian structures on the categories $\Mod^c_A$ arise from $k$-modules with genuine involution $(M,N,\alpha)$ over $A$  \cite[Definition 3.2.3]{9authI}: The first entry $M$ plays the same role as $M$ did above; it is an $A \otimes_k A$-module, equipped with the structure of a homotopy fixed point in the category $\Mod_{A \otimes_k A}$ under the $\Ct$-action flipping the two factors.

The additional data of a module with genuine involution consists of an $A$-module spectrum $N$, and an $A$-linear map $\alpha \colon N \rightarrow M^\tC$; to make sense of the latter term, note that upon forgetting the $A \otimes_k A$-action, the involution equips $M$ with the structure of a (naive) $\Ct$-$k$-module spectrum. The $k$-module spectrum $M^\tC$ then becomes an $(A \otimes_k A)^\tC$-module via the lax monoidality of the Tate construction and from here obtains an $A$-module structure on $M^\tC$ by pullback along the Tate diagonal $A \rightarrow (A \otimes_k A)^\tC$, which is a map of $\Eone$-ring spectra, see \cite[Chapter III.1]{NS} for an exposition of the Tate diagonal in the present language. Let us immediately warn the reader that the Tate diagonal is not generally $k$-linear for the $k$-module structure on $(A \otimes_k A)^\tC$ arising from the unit map $k \rightarrow (k \otimes_k k)^\tC = k^\tC$, as this map is usually different from the Tate-diagonal of $k$ (in particular, this is the case for $k = \mathrm{H}\mathbb Z$ by \cite[Theorem III.1.10]{NS}). 

Even if only interested in discrete $R$, one therefore has to leave not only the realm of discrete $R$-modules to form the Tate construction, but even the realm of derived categories, as no replacement for the Tate diagonal can exist in that regime. In fact, $k$ really only plays an auxiliary role in the entire discussion above and for the purposes of the present paper it is best to take $k = \mathbb S$, as in that case the following procedure captures all 
hermitian structures on $\Mod^c_A$. \\

The hermitian structure associated to a module with genuine involution $(M,N,\alpha)$ is given by the pullback
\[\xymatrix{
\QF_M^\alpha(X) \ar[rr] \ar[d] & & \hom_A(X,N) \ar[d]^{\alpha_*} \\
  \hom_{A \otimes_k A}(X \otimes_k X,M)^\hC \ar[r] & \hom_{A \otimes_k A}(X \otimes_k X,M)^\tC \ar[r] & \hom_A(X,M^\tC)
}\]
where the $\Ct$-action on $\hom_{A \otimes_k A}(X \otimes_k X,M)$ is given by flipping the factors in the source and the involution on $M$. It is a Poincaré structure on $\Modp_A$ if $M$ restricts to an object of $\Modp_{A}$ under either inclusion $A \rightarrow A \otimes_k A$, and furthermore $M$ is invertible, i.e.\ the natural map
\[
A \rightarrow \hom_A(M,M)
\]
is an equivalence; again we shall tacitly assume that this condition is satisfied whenever we speak of modules with involution over $A$ in the sequel. In this case the associated duality is given by $X \mapsto \hom_A(X,M)$ regarded as an $A$-module via the extraneous $A$-module structure on $M$. Given a subgroup $c \in \K_0(A)$ one obtains a Poincar\'e structure on $\Mod^c_A$ if in addition $c$ is closed under the duality on $\K_0(A)$ induced by $M$. In case $c = \K_0(A)$ or $c=0$ this condition is of course vacuous, and for $c = \langle A \rangle$ it simply translates to $[M]\in \langle A \rangle$.  \\

The hermitian structure $\QF^{\g\fpm}_M$ for $\fpm = (Q,\rho,\tau,\sigma)$ a form parameter with underlying module $M$ then corresponds to the $\mathrm H\mathbb Z$-module with genuine involution over $\mathrm HR$
\[\g\lambda = (M,\mathrm{cof}(M_\hC \rightarrow M_\Ct \xrightarrow{\tau} Q), \cof(M_\hC \rightarrow Q) \xrightarrow{\rho} \cof(M_\hC \rightarrow M^\Ct) \rightarrow M^\tC);\]
the underlying spectra are easily determined by evaluating the fibre sequence
\[\Bil_\QF(X,X)_\hC \longrightarrow \QF(X) \longrightarrow \Lin_\QF(X)\]
for $X = R \in \Dperf(R)$, but it is somewhat tricky to lift $\mathrm{cof}(M_\hC \xrightarrow{\tau} Q)$ and its map to $M^\tC$ from $\Spa$ to $\D(R)$ a priori as their constituents are not themselves $R$-linear, see \cite[Section 4.2]{9authI}.

From this formula we for example find $\QF^{\g\q}_M$ and $\QF^{\g\s}_M$ associated to
\[(M,\tau_{\geq 2}M^\tC, \mathrm{can}) \quad \text{and} \quad (M,\tau_{\geq 0}M^\tC, \mathrm{can}),\]
respectively, where $\mathrm{can}$ is the canonical map from the connective cover. These structures extend to an arbitrary connective $\Eone$-ring $A$, where we define $\QF^{\geq i}_M$ associated to $(M,\tau_{\geq i} M^\tC,\mathrm{can})$. This includes the cases $i \in \{\pm \infty\}$, where we obtain the homotopy symmetric and quadratic Poincar\'e structures
\[\QF^\q_M(X) = \hom_{A \otimes A}(X \otimes X,M)_\hC  \quad \text{and} \quad \QF^\sym_M = \hom_{A \otimes A}(X \otimes X,M)^\hC\]
corresponding to 
\[(M,0,0) \quad \text{and} \quad (M,M^\tC,\id),\]
respectively. Their Poincar\'e objects (in the case of a discrete ring) where originally introduced into surgery theory by Ranicki as the base objects for higher $\L$- or Witt groups, compare \cite{Ranickiblue}, and then reworked into the present language by Lurie in \cite{Lurie-L}. Let us warn immediately, that if $R$ is a discrete ring and $2$ is not assumed invertible in $R$, the functor $\QF^\sym_M$ will essentially never satisfy the connectivity conditions that we have to assume in our main results, whereas our methods give information about each of the finite stages. \\

Another very interesting example, that our results apply to, is the universal Poincar\'e structure $\QF^\mathrm u$ on the category of finite spectra, corresponding to the module with genuine involution $(\mathbb S, \mathbb S, u)$, where $u$ is the unit of the $\Einf$-ring $\mathbb S^\tC \simeq \mathbb S_2^\wedge$, see \cite[Section 4.1]{9authI}. This Poincar\'e structure was originally considered by Weiss and Williams in \cite{WW3} (in different language of course) during their pursuit of a unification of surgery and concordance theory, and generalises to arbitrary $B \in \Sps$, as in the case of Waldhausen's $\mathrm A$-spectra, as follows: Consider the categories $\Spa/B = \Fun(B,\Spa)$ for some $B \in \Sps$. Entirely parallel to the discussion above, one can derive hermitian structures on the compact objects of $\Spa/B$ from triples $(M,N, \alpha)$ with $M \in (\Spa/B \times B)^\hC$ and $\alpha \colon N \rightarrow (\Delta^*M)^\tC$ a map in $\Spa/B$, where $\Delta \colon B \rightarrow B \times B$ is the diagonal, compare \cite[Section 4.4]{9authI}. The most important examples of such functors are the visible Poincaré structures $\QF_\xi^\vis$ given by the triples 
\[
(\Delta_! \xi, \xi, u \colon \xi \rightarrow (\Delta^*\Delta_!\xi)^\tC),
\]
where $\xi \colon B \rightarrow \mathrm{Pic}(\mathbb S) \subset \Spa$ is some stable spherical fibration over $B$, $\Delta_! \colon \Spa/B \rightarrow \Spa/(B \times B)$ is the left adjoint to $\Delta^*$ and $u$ is the unit of this adjunction (which factors through $\xi \rightarrow (\Delta^*\Delta_!\xi)^\hC$ since $\Delta$ is invariant under the $\Ct$-action on $B \times B$). These hermitian structures are automatically Poincaré with associated duality given by 
\[
X \longmapsto \hom_B(X,\Delta_!\xi),
\]
the Costenoble-Waner duality functor twisted by $\xi$. The case of the universal Poincar\'e category is contained in this discussion as the case $B = \ast$ and $\xi = \mathbb S$. In \cite[Section 4.7]{9authII} we produced equivalences
\[\GW((\Spa/B)^\omega, \QF^\lambda_\xi) \simeq \mathrm{LA}^\lambda(B,\xi)\]
extending Waldhausen's equivalence
\[\K((\Spa/B)^\omega) \simeq \mathrm{A}(B)\]
for $\lambda \in \{\q,\s,\vis\}$, where the top right hand side denotes the spectra constructed by Weiss and Williams and the lower right hand side Waldhausen's algebraic $\K$-spectrum of the space $B$. \\

Finally, we mention once more the hyperbolic categories $\Hyp(\C) = \C\op \times \C$ for $\C$ a stable $\infty$-category, with Poincar\'e structure given by $\hom_\C \colon \C\op \times \C \rightarrow \Spa$. The duality in this case is given by $(X,Y) \mapsto (Y,X)$ and the diagonal induces an equivalence
\[\Poinc(\Hyp(\C)) \simeq \grpcr(\C),\]
through which all of our results on Grothendieck-Witt spaces contain their $K$-theoretic counterparts, since essentially per construction 
\[\Cob(\Hyp(\C)) \simeq \Span(\C)\]
and thus $\GW(\Hyp(\C)) \simeq \K(\C)$, compare \cite[Section 2.2]{9authI}.

\section{Weight structures on Poincar\'e categories}\label{sec:weight}

In this section we briefly recall the definition of a weight structure on a stable $\infty$-category. We then consider weight structures on Poincar\'e categories, discuss the examples mentioned in the introduction, and record a few basic interactions.

\subsection{Weight structures on stable $\infty$-categories}\label{subsec:weightstable}

Weight structures were originally introduced for triangulated categories by Bondarko in \cite{Bondarko2010weights}, and a detailed discussion on stable $\infty$-categories is for example given in \cite{Sosniloweight}. We recall the definition and basic properties for the reader's convenience.

\begin{definition}\label{def:weightstr}
A \emph{weight structure} on a stable $\infty$-category consists of two full subcategories $\C_{[0,\infty]}$ and $\C_{[-\infty,0]}$ of $\C$ such that the following conditions hold:
\begin{enumerate}
\item[i)] $\C_{[0,\infty]}$ and $\C_{[-\infty,0]}$ are both closed under retracts in $\C$,
\item[ii)] for $X \in \C_{[-\infty,0]}$ and $Y \in \C_{[0,\infty]}$ the mapping spectrum $\hom_\C(X,Y)$ is connective (i.e. has trivial negative homotopy groups) and
\item[iii)] for every $X \in \C$, there is a fibre sequence
\[Y \rightarrow X \rightarrow Z\]
with $Y \in \C_{[-\infty,0]}$ and $Z^{[-1]} \in \C_{[0,\infty]}$.
\end{enumerate}
For $a \leq b \in \mathbb{Z}$ we set $\C_{[a,\infty]}$ and $\C_{[-\infty,b]}$ to be the full subcategories of $\C$ spanned by all objects $X$ with $X^{[-a]} \in \C_{[0,\infty]}$ or $X^{[-b]} \in \C_{[-\infty,0]}$, respectively, and put
\[\C_{[a,b]} = \C_{[a,\infty]} \cap\C_{[-\infty,b]}.\]
The category $\C_{[0,0]} = \C^\heart$ is the \emph{heart} of the weight structure. Objects of $\C_{]-\infty,\infty[} = \bigcup_{n \in \mathbb Z} \C_{[-n,n]}$ are said to be \emph{bounded}, and the weight structure is \emph{exhaustive} if $\C = \C_{]-\infty,\infty[}$.

We call a stable $\infty$-category equipped with a weight structure a \emph{weighted category}. An exact functor $F \colon \C \rightarrow \D$ between weighted categories we call \emph{bounded below} or \emph{above} by $d \in \mathbb Z$, if 
\[F(\C_{[0,\infty]}) \subseteq \D_{[d,\infty]} \quad \text{or} \quad F(\C_{[-\infty,0]}) \subseteq \D_{[-\infty,d]},\]
respectively. We say $F$ has \emph{degree} $d$ if it is both bounded below and above by $d$.
\end{definition}

Note the order reversal of the terms in item iii) as compared to a $t$-structure. As will be visible in the examples below, this has the effect that weight decompositions as in iii) are not unique.

If $X \in \C_{[-\infty,b]}$ and $Y \in \C_{[a,\infty]}$, then the spectrum $\hom_\C(X,Y)$ is $(a-b)$-connective. In particular,  
\[\C_{[-\infty, a]} \cap \C_{[b, \infty]}=0 \quad \mathrm{if}\quad a<b.\]

\begin{lemma}
We have
\[\{X \in \C \mid \hom_\C(X,Y) \text{ is connective for all } Y \in \C_{[b,\infty]}\} = \C_{[-\infty,b]}\]
and dually
\[\{X \in \C \mid \hom_\C(Y,X) \text{ is connective for all } Y \in \C_{[-\infty,a]}\} = \C_{[a,\infty]}.\]
In particular, $\C_{[a,\infty]}$ is closed under pushouts and $\C_{[-\infty,b]}$ under pullbacks, and all $\C_{[a,b]}$ are closed under extensions and thus additive. 
\end{lemma}

\begin{proof}
 The inclusion of the right hand sides into the left holds by axiom ii), and for the converse in the first statement pick a decomposition  $Y \rightarrow X \rightarrow Z$ with $Y \in \C_{[-\infty,b]}$ and $Z \in \C_{[b+1,\infty]}$ using iii), and note that per assumption on $X$ its second map is null-homotopic. Thus the rotated fibre sequence
\[Z\qshift{-1} \longrightarrow Y \longrightarrow X\]
splits, displaying $X$ as a retract of $Y$, which suffices for the claim by axiom i). 
\end{proof}

We will refer to objects in $\C_{[n,\infty]}$ as $n$-connective as this will fit our examples, and dually we shall call objects in $\C_{[-\infty,n]}$ $n$-coconnective. We warn the reader explicitly, that in most examples $n$-coconnectivity is neither directly related to $n$-connectivity in a $t$-structure with the same connective part (which may or may not exist) nor to the notion of $n$-truncated objects in $\C$. Furthermore, we call a map $X \rightarrow Y$ an $n$-connective if its fibre lies in $\C_{[n,\infty]}$, and $n$-coconnective if its cofibre lies in $\C_{[-\infty,n]}$.

Finally, let us mention that an exhaustive weight structure is entirely determined by its heart, which can be an arbitrary full subcategory $H \subseteq \C$ that is closed under taking finite sums and retracts, has $\C$ as its stable envelope and whose morphism spectra are connective: One inductively sets $H_{0,k+1}$ to consist of all objects $X$ which can be written in an extension 
\[Y \rightarrow X \rightarrow Z\]
with $Y \in H_{0,k}$ and $X^{[-k-1]} \in H$, and then 
\[\C_{[0,\infty]} = \bigcup_{k \in \mathbb N} H_{0,k},\]
and similarly for $\C_{[-\infty,0]}$. It is then an elementary exercise to check that these definitions give an exhaustive weight structure, that the construction really is inverse to taking hearts and that $\C_{[0,k]} = H_{0,k}$. 

In fact, even the category $\C$ is determined by the heart: Sosnilo shows in \cite[Section 3]{Sosniloweight} that taking hearts sets up a fully faithful functor $(-)^\heart \colon \Catw \rightarrow \Cata$ from the category of stable $\infty$-categories equipped with an exhaustive weight structure and functors of degree $0$ to the category of additive ($\infty$-)categories and product preserving functors. Its essential image is $\Catwic$, the full subcategory spanned by the weakly idempotent complete additive categories, \emph{i.e.} those for which every map that admits a retraction is the inclusion of a direct summand. We shall also refer to these as $\flat$-additive for brevity ($\flat$ being an accidental lower than $\natural$, the standard symbol for idempotent completion). A typical example of an additive category that is not weakly idempotent complete is given by the finitely generated free modules over certain rings, e.g.\ the sections in the tangent bundles of $\mathrm S^2$ as a module over $\mathrm{C}(\mathrm S^2,\mathbb R)$. Its weak idempotent closure is the category of stably free modules. 

The functor $(-)^\flat$ admits a right adjoint $\Stab \colon \Cata \rightarrow \Catw$, which takes any additive category $\mathcal A$ to its stable envelope $\Stab(\mathcal A)$, \emph{i.e.} the smallest stable subcategory of $\Fun^\times(\mathcal A\op,\Spa)$ containing the (connective) deloopings of the representable functors $\mathcal A\op \rightarrow \mathrm{Grp}_{\Einf}(\Sps)$, where the superscript denotes product preserving functors. It is easily checked invariant under weak idempotent completion and necessarily gives the inverse to $(-)^\heart$ on weakly idempotent complete additive categories. The invariance of stablisation under this operation will allow us to also treat general additive categories in many cases. For example, due to the group completion theorem forming weak idempotent completions has no effect on direct sum $\K$-theory, \emph{i.e.} the association $\mathcal A \mapsto \core(\mathcal A)^\grp$, see for example \cite[Lemma 7.5]{HLS}, and an analogous statement holds in the hermitian context.

\begin{example}\label{examplesofweight}
\begin{enumerate}
\item\label{CCvsStab} Let $\mathcal A$ be an ordinary additive category, and denote by $\CCbdd_\infty(\mathcal A)$ the $\infty$-category of bounded chain complexes in $\mathcal A$, defined for example as the differential graded nerve of the ordinary category $\CCbdd(\mathcal A)$ of bounded chain complexes of $\mathcal A$ with its canonical enrichment in chain complexes over the integers. Then $\CCbdd_\infty(\mathcal A)$ is stable as a consequence of \cite[Proposition 1.3.2.10]{HA} and the full subcategories spanned by complexes concentrated in non-negative and non-positive degrees define a weight structure on $\CCbdd_\infty(\mathcal A)$ whose heart identifies with the weak idempotent completion $\mathcal A^\flat$ (given by all those objects $X$ in the idempotent completion $\mathcal A^\natural$ for which there exists an object $Y \in \mathcal A$ with $X \oplus Y \in \mathcal A$): 

The only non-trivial step in verifying the existence part is that both subcategories are closed under retracts (the fibre sequences in iii) are obtained by naively chopping a chain complex into its non-negative and negative parts). Let us explain the closure under retracts in case of connective objects, the other one being dual. If for example $Y \in \CCbdd_\infty(\mathcal A)$ is a retract of $X$ and $X_j = 0$ for $j \leq i$ but  $Y_i\neq 0$, then for the least such $i$, the $i$-th part $s \colon Y_i \rightarrow Y_{i+1}$ of a witnessing chain homotopy splits the differential $d_{i+1} \colon Y_{i+1} \rightarrow Y_{i}$. 

But it is then not difficult to check that $Y$ is chain homotopy equivalent to a complex terminating in degree $i+1$, namely to
\[Z = \big[\dots \longrightarrow Y_{i+3} \xrightarrow{(d_{i+3},0)} Y_{i+2} \oplus Y_i \xrightarrow{(d_{i+2},s)} Y_{i+1} \longrightarrow 0 \longrightarrow 0 \longrightarrow \dots\big]\]
via the maps $f \colon Y \rightarrow Z$ and $g \colon Z \rightarrow Y$  which are the identity above degree $i+2$ and below $i$ and in that range given by
\[f_{i} = 0,\ f_{i+1} = \id,\ f_{i+2} = \mathrm{incl} \quad \text{and} \quad g_i = 0,\ g_{i+1} = \id-sd_{i+1},\ g_{i+2} = \pr.\]
The witnessing homotopies vanish except for $(0,d_{i+1}) \colon Y_{i+1} \rightarrow Y_{i+1} \oplus Y_i$ and $s \colon Y_i \rightarrow Y_{i+1}$ as the $(i+1)$st and $i$th components, respectively. (If $\mathcal A$ happens to be weakly idempotent complete, then $d_{i+1}$ admits a kernel, on account of the split $s$, and one may take instead of $Z$ the slightly simpler complex
\[\dots \longrightarrow Y_{i+3} \xrightarrow{d_{i+3}} Y_{i+2} \xrightarrow{d_{i+2}} \mathrm{ker}(d_{i+1}) \longrightarrow 0 \longrightarrow 0 \longrightarrow \dots\]
which differs from each of $Z$ and $Y$ by a contractible summand with two consecutive entries $Y_i$; the above are then simply the composite equivalences.)

Applied to a chain complex $Y$ that is already concentrated in non-positive degrees, the argument also concentrates a chain complex which lies in the heart of the weight structure into degrees $1$ and $0$ with split differential. But for two such complexes $X$ and $Y$ one finds 
\[\Hom_{\CCbdd_\infty(\mathcal A)}(X,Y) \simeq \Hom_{\mathcal A^\natural}(\mathrm{H}_0(X),\mathrm{H}_0(Y)),\]
where $\mathcal A^\natural$ is the idempotent completion of $\mathcal A$. It follows that the heart of the weight structure is indeed precisely the weak idempotent completion of $\mathcal A$.

Sosnilo's result mentioned above implies that $\CCbdd_\infty(\mathcal A) \simeq \Stab(\mathcal A^\flat) \simeq \Stab(\mathcal A)$ for every ordinary additive category. Stabilisation can hence be thought of as a version of taking finite complexes in non-ordinary additive categories. Let us also remind the reader that $\CCbdd_\infty(\mathcal A)$ is generally different from the derived category of $\mathcal A$ (for $\mathcal A$ abelian with sufficiently many projective objects, say): For example, in the stabilisation of the category of finite generated abelian groups the null-composite sequence
\[\mathbb Z \xrightarrow{n} \mathbb Z \longrightarrow \mathbb Z/n\]
is not a fibre sequence (since all fibre sequences in the heart of a weight structure necessarily split).

\item\label{3.1.3 (2)} Given a weight structure on $\C$ one can put \[(\C\op)_{[0,\infty]} = (\C_{[-\infty,0]})\op \quad \text{and} \quad (\C\op)_{[-\infty,0]} = (\C_{[0,\infty]})\op\]
to obtain a weight structure on $\C\op$, and unless explicitly stated we shall always equip opposite categories with this weight structure. Similarly, one can define a new weighted category $\C^{\langle p \rangle}$ with underlying category $\C$ by
\[(\C^{\langle p \rangle})_{[0,\infty]} = \C_{[p,\infty]} \quad \text{and} \quad (\C^{\langle p \rangle})_{[-\infty,0]} = \C_{[-\infty,p]}.\]
Clearly, exhaustiveness is preserved by either construction. Furthermore, the product and direct sum (i.e.\ the coproduct in $\Catx$) of any family of  stable $\infty$-categories inherit a weight structure from their constituents, where both connectivity and coconnectivity are measured pointwise. In the case of the direct sum also exhaustiveness is preserved, but this fails for infinite products.

\item Given a connective $\Eone$-ring $A$, there is a canonical weight  structure on the category $\Modp_A$ of perfect (i.e., compact) $A$-modules. More generally, for a subgroup $c \subseteq \K_0(A)$, with $[A] \in c$, the same is true for the full subcategory $\Mod^c_A\subset \Modp_A$ of $A$-modules whose $\K_0$-class lies in $c$. Namely, set $(\Mod^c_A)_{[0,\infty]}$ to consist of all  $A$-modules $Y \in\Mod^c_A$ whose homotopy groups vanish in negative degrees and $(\Mod^c_A)_{[-\infty,0]}$ to consist of all $A$-modules $X \in \Mod^c_A$ with $\hom_A(X,A)$ having vanishing negative homotopy groups. Note that connectivity for the weight structure is equivalent to connectivity in the usual $t$-structure of $\Mod_A$, but coconnectivity in the weight structure is quite different from being truncated in the $t$-structure, i.e.\ having vanishing positive homotopy groups. This latter condition is classically also referred to as ``coconnectivity'', but we emphasize that this is not our meaning of the term and that $t$-structures do not appear in this paper.

Condition i) is clear, ii) holds because $\hom_A(X,Y) \simeq \hom_A(X,A) \otimes_A Y$ and iii) is most easily proven by induction on $d$ for the statement that there is a $Y \in (\Mod^c_A)_{[-\infty,d]}$ together with a $d$-connective map $f \colon Y \rightarrow X$: Since perfect modules are bounded, we can take $Y=0$ for very negative $d$ and to go from $d$ to $d+1$ pick a map $g \colon \bigoplus A^{[d]} \rightarrow \fib(f)$ that is surjective on $\pi_d$; the sum can be chosen finite since the lowest homology group of any compact module is finitely generated over $\pi_0(A)$ (since the modules for which this is the case form a stable subcategory of $\Mod_A$ that contains $A$ and is closed under retracts). Thus $g$ is $d$-connective and we let $P$ be the cofibre of the composite $\bigoplus A^{[d]} \rightarrow \fib(f) \rightarrow Y$. Then per construction there is on the one hand a fibre sequence
\[Y \longrightarrow P \longrightarrow \bigoplus A^{[d+1]}\]
so $P \in (\Mod^c_A)_{[-\infty,d+1]}$ using the assumption on $c$, and on the other a factorisation $Y \rightarrow P \xrightarrow{h} X$ of $f$, whose second constituent resides in a fibre sequence
\[\cof(g) \longrightarrow P \xrightarrow{h} X.\]
This witnesses that the map $h \colon P \rightarrow X$ is as desired.

Since $\hom_A(X,A)$ is again perfect (albeit over $A\op$), and thus bounded, this weight structure is exhaustive and $(\Mod^c_A)^\heart$ is given by the finitely generated projective $A$-modules $X$, i.e. the retracts of finite direct sums of $A$, with $[X] \in c$: If $X$ is connective, we can pick a map $g \colon A^n \rightarrow X$, that is surjective on $\pi_0$. Then 
\[g_* \colon \hom_A(X,A^n) \longrightarrow \hom_A(X,X)\]
is isomorphic to 
\[(\mathrm{id} \otimes g)_* \colon \hom_A(X,A) \otimes_A A^n \longrightarrow \hom_A(X,A) \otimes_A X\]
and therefore surjective on $\pi_0$ if also $\Hom_A(X,A)$ is connective. A preimage of the identity of $X$ displays $X$ as a retract of $A^n$.

The assumption that $[A] \in c$ is not optimal but one cannot remove it entirely: For example on $\Mod^0_{\mathrm H\mathbb Z} \simeq \D^{\mathrm{min}}(\mathbb Z)$ a weight structure as above would have vanishing heart, which is absurd. It is easily verified directly from the construction above or indeed implied by \cite[Theorem 4.3.2]{Bondarko2010weights} that a necessary and sufficient condition is that the stable hull of $\Mod^c_A \cap (\Modp_A)^\heart$ in $\Modp_A$ is $\Mod^c_A$. For example this is also the case if $n \cdot [A] \in c$ for some $n \geq 2$.

\item Specialised to an Eilenberg-Mac Lane spectrum $\mathrm H R$ for some ordinary ring $R$ the weight structure constructed in the previous point corresponds by inspection to the one on $\Dperf(R) = \CCbdd_\infty(\mathcal P(R))$ from the first example. Complexes in $\Dperf(R)_{[a,b]}$ are often said to have projective (or $\mathrm{Ext}$-)amplitude $[-b,-a]$, but we shall refrain from engaging with this slightly arcane terminology.

\item The same construction as in (\ref{3.1.3 (2)}) can be used to give a weight structure to categories of parametrised spectra, i.e.\ $\Spa/B = \Fun(B,\Spa)$, where the connective objects are measured pointwise; their existence can also directly be deduced from (\ref{3.1.3 (2)}) through the equivalences
\[\Spa/B \simeq \prod_{X \in \pi_0(B)} \Spa/X \simeq  \prod_{b \in \pi_0(B)} \Mod_{\mathbb S[\Omega_b B]}\]
determined by a choice of basepoint $b$ in each component $X$ of $B$: They restrict to equivalences
\[(\Spa/B)^\omega \simeq \bigoplus_{X \in \pi_0(B)} (\Spa/X)^\omega \simeq  \bigoplus_{b \in \pi_0(B)} \Modp_{\mathbb S[\Omega_b B]},\]
since a compact object of an infinite product of stable $\infty$-categories vanishes in all but finitely many components (test it against the filtered colimit formed by the parts living in finitely many components), whence the claim follows from (\ref{3.1.3 (2)}).

Unwinding definitions, the coconnective objects are given by those $X \in (\Spa/B)^\omega$ whose Costenoble-Waner dual $\hom_B(X,\Delta_!\mathbb S_B) \in \Spa/B$ is connective, and the heart is then spanned by spectra of the form $i_!(\mathbb S)$, for $i \colon * \rightarrow B$ the inclusion of a point, under direct sums and retracts.
\end{enumerate}
\end{example}

For easier reference later we record the following simple statements:

\begin{lemma}\label{connlift}
If $f \colon X \rightarrow Y$ is $n$-connective then for any $Z \in \C_{[-\infty,n]}$ the map
\[\pi_0\Hom_\C(Z,X) \xrightarrow{f_*} \pi_0\Hom_\C(Z,Y)\]
is surjective, and even bijective if $Z \in \C_{[-\infty,n-1]}$.
\end{lemma}

\begin{proof}
There is a fibre sequence of spectra
\[\hom_\C(Z,\fib(f)) \longrightarrow \hom_\C(Z,X) \longrightarrow \hom_\C(Z,Y)\]
and by assumption $\fib(f) \in \C_{[n,\infty]}$. The claim follows from axiom (2). 
\end{proof}

\begin{lemma}\label{connfib}
Given $a \leq c < b$ and $X \in \C_{[a,b]}$ then any fibre sequence $Y \rightarrow X \rightarrow Z$ with $Y \in \C_{[-\infty,c]}$ and $Z \in \C_{[c+1,\infty]}$ automatically has $Y \in \C_{[a,c]}$ and $Z \in \C_{[c+1,b]}$.
\end{lemma}

\begin{proof}
The claim follows by inspecting the rotated fibre sequences
\[Z^{[-1]} \longrightarrow Y \longrightarrow X \quad \text{and} \quad X \longrightarrow Z \longrightarrow Y^{[1]}\]
using the fact that the categories $\C_{[a,\infty]}$ and $\C_{[-\infty,b]}$ are closed under extensions.
\end{proof}

\begin{lemma}\label{heartsplit}
Let $\C$ be a weighted category. Then any fibre sequence $X \rightarrow Y \rightarrow Z$ with $Z \in \C_{[-\infty,0]}$ and $X \in \C_{[0,\infty]}$ splits.

In particular, this is the case for for fibre sequence of $\C$ lying in $\C^\heart$.
\end{lemma}

\begin{proof}
In the rotated sequence
\[Y \longrightarrow Z \xrightarrow{\partial} X^{[1]}\]
we find $\partial$ nullhomotopic since $\pi_0\Hom_\C(Z,X^{(1)}) = \pi_{-1}\hom_\C(Z,X)= 0$ by assumption.
\end{proof}

\subsection{Weight structures on Poincar\'e categories}

\begin{definition}
A \emph{Poincar\'e category of dimension at least $d$} is a Poincar\'e category  $(\C,\QF)$ equipped with an exhaustive weight structure, such that the duality functor $\Dual_\QF \colon \C\op \rightarrow \C$ is bounded below by $d$, explicitly 
\[\Dual_\QF(X) \subseteq \C_{[d,\infty]} \quad \text{for all } X \in \C_{[-\infty,0]}.\]  
Similarly, it is of \emph{dimension at most $d$} if it is bounded above by $d$, i.e.\
\[\Dual_\QF(X) \subseteq \C_{[-\infty,d]} \quad \text{for all } X \in \C_{[0,\infty]}\]  
and of dimension exactly $d$ if both are true. 
In either case, we shall write $\Poinc(\C,\QF)_{[a,b]}$ for the collection of path components in $\Poinc(\C,\QF)$ spanned by those objects $(X,q)$ with $X \in \C_{[a,b]}$.
\end{definition}

It is easy to see from exhaustiveness that it in fact suffices to check $\Dual_\QF(X) \subseteq \C_{[d,\infty]}$ for all $X \in \C^\heart$ to conclude that $(\C,\QF)$ has dimension at least $d$, and that in this case automatically
\[\Dual_\QF(X) \subseteq \C_{[d-l,\infty]} \quad \text{for all } X \in \C_{[-\infty,l]}.\]  
A similar statement of course applies to upper dimension bounds. In particular, $(\C,\QF)$ has dimension exactly $d$ if and only if 
\[\Dual_\QF(X) \subseteq \C_{[d,d]} \quad \text{for all } X \in \C^\heart.\]

\begin{remark}
Let us warn the reader, that the definition is slightly abusive in the sense, that a Poincar\'e category of dimension at least $d$ may not have an exact dimension at all, consider e.g.\ $\C\times \D$ where $\C$ and $\D$ are of different exact dimensions, or more explicitly $(\Dperf(K \times L), \QF^\sym_{M})$ with $M = K[d] \oplus L[d+1]$ with its standard weight structure.
\end{remark}

\begin{lemma}\label{weight:dualvsbil}
An exhaustive weight structure on a stable $\infty$-category $\C$ yields a Poincar\'e category $(\C,\QF)$ of dimension at least $d$ if and only if for all $X,Y \in \C^\heart$ the spectrum $\Bil_\QF(X,Y)$ is $d$-connective. In particular, this is the case if $\QF$ takes $d$-connective values on the heart of $\C$.
\end{lemma}

\begin{proof}
For $X \in \C$ we have $\Dual_\QF X \in \C_{[d,\infty]}$ if and only if $\hom_\C(Y,\Dual_\QF X)$ is connective for all $Y \in \C_{[-\infty,d]}$, or in other words if and only if $\Bil_\QF(Y,X)$ is $d$-connective for all $Y \in \C_{[-\infty,0]}$. But if this is true for $Y \in \C^\heart$ it follows for all $Y \in \C_{[-n,0]}$ by induction using the fibre sequences from the definition of weight structures together with \ref{connfib} above.
Exhaustiveness thus implies that $\Dual_\QF X \in \C_{[d,\infty]}$ if and only if $\Bil_\QF(X,Y)$ is $d$-connective for all $Y \in \C^\heart$. 

For the final statement simply note that $\Bil_\QF(X,Y)$ is (per construction) a retract summand of $\QF(X \oplus Y)$.
\end{proof}

\begin{example}\label{examplesofpoincweight}
\begin{enumerate}
\item\label{item:Poincop} Given a Poincar\'e category $(\C,\QF)$, the duality provides an equivalence $\C\op \simeq \C$, so in particular the composition
\[\C \xrightarrow{\Dual_\QF\op} \C\op \xrightarrow{\QF} \Spa\]
makes $\C\op$ into a Poincar\'e category with duality equivalence $\Dual_\QF\op \colon \C\to \C\op$. If $(\C,\QF)$ is equipped with a weight structure of dimension at least $d$, then equipping $\C\op$ with the opposite weight structure provided by \ref{examplesofweight} makes it have dimension at most $-d$, and similarly with upper and lower bounds flipped. 

We can therefore restrict attention to one of the two cases, and to conform better with the intuition from geometric topology, where one is forced to consider surgery `below the middle dimension', we choose that of lower dimension bounds.
\item If $(\C,\QF)$ is a Poincar\'e category of dimension at least/most $d$, then $(\C,\QF\qshift{p})$ has dimension at least/most $d+p$ and $(\C^{\langle p \rangle},\QF)$ has dimension at least/most $d-2p$.
\item For an invertible form parameter $\fpm$ over some discrete ring $R$ with underlying module with involution $M$, the associated Poincar\'e structure $\QF^{\g\fpm}_M$ on $\Dperf(R)$ has dimension $0$, for the weight structure considered in Example \ref{examplesofweight} above: The duality is given by \[\Dual_{\QF^{\g\fpm}_M}(X) \simeq \mathbb R\Hom_R(X,M)\]
which clearly preserves $\PProj(R) = \Dperf(R)^\heart$; recall our convention that $\lambda$ being invertible in particular means that $M$ is projective over $R$. 
\item Let us discuss the projectivity assumption from the previous point in more detail: Recall that for a general form parameter $\lambda$ the hermitian structure $\QF^{\g\fpm}_M$ is Poincar\'e if and only if $M$ is perfect as an $R$-module and $\mathbb R\Hom_R(M,M) \simeq R[0]$ via the action through the right factor, say. In the most common case, namely when $R$ is commutative and the $R \otimes_\mathbb Z R$-module structure on $M$ arises from an $R$-module structure by pullback along the multiplication $R \otimes_\mathbb Z R \rightarrow R$, this condition is in fact equivalent to $M$ being an $\otimes_R$-invertible $R$-module, which in particular means that $M$ is projective, see e.g. \cite[Tag 0FNP]{StacksProj}. 

In general, however, tilting theory provides examples of form parameters $\lambda$ for which $\QF^{\g\fpm}_M$ is Poincar\'e without $M$ being projective (and in this case $\QF^{\g\fpm}_M$ is of dimension bounded between $0$ and the projective amplitude of $M$, but need not have an exact dimension). For example let $K$ be a field and $A$ a finite dimensional $K$-algebra with skew involution $w$ that is Gorenstein (i.e. the injective dimension of $A$ as a left (or because of the involution equivalently as a right) module is finite). Then the functor $\Hom_K(-,K) \colon \D(A)\op \rightarrow \D(A\op)$ restricts to perfect objects and gives an equivalence 
\[w^* \circ \Hom_K(-,K) \colon \Dperf(A)\op \longrightarrow \Dperf(A)\]
which is represented by $\mathrm D_{A/K} = \Hom_K(A,K)$ regarded as a $A \otimes_K A$-module via the its natural $A\op \otimes_K A$-module structure and the involution, see e.g. \cite[Proposition 9.2.17]{KrauseH}. Together with the map $w^* \colon \mathrm D_{A/K} \rightarrow \mathrm D_{A/K}$ this turns $\mathrm D_{A/K}$ into a $K$-module with involution over $A$, which gives rise to a Poincar\'e structure $\QF^{\g\sym}_{\mathrm D_{A/K}}$. The corresponding form parameter
\[(\mathrm D_{A/K})_{\Ct} \xrightarrow{\id + w^*} (\mathrm D_{A/K})^{\Ct} \xrightarrow{\id} (\mathrm D_{A/K})^\Ct\]
is invertible in our convention if and only if $\mathrm D_{A/K}$ is projective over $A$ (in either of its module structures). Since for checking both injectivity and projectivity of a finitely generated module over $A$, it suffices to consider lifting problems among finitely generated modules (for injectivity use Baer's criterion and the fact that $A$ is noetherian), on which $w^*\Hom_K(-,K)$ restricts to an equivalence, it follows that $\mathrm D_{A/K}$ is projective over $A$ if and only if $A$ is injective over itself. This, however, need not generally be the case: For example, consider $A=\mathrm{T}_2(K)$, the ring of upper triangular $(2 \times 2)$-matrices over $K$. This ring carries a skew-involution $w \colon \mathrm{T}_2(K) \rightarrow \mathrm{T}_2(K)\op$ by flipping the two diagonal entries. Let us set $\mathrm D = \mathrm D_{\mathrm T_2(K)/K}$ for brevity. As a module over $\mathrm T_2(K)$, one easily finds $D \cong B \oplus S$, where $B$ is the $\mathrm T_2(K)$-module defined by the ring homomorphism $\mathrm T_2(K) \rightarrow K$ extracting the lower diagonal entry (in particular the underlying vector space of $B$ is $K$), and $S$ is the standard representation of $\mathrm T_2(K)$ on $K^2$. But $B$ cannot be projective, since there is a non-split short exact sequence $T \rightarrow S \rightarrow B$, where $T$ is defined analogously to $B$ using the upper diagonal entry. Instead, this sequence is a projective resolution of $B$, since as $\mathrm T_2(K)$-modules one has $\mathrm T_2(K) \cong T \oplus S$, making both $T$ and $S$ projective. Note also that $S$ carries a symmetric form $q$ with respect to $\mathrm D$ induced by the map
\[w^*S \otimes_{\mathrm T_2(K)} S \longrightarrow K, \quad (x,y) \otimes (z,w) \longmapsto xw + yz.\]
with Lagrangian $T$. It is easily checked to be unimodular, so induces equivalences $S \cong w^* \Hom_K(S,K)$ and also $B \cong w^*\Hom_K(T,K)$. This implies that $S$ and $B$ are injective. In total, the above sequence thus also functions as an injective resolution of $T$ showing that $\mathrm T_2(K)$ is indeed Gorenstein of dimension $1$.

In the language of representation theory, $\mathrm D$ is the simplest example of a tilting complex which is not a tilting module and we thank Henning Krause for pointing it out to us.

We finally mention that there is an equivalence
\[\Met(\Dperf(K),\QF_K^{\g\sym}) \simeq (\Dperf(\mathrm T_2(K)),\QF_\mathrm D^{\g\sym})\]
with underlying functor
\[\hom_{\Ar(\Dperf(K))}(K \xrightarrow{\iota_1} K^2, - ) \colon \Ar(\Dperf(K)) \rightarrow \Dperf(\mathrm T_2(K)),\]
and Poincar\'e structure induced by the unimodular form $q$. In particular, the Grothendieck-Witt spectrum of $(\mathrm T_2(K),\QF_\mathrm D^{\g\sym})$ is simply $\K(K)$.

\item More generally, if $(M,N,\alpha)$ is an invertible module with genuine involution over some $\Eone$-ring $A$, then $(\Modp_A,\QF^\alpha_M)$ has dimension at least or most $d$ if and only if $M \in (\Modp_A)_{[d,\infty]}$ or $M \in (\Modp_A)_{[-\infty,d]}$, respectively. In particular, also $(\Dperf(R),\QF^\q_M)$ and $(\Dperf(R),\QF^\s_M)$ have dimension $0$ for discrete, projective $M$.

\item A completely analogous statement holds for parametrised spectra. For example, for $\lambda \in \{\q,\vis,\s\}$ the Poincar\'e category $((\Spa/B)^\omega,\QF^\lambda_\xi)$ has dimension  $d$ if and only if the stable sperical fibration $\xi \colon B \rightarrow \mathrm{Pic}(\mathbb S)$ takes values in spheres of dimension $d$:
We need to check that any object of $(\Spa/B)^\omega$ is taken into $(\Spa/B)^\omega_{[d,d]}$ by the $\xi$-twisted Constenoble-Waner duality
\[\Dual_{\QF_\xi^{\g\lambda}} \simeq \hom_B(-,\Delta_!\xi).\] 
By the description of the heart in \ref{examplesofweight} (4), it suffices to check this for the various $i_!(\mathbb S)$, where $i \colon * \rightarrow B$ is the inclusion of a point. We compute
\begin{align*}
\hom_B(i_!\mathbb S,\Delta_!\xi) & \simeq (\pr_2)_* \mathrm{F}(\pr_1^*i_!\mathbb S,\Delta_!\xi) \\
&\simeq (\pr_2)_* \mathrm{F}((i \times \id_B)_!\mathbb S_B,\Delta_!\xi) \\
&\simeq (\pr_2)_*(i \times \id_B)_* \mathrm{F}(\mathbb S_B,(i \times \id_B)^*\Delta_!\xi) \\
&\simeq (i \times \id_B)^*\Delta_!\xi \\
&\simeq i_!i^*\xi,
\end{align*}
where the first line is the definition (with $\mathrm{F}(X,Y) \in \Spa/A$ for $X,Y \in \Spa/A$ denoting the internal morphism object), the second and last line are obtained from the Beck-Chevalley formula for the cartesian squares
\[\xymatrix{
B \ar[r] \ar[d]^{i \times \id_B} & \ast \ar[d]^i && \ast \ar[r]^-i \ar[d]^i & B \ar[d]^\Delta 
\\
B \times B \ar[r]^-{\pr_1} & B && B \ar[r]^-{i \times \id_B} & B \times B,
}\]
the third line arises from one of the projection formulae, and the fourth is hopefully obvious. Since $i^*\xi\simeq \mathbb S^d$, we conclude that $i_!i^*\xi \in (\Spa/B)^\omega_{[d,d]}$.

\item For a stable $\infty$-category $\C$ equipped with an exhaustive weight structure, the hyperbolic Poincar\'e category $\Hyp(\C)$ is most naturally seen as having dimension $0$ using the product weight structure on $\C \times \C\op$, since the duality $(X,Y) \mapsto (Y,X)$ preserves $\Hyp(\C)^\heart = \C^\heart \times (\C^\heart)\op$. However, we can also equip the underlying category $\C \times \C\op$ with the weight structure $\C^{\langle p \rangle} \times \C\op$, which has dimension $-p$. As mentioned in the introduction, we will exploit this phenomenon to give a new proof of Quillen's `+-equals-$\Q$'-theorem.
\end{enumerate}
\end{example}

To describe the interaction of passing to hearts with Poincar\'e structures, we have to introduce a bit terminology: We shall call a functor $F \colon \mathcal A \rightarrow \mathcal B$ \emph{additively quadratic} if its cross-effect 
\[\Bil_F \colon \mathcal A \times \mathcal A \longrightarrow \mathcal B^\natural, \quad (a,b) \longmapsto \mathrm{ker}(F(a \oplus b) \rightarrow F(a) \oplus F(b))\]
preserves direct sums in each variable; here $\mathcal B^\natural$ denotes the idempotent completion of $\mathcal B$. Let us immediately say, that for two stable $\infty$-categories this condition is weaker than being quadratic in the sense required in the definition of a Poincar\'e structure (for example every product preserving functor between stable $\infty$-categories is additively quadratic, but quadratic if and only if it is exact). The distinction is discussed in detail in \cite[Section 4.2]{9authI} and \cite[Section 2]{polynomial}. One particular upshot is that additively quadratic functors $\mathcal A\op \rightarrow \Spa$ uniquely extend to quadratic functors $\Stab(\mathcal A)\op \rightarrow \Spa$, see \cite[Propostion 4.2.15]{9authI} (for module categories, but the proof carries over verbatim) or \cite[Theorem 2.19]{polynomial}. Combined with Sosnilo's description of exhaustive weight structures in terms of their heart (as discussed in the previous section), we therefore deduce a characterisation of Poincar\'e categories of dimension $0$: Let $\Catpw$ denote the category spanned by these with degree $0$ Poincar\'e functors as morphisms, and let $\Cataq$ denote the category of additive $\infty$-categories equipped with additively quadratic functors, that is the cartesian unstraightening of the functor 
\[(\Cata)\op \longrightarrow \Cat, \quad \mathcal A \longmapsto \Fun^{\mathrm{aq}}(\mathcal A\op, \Spa)\]
where the superscript denotes the full subcategory of additively quadratic functors. We obtain:

\begin{proposition}\label{prop:poinconstab}
The functor $(-)^\heart \colon \Catpw \rightarrow \Cataq$ is an equivalence onto the following subcategory of the target: Its objects are those pairs $(\mathcal A, \QF \colon \mathcal A\op \rightarrow \Spa)$ such that
\begin{enumerate}
\item $\mathcal A$ is weakly idempotent complete,
\item the cross-effect $\Bil_\QF$ of $\QF$ takes values in connective spectra, and
\item $\Omega^\infty \Bil_\QF(a,b) \simeq \Hom_\mathcal A(a,\Dual b)$ naturally in $a,b \in \mathcal A$ for an equivalence $\Dual \colon \mathcal A\op \rightarrow \mathcal A$.
\end{enumerate}
The functor $\Dual$ is then uniquely determined by $\QF$, and we denote it by $\Dual_\QF$, and the morphisms in the image of $(-)^\heart$ are given by those pairs $(F \colon \mathcal A \rightarrow \mathcal B, \eta \colon \QF \Rightarrow \Phi \circ F\op)$ such that the natural map
\[\eta_\sharp \colon F(\Dual_\QF(a)) \longrightarrow \Dual_\Phi(F(a))\]
is an equivalence for every $a \in \mathcal A$.

The inverse process takes $(\mathcal A,\QF)$ to $\Stab(\mathcal A)$ equipped with the unique weight structure with $\Stab(\mathcal A)^\heart \simeq \mathcal A$ and the unique quadratic extension of $\QF$.
\end{proposition}

We shall denote the subcategory of $\Cataq$ described in the Proposition by $\Catpad$.

\begin{example}
\begin{enumerate}
\item The category $\Catpad$ contains, in particular, Schlichting's additive form categories with strong duality from \cite[Definition 2.1]{Schlichting2019}: They consist essentially of those pairs $(\mathcal A,\QF)$, where $\mathcal A$ is an ordinary additive category (which we take to be weakly idempotent complete), and $\QF$ factors over the Eilenberg-Mac Lane functor $\mathrm{H} \colon \Ab \rightarrow \Spa$, see \cite[Remark 2.4]{Schlichting2019}. By stabilisation we therefore obtain Poincar\'e categories of dimension $0$ from such examples (even non weakly idempotent complete ones). We will explore this relation further in Section \ref{subsec:applQ} below.
\item Many examples of additively quadratic functors $\mathcal A\op \rightarrow \Spa$, that do not necessarily factor through the inclusion $\mathrm{Ab} \subset \Spa$, can be obtained from Ranicki's additive categories with chain duality $(\mathcal A,T,e)$, see \cite[Definition 1.1]{Ranickiblue}. By definition $\mathcal A$ is an ordinary additive category, $T$ is a product preserving functor from $\mathcal A\op$ to the ordinary category $\CCbdd(\mathcal A)$ of bounded chain complexes in $\mathcal A$, and $e \colon T^2 \Rightarrow \mathrm{incl}$ a natural transformation such that
\begin{enumerate}
\item[a)] the composite 
\[TX \xrightarrow{T(e_X)} T^3X \xrightarrow{e_{TX}} TX\]
is the identity, and
\item[b)] $e_X \colon T^2X \rightarrow X[0]$ is a chain homotopy equivalence for all $X \in \mathcal A$;
\end{enumerate}
here $T$ is extended to a functor $\CCbdd(\mathcal A)\op \rightarrow \CCbdd(\mathcal A)$ by following the functor $\CCbdd(\mathcal A) \rightarrow \CCbdd(\CCbdd(\mathcal A))$ it induces with the formation of the total complex from a double complex. 

Given this data Ranicki defines 
\[X \otimes_T Y = \mathcal{H}\mathrm{om}_{\CCbdd(\mathcal A)}(TX,Y) \in \CCbdd(\mathbb Z)\]
for $X,Y \in \mathcal A$ and observes that this canonically gives rise to a functor
\[\Bil^T \colon \mathcal A \rightarrow \CCbdd(\mathbb Z)^\hC, \quad X \longmapsto X \otimes_T X,\]
where the target is nothing but the category of bounded $\Ct$-chain complexes. The composite
\[\mathcal A \xrightarrow{\Bil^T} \CCbdd(\mathbb Z)^\hC \longrightarrow \CCbdd_\infty(\mathbb Z)^\hC \xrightarrow{\mathbb L} \D(\mathbb Z)^\hC \xrightarrow{\mathrm H} \Spa^\hC\]
can then be followed by either homotopy fixed points or orbits to obtain functors $\widetilde{\QF}^\sym_T$ and $\widetilde{\QF}^\q_T$, that record the spectra underlying Ranicki's functors $\mathrm W^{\%}$ and $W_{\%}$ to $\mathrm{Ch}(\mathbb Z)$.
By direct inspection both are additively quadratic with cross-effect taking $(X,Y)$ to $\mathrm{H}\mathbb L(X \otimes_T Y)$.

By the discussion above we can therefore extend $\widetilde{\QF}^\q_T$ and $\widetilde{\QF}^\sym_T$ uniquely to hermitian structures on $\CCbdd_\infty(\mathcal A)\op$ using the equivalence $\CCbdd_\infty(\mathcal A) \simeq \Stab(\mathcal A)$ from \ref{examplesofweight} \eqref{CCvsStab}. It is a simple computation that this hermitian structure is in fact Poincar\'e on account of the second condition above (which went unused so far) with duality given by the inverse of the extension of $T$ to a functor
\[T \colon \CCbdd_\infty(\mathcal A)\op \longrightarrow \CCbdd_\infty(\mathcal A);\]
we note in passing that this extension can indeed be obtained from $T \colon \CCbdd(\mathcal A)\op \rightarrow \CCbdd(\mathcal A)$, since this functor is canonically enriched over $\mathrm{Ch}(\mathbb Z)$.

For $\mathcal A = \mathcal P(R)$, the category of finitely generated projective $R$-modules, we obtain Poincar\'e structures on $\CCbdd_\infty(\mathcal P(R))\op \simeq \Dperf(R)\op$. To better match up Ranicki's with our conventions one can define Poincar\'e structures $\QF^\sym_T$ and $\QF^\q_T$ on $\CCbdd_\infty(\mathcal A)$ by
\[\QF_T^\sym = \widetilde{\QF}^\sym_T \circ T \quad \text{and} \quad \QF_T^\q = \widetilde{\QF}^\q_T \circ T,\]
compare \ref{examplesofpoincweight} \eqref{item:Poincop}.
Using the weight structure from \ref{examplesofweight} \eqref{CCvsStab} we then find that these latter Poincar\'e structures are of dimension $i$ precisely if $T(X)$ is (equivalent to) a complex concentrated in degree $-i$ for every $X \in \mathcal A$.
\end{enumerate}
\end{example}

Finally, we record two elementary statements for easier reference.  Note that for a Poincar\'e category of dimension at least $d$, we have $\D_\QF(X) \in \C_{[d-b,\infty]}$, whenever $X \in \C_{[-\infty,b]}$, but that no statement about the duals of connective objects can be inferred. The dual statement, however, applies to Poincar\'e categories of dimension bounded above, so if $(\C,\QF)$ is of exact dimension $d$ then $\D_\QF(X) \in \C_{[d-b,d-a]}$ if $X \in \C_{[a,b]}$. The individual statements imply:

\begin{lemma}\label{highlyconnectedvanishes}
If $(\C,\QF)$ is a Poincar\'e category of dimension at most $d$ and $X \in \Poinc(\C,\QF^{[1]})_{[n,\infty]}$ with $2n>d+1$, then $X = 0$, and similarly, if $(\C,\QF)$ has dimension at least $d$ and $X \in \Poinc(\C,\QF^{[1]})_{[-\infty,m]}$ and $d+1>2m$.
\end{lemma}

Here, and also below, we shall state connectivity estimates for Poincar\'e objects in $(\C,\QF^{[1]})$, since the objects of the cobordism category $\Cob(\C,\QF)$ are given by $\Poinc(\C,\QF^{[1]})$.

\begin{proof}
In the former case we find $X \simeq \Dual_\QF X^{[1]} \in \C_{[-\infty,d+1-n]}$ using the calculation above, so $X \in \C_{[n,d+1-n]}$, which consists only of the null objects if $2n>d+1$. The second case is dual. 
\end{proof}

\begin{lemma}\label{leftarrightar}
If $(\C,\QF)$ is a Poincar\'e category of dimension at least $d$, then for any Poincar\'e cobordism $C \in \Poinc(\Q_1(\C,\QF^{[1]})$, whose underlying span $X \longleftarrow W \longrightarrow Y$ has the map $W \rightarrow Y$ $(d-\mc+1)$-coconnective, also the map $W \rightarrow X$ is $\mc$-connective.

If $(\C,\QF)$ is of dimension exactly $d$ the converse holds as well.
\end{lemma}

\begin{proof}
By Lefschetz duality we find an equivalence $\fib(W \rightarrow X) \simeq \Dual_\QF(\fib(W \rightarrow Y))$, from which the claim is immediate.
\end{proof}

\section{Parametrised algebraic surgery}\label{sec:para}

In this section we explain how to do surgery on the entire cobordism category and not just a single of its objects or morphisms. The method is strongly inspired by the surgery moves Galatius and Randal-Williams used to manipulate the geometric cobordism category in \cite{GRW}; we replace the geometric surgery by  algebraic surgery, originally developed by Ranicki \cite{rsurg2, rsurg1}, reworked into the language of Poincar\'e categories by Lurie in \cite{Lurie-L} and finally upgraded into the present form in \cite[Section 2.3]{9authII}.

We start by recalling from \cite[Section 2.3]{9authII} the basic results on algebraic surgery, bringing them into the form that will be used later, and by deriving some immediate consequences, mostly well-known at least in the classical case of module categories. Our main technique, forwards and backwards surgery, is introduced in subsection \ref{subsec:forwardssurgery}; it will serve as the basis for all of our later arguments. The final subsection analyses the behaviour of doubling and reflecting cobordisms at the level of detail required later.

\subsection{Review of algebraic surgery}\label{sec:review_of_surgery}

Just as in geometric topology we will use algebraic surgery to modify objects up to cobordism. Given a Poincar\'e object $(X,q) \in \Poinc(\C,\QF)$, a \emph{surgery datum} on $(X,q)$ consists of a map $f \colon T \rightarrow X$ in $\C$, together with a null homotopy $r$ of $f^*(q) \in \Omega^\infty\QF(T)$. Associated to such a datum is a cobordism
\[(X,q) \longleftarrow \chi(f) \longrightarrow X_f,\]
the \emph{trace} of the surgery, whose right end we call the \emph{result of surgery}; the underlying span of this cobordism is given as follows: $\chi(f)$ is the fibre of the composite 
\[X \xrightarrow{q_\sharp} \Dual_\QF(X) \xrightarrow{\Dual_\QF f} \Dual_\QF(T).\] 
It receives a canonical map from $T$, since its precomposition with $f$ is canonically identified with $f^*(q)_\sharp$, which $r$ provide a null-homotopy of. The object underlying $X_f$ is then the cofibre of this map $T \rightarrow \chi(f)$. 

Let us also recall that surgery is a functorial equivalence in the following sense: 
A surgery datum on $(X,q)$ amounts to the same thing as a hermitian form in $\Met(\C,\QF)$ restricting to $(X,q)$ along $\mathrm{met} \colon \catforms(\Met(\C,\QF)) \rightarrow \catforms(\C,\QF)$ and we define the \emph{category of surgery data} in $(\C, \QF)$ by
\[\Surg(\C, \QF) =\Poinc(\C, \QF)\times_{\catforms(\C, \QF)} \catforms(\Met(\C, \QF)).\]
If we denote for $\mathcal A \in \Cat$ by $\dec(\mathcal A)$ the pull-back $\Ar(\mathcal A)\times_{\mathcal A} \core \mathcal A$ along the source-projection, the results of \cite[Section 2.3]{9authII} upgrade the extraction of the trace as above to a commutative diagram
\[\xymatrix{
\Surg(\C, \QF\qshift 1)\ar[r]^\chi \ar[d] & \dec\Cob(\C, \QF)\ar[d]^\pr \\
\Poinc(\C, \QF\qshift 1) \ar[r]^-{\mathrm{cyl}} & \core\Cob(\C,\QF)
}\]
natural in the Poincar\'e category $(\C,\QF)$, whose top horizontal map is an equivalence, the \emph{surgery equivalence}, and whose bottom horizontal map is the equivalence taking a map $b \colon X \rightarrow Y$ in $\Poinc(\C,\QF\qshift{1})$ to the span
\[X \xleftarrow{\id} X \xrightarrow{b} Y\]
with its canonical structure as a cobordism $\xymatrix{X \ar@{~>}[r] & Y}$.  The symbol $\dec$ originates from the equivalence $\nerv_n\dec(\mathcal A) \simeq \nerv_{1+n}(\mathcal A)$ which identifies the Rezk nerve of $\dec\mathcal A$ with the décalage of $\nerv(A)$.

The surgery equivalence takes a morphism in $\Surg(\C, \QF\qshift 1)$ from $f$ to $g$ with underlying square
 \[\xymatrix{
 S \ar[r]^a\ar[d]_f & T \ar[d]^g  \\
 X \ar[r]^b & Y 
 }\]
 and $b$ an equivalence, to a commutative square in $\Cob(\C, \QF)$ 
 \[\xymatrix{
  X \ar@{~>}[r]^{\mathrm{cyl}(b)} \ar@{~>}[d]^{\chi(f)} & Y \ar@{~>}[d]^{\chi(g)}\\
  X_f \ar@{~>}[r]^{\chi(a)} & Y_g.
 }\]
 The lower cobordism $\chi(a)$ arises (necessarily) by surgery along some surgery datum on $X_g$; its underlying map is of the form $\cof(a) \to X_g$, and since this will suffice for our purposes let us refrain from spelling out the hermitian refinement.
 
 The inverse of $\chi$ takes a cobordism with underlying span $X\leftarrow Z \rightarrow Y$ to the map $\fib(Z\to Y) \to X$, along which the Poincar\'e form on $X$ canonically inherits a null-homotopy.

\begin{convention}
\begin{enumerate}
 \item When writing a Poincar\'e cobordism $C$ as $X\leftarrow W \to Y$ from left to right, we view $X$ as the source and $Y$ as the target of the cobordism; in other words $X = d_1(C)$ and $Y = d_0(C)$ for the simplicial structure on $\Q(\C,\QF)$. As already apparent here, we will often suppress hermitian forms from notation.  
 \item By default we shall view the trace of a surgery as a cobordism \emph{from} the original object \emph{to} the result  of surgery.
\end{enumerate}
\end{convention}

Now let $K$ be a small $\infty$-category. From the surgery equivalence, we obtain an induced equivalence
\[\Hom_{\Cat}(K, \Surg(\C, \QF\qshift 1)) \simeq \Hom_{\Cat}(K, \dec \Cob(\C, \QF)) \simeq \Hom_{\Cat}(\pcone K, \Cob(\C, \QF))\]
where $\pcone K$, the left pre-cone, is defined by the pushout
\[\xymatrix{K \ar[r]^-{(\mathrm{id},0)} \ar[d] & K \times [1] \ar[d] \\
            |K| \ar[r] &\pcone K}\]
in $\Cat$, which makes it adjoint to the décalage construction. We are more interested in describing functors out of the actual left cone $\lcone K$, rather than the pre-cone, so we make the following definition:

\begin{definition}
The \emph{space of $K$-fold surgery data} in a Poincar\'e category $(\C, \QF)$ is the space $\Pmsd{K}(\C, \QF)$ defined by the pull-back
\[\xymatrix{
 \Pmsd{K}(\C, \QF) \ar[r] \ar[d] & \Hom_{\Cat}(K, \Surg(\C, \QF)) \ar[d]^{\fgt}\\
 \Poinc(\C, \QF) \ar[r]^(.35){\mathrm{const}} & \Hom_{\Cat}(K, \Poinc(\C, \QF)).
}\]
\end{definition}

For $K=*$, in which case we abbreviate the above to $\Psd(\C, \QF)$, one of course simply obtains the core of $\Surg(\C,\QF)$, i.e.\ the space of surgery data, and more generally, if $|K|\simeq *$ one simply has
\[\Pmsd{K}(\C,\QF) \simeq \Hom_{\Cat}(K,\Surg(\C,\QF)).\]

From the considerations above we immediately conclude:

\begin{proposition}\label{prop:multiple_surgery_equivalence}
The surgery equivalence induces an equivalence 
\[\chi \colon \Pmsd K (\C, \QF\qshift 1) \longrightarrow \Hom_{\Cat}(\lcone K, \Cob(\C, \QF)),\]
such that the extraction of the underlying Poincar\'e object on the left corresponds to evaluation at the cone point on the right.
\end{proposition}

\begin{remark}\label{rem:alternativedefsofmultsurg}
\begin{enumerate}
\item Unwinding definitions one finds a cartesian square
\[\xymatrix{\Pmsd{K}(\C,\QF) \ar[r] \ar[d] & \Hom_{\Cat}(K,\catforms(\Met(\C,\QF))) \ar[d]^{\mathrm{met}} \\
            \Poinc(\C,\QF) \ar[r]^-{\const} & \Hom_{\Cat}(K,\catforms(\C,\QF)),}\] that we will occansionally use below.
\item The space $\Pmsd K(\C,\QF)$ can also be regarded as the (cocartesian) unstraightening of the functor 
\[\Hom_{\Cat}(K,\Surg_-(\C,\QF)) \colon \Poinc(\C,\QF) \longrightarrow \Sps\]
where $\Surg_{(X,q)}(\C, \QF)$ is the fibre of the forgetful map $\Surg(\C, \QF) \to \Poinc(\C, \QF)$ over $(X,q)$. In formulae this gives
 \[\Pmsd K (\C, \QF) \simeq \int_{(X,q)} \Hom_{\Cat}(K, \Surg_{(X,q)}(\C, \QF))\]
whereas
 \[\Hom_{\Cat}(K, \Surg(\C, \QF) \simeq \Hom_{\Cat}\biggl(K, \int_{(X, q)} \Surg_{(X,q)}(\C, \QF)\biggr).\]
\end{enumerate}
\end{remark}

Besides the map 
\[\Pmsd K (\C, \QF)\to  \Poinc(\C,\QF)\] 
forgetting the multiple surgery datum, there is also a map
\[\Pmsd K(\C, \QF) \to \Hom_{\Cat}(\rcone K, \C)\]
extracting the diagram underlying the surgery datum: Per construction we have a map $\Surg(\C,\QF) \rightarrow \Ar(\C)$ forgetting hermitian structures.
It induces a map
\[\Hom_{\Cat}(K,\Surg(\C,\QF)) \longrightarrow \Hom_{\Cat}(K \times [1], \C)\]
which (by definition of $\Surg(\C,\QF)$) restricts as desired. 
Somewhat informally, we will often write this underlying diagram as $f \colon T\to X$, where $X$ is the value at the cone point, and $T$ is the restriction $K\to \C$. In generalisation of the case $K=*$ the result of surgery in this generality is the diagram $X_f \colon K \rightarrow \Cob(\C,\QF)$, or equivalently an object of $X_f \in \Poinc\Q_K(\C,\QF\qshift{1})$.

\begin{lemma}\label{lem:cofinal_surgery_data}
If $K\to L$ is cofinal, then the square
\[\xymatrix{
 \Pmsd L (\C, \QF) \ar[r] \ar[d] & \Pmsd K(\C, \QF) \ar[d]\\
 \Hom_{\Cat}(\rcone L, \C) \ar[r] & \Hom_{\Cat}(\rcone K, \C)
}\]
is cartesian.
\end{lemma}

\begin{proof}
The map $\Hom_{\Cat}(K, \catforms(\C, \QF))\to \Hom_{\Cat}(K, \C)$ identifies, via the equivalence $\catforms((\C,\QF)^K) \simeq \Fun(K,\catforms(\C,\QF))$ of \cite[6.3.15 Corollary]{9authI}, with the right fibration
\[\spcforms((\C, \QF)^K) \to \grpcr (\C^K)\]
classified by the infinite loop space of the hermitian structure $\QF^K$ on $\C^K = \Fun(K,\C)$. From the explicit formula
\[\QF^K(F) = \lim_{K\op} \QF\circ F\op\]
we conclude that the square
\[\xymatrix{
 \Hom_{\Cat}(L, \catforms(\C, \QF)) \ar[r] \ar[d] & \Hom_{\Cat}(K, \catforms(\C, \QF)) \ar[d]\\
 \Hom_{\Cat}(L, \C) \ar[r] & \Hom_{\Cat}(K, \C)
}\]
is cartesian if $K\to L$ is cofinal. 

By using the description of $\Pmsd{L}(\C,\QF)$ from \ref{rem:alternativedefsofmultsurg} (1), we therefore find $\Pmsd L (\C, \QF)$ given as the limit of
\[\xymatrix{
\Hom_{\Cat}(K, \catforms(\Met(\C, \QF))) \ar[r]\ar[d]
  & \Hom_{\Cat}(K, \catforms(\C, \QF)) \ar[d]
  & \Poinc(\C,\QF) \ar[l]\ar[d]
\\
\Hom_{\Cat}(K,\Ar(\C)) \ar[r]
  & \Hom_{\Cat}(K,\C)
  & \grpcr(\C) \ar[l]
\\
\Hom_{\Cat}(L,\Ar\C) \ar[r] \ar[u]
  & \Hom_{\Cat}(L,\C)\ar[u]
  & \grpcr(\C) \ar[l]\ar[u],
}\]
computing vertical limits first. Computing horizontal pullbacks first gives the claim.
\end{proof}

We conclude this section by recording some easy general facts about surgery for later reference; most of them are of course well-known at least in the case of module categories. We start with the general connectivity estimate for  the trace.

\begin{lemma}\label{lem:basic_surgery_estimate}
Let $(\C, \QF)$ be a Poincar\'e category of dimension at least $d$, $f \colon T\to X$ a surgery datum on $X \in \Poinc(\C, \QF^{[1]})$ and consider the associated cobordism
\[X \longleftarrow \chi(f) \longrightarrow  X_f.\]
If $T$ is $k$-connective, then $\chi(f) \rightarrow X_f$ is $k$-connective; and if $T$ is $k$-coconnective then $\chi(f) \rightarrow X$ is $d-k$-connective. Furthermore, if $T\in \C_{[k,k]}$ and $(\C, \QF)$ is of exact dimension $d$, then the fibres of these two maps are concentrated in degrees $k$ and $d-k$, respectively.
\end{lemma}

\begin{proof}
This is immediate from the fibre sequences
\[\Dual_{\QF}(T) \longrightarrow \chi(f) \longrightarrow X \quad \text{and} \quad T \longrightarrow \chi(f) \longrightarrow X_f\]
in the definition of $\chi(f)$ and $X_f$.
\end{proof}

To improve connectivity of objects and morphisms, we will use the following two estimates.

\begin{lemma}\label{lem:surgery_below_the_middle_dimension}
Let $(\C, \QF)$ be a Poincar\'e category of dimension at least $d$ and $f \colon T\to X$ a surgery datum on $X \in \Poinc(\C, \QF\qshift 1)$, such that $T$ is $l$-coconnective. If  $T \rightarrow X$ is $\oc$-connective, then  $X_f$ is $\oc+1$-connective  as long as $l+\oc<d$. The converse statement holds as long as $l+\oc\leq d$.  
\end{lemma}

In particular, if $T \in \C_{[\oc,\oc]}$ we find $X_f$ is $\oc+1$-connective if $f$ is $\oc$-connective, as long as $2\oc<d$, i.e. as long as $T$ is a surgery datum below the middle dimension (note that in this case $X \in \C_{[\oc,\infty]}$ is a definite requirement for $X/T \in \C_{[\oc+1,\infty]}$).

\begin{proof}
Checking vertical fibres one finds the square
\[\xymatrix{\chi(f) \ar[r] \ar[d] & X \ar[d] \\
X_f \ar[r] & X/T}\]
cartesian. Considering horizontal fibres one then finds a fibre sequence
\[\Dual_\QF(T) \longrightarrow X_f \longrightarrow X/T\]
and by assumption the left hand term is in $\C_{[d-l,\infty]}$, so closure of $\C_{[\oc+1,\infty]}$ under extensions gives the claim once $d-l \geq \oc+1$. The converse is similar.
\end{proof}

\begin{lemma}\label{extraconnectivitycalc}
Given a surgery datum in $\Q_1(\C, \QF\qshift 1)$ with underlying diagram
\[\xymatrix{S \ar[d]^g & T \ar[r]^{\id}\ar[l]\ar[d] &T \ar[d]^f\\
            X & \ar[r] W \ar[l] & Y}\]
such that the total fibre of the left square is $\mc$-connective, then the result of surgery has left pointing map $\mc+1$-connective.
\end{lemma}

\begin{proof}
We need to show that in the result of surgery
\[X_{g} \longleftarrow \fib(W\to \Dual_{\QF\qshift 1} S)/T \longrightarrow Y_f\]
the lower left arrow is $\mc+1$-connective. But its fibre is the same as that of $W/T \rightarrow X/S$, and thus the suspension of the totel fibre of the square in question.
\end{proof}

To take care that the connectivity of morphisms is not destroyed by surgery, we use the following estimate.

\begin{lemma}\label{connsurgmor}
Let $(\C,\QF)$ be a Poincar\'e category of dimension at least $d$ and let
\[\xymatrix{R \ar[d]^f & S \ar[l]\ar[r] \ar[d]^g & T \ar[d]^h \\
            X & W \ar[r] \ar[l]& Y}\]
be a surgery datum on an object $\Poinc\Q_1(\C,\QF^{[1]})$ with $W \rightarrow X$ is $\mc$-connective. Then if $S \rightarrow T$ is $l$-coconnective and $S \rightarrow R$ $k$-connective, also the left pointing map in the result of the surgery is $\mc$-connective as long as $\mc \leq k+1, d-l$. 
\end{lemma}

\begin{proof}
Set $M = R \cup_ST$ for legibility. Then the trace of the surgery yields a commutative diagram
\[\xymatrix{ X & \ar[r]\ar[l]W & Y \\
             \chi(f) \ar[u]\ar[d]& \ar[r]\ar[d]\ar[l]\ar[u]\fib(W \rightarrow \Dual_\QF(M)^{[1]}) & \ar[u]\ar[d]\chi(h) \\
             X_f & \ar[r]\ar[l]\fib(W \rightarrow \Dual_\QF(M)^{[1]})/S & Y_h}\]
whose bottom row is the result of surgery. Let $F_b$, $F_m$ and $F_u$ denote the fibres of the bottom, middle and upper left pointing arrows, respectively. Then computing the total fibre of the lower left square first vertically and then horizontally results in a (rotation of the) fibre sequence
\[F_m \longrightarrow F_b \longrightarrow R/S\]
and doing the same for the upper left square gives
\[\fib(\Dual_\QF(M) \rightarrow \Dual_\QF(R)) \longrightarrow F_m \longrightarrow F_u.\]
But the source in the second term is equivalent to $\Dual_\QF(M/R) \simeq \Dual_\QF(T/S)$. Thus if $T/S \in \C_{[-\infty,l]}$ this terms lies in $\C_{[d-l,\infty]}$ and the result follows.
\end{proof}

We shall also need the following \emph{automatic disjointness}, which reflects the geometric fact that any two embedded submanifolds of less than complementary dimension can be isotoped apart:

\begin{lemma}\label{autodisj}
Let $(\C,\QF)$ be a Poincar\'e category of dimension at least $d$ and suppose given pieces of surgery data $f \colon S \rightarrow X$ and $g \colon T \rightarrow X$ on $X \in \Poinc(\C,\QF^{[1]})$. Assume that $S \in \C_{[-\infty,k]}$ and $T \in \C_{[-\infty,l]}$. Then as long as $l+k \leq d$, there is a refinement of $f + g \colon S \oplus T \rightarrow X$ to a surgery datum restricting to the original two pieces under precomposition.
\end{lemma}

\begin{proof}
We have to show that nullhomotopies of the restrictions to $S$ and $T$ of the Poincar\'e form $q \in \Omega^{\infty-1}\QF(X)$ extend to a nullhomotopy of the restriction of $q$ in $\Omega^{\infty-1}\QF(S \oplus T)$. To see this simply recall that 
\[\QF(S \oplus T) = \QF(S) \oplus \QF(T) \oplus \Bil_{\QF}(S,T),\]
so it will suffice to show that $\Bil_{\QF}(S,T)$ is connective, but this is guaranteed by the assumption on $k$ and $l$ together with Lemma \ref{weight:dualvsbil}.
\end{proof}

\subsection{Forwards and backwards surgery}
\label{subsec:forwardssurgery}

Recall that the cobordism category $\Cob(\C,\QF)$ for a Poincar\'e category $(\C,\QF)$ is the $\infty$-category associated to the complete Segal space $\Poinc(\Q(\C, \QF\qshift{1}))$. Adding the space of surgery data to the mix yields another simplicial space $\Psd{(\Q(\C, \QF\qshift{1}))}$, whose associated category we will denote $\Cobsd(\C, \QF)$, the \emph{cobordism category with surgery data}. Since the functor $\Surg\colon \Catp\to \Cat$ preserves limits (each of its constituents $\Poinc$, $\catforms$, and $\Met$ does), and therefore $\Psd \colon \Catp \rightarrow \Sps$ also preserves limits, \cite[Lemmas 2.2.4]{9authII} implies that $\Psd{(\Q(\C, \QF\qshift{1}))}$ is also a complete Segal space, giving us good control over $\Cobsd(\C,\QF)$.

There is an evident forgetful functor $\fgt \colon \Cobsd(\C,\QF) \rightarrow \Cob(\C,\QF)$ induced by the natural projection $\Psd\to \Poinc$, which may be rewritten as
\[\Psd \xrightarrow{\chi} \Poinc\Q_1 \xrightarrow{d_1} \Poinc\Q_0 = \Poinc.\] 
But similarly the transformation 
\[\Psd \xrightarrow{\chi} \Poinc\Q_1 \xrightarrow{d_0} \Poinc\] 
induces another functor $\srg \colon \Cobsd(\C,\QF) \rightarrow \Cob(\C,\QF)$ that extracts the result of surgery.

The key insight allowing for parametrised surgery arguments is that these two functors should be related by a natural transformation, and hence become homotopic upon realisation, once restricted appropriately.

To explain this, we note that by construction the functors are related by a cobordism, namely by the trace of the surgery, i.e. $\chi \colon \Psd\Q(\C,\QF^{[1]}) \rightarrow \Poinc\Q_1\Q(\C,\QF^{[1]})$. In more categorical language this can be regarded as a functor into the category of vertical morphisms in the double category $\DCob(\C,\QF^{[-1]})$ represented by the double Segal space $\Poinc\Q^{(2)}(\C,\QF^{[1]})$. Such a cobordism, however, is not the same thing as a natural transformation between its source and target, and generally does not induce a homotopy at the level of classifying spaces. 

It turns out, however, that the double category $\DCob(\C,\QF^{[-1]})$ contains the double category  $\Sq(\Cob(\C,\QF))$ of commutative squares in $\Cob(\C,\QF)$: A square in the double category $\DCob(\C,\QF^{[-1]})$ is by definition a  $(1,1)$-simplex in $\Poinc\Q^{(2)}(\C,\QF^{[1]})$, i.e.\ a cobordism between cobordisms, which may be written as a diagram
\[\xymatrix{X & \ar[r]\ar[l]W & Y \\
            V \ar[u]\ar[d]& \ar[r]\ar[l]\ar[u]\ar[d]C &\ar[u]\ar[d] U \\
            X'&\ar[r]\ar[l] W'& Y'}\]
(discarding the hermitian structures from the notation) and we will argue that such a datum can be regarded as a commutative square 
\[\xymatrix{X \ar@{~>}[r]\ar@{~>}[d]& Y \ar@{~>}[d] \\
            X' \ar@{~>}[r] & Y'}\]
in $\Cob(\C,\QF)$ just in case its top right and lower left corner are cartesian. Restricting the domains of $\fgt$ and $\srg$ to the wide subcategory $\Cobfw(\C,\QF)$ of $\Cobsd(\C,\QF)$ consisting of surgery data that produce such diagrams one then obtains the desired transformation $\fw \colon \fgt \Rightarrow \srg$, since functors into the category of vertical morphisms in $\Sq(\D)$ are the same thing as natural transformations. We shall also need a version of this transformation for multiple surgery, and dual statements obtained from the symmetry of the cobordism category.\\

\begin{remark}
We have chosen to distinguish between a complete double Segal space such as $\Poinc\Q^{(2)}(\C,\QF\qshift{2})$ and its associated double category $\Cob^{(2)}(\C,\QF)$ for conceptual clarity (in particular, we shall refer to the former as the nerve of the latter). We therefore invite the reader to choose their favourite model for double categories in the following; of course, one perfectly valid such choice are complete double Segal spaces themselves, i.e.\ Segal objects in complete Segal spaces (and we require the reader's choice to be equivalent to this one), in which case the distinction above is merely syntactic. 
\end{remark}

To start the analysis recall that the double category of commutative squares $\Sq(\D)$ in a given category $\D$ is associated to the two-fold complete Segal space
\[(i,j) \mapsto \Hom_{\Cat}([i] \times [j], \D) \simeq \Hom_{\sSps}(\Delta^i \times \Delta^j, \nerv(\D)),\] 
(the two terms are equivalent because the Rezk nerve $\nerv$ is fully faithful). Our first order of business will be to verify that the sub-Segal-space of $\Poinc\Q^{(2)}(\C,\QF^{[1]})$ spanned by the squares described above indeed agrees with the nerve of $\Sq(\Cob(\C,\QF))$.

From the natural equivalence
\[\Hom_{\Cat}(K,\Cob(\C,\QF)) \simeq \Poinc\Q_K(\C,\QF\qshift{1}),\]
we obtain: 

\begin{observation}\label{squaresindouble}
The double category $\Sq(\Cob(\C,\QF))$ is associated to the two-fold complete Segal space
\[(i,j) \longmapsto \Poinc\Q_{[i] \times [j]}(\C,\QF^{[1]})\]
for every Poincar\'e category $(\C,\QF)$.
\end{observation}

This makes it very easy to construct the embedding of $\Sq(\Cob(\C,\QF))$ into $\DCob(\C,\QF)$.

\begin{construction}\label{squaresindoubleII}
Both $\Q_{[i] \times [j]}(\C)$ and $\Q^{(2)}_{i,j}(\C)=\Q_i\Q_j(\C)$ are (essentially by definition) full subcategories of 
\[\Fun(\Twar([i] \times [j]),\C) \simeq \Fun(\Twar([i]) \times \Twar([j]), \C)\]
and their hermitian structure are the restricted ones. In the former case all squares associated to maps $[3] \rightarrow [i] \times [j]$ are required to be pullbacks, whereas in the latter case one only requires this for maps $[3] \rightarrow [i] \times [j]$ that are of the form $k=k \rightarrow l=l$ in one of the two factors (since limits in functor categories are computed pointwise). Therefore, we clearly find a fully faithful functor
\[\Q_{[i] \times [j]}(\C) \longrightarrow \Q^{(2)}_{i,j}(\C)\]
natural in $(i,j) \in \bbDelta \times \bbDelta$. 
\end{construction}

One easily decodes that it is an equivalence for $i=0$ or $j=0$, but not in general as the following example shows:

\begin{observation}\label{ex:Q_1Q_1}
The category $Q^{(2)}_{1,1}(\C)$ is the full functor category $\Fun(\Twar([1]\times [1]), \C)$, whereas directly decoding definitions shows that a (solid) diagram
\[\xymatrix{
 X_{0\leq 0,0 \leq 0}\ar@{--}[rd] &\ar[r]\ar[l] X_{0\leq 1,0 \leq 0}  & X_{1 \leq 1, 0 \leq 0} 
 \\
 X_{0\leq 0,0 \leq 1} \ar[u]\ar[d]&\ar[r]\ar[l]\ar[u]\ar[d] X_{0\leq 1,0 \leq 1} \ar@{--}[rd] & X_{1 \leq 1, 0 \leq 1} \ar[u]\ar[d]
 \\
 X_{0\leq 0,1 \leq 1} &\ar[r]\ar[l] X_{0\leq 1,1 \leq 1}& X_{1 \leq 1, 1 \leq 1} }\]
defines an element in $\Q_{[1] \times [1]}(\C)$ if and only if the lower left and the upper right squares are cartesian. Alternative to a direct check, one can use the decomposition 
\[\Q_{[1]\times[1]}(\C)\simeq \Q_2(\C)\times_{\Q_1(\C)} \Q_2(\C)\]
arising from $[1]\times [1]=[2]\cup_{[1]} [2]$, indicated by the dashed line in the diagram, to reduce this to the characterisation of $\Q_2(\C)$ inside $\Fun(\Twar[2], \C)$.
We further note that if the lower left square is cartesian, then so is the upper right square in the diagram obtained by the duality of $\Q^{(2)}_{1,1}(\C,\QF)$, and vice versa: To see this, apply the  commutativity of 
\[\xymatrix{
 \Q_1(\C) \ar[rr]^{\Dual_{\QF^{\Twar[1]}}} \ar[d]^{\fib(d_0)}
 && \Q_1(\C)\op \ar[d]^{\cof(d_1)\op}
 \\
 \C \ar[rr]^{\Dual_\QF} 
 && \C\op
}\]
both to $(\C, \QF)$ and to $\Q_1(\C, \QF)$. It follows $\Q_{[1]\times [1]}(\C)$ is closed under the duality of $\Q^{(2)}_{1,1}(\C, \QF)$, and therefore the inclusion $\Q_{[1] \times [1]}(\C,\QF) \rightarrow \Q^{(2)}_{1,1}(\C,\QF)$ constructed above is Poincar\'e. Furthermore, it follows that for a Poincar\'e object of $\Q^{(2)}_{1,1}(\C, \QF)$ cartesianness of one of the two squares implies cartesianness of the other.
\end{observation}

From this example and the Segal conditions on both sides, we see that $\Q_{[i]\times [j]}(\C, \QF)\subset \Q^{(2)}_{i,j}(\C, \QF)$ is a Poincar\'e subcategory for all values of $i,j$.  
Upon applying Poincar\'e objects we thus find an inclusion of path components
\[\Poinc\Q_{[i] \times [j]}(\C,\QF\qshift{1}) \longrightarrow \Poinc\Q^{(2)}_{i,j}(\C,\QF\qshift{1}),\]
and investing \ref{squaresindouble}, we obtain the desired faithful inclusion
\[\fw \colon \Sq(\Cob(\C,\QF)) \longrightarrow \DCob(\C,\QF^{[-1]}).\]

From $\Twar(\D\op) \simeq \Twar(\D)$ one furthermore obtains an equivalence $\refl \colon \Q(\C,\QF) \rightarrow \Q(\C,\QF)^\rev$ of simplicial Poincar\'e categories, the \emph{reflection}, where the superscript on the right denotes the order reversal of the simplicial object. Postcomposing the reverse of $\fw$ in, say, the first simplicial direction, with this reflection gives the other embedding
\[\bw \colon \Sq(\Cob(\C,\QF))^{\mathrm{op_1}} \longrightarrow \DCob(\C,\QF^{[-1]}),\]
as promised, which is characterized by the upper left and lower right square in \ref{ex:Q_1Q_1} being cartesian.

Let us note immediately that both 
\[\refl_1 \colon \DCob(\C,\QF^{[-1]}) \longrightarrow \DCob(\C,\QF^{[-1]})^{\mathrm{op_1}} \quad \text{and} \quad \refl_2 \colon \DCob(\C,\QF^{[-1]}) \longrightarrow \DCob(\C,\QF^{[-1]})^{\mathrm{op_2}}\] 
exchange the essential images of the (appropriate opposites of the) embeddings $\fw$ and $\bw$; this holds essentially by definition for $\refl_1$ and for $\refl_2$ it follows from the characterisation of the image in the observation above.

Now for a category $K$ and a double category $\D$ let us denote by $\Hom_{\Cat}^1(K,\D)$ the category with objects functors from $K$ into the first direction of $\D$ and morphisms given by the second direction; at the level of nerves this corresponds to forming the simplicial space
\[[n] \longmapsto \Hom_{\sSps}(\nerv K,\nerv_{-, n}\D).\]

Per construction then, the transformation $\chi \colon \Psd\Q(\C,\QF^{[1]}) \rightarrow \Poinc\Q_1\Q(\C,\QF^{[1]})$ is the same as a functor
\[\chi \colon \Cobsd(\C,\QF) \longrightarrow \Hom_{\Cat}^1([1],\Cob^{(2)}(\C,\QF\qshift{-1})).\]
Furthermore, there is a tautological equivalence
\[\Hom_{\Cat}^1(K,\Sq(\C)) \simeq \Fun(K,\C)\]
for any category $\D$.

\begin{definition}\label{def:forwards}
We call a surgery datum on a Poincar\'e cobordism (i.e.\ a morphism in $\Cobsd(\C,\QF)$) \emph{forwards} if $\chi$ takes  it into the image of the subcategory embedding
\[\fw \colon \Hom_{\Cat}^1([1],\Sq\Cob(\C,\QF)) \longrightarrow \Hom_{\Cat}^1([1],\Cob^{(2)}(\C,\QF\qshift{-1})).\]
The wide subcategory of $\Cobsd(\C,\QF)$ spanned by these pieces of surgery data we denote $\Cobfw(\C,\QF)$, and there is an entirely analogous definition of $\Cobbw(\C,\QF)$. We call these  the \emph{cobordism categories with forwards (or backwards) surgery data}.
\end{definition}

Per construction $\chi$ restricts to a functor
\[\chi \colon \Cobfw(\C,\QF) \longrightarrow \Ar(\Cob(\C,\QF))\]
that adjoins to the desired transformation
\[\chifw \colon \fgt \Longrightarrow \srg\]
of functors $\Cobfw(\C,\QF) \rightarrow \Cob(\C,\QF)$; applying the same analysis to the backwards inclusion we obtain a transformation
\[\chibw \colon\srg \Longrightarrow \fgt\]
of functors $\Cobbw(\C,\QF) \rightarrow \Cob(\C,\QF)$.

Forwards surgery data can be explicitly characterised as follows:

\begin{proposition}
A surgery datum 
\[\xymatrix{R \ar[d]^f & S \ar[r]\ar[d]^g\ar[l] & T\ar[d]^h \\
            X &  W\ar[r]\ar[l] & Y,}\]
in $\Q_1(\C, \QF\qshift 1)$ is forwards if and only if the map $S\to R$ is an equivalence and backwards if the map $S\to T$ is an equivalence.
\end{proposition}

Recall our convention that the diagram in the proposition is to be read as a morphism from $f$ to $h$ in $\Cobsd(\C,\QF)$.

\begin{proof}
Writing the associated cobordism as
\[\xymatrix{ 
X & \ar[r]\ar[l]W & Y \\
\chi(f) \ar[u]\ar[d]& \ar[r]\ar[d]\ar[l]\ar[u]W' & \ar[u]\ar[d]\chi(h) \\
X_f & \ar[r]\ar[l]W'' & Y_h
}\]
the lower line is obtained from the middle one by taking the quotient of $R\leftarrow S \to T$, so the left lower square is cartesian if and only if $S\to R$ is an equivalence. The claim in the forwards case then follows from \ref{ex:Q_1Q_1} and the fact that the cobordism is Poincar\'e. The backwards case follows by reflection.  
\end{proof}

\begin{remark}
This proposition makes our forwards and backwards surgery moves correspond to those of Galatius and Randal-Williams in a very close manner: In their implementation one first replaces $|\Cob_d|$ by what has become known as the space of long manifolds, i.e.\ the space of $d$-dimensional submanifolds  in $\mathbb R \times \mathbb R^{(\infty)}$, such that the projection to the first coordinate is proper is and eventually displays $M$ as a cylinder of some slice in both the positive and negative directions. Roughly speaking, the equivalence from $|\Cob_d|$ takes an object $N$ to $\mathbb R \times N$ (for some embedding $N \rightarrow \mathbb R^{(\infty)}$ of which there is a contractible choice), and a cobordism $W$ represents the path given by sliding the cylindrical extension of a similarly embedded $W$ from positive to negative infinity; such a path indeed starts at $\mathbb R \times \partial_0 W$ and ends at $\mathbb R \times \partial_1 W$.

In this language Galatius and Randal-Williams perform a surgery move on some geometric surgery datum $\varphi \colon \mathrm S^k \times \mathrm D^{d-1-k} \rightarrow N$, with $N$ an object of $\Cob_d$, by extending it to an embedding $\id \times\varphi \colon\mathbb R_{\geq 0} \times \mathrm S^k \times \mathrm D^{d-1-k} \rightarrow \mathbb R \times N$, and then sliding a handle (i.e. the trace of the surgery) in from positive infinity; compare \cite[Figure 8]{GRW}. In particular, the slide  connects the result of surgery to the original manifold by a canonical path in the space of long manifolds. If more generally $N$ is a vertex of some higher simplex, whose centre in the geometric realisation of $\Cob_d$ gives rise to a long manifold $M$, say, which contains $N$ as a slice, a similar move requires the extension of $\varphi$ to an embedding $\mathbb R_{\geq 0} \times \mathrm S^k \times \mathrm D^{d-1-k} \rightarrow M$ compatible with the projections to the first coordinate. This is roughly the data Galatius and Randal-Williams start their surgery moves with. Note that this process performs a surgery parallel to $\varphi$ on all slices of $M$ to the right of $N$ (i.e.\ at higher values of the first coordinate) as well, and that this seems necessary to produce a canonical  path to the original situation.

The proposition above can be read as saying that in order for our transformation $\chifw$ to exist, our surgery data needs to satisfy a similar condition: We can perform surgery on some datum $\varphi \colon T \rightarrow C_{i \leq i}$ that is part of a higher simplex $C \in \Poinc\Q_n(\C,\QF\qshift{1})$ of $\Cob(\C,\QF)$ in a forwards manner once we have extended it to a diagram of the form
\[\xymatrix{\dots \ar[r]& 0\ar[d] & 0 \ar[d]\ar[r]\ar[l]& T \ar[d]^{\varphi}& T \ar[r]^\id\ar[l]_\id\ar[d]& T \ar[d]&\ar[l] \dots \\
            \dots \ar[r]& C_{i-1\leq i-1} & C_{i-1\leq i}\ar[r]\ar[l]& C_{i\leq i}& C_{i \leq i+1} \ar[r]\ar[l]& C_{i+1 \leq i+1} &\ar[l] \dots}\]
which clearly has the same behaviour as the kind of data used by Galatius and Randal-Williams. In fact, unwinding definitions one finds that $\chifw$ is simply induced by regarding the trace of the surgery datum $\varphi$ as a morphism in $\Cob(\C,\QF)$, which in turn seems to correspond exactly to the slide performed by Galatius and Randal-Williams. 

We thus, in particular, expect our analysis to be portable to the geometric set-up, where it would lead to variant of parametrised surgery that dispenses with the need for the space of long manifolds.
\end{remark}

Let us now upgrade the discussion to multiple surgery. Applying the category of surgery data to each entry in the $\Q$-construction yields a two-fold complete Segal space
\[ (m,n)\longmapsto N_m\Surg(\Q_n(\C, \QF\qshift 1))\]
whose associated double category we shall denote $\SurgCob(\C, \QF)$. The surgery equivalence then induces an equivalence of double categories
\[\chi \colon \SurgCob(\C, \QF) \simeq \dec_1 \DCob(\C, \QF\qshift{-1}),\]
where $\dec_1$ indicates the décalage construction in the first direction. For ease of notation we shall simply refer to the first direction of $\SurgCob(\C,\QF)$ as the surgery direction, and to the second one as the cobordism direction. 
We denote by 
\[\SurgCobfw(\C,\QF)\quad \mathrm{and}\quad \SurgCobbw(\C,\QF)\]
the subcategories of $\SurgCob(\C,\QF)$ spanned in the cobordism direction by the forwards and backwards morphisms, respectively, and all morphisms in the surgery direction. In other words, $N_{m,n}\SurgCobfw(\C, \QF)\subseteq N_{m,n}\SurgCob(\C, \QF)$ is the full subspace on those elements such that for each vertex in the first direction and each edge in the second direction, the restricted element of $N_{0,1}\SurgCob(\C, \QF) = \Psd \Q_1(\C, \QF)$ is a forwards surgery datum.

Straight from the definitions we then find the following \emph{forwards and backwards surgery equivalences}:

\begin{proposition}
\label{forwards_surgery_equivalence}
The equivalence $\chi$ restricts to an equivalence
\[\SurgCobfw(\C, \QF) \simeq \dec_1 (\Sq(\Cob(\C, \QF)))\]
along the forwards embedding, and similarly to an equivalence
\[\SurgCobbw(\C, \QF)\simeq \dec_1 (\Sq(\Cob(\C, \QF))^{\rev_1}).\]
\end{proposition}

To extract the case of $K$-fold surgery data let $\Cobmsd{K}(\C, \QF)$ denote the category associated to the complete Segal space $\Pmsd{K}(\Q(\C, \QF\qshift 1))$; in categorical language this is described by the cartesian square
\[\xymatrix{\Cobmsd{K}(\C, \QF) \ar[r] \ar[d] & \Hom_{\Cat}^1(K,\SurgCob(\C,\QF)) \ar[d]^\fgt 
\\
\Cob(\C, \QF) \ar[r]^-\const & \Fun(K,\Cob(\C,\QF)).
}\]
Using either description  we obtain an equivalence
\[\Cobmsd{K}(\C, \QF) \simeq \Hom_{\Cat}^1(\lcone{K}, \DCob(\C, \QF\qshift{-1}))\]
from \ref{prop:multiple_surgery_equivalence}. Inside $\Cobmsd K(\C, \QF)$ we find subcategories
\[\Cobmfw K(\C, \QF) \quad \mathrm{and} \quad \Cobmbw K(\C, \QF)\]
consisting of those morphisms that lie in $\SurgCobfw(\C, \QF)$ and $\SurgCobbw(\C, \QF)$, respectively. We thus find:

\begin{corollary}\label{prop:forwardssurgery}\label{corhomotopy}
The surgery equivalence induces equivalences 
\[\Cobmfw{K}(\C, \QF) \xrightarrow{\chifw} \Fun(\lcone K, \Cob(\C, \QF)) \quad \text{and} \quad \Cobmbw{K}(\C, \QF) \xrightarrow{\chibw} \Fun((\lcone K)\op, \Cob(\C, \QF))\]
natural in $K \in \Cat$ and $(\C, \QF) \in \Catp$. 
\end{corollary}

From the case $K=\emptyset$ we learn that the extraction of underlying objects on the left hand side corresponds to the evaluation at the cone point on the right hand side, and for $K=\ast$ we recover the transformation
\[\chifw \colon \fgt \Longrightarrow \srg\]
on $\Cobfw(\C,\QF)$ constructed above, along with its backwards companion.

\begin{proof}
From the forwards surgery equivalence we obtain an equivalence
\[\Cobmfw{K}(\C, \QF\qshift 1)\simeq \Hom^1_{\Cat}(\lcone K, \Sq(\Cob(\C, \QF))) \simeq \Fun(\lcone K,\Cob(\C,\QF));\] the backwards case is analogous. 
\end{proof}

\subsection{Reflection of cubes in the cobordism category}
\label{sec:reflection_and_doubling}

In the present section we study the process of reflecting commutative squares in the cobordism category. More specifically, we analyse under what circumstances a commutative square in $\Cob(\C, \QF)$ as in the middle of 
\[\xymatrix{X \ar@{~>}[d]^U & Y \ar@{~>}[l]_{\refl(W)} \ar@{~>}[d]^V && X \ar@{~>}[r]^W \ar@{~>}[d]^U & Y \ar@{~>}[d]^V && X \ar@{~>}[r]^W & Y \\
            X' & Y' \ar@{~>}[l]_{\refl(T)} && X' \ar@{~>}[r]^T & Y' && X' \ar@{~>}[r]^T\ar@{~>}[u]_{\refl(U)} & Y' \ar@{~>}[u]_{\refl(V)}}\]
induces commutative squares as on the left and the right, and how this can be read off from surgery data defining the original square. To this end view the middle square as an element in $\Sq(\Cob(\C,\Q)) \subseteq \Cob^{(2)}(\C,\QF\qshift{-1})$. Then the outer two diagrams canonically upgrade to squares in the double category $\Cob^{(2)}(\C,\QF\qshift{-1})$ as well and we are asking that they lie in $\Sq(\Cob(\C,\QF))$. Whenever this is correct the squares above combine into a larger commutative square
\[\xymatrix{X \ar@{~>}[r]^W \ar@{~>}[d]^U & Y \ar@{~>}[d]^V \ar@{~>}[r]^{\refl(W)} & X  \ar@{~>}[d]^U \\
              X' \ar@{~>}[r]^T\ar@{~>}[d]^{\refl(U)} & Y' \ar@{~>}[r]^{\refl(T)} \ar@{~>}[d]^{\refl(V)}& X' \ar@{~>}[d]^{\refl(U)}\\   
              X \ar@{~>}[r]^W & Y \ar@{~>}[r]^{\refl(W)} & X}\]
which we call the double of the original square. This doubling construction will be important while doing surgery in the middle dimension; roughly, it will allows us to undo the effect of certain surgeries while keeping other surgeries done. We will also need the variants for higher dimensional and longer cubes later, so treat diagrams of the form $[a]^r \rightarrow \Cob(\C,\QF)$ generally for natural numbers $a,r>0$.

To get started recall the two embeddings of Poincar\'e categories
\[\fw\colon \Q_{[m] \times [n]}(\C, \QF) \to \Q_m \Q_n (\C, \QF) \leftarrow \Q_{[m] \times [n]}(\C, \QF)\colon \bw\]
where the first one is bisimplicial and the second one simplicial in the second and reverse simplicial in the first variable. The reflection equivalence $\Q(\C,\QF)^\mathrm{rev} \simeq \Q(\C,\QF)$ in the first or the second variable exchanges their images.

\begin{definition}
We denote by
\[\Qsq_{m,n}(\C, \QF)\subseteq \Q_m\Q_n(\C, \QF)\] 
the Poincar\'e category given by the intersection of the essential images of $\fw$ and $\bw$. A diagram in the underlying stable $\infty$-category $\Q^\square_{1,1}(\C)$ will be called \emph{independent square}.
\end{definition}

Explicitly, a diagram $X\colon \sd([1]\times[1])\to \C$ is an independent square if in
\[\xymatrix{
 X_{0\leq 0,0 \leq 0} &\ar[r]\ar[l] X_{0\leq 1,0 \leq 0} & X_{1 \leq 1, 0 \leq 0} \\
 X_{0\leq 0,0 \leq 1} \ar[u]\ar[d]&\ar[r]\ar[l]\ar[u]\ar[d] X_{0\leq 1,0 \leq 1} & X_{1 \leq 1, 0 \leq 1} \ar[u]\ar[d]\\
 X_{0\leq 0,1 \leq 1} &\ar[r]\ar[l] X_{0\leq 1,1 \leq 1} & X_{1 \leq 1, 1 \leq 1} }\]
all four displayed squares are bicartesian.

Varying $m$ and $n$ we obtain a two-fold complete Segal object $\Qsq(\C, \QF)$ of $\Catp$, and applying Poincar\'e objects we obtain a two-fold complete Segal space $\Poinc(\Qsq(\C, \QF))$ which is a bisimplicial subspace of the nerve of $\DCob(\C, \QF\qshift 2)$ and, through the forwards embedding, also of $\Sq(\Cob(\C, \QF\qshift 1))$; as we will have to iterate this construction in a moment, we refrain from translating this definition to the multi-categorical perspective. The commutative squares in $\Cob(\C, \QF\qshift 1)$ which are independent are easily recognised by means of their surgery data as follows:

\begin{proposition}\label{lem:bicartesian_square_surgeres_to_independent}
Let $T$ be a $([1]\times [1])$-fold surgery datum on some $(C,q)\in \Poinc(\C, \QF)$, giving rise by surgery to an object $C_T\in \Poinc(\Q_{[1]\times [1]}(\C, \QF))\subseteq \Poinc(\Q_1\Q_1(\C, \QF))$. Then, the following are equivalent:
\begin{enumerate}
 \item $C_T$ is independent.
 \item The underlying diagram $T\colon [1]\times [1]\to \C$ is bicartesian. 
\end{enumerate}
\end{proposition}

The proof is essentially by inspection. It is made most transparent by the following observation:

\begin{lemma}\label{lem:2cells}
If $K$ has a terminal object $t$ then the composite
\[\Hom_{\Cat}(\lcone K,\Cob(\C,\QF)) \simeq \Pmsd K(\C,\QF\qshift 1) \longrightarrow \Hom_{\Cat}(K,\C)\]
extracting the underlying surgery data is given by
\[\big[F \colon \lcone K \rightarrow \Cob(\C,\QF) \big] \longmapsto \big[K \rightarrow \C, \quad k \mapsto \fib(\widetilde F_{* \rightarrow t} \rightarrow \widetilde F_{k \rightarrow t})\big],\]
where $\widetilde F \colon \Twar(\lcone K) \rightarrow \C$ is the diagram adjoint to $F$ and $*$ denotes the cone point.
\end{lemma}

Let us unwind the statement for a moment: Per definition the image of $F$ on the right hand side takes some $k \in K$ to the surgery datum associated to the cobordism $\xymatrix{F(*) \ar@{~>}[r] & F(k)}$ and per construction this is given by $\fib(\widetilde F_{* \rightarrow k} \rightarrow \widetilde F_{k \rightarrow k})$. But the square
\[\xymatrix{\widetilde F_{* \rightarrow t} \ar[r] \ar[d] & \widetilde F_{k \rightarrow t} \ar[d] \\
            \widetilde F_{* \rightarrow k} \ar[r] & \widetilde F_{k \rightarrow k}}\]
in $\Twar(\lcone K)$ is cartesian by definition of the $\Q$-construction and identifies the fibre above with the fibre in the statement of the lemma. The trouble is now that the displayed square in $\C$ is not natural in $k$ in any sense. But for \ref{lem:bicartesian_square_surgeres_to_independent} we shall require a natural identification and the point of the lemma that this can indeed be achieved.

\begin{proof}
As explained in Section \ref{sec:review_of_surgery} the extraction of surgery data gives a map
\[\Hom_{\Cat}(\pcone(K),\Cob(\C,\QF)) \longrightarrow \Hom_{\Cat}(K,\C)\]
for an arbitrary category $K$ via. According to the constructions in \cite[Proposition 2.4.8 \& Corollary 2.4.9]{9authII} it is given for $K=[n]$ by
\[\Hom_{\Cat}(\pcone([n]),\Cob(\C,\QF)) = \Poinc(\Q_{1+n}(\C,\QF)) \xrightarrow{s_n} \Poinc(\mathrm{Pair}(\Met((\C,\QF)^{[n]})))\]
\[\xrightarrow{\fgt} \core(\Ar(\C^{[n]})) = \Hom_{\Cat}([n] \times [1],\C) \xrightarrow{\ev_0} \Hom_{\Cat}([n],\C);\]
the evaluation at $1 \in [1]$ in the last step extracts the constant diagram on the underlying object (i.e.\ the evaluation at the cone point of $\lcone K$). Unwinding the definition of $s_n$ the above composite is also given (still for $K=[n]$) by the composite
\begin{align*}
\Hom_{\Cat}(\pcone(K),\Cob(\C,\QF)) \subseteq& \Hom_{\Cat}(K \times [1], \Cob(\C,\QF)) \\
\xrightarrow{\fgt}& \Hom_{\Cat}(K \times [1],\Span(\C)) \\
\subseteq &\Hom_{\Cat}(\Twar(K \times [1]),\C) \\
\xrightarrow{\iota^\ast}& \Hom_{\Cat}(\Twar(K) \times [1], \C) \\
= &\Hom_{\Cat}(\Twar(K), \Ar(\C)) \\
\xrightarrow{\fib}& \Hom_{\Cat}(\Twar(K), \C)\\
\xrightarrow{s_!}&\Hom_{\Cat}(K, \C)
\end{align*}
with $\iota \colon \Twar(K) \times [1] \rightarrow \Twar(K) \times \Twar([1]) = \Twar(K \times [1])$ induced by the identity in the first coordinate and the arrow $(0 \leq 1 ) \rightarrow (1 \leq 1)$ in the second, and $s \colon \Twar(K) \rightarrow K$ is the source projection. But the above composite constitutes a natural transformation in arbitrary $K\in \Cat\op$: This is evident for all steps except the left Kan extension in the very last one, which is only lax natural via the Beck-Chevalley transformations of the adjunction between $s_!$ and the (evidently natural) $s^*$; see \cite{HHLN1}. But this Kan extension is an optical illusion: Since $s \colon \Twar(K) \rightarrow K$ is a cartesian fibration with contractible fibres it is a localisation \cite[Lemma 5.5]{HHLN2}, so that $s^* \colon \Hom_{\Cat}(K,\C) \rightarrow \Hom_{\Cat}(\Twar(K),\C)$ is an equivalence oonto those components $\Hom_{\Cat}^s(\Twar(K),\C)$ in the target consisting of all functors $\Twar(K) \rightarrow \C$ inverting all maps in $\Twar(K)$ whose source component is an equivalence. It follows that on this part the requisite Beck-Chevalley transformations are invertible, and $s_! \colon \Hom_{\Cat}^s(\Twar(K),\C) \rightarrow \Hom_{\Cat}(\Twar(K),\C)$ is natural in $K$. It remains to check that composite of the first six steps takes values in $\Hom_{\Cat}^s(\Twar(K),\C)$. But some $F \colon K \times[1] \rightarrow \Cob(\C,\QF)$ and a map $(k \rightarrow l) \rightarrow (k' \rightarrow l')$ in $\Twar(K)$ are indeed taken to the induced map on horizontal fibres in
\[\xymatrix{
\widetilde F(k\rightarrow l,0 \leq 1) \ar[r] \ar[d] & \widetilde F(k \rightarrow l,1 \leq 1) \ar[d] \\
\widetilde F(k'\rightarrow l',0 \leq 1) \ar[r]  & \widetilde F(k' \rightarrow l',1 \leq 1) \\
}\] 
and in case the map $k \rightarrow k'$ is an equivalence, this diagram is in the image of $\Hom_{\Cat}([3],\Cob(\C,\QF)) = \Poinc\Q_3(\C,\QF\qshift{-1})$ under the map 
\[[3] \longrightarrow K \times[1], \quad (k,0) \rightarrow (k,1) \rightarrow (l',1) \rightarrow (l,1).\]
It is therefore cartesian by definition of the $\Q$-construction. 


The proof of the actual statement is now short: If $K$ has a terminal object $t$, then $|K| \simeq \ast$, so $\pcone(K) = {\lcone K}$ and the source functor $s \colon \Twar(K) \rightarrow K$ admits $k \mapsto (k \rightarrow t)$ as a fully faithful adjoint, so that the resulting functor $\Hom_{\Cat}(\Twar(K),\Ar(\C)) \rightarrow \Hom_{\Cat}(K,\Ar(\C))$ agrees with $s_!$ and can thus be used to extract the the diagram $K \rightarrow \C$ in the last step. This results in the functor from the statement being induced by the functoriality of
\[k \longmapsto \fib(\widetilde F((k,0) \rightarrow (t,1)) \rightarrow \widetilde F((k,1) \rightarrow (t,1))\]
and for functors pulled back along $K \times[1] \rightarrow {\lcone K}$ this is equivalent to the formula we want to show. 
\end{proof}

\begin{proof}[Proof of Proposition \ref{lem:bicartesian_square_surgeres_to_independent}]
Let 
 \[\xymatrix@-1pc{
  a \ar[r] \ar[d] & b \ar[d] \\ c \ar[r] & d
 }\]
 be the underlying diagram in $\C$ of a $[1] \times [1]$-fold surgery datum on a Poincar\'e object $C$ of $\C$. The result of surgery is then given by:
 \[\xymatrix{
 \fib(C\to \Dual_\QF(a))/a  & \fib(C\to \Dual_\QF(b))/a \ar[l]  \ar[r] & \fib(C \rightarrow \Dual_\QF(b))/b \\
 \fib(C\to \Dual_\QF(c))/a \ar[u]\ar[d] & \fib(C\to \Dual_\QF(d))/a \ar[l]\ar[u]\ar[r]\ar[d]  & \fib(C \rightarrow \Dual_\QF(a))/b \ar[d]\ar[u]\\
 \fib(C \rightarrow \Dual_\QF(c))/c & \ar[l] \fib(C \rightarrow \Dual_\QF(d))/c \ar[r] & \fib(C \rightarrow \Dual_\QF(d))/d }\]
with maps induced by the given diagram in $\C$. The claim now essentially follows by inspection, once we have also identified the homotopies in the squares above with the ones induced by the original square. In fact, by the final comment in \ref{ex:Q_1Q_1} it suffices to do this in the lower right square, where it is a consequence of \ref{lem:2cells}: As part of the diagram on $\lcone([1] \times [1])$ this square extends to
 \[\xymatrix{
 \fib(C\to \Dual_\QF(d)) \ar@/_0.5pc/[rd] \ar@/^1pc/[rrd] \ar@/_1pc/[rdd] \ar@{..>}@/^1.5pc/[rrdd]& \\
 & \fib(C\to \Dual_\QF(d))/a \ar[r]\ar[d]  & \fib(C \rightarrow \Dual_\QF(d))/b \ar[d]\\
  &  \fib(C \rightarrow \Dual_\QF(d))/c \ar[r] & \fib(C \rightarrow \Dual_\QF(d))/d }\]
and unwinding definitions \ref{lem:2cells} precisely identifies the diagram formed by the fibres of the slanted maps as the original original surgery datum, which implies the claim.
\end{proof}

We will also need a version of these considerations for higher cubes. We therefore put:

\begin{definition}
We define $\Qsq_{n_1,\dots,n_r}(\C,\QF) \subseteq \Q_{n_1}\dots\Q_{n_r}(\C,\QF)$ as the full subcategory spanned by the intersection of $\Q_{[n_1] \times \dots \times [n_r]}(\C,\QF)$ with its $r$ reflections.
\end{definition}

Since the forwards embedding and the reflections are Poincar\'e functors, this subcategory is closed under the duality of $\Q_{n_1}\dots\Q_{n_r}(\C,\QF)$ and thus Poincar\'e itself. Similarly, one checks that it really defines an $r$-fold complete Segal object of $\Catp$. The $r$-fold category associated to its Poincar\'e objects is thus a subcategory of $\Cob^{(r)}(\C,\QF\qshift{-r})$, the same objects and morphisms (in all directions), but severely restricted higher squares. It can be regarded as the largest subcategory of the $r$-fold category of commutative $r$-cubes in $\Cob(\C,\QF\qshift{-1})$ that is closed under all $r$ reflections.

Our next goal is to characterise the pieces of $[1]^r$-fold surgery data that give rise to objects in $\Poinc\Qsq_{1,\dots,1}(\C,\QF)$ in generalisation of \ref{lem:bicartesian_square_surgeres_to_independent} above. 
It will be convenient to do this directly for $[a]^r$-fold cubes (and of course the case of an arbitrary product is in the end only notationally more involved).

We start out by characterising the image of $\Q_{[1]^r}(\C,\QF) \subseteq \Q_1^{(r)}(\C,\QF)$. To this end abbreviate $\Twar[1]$ by $H$ (for horn). Then the underlying category of the target consists of diagrams $H^r \rightarrow \C$. Now, the category $H^r$ is glued from $2^r$ many $r$-dimensional cubes $[1]^r$, that we index as follows: For some $v \in [1]^r$ set $H^r_v$ to be the full subcategories spanned by all $u \leq w \in \H^r$ with $u \leq v \leq w$. We will refer to the to $r$-cubes $H_0^r$ and $H_1^r$ as the initial and terminal cube of $H^r$, respectively.

\begin{lemma}\label{lem:Q1r}
A diagram $X \colon H^r \rightarrow \C$ lies in $\Q_{[1]^r}(\C,\QF)$ if and only if all of the $H^r_v$ but possibly the initial and terminal ones are taken to strongly cocartesian cubes in $\C$, and $X \in \Qsq_{1,\dots,1}(\C,\QF)$ if all $H_v^r$ are taken to strongly cocartesian cubes.
\end{lemma}

In particular, we see that a Poincar\'e object in $\Q_{[1]^r}(\C,\QF)$ lies in $\Poinc\Qsq_{1,\dots,1}(\C,\QF)$ if its terminal cube is strongly cocartesian.

\begin{proof}
We begin with the first statement. So let $X$ satisfy the condition on its subcubes. Given then a chain $u \leq v \leq w$ in $[1]^r$ we have to check that the square 
\[\xymatrix{u \leq w \ar[r]\ar[d] & u \leq v \ar[d] \\  
            v \leq w \ar[r]       & v \leq v}\]
goes to a cartesian square. If one of the inequalities is actually an equality this is trivial and otherwise, this square can be pasted from $2$-dimensional faces of $H_v^r$ (which are assumed cocartesian): Namely, choose a sequence 
\[u = u_0 < u_1 < \dots < u_n < v =v_0 < v_1 \dots < v_m = w\] such that the difference in each step is only a single coordinate and consider the squares
\[\xymatrix{u_i \leq v_{j+1} \ar[r] \ar[d] & u_i \leq v_j \ar[d] \\
            u_{i+1} \leq v_{j+1} \ar[r] & u_{i+1} \leq v_j.}\]
We will prove the converse by induction on $r$. The cases $r=0,1$ are vacuous and $r=2$ is precisely \ref{ex:Q_1Q_1}. Now per construction 
\[X \in \Q_{[1]^{r+1}}(\C,\QF) \subseteq \Q_{[1]^r}\Q_1(\C,\QF),\]
so we learn from the induction hypothesis that for any of the $r$-cube $H^r_v$ in $H^r$, that is neither the initial nor terminal one, the three $r$-cubes
\[H^r_v \times \{0 \leq 0\}, \quad H^r_v \times \{0 \leq 1 \}, \quad \text{and} \quad H^r_v \times \{1 \leq 1\}\]
of $H^{r+1}$ go to strongly cocartesian diagrams under $X$. Since an $r+1$-cube is strongly cocartesian if two of its opposing faces are, this implies  that all but potentially the four cubes
\[H_0^{r+1} = H^r_0 \times \{(0 \leq 1) \rightarrow (0 \leq 0)\}, \quad H_0^r \times \{(0 \leq 1) \rightarrow (1 \leq 1)\}, \quad H_1^r \times \{(0 \leq 1) \rightarrow  (0\leq 0)\},\] \[\text{and} \quad H_1^{r+1} = H_1^r \times \{(0 \leq 1) \rightarrow (1 \leq 1)\}\]
are cocartesian. But splitting off the first instead of the last coordinate of $[1]^{r+1}$ covers the middle two.

To obtain the second statement simply observe that any reflection exchanges the initial and terminal cubes for one of the others.
\end{proof}

From here it is easy to characterise the associated surgery data. To formulate the result we need:

\begin{lemma}\label{lem:strongly_cocartesian_of_sidelength_a}
For a finite set $S$, a number $a \in \mathbb N$, and a diagram $X\colon [a]^S\to \C$ in an $\infty$-category $\C$ with finite colimits, the following are equivalent:
\begin{enumerate}
 \item\label{item:strongly_cocartesian_1} $X$ is left Kan extended from the full subposet of functions supported in at most one element, 
 \item\label{item:strongly_cocartesian_2} for each collection of strictly monotone maps
 \[f_s\colon [1]\longrightarrow [a]\quad\mathrm{with}\quad s(0)=0 \quad (s\in S),\]
 the restricted cubical diagram of side-length $1$
 \[[1]^S \xrightarrow{\prod_{s\in S} f_s} [a]^S \xrightarrow{X} \C\]
 is strongly cocartesian, 
 \item\label{item:strongly_cocartesian_3} for each collection of strictly monotone maps $f_s\colon [1]\longrightarrow [a]$, with $s\in S$, the restricted cubical diagram of side-length $1$
 \[[1]^S \xrightarrow{\prod_{s\in S} f_s} [a]^S \xrightarrow{X} \C\]
 is strongly cocartesian.
\end{enumerate}
\end{lemma}

\begin{definition}\label{def:strongly_cocartesian}
We call a diagram $X\colon [a]^S\to \C$ \emph{strongly cocartesian} if it satisfies the equivalent conditions of Lemma \ref{lem:strongly_cocartesian_of_sidelength_a}.
\end{definition}

\begin{proof}[Proof of Lemma \ref{lem:strongly_cocartesian_of_sidelength_a}]
For a map $g\colon S\to [a]$, we consider the poset $P$ of all $h\colon S\to [a]$ with $h\leq g$ (pointwise) and $h$ supported in at most a single element, and the subposet $P'\subset P$ where $h\in P$ is discarded if there exists some $s$ such that $0<h(s)<g(s)$. Then we have maps 
\[\colim_{h\in P'} X(h)\longrightarrow \colim_{h\in P} X(h) \longrightarrow X(g)\]
where the first is an equivalence since $P'\subset P$ is cofinal. By the pointwise formula for Kan extensions, condition (\ref{item:strongly_cocartesian_1}) is equivalent to the right map being an equivalence, for all $g$, while condition (\ref{item:strongly_cocartesian_2}) is equivalent to the composite map being an equivalence.

This proves the equivalence of (\ref{item:strongly_cocartesian_1}) and (\ref{item:strongly_cocartesian_2}). Clearly (\ref{item:strongly_cocartesian_3}) implies (\ref{item:strongly_cocartesian_2}), and the converse is proven by induction in the number $n$ of elements $s$ such that $f_s(0)>0$: The induction beginning $n=0$ is precisely condition (\ref{item:strongly_cocartesian_2}), and the inductive step is proven using the pasting law for strongly cocartesian cubes (which easily follows from the pasting law for cocartesian cubes, using \cite[6.1.1.15]{HA}). 
\end{proof}

\begin{proposition}\label{lem:strongly_cocartesian_surgeres_to_independent}
Let $T$ be a $\powerar$-fold surgery datum on some $C\in \Poinc(\C, \QF)$, giving rise by surgery to  $C_T\in \Poinc(\Q_{\powerar}(\C, \QF))\subset \Poinc(\Q_a\dots \Q_a(\C, \QF))$. Then, the following are equivalent:
\begin{enumerate}
 \item $C_T\in \Poinc(\Qsq_{a,\dots, a}(\C, \QF))$.
 \item The underlying diagram $T\colon \powerar\to \C$ is strongly cocartesian.
\end{enumerate}
\end{proposition}

\begin{proof}
By the Segal condition and Lemma \ref{lem:strongly_cocartesian_of_sidelength_a}, it is enough to consider the case $a=1$ (the extra condition for non-injective maps $[1]\to [a]$ is vacuous).

Using \ref{lem:Q1r} in place of \ref{ex:Q_1Q_1} the proof of \ref{lem:bicartesian_square_surgeres_to_independent} applies essentially verbatim to this case upon replacing the lower right square occuring there by the terminal cube in $H^{r}$; alternatively one can also deduce the statement from \ref{lem:Q1r} and \ref{lem:bicartesian_square_surgeres_to_independent} by induction on $r$: Lemma \ref{lem:Q1r} is the induction start. Using that a cube is strongly cocartesian if two of its opposing faces are, we find by induction hypothesis that (2) is equivalent to two opposing boundary faces of $C_T \colon H^{r+1} \rightarrow \C$ lying in $\Qsq_{1,\dots,1}(\C,\QF)$; let us call $F$ the one intersecting the initial cube of $H^{r+1}$ and $G$ the other one, which intersects the terminal cube. By \ref{lem:bicartesian_square_surgeres_to_independent} the condition $F,G \in \Qsq_{1,\dots,1}(\C,\QF)$ means precisely that all constituent $r$-cubes of $F$ and $G$ are strongly cocartesian. We now claim that this is the case if and only if $C_T$ takes the terminal and initial cubes of $H^{r+1}$ to strongly cocartesian cubes in $\C$.

To see the claim, let us take stock of the cubes that are send to cocartesian cubes by $C_T$ regardless of the assumptions in (1) and (2): Since $\C_T \in \Poinc(\Q_{[1]^{r+1}}(\C,\QF))$ \ref{lem:bicartesian_square_surgeres_to_independent} implies that all but the initial and terminal $r+1$-cubes of $H^{r+1}$ are taken to strongly cocartesian cubes by $C_T$. This implies that, on the one hand, all faces of the initial and terminal cubes of $H^{r+1}$ that do not contain the 'outer corners' (i.e. the vertices $(0 \leq 0, \dots, 0 \leq 0)$ and $(1 \leq 1, \dots, 1 \leq 1)$) are strongly cocartesian. On the other, also all $r$-cubes in $F$ but its initial cube and all cubes in $G$ but the terminal one are strongly cocartesian as a consequence.

The claim we made above is thus equivalent to the statement that the initial cube of $F$ and the terminal cube of $G$ are strongly cocartesian if and only if the initial and terminal cube of $H^{r+1}$ are. With half the faces of these $r+1$-cubes already strongly cocartesian this follows from another application of the fact that a cube is strongly cocartesian if two opposing faces are.

In total, we have thus shown that (2) is equivalent to all constituent cubes of $C_T \colon H^{r+1} \rightarrow \C$ being strongly cocartesian. Another application of \ref{lem:bicartesian_square_surgeres_to_independent} thus finishes the proof.
\end{proof}

We now proceed to construct the doubles of cubes in $\Cob(\C,\QF)$. This is facilitated by the following \emph{reflection lemma}. For its formulation recall that a complete Segal object uniquely extends to a limit preserving functor on $\Cat\op$. In particular, we obtain 
\[\Qsq \colon \Cat^r \longrightarrow \Catp\]
as such an extension for each $r$.
Since each $\Qsq_{n_1, \dots, n_r}(\C, \QF)\subseteq \Q_{n_1}\dots \Q_{n_r}(\C, \QF)$ is a full subcategory, the same is true for $\Qsq_{K_1, \dots, K_r}(\C, \QF)\subseteq \Q_{K_1}\dots \Q_{K_r}(\C, \QF)$, and the inclusion functor is Poincar\'e, $\Catp\subset \Cath$ being closed under limits.

\begin{lemma}\label{lem:reflection_lemma}
For any Poincar\'e category $(\C, \QF)$, reflection in the second summand and the identity in the first summand of each entry combine to a canonical equivalence of Poincar\'e categories for any $a\in \NN$
\[\Qsq_{[a]\cup_{\{a\}} [a], \dots, [a]\cup_{\{a\}} [a]}(\C, \QF) 
\longrightarrow 
\Qsq_{[2a], \dots , [2a]}(\C, \QF).\]
\end{lemma}

Both sides in fact exhibit semi-simplicial functoriality in the number of their indices: The values at $[r]$ are given by the expressions with $|[r]| = r+1$ indices. Informally, the face map $d_i$ is given by evaluation at $0$ in $[2a]$ or the first summand of $[a]\cup_{\{a\}} [a]$ in the coordinate labelled by $i\in [r]$, and more formally, arises on the right from the inclusion 
\[\Qsq_{[2a], \dots , [2a]}(\C, \QF) \subseteq \Q_{\prod_{[r]} [2a]}(\C,\QF)\]
which one checks closed under the face map of the co-semi-simplicial structure of $r \mapsto \prod_{[r]} [2a]$ given by the prescription above ($r \mapsto \Q_{\prod_{[r]} [2a]}(\C,\QF)$ in fact canonically extends to a simplicial object by doubling entries, but this does not persist to the left hand side). The left hand side in \ref{lem:reflection_lemma} is entirely analogous. The proof of \ref{lem:reflection_lemma} will also show:

\begin{add}\label{rem:reflection_lemma_natural}
The equivalence of the reflection lemma is natural in the semi-simplicial structure just described.
\end{add}

\begin{proof}[Proof of \ref{lem:reflection_lemma}]
We have maps of posets
\[\Twar([a]\cup_{\{a\}} [a])^{r+1}\to (\Twar[2a])^{r+1}\]
that map, in each factor, the first summand to $\Twar([a]) \subseteq \Twar[2a]$ by the identity and the second summand to $\Twar([a,2a])$ via the map induced by $n \mapsto 2a-n$. Identifying the $r+1$-fold product with the product indexed over the set  $[r]$, both sides inherit co-semi-simplicial structures from the construction above and it is easily checked that the maps above are part of a natural transformation. Taking hermitian functor categories, we obtain a natural transformation
\begin{equation*}\tag{$\ast$}
\Fun((\Twar[2a])^{r+1}, (\C, \QF)) \to \Fun((\Twar([a]\cup_{\{a\}} [a])^{r+1}), (\C, \QF)) 
\end{equation*}
of semi-simplicial objects in $\Cath$. We claim that it restricts to an equivalence
\[\Qsq_{[2a], \dots , [2a]}(\C, \QF)\to \Qsq_{[a]\cup_{\{a\}} [a], \dots , [a]\cup_{\{a\}} [a]}(\C, \QF)\] 
in $\Catp$. 

To see this, we first note that the reflection map in one variable 
\[\Q_{[2a]}(\C, \QF)\to \Q_{[a]\cup_{\{a\}}[a]}(\C, \QF)\]
is equivalent, by the Segal condition for the left hand side and compatibility of the $\Q$-construction with colimits on the right, to the map
\[\Q_{[a]}(\C, \QF)\times_{\Q_{\{a\}}(\C, \QF)} \Q_{[a,2a]}(\C, \QF)\to \Q_{[a]}(\C, \QF)\times_{\Q_{\{a\}}(\C, \QF)} \Q_{[a]}(\C, \QF)\]
which is the identity on the first factor and reflection on the second and therefore an equivalence. 

We deduce iteratively that ($\ast$) restricts to an equivalence
\[\Q_{[2a]} \dots Q_{[2a]}(\C, \QF) \to \Q_{[a]\cup _{\{a\}}[a]} \dots \Q_{[a]\cup _{\{a\}}[a]} (\C, \QF)\]
and we are left to see that the full subcategory $\Qsq_{[2a], \dots, [2a]}(\C, \QF)$ on the left corresponds to the full subcategory $\Qsq_{[a]\cup_{\{a\}} [a], \dots , [a]\cup_{\{a\}} [a]}(\C, \QF)$ on the right. By the Segal conditions it suffices to consider restriction of the diagrams along neighboring edges of the form $i\leq i+1$ of $[2a]$ in each direction. Then, the map ($\ast$) is given by reflection in some of the entries and identity in the others and the claim follows from the fact that $\Qsq(\C, \QF)\subset \Q^{(r+1)}(\C, \QF)$ is closed under each of the reflections. 
\end{proof}

To obtain the doubled cube we now simply observe that an independent cube \[(C,q) \in \Poinc\Qsq_{1,1, \dots,1}(\C,\QF\qshift{1}) \subseteq \Hom_{\Cat}(\lcube 1 {r+1}, \Cob(\C,\QF))\] 
can be pulled back to an object of $\Poinc\Qsq_{[1]\cup_1 [1], \dots, [1]\cup_1 [1]}(\C, \QF)$ via the fold map $[1] \cup_{\{1\}}[1] \rightarrow [1]$, to which we can the apply the reflection lemma. There results a map
\[\Db \colon \Poinc\Qsq_{1,1, \dots,1}(\C,\QF\qshift{1}) \longrightarrow \Poinc\Qsq_{2,2, \dots,2}(\C,\QF\qshift{1}) \subseteq \Hom_{\Cat}(\lcube 2 {r+1},\Cob(\C,\QF))\]
natural in both $(\C,\QF)$ and also when regarded as a map of semi-simplicial objects in the number of indices.

\bigskip \noindent

We finally translate the doubling construction along the surgery equivalence. 

\begin{definition}\label{def:cocartsurgdata}
For $a,r\in \NN$, a \emph{disjoint $[r]$-tuple of surgery data of length $a$} in $(\C, \QF)$ is an $\lcube a {r+1}$-fold surgery datum in $(\C, \QF)$ whose underlying diagram $T\to X$ satisfies:
\begin{enumerate}
\item $T(0)=0$, and 
\item $T\colon \lcube a {r+1}\to \C$ is a strongly cocartesian cubical diagram.
\end{enumerate}
We define $\Pccsd a {r}(\C, \QF)\subseteq \Pmsd{\lcube a {r+1}}(\C, \QF)$ to be the full subspace spanned by such surgery data, and we let  $\Cobccsd a{r}(\C,\QF)$ to be the category associated to the simplicial space $\Pccsd a {r}(\Q(\C, \QF\qshift 1))$.
\end{definition} 

\begin{remark}
By definition, the underlying diagram of such a disjoint $[r]$-tuple is determined by the $(r+1)$-many $\{1\leq\dots\leq a\}$-fold pieces of surgery data that make up the edges of the cube emanating from $0$. However, we warn the reader that the similar statement is not true for the hermitian structures: The functor $\QF$ is not linear and will thus generally not send strongly cocartesian cubical diagrams in $\C$ to strongly cartesian cubical diagrams in $\Spa$. Thus the word ``tuple'' in the definition above should not be taken too literal; informally, disjointness of surgery data is an extra datum rather than just a property. We will see below, that disjointness is, however, still a pairwise datum, in the sense that the hermitian structure is determined by its restriction to all $2$-dimensional faces emanating from $0$ in the given cube since $\QF$ is quadratic, see subsection~\ref{coskeletalsection} below. 
\end{remark}

It is easily seen that $\Pccsd a {r} \colon \Catp \rightarrow \Sps$ preserves pullbacks, so that $\Pccsd a {r}(\Q(\C, \QF))$ is a complete Segal space and thus $\Cobccsd a {r}(\C, \QF)$ a (non-full) subcategory of $\Cobmsd{\lcube a {r+1}}(\C, \QF)$. By the discussion preceeding \ref{rem:reflection_lemma_natural} the association 
\[[r] \longmapsto \Cobccsd a {r}(\C, \QF)\]
defines a semi-simplicial category $\Cobccsd a {-}(\C, \QF)$. The notation is meant to suggest that an element of $\Pccsd a r(\C, \QF)$ consists of $[r]$ many disjoint $[a]$-fold surgery data in $\Pmsd{[a]}(\C, \QF)$.

Combining \ref{lem:strongly_cocartesian_surgeres_to_independent} and \ref{lem:reflection_lemma}, we obtain:

\begin{corollary}
Reflection of surgery data induces a map 	
\[\Db \colon \Cobccsd 1{r}(\C,\QF)  \longrightarrow \Cobccsd 2 {r}(\C,\QF),\]
of semi-simplicial categories (in the $r$-variable), such that the diagram
\[\xymatrix{
\nerv_n \Cobccsd 1{r}(\C,\QF) \ar[r]^\Db \ar[d]_\chi 
  & \nerv_n \Cobccsd 2{r}(\C,\QF) \ar[d]_\chi 
\\
 \Poinc\Qsq_{1,1, \dots,1}(\Q_n(\C,\QF\qshift{1})) \ar[r]^\Db 
   & \Poinc\Qsq_{2,2, \dots,2}(\Q_n(\C,\QF\qshift{1}))
}\]
commutes.
\end{corollary}

We finally remark that we shall also denote the surgery datum of the composite
\[\xymatrix{C \ar@{~>}[r]^{\chi_f} & C_f \ar@{~>}[r]^{\refl(\chi_f)} & C}\]
in $\Cob(\C,\QF)$ for some $\big[(f,\eta) \colon S \rightarrow (C,q)\big] \in \Surg(\C,\QF\qshift{1})$ by $\Db(f)$. It is easy to check that its underlying map is given by $(f,0) \colon S \oplus \Dual_\QF(S) \rightarrow C$, whence
\[(f^*q,0,0) \simeq \Db(f)^*q \in \pi_0\QF\qshift{1}(S \oplus \Dual_\QF(S)) \simeq \pi_0\QF\qshift{1}(S) \oplus \pi_{-1}\hom_\C(S,S) \oplus \pi_0\QF\qshift{1}(\Dual_\QF(S)).\] It is slightly more work to check that the induced nullhomotopy of this form is given by $(\eta,\id_S,0)$, but conveniently we shall not have to make use of this latter fact, so leave details to the reader.

\section{Surgery complexes}\label{sec:cplx}

Just as in \cite{GRW} the proof of Theorem \ref{thm:main} consists of three steps: surgery on objects below and then in the middle dimension and surgery on morphisms. These are the contents of sections \ref{sec:objects} and \ref{sec:morph}, respectively. The goal of the present section is to expound the common structure of the arguments in \ref{sketch}, and establish a key property of what we call multiple disjoint surgery data in \ref{coskeletalsection}, that is central to all three cases.

\subsection{The surgery argument}\label{sketch}

We start out by giving an outline of the surgery arguments. Let us first define the filtration of the cobordism category we will study in the remainder of this paper.

\begin{definition}
Let $(\C,\QF)$ be a Poincar\'e category, whose underlying stable $\infty$-category is equipped with a weight structure. We call a cobordism $C_{0\leq 0} \longleftarrow C_{0\leq 1} \longrightarrow C_{1\leq 1}$
\emph{$\mc$-connective} if the left pointing map is $\mc$-connective. We denote by 
\[\Cob^\mc(\C, \QF)\subset \Cob(\C, \QF)\]
the subcategory containing all objects and the $\mc$-connective cobordisms; and by 
\[\Cob^{\mc,\oc}(\C, \QF)\subset \Cob^\mc(\C, \QF)\]
the further full subcategory spanned by the $\oc$-connective objects only. 
\end{definition}

Let us now describe the strategy for surgery on objects below the middle dimension, which is the simplest of the cases, in more detail. For a Poincar\'e category of some dimension $d$, and numbers $\oc$ and $\mc$ satisfying certain numerical assumptions, we wish to show that the inclusion
\[\Cob^{\mc,\oc+1}(\C,\QF) \longrightarrow \Cob^{\mc,\oc}(\C,\QF)\]
becomes an equivalence upon realisation. This is achieved by constructing an inverse, via choosing appropriate surgery data and performing the parametrised surgery developed in the previous section. To be of use pieces of surgery data have to meet two design criteria: They have to improve the connectivity of objects while not lowering the connectivity of morphisms. To formalise this, recall the cobordism category with forwards surgery data $\Cobfw(\C, \QF)\subset \Cobsd(\C, \QF)$ from \ref{def:forwards}, and denote
\[\Cobfw^{\mc,\oc}(\C,\QF)\subset \Cobfw(\C,\QF)\]
the subcategory spanned by those objects, for which
\begin{enumerate}
\item the arrow $[1] \rightarrow \Cob(\C,\QF)$ resulting from surgery lies in $\Cob^{\mc,\oc}(\C,\QF)$, and
\item its endpoint $\{1\} \rightarrow \Cob(\C,\QF)$ even lies in $\Cob^{\mc,\oc+1}(\C,\QF)$.
\end{enumerate}
and similarly morphisms such that
\begin{enumerate}
\item the square $[1] \times [1] \rightarrow \Cob(\C,\QF)$ resulting from surgery lies in $\Cob^{\mc,\oc}(\C,\QF)$, and
\item its restriction to the end $[1] \times \{1\} \rightarrow \Cob(\C,\QF)$ (i.e. the surgered parts) even lies in $\Cob^{\mc,\oc+1}(\C,\QF)$.
\end{enumerate}

Now consider the following (solid) commutative diagram
\[\xymatrix{
 \Cob^{\mc,\oc+1}(\C,\QF) \ar[d] 
   & \Ar(\Cob^{\mc,\oc}(\C,\QF)) \times_{\Cob^{\mc,\oc}(\C,\QF)} \Cob^{\mc,\oc+1}(\C,\QF) \ar[l]_-{\ev_1} \ar[dl]_-{\ev_0} 
 \\
 \Cob^{\mc,\oc}(\C,\QF) \ar@/_{2ex}/@{.>}[r]_s
   & \Cobfw^{\mc,\oc}(\C,\QF) \ar[u]^{\chifw} \ar[l]_{\fgt}
}
\]
where the pullback in the top right corner is formed using the evaluation of an arrow at its target.

Suppose one could find a section $s \colon |\Cob^{\mc,\oc}(\C,\QF)| \rightarrow |\Cobfw^{\mc,\oc}(\C,\QF)|$ as indicated which makes the diagram 
\[\xymatrix{|\Cob^{\mc,\oc+1}(\C,\QF)| \ar[rd] \ar[rr]^{0 \Rightarrow \mathrm{id}} && |\Cobfw^{\mc,\oc}(\C,\QF)| \\
& |\Cob^{\mc,\oc}(\C,\QF)| \ar[ru]^s}\]
commute, where the top horizontal map adds the surgery datum $0$. Then a straight-forward chase in the square above shows that the composite along the right side of the rectangle (using $s$ as the first map) gives an inverse to the inclusion $|\Cob^{\mc,\oc+1}(\C,\QF)| \rightarrow |\Cob^{\mc,\oc}(\C,\QF)|$; in fact, the square identifies the inclusion $|\Cob^{\mc,\oc+1}(\C,\QF)| \rightarrow |\Cob^{\mc,\oc}(\C,\QF)|$ with the forgetful map $|\Cobfw^{\mc,\oc}(\C,\QF)| \rightarrow |\Cob^{\mc,\oc}(\C,\QF)|$, and the two condition imposed on $s$ make it a both sided inverse to the latter (one does not really need to invest that $\chifw$ is an equivalence, however, and such a statement would indeed not be true in the analogous situation in \cite{GRW}).

\bigskip\noindent

However, we do not know how to directly produce such a split $s$ and in fact, given the arbitrariness of the choices of surgery data it may seem surprising that it even exists (presumably it can in particular not be obtained from a section at the level of categories). Once translated to our language, the key insight originally due to Madsen and Weiss, and more explicitly articulated by Galatius and Randal-Williams (see \cite[Section 3.3]{GRW}), is that this obstacle can be circumvented by introducing another degree of freedom and allowing multiple surgeries on disjoint surgery data in the above argument; we now set out to explain our implementation of this idea.

Roughly speaking, we replace the category $\Cobfw^{\mc,\oc}(\C,\QF)$ by the (semi)simplicial subcategory $[r] \mapsto \Cobccsd 1 {r}^{\mc,\oc}(\C,\QF)$ of $\Cobmsd{\lcube 1 {r+1}}(\C,\QF)$ (compare \ref{def:cocartsurgdata}) spanned by disjoint $[r]$-tuples of surgery data of length $1$ that  on $[1]^{r+1}\setminus \{0\}$ satisfy a connectivity assumption similar to that imposed on $\Cobfw^{\mc,\oc}(\C,\QF)$, see \ref{connectedsurgery} below.

The first two conditions can be read as the cube of surgery data consisting of $r+1$ individual pieces of surgery data (the evaluations at the singleton sets), that are `disjoint': In situation of surgery on manifolds the union of two sets of embedded spheres with trivial normal bundle can be regarded as a single surgery datum if and only if they are disjoint from one another. Let us also note again, that this condition is preserved by the boundary, but not the degeneracy maps (doubling data clearly does not preserve disjointness), so we only obtain a semisimplicial category this way. The map $\srg \colon |\Cobfw^{\mc,\oc}(\C,\QF)| \rightarrow |\Cob^{\mc,\oc+1}(\C,\QF)|$ then extends to a map
\[|\Cobccsd 1 {-}^{\mc,\oc}(\C,\QF)| \longrightarrow |\Cob^{\mc,\oc+1}(\C,\QF)|,\]
whence it suffices to split
\[|\Cobccsd 1 {-}^{\mc,\oc}(\C,\QF)| \longrightarrow |\Cob^{\mc,\oc}(\C,\QF)|\]
suitably, see \ref{generalmethodobjects} below, for a precise statement. 

The reason one has a chance at finding this kind of split more easily, is that the particular type of surgery data we shall consider, while not unique per se, are `unique up to disjointness' under appropriate dimension assumptions. More precisely, we shall consider a subcategory $\SC_r^{\mc,\oc}$ inside $\Cobccsd 1 {r}^{\mc,\oc}(\C,\QF)$, which, following Galatius and Randal-Williams, we dub the surgery complex, such that the forgetful map 
\[\nerv_n(\SC^{\mc,\oc}(\C,\QF)) \longrightarrow \const\ \nerv_n\Cob^{\mc,\oc}(\C,\QF)\]
is a trivial fibration of semi-simplicial spaces for every fixed $n \in \bbDelta$; see \cite[section A.5]{SAG} for a discussion of Kan and trivial fibrations in this context. This implies that $|\SC^{\mc,\oc}(\C,\QF)| \rightarrow |\Cob^{\mc,\oc}(\C,\QF)|$ is an equivalence, and composing an inverse with the inclusion of the source into $|\Cobccsd 1 {-}^{\mc,\oc}(\C,\QF)|$ provides the desired split.

To give some intuition for the assertion concerning the trivial fibration, let us consider for a moment the analogous statement at the level of components. Any semi-simplicial space induces an equivalence relation of its zero simplices by stipulating $d_0x \sim d_1x$ for all $1$-simplices $x$ (and then letting this generate the relation). In the case at hand the generating relation for two pieces of surgery data $S \rightarrow X$ and $T \rightarrow X$ in $\nerv_n\SC_0^{\mc,\oc}(\C,\QF) \subseteq \nerv_n\Cobsd(\C,\QF)$ is that the map $S \oplus T \rightarrow X$ can be refined to a surgery datum, that restricts to the original two pieces. Such a refinement provides a path between the results of surgery on $S$ and $T$, thereby witnessing that two pieces of surgery data related by disjointness may be interchangeably used to promote an element from $\pi_0|\Cob^{\mc,\oc}(\C,\QF)|$ to one of $\pi_0|\Cob^{\mc,\oc+1}(\C,\QF)|$. 

To understand why it is reasonable to expect surgery data to be unique up to this relation, consider again the geometric context: By the transversality theorem any two sets of embedded spheres below the middle dimension are automatically disjoint up to isotopy (we already recorded the algebraic analogue in \ref{autodisj}). The full assertion on the level of (semi)simplicial spaces then follows the general paradigm that equivalence relations of the kind above are shadows of underlying (semi)simplicial objects, see for example \cite[Sections 6.1.1 \& 6.1.2]{HTT} (the question about semisimplicial versus simplicial objects corresponds to whether the generating relation is already reflexive or not, and is immaterial for the discussion here). In fact, it turns out that there is comparatively little extra work to be undertaken: Essentially since disjointness is a pairwise condition, the forgetful map
\[\nerv_n(\Cobccsd 1 {-}^{\mc,\oc}(\C,\QF)) \longrightarrow \const\ \nerv_n\Cob^{\mc,\oc}(\C,\QF)\]
is $1$-coskeletal for every fixed $n \in \bbDelta$, see Section \ref{coskeletalsection} below, and the definitions are geared for this to persist to the surgery complexes.

From this knowledge, it suffices to check the filling condition for a trivial fibration for $0$- and $1$-simplices to deduce that the forgetful map from the surgery complex is a trivial fibration. Since this low-degree statement already suffices to prove the $\pi_0$-statement of Theorem \ref{thm:main}, we have separated it into its own section \ref{subsec:surgery_on_component}, before putting all ingredients together in \ref{subsec:surgbelow}. \\

This concludes our overview of the argument for surgery on objects below the middle dimension. Let us formalise the ingredients as follows:

\begin{definition}\label{connectedsurgery}
For a category $I$ let $\Cobmfw I^{\mc,\oc}(\C,\QF)$ denote the subcategory of $\Cobmfw I(\C,\QF)$ spanned those objects such that 
\begin{enumerate}
\item the diagram $\lcone I \rightarrow \Cob(\C,\QF)$ resulting from surgery lies in $\Cob^{\mc,\oc}(\C,\QF)$, and
\item its restriction to $I \rightarrow \Cob(\C,\QF)$ even lies in $\Cob^{\mc,\oc+1}(\C,\QF)$.
\end{enumerate}
and similarly morphisms such that
\begin{enumerate}
\item the diagram $[1] \times (\lcone{I}) \rightarrow \Cob(\C,\QF)$ resulting from surgery lies in $\Cob^{\mc,\oc}(\C,\QF)$, and
\item its restriction to $[1] \times I \rightarrow \Cob(\C,\QF)$ even lies in $\Cob^{\mc,\oc+1}(\C,\QF)$.
\end{enumerate}
\end{definition}

\begin{proposition}\label{generalmethodobjects} 
Let $I \colon \Deltainj \rightarrow \Cat$ be a co-semi-simplicial category with $|I_r| \simeq *$ for all $r \in \bbDelta$ and let $\coll \subseteq |\Cob^{\mc,\oc}(\C,\QF)|$ be a collection of path components such that the forgetful map
\[|\Cobmfw{I}^{\mc,\oc}(\C,\QF)|  \longrightarrow |\Cob^{\mc,\oc}(\C,\QF)|\]
admits a section $s$ over $\coll$ for which furthermore the diagram 
\[\xymatrix{|\Cob^{\mc,\oc+1}(\C,\QF)| \times_{|\Cob^{\mc,\oc}(\C,\QF)|} \coll \ar[rd]_{\mathrm{pr}} \ar[rr]^{0 \Rightarrow \mathrm{id}} && |\Cobmfw{I}^{\mc,\oc}(\C,\QF)| \times_{|\Cob^{\mc,\oc}(\C,\QF)|} \coll \\
& \coll \ar[ru]_{(s,\mathrm{id})}}\]
commutes. Then
\[|\Cob^{\mc,\oc+1}(\C,\QF) \times_{|\Cob^{\mc,\oc}(\C,\QF)|} \coll \longrightarrow \coll\]
is an equivalence.
\end{proposition}

\begin{proof}
Consider the following solid commutative diagram (natural in $r \in \Deltainj$)
\[\xymatrix{
|\Cob^{\mc,\oc+1}(\C,\QF)| \ar[dd] & \Fun([1],|\Cob^{\mc,\oc}(\C,\QF)|) \times_{|\Cob^{\mc,\oc}(\C,\QF)|} |\Cob^{\mc,\oc+1}(\C,\QF)| \ar[d]\ar[l]_-{\ev_1} \ar@/_3pc/[ddl]_-{\ev_0} 
\\
& \Fun(\lcone I_r,|\Cob^{\mc,\oc}(\C,\QF)|) \times_{\Fun(I_r,|\Cob^{\mc,\oc}(\C,\QF)|)} \Fun(I_r,|\Cob^{\mc,\oc+1}(\C,\QF)|) \ar@/_1pc/[ld]^-{\ev_*} \ar@/_1pc/@{-->}[u]
\\
|\Cob^{\mc,\oc}(\C,\QF)| \ar@/_2pc/@{-->}[r]_s & |\Cobmfw{I_r}^{\mc,\oc}(\C,\QF)| \ar[u]^{\chifw} \ar[l]_{\fgt}
}\]
where the top right map is restriction along the evident projection $\lcone I \rightarrow [1]$. As indicated with the upper dashed arrow it is an equivalence since $|I_r| \simeq *$ by assumption. After pullback to $\coll$ and taking the realisation in the $r$-direction, the lower map per assumption admits a split as indicated using the lower dashed arrow. We wish to show that taking the long way around the right side of the diagram gives an inverse to its left vertical map. That the composite starting and ending in the lower left corner is the identity is immediate from the commutativity of the diagram using that $s$ is assumed a section of the forgetful map. For the composite starting and ending in the upper left corner this follows immediately from the addititional assumption on $s$.
\end{proof}

\begin{remark}\label{simptechnique}
Let us briefly remark on two details in the comparison of our strategy to that of Galatius and Randal-Williams:
\begin{enumerate}
\item Up to some point-set noise, $1$-coskeletality for fixed $n$ makes 
\[\nerv_n(\Cobccsd 1 {-}(\C,\QF)) \longrightarrow \const \nerv_n(\Cob(\C,\QF))\]
into a topological flag complex as in \cite[Section 6.2]{GRW}. In our set-up the fact that a $1$-coskeletal map is a trivial fibration if it behaves appropriately on $0$- and $1$-simplices replaces the `simplicial technique' that Galatius and Randal-Williams develop to show that their flag complexes give equivalences upon realisation, see \cite[Theorem 6.2]{GRW}. 

This complication in their proof is caused by the $1$-categorical nature of their set-up: The map of semisimplicial \emph{topological spaces} that they consider only becomes a trivial fibration after passing to the underlying semisimplicial object in the $\infty$-category $\Sps$; below the middle dimension any set of embedded spheres with trivial normal bundles is automatically disjoint up to isotopy from any other, but certainly this is not true on the nose. The price to pay is that in the higher categorical approach `disjointness' is a structure on a set of surgery data (namely an isotopy moving the pieces apart) and not a property. See also \cite{GRWer} for the correction of a small gap in \cite{GRW} caused by this subtlety.
\item Through our direct use of the above propositions we skip the step of recombining the multiple surgery data into a single piece as Galatius and Randal-Williams do in \cite[Section 6.1]{GRW}: They show that the geometric analogues of the inclusions $|\SC_0^{\mc,\oc}(\C,\QF)| \longrightarrow |\SC^{\mc,\oc}|$ induce equivalences as well. The analogous argument can be used to in the end produce a section of 
\[|\Cobfw^{\mc,\oc}(\C,\QF)| \longrightarrow |\Cob^{\mc,\oc}(\C,\QF)|,\]
with which we started our sketch above, but we found it rather more simple to directly do multiple surgery. A similar modification should also be possible in the geometric case.
\end{enumerate}
\end{remark}

Finally note, that it is not actually necessary in the above sketch that the surgery complex inside $\Cobmsd{\lcube 1 {r+1}}^{\mc,\oc}(\C,\QF)$ already consists of the data that one eventually performs surgery on (or in fact that it is really a subcategory): It is perfectly possible to interject a map that massages the surgery data in the surgery complex into a form suitable for the second step. Indeed, such an intermediate step is required in the middle dimension, and for reasons we will explain in sections \ref{sec:objects} and \ref{sec:morph} we will, furthermore, have to contend with using surgery data of shape $\lcube 2 {r+1}$ instead of $\lcube 1{r+1}$.

In total then, we will produce under appropriate assumptions three \emph{surgery complexes} $\SC_r^{\mc,\oc}(\C,\QF)$ (for surgery on objects below the middle dimension, \ref{def:surgerycomplexbelow}), $\SSC_r^{\oc}(\C,\QF)$ (for surgery on objects in the middle dimension, \ref{def:surgerycomplexmiddle}) and $\SC_r^{\mc}(\C,\QF)$ (for surgery on morphisms, \ref{def:surgerycomplexmorph}) semisimplicial in $r$, equipped trivial fibrations to the nerves of $\Cob^{\mc,\oc}(\C,\QF)$ and $\Cob^{\mc}(\C,\QF)$, to produce the requisite splits.

\subsection{Forgetting surgery data}\label{coskeletalsection}            
          
Recall that $\Cobccsd a k(\C, \QF)$ denotes the cobordism category with disjoint $[r]$-tuples of surgery data of length $a$ from Definition~\ref{def:cocartsurgdata}. Forgetting surgery data provides us with a map 
\[\nerv_n(\Cobccsd a {k}(\C,\QF)) \longrightarrow \nerv_n(\Cob(\C,\QF)),\]
and our goal in this section is to establish that this map is, for fixed $n\in \bbDelta\op$ and $a\in \NN$, a $1$-coskeletal map of semi-simplicial spaces in $k \in \bbDelta\op_\inj$ (where the target is regarded as a constant semi-simplicial space). This result replaces the `simplicial technique' for topological flag complexes developed by Galatius and Randal-Williams in their study of geometric surgery complexes, compare \ref{simptechnique}.

We remind the reader that a surgery datum in $\Pmsd{\lcube a {r+1}}(\C,\QF\qshift{1})$ defines an element of $\Cobccsd a {r}(\C,\QF)$ if it is strongly cocartesian, i.e.\ if it is left Kan extended from $\axesar$ (the union of the axes through $0$) and takes the value $0$ at the basepoint $(0,\dots, 0)$; we recorded different characterisations of being strongly cocartesian in \ref{lem:strongly_cocartesian_of_sidelength_a}. \\ 

Recall then, that for a complete $\infty$-category $\C$, a map of semi-simplicial objects in $\C$, $X \rightarrow Y$ is called $p$-coskeletal if when regarded as a functor $\Delta_\inj\op \rightarrow \Ar(\C)$ it is right Kan-extended from $\bbDelta_{\leq p, \inj}\op$ relative to the target projection $\Ar(\C) \rightarrow \C$. Using (the dual of) \cite[Propositions 4.3.1.9 \& 4.3.1.10]{HTT} this translates to 
\[X(\Delta^n) \longrightarrow X(\mathrm{sk}_p(\Delta^n)) \times_{Y(\mathrm{sk}_p(\Delta^n))}Y(\Delta^n)\]
being an equivalence for every $n\in \Dinj$; here $X(-)\colon \Fun(\Dinj\op,\Sps)\op \to \C$ denotes the limit preserving extension of $X$, $\Delta^n=\Hom_{\Dinj}(-, [n])$ is the semi-simplicial $n$-simplex, and $\mathrm{sk}_p(\Delta^n)$ is its $p$-skeleton, so that in particular 
\[X(\mathrm{sk}_p(\Delta^n)) = \lim_{[i] \in (\Dinj^{\leq p}/n)\op} X_i.\]
 An easy skeletal induction shows that this is equivalent to
\[ X(\Delta^n) \longrightarrow X(\partial \Delta^n) \times_{Y(\partial \Delta^n)} Y(\Delta^n)\]
being an equivalence for every $n > p$, where $\partial\Delta^k$ is the boundary of the $n$-simplex. Note that for $\C = \Sps$ the space $X(T)$ is just the morphism space from $T$ to $X$ in the category of semi-simplicial spaces.

In this most important special case, a $p$-coskeletal map thus automatically satisfies the filling condition for a trivial fibration, i.e.\ that 
\[\pi_0X(\Delta^n) \longrightarrow \pi_0\left(X(\partial\Delta^n) \times_{Y(\partial\Delta^n)} Y(\Delta^n)\right)\]
is surjective, for all $n>p$. Thus, if a $1$-coskeletal map $f \colon X \rightarrow Y$ of semi-simplicial spaces has both 
\[X_0 \xrightarrow{f_0} Y_0 \quad \text{and} \quad X_1 \xrightarrow{(d_0,d_1,f_1)} X_0 \times X_0 \times_{Y_0 \times Y_0} Y_1\]
surjective on $\pi_0$, then it is a trivial fibration and therefore $|X| \rightarrow |Y|$ an equivalence whenever $Y$ is Kan \cite[Lemma A.5.3.7]{SAG}. We will only have to apply this statement when $Y$ is constant, and constant semi-simplicial spaces are certainly Kan, since horns have contractible realisation. 

\begin{theorem}\label{coskeletality}
The forgetful map
\[\fgt \colon \nerv_n\Cobccsd a {-}(\C,\QF) \longrightarrow \const \nerv_n\Cob(\C,\QF)\]
is a $1$-coskeletal map of semisimplicial spaces for every $n \in \bbDelta$, every $a\in \NN$, and every Poincar\'e category $(\C,\QF)$.
\end{theorem}
          
The remainder of this section is devoted to the proof of this result, in fact we will show a categorical upgrade. To this end denote by
\[\Fun^{\cc}(\powerar, \C)\subseteq \Fun(\powerar, \C)\]
the full subcategory on the strongly cocartesian diagrams and regard it as a hermitian category with hermitian structure restricted from $(\C, \QF)^{\powerar}$, i.e.\ with $\QF^{\powerar}(F) = \QF(F(a,\dots,a))$. For varying $r$, we can view this as a semi-simplicial object in the category $\Cath$ (compare the discussion after \ref{lem:reflection_lemma});  evaluation at $0$ defines a map
\[\ev_0 \colon (\Fun^{\cc}(\powerar, \C),\QF^{\powerar}) \longrightarrow (\C,\QF)\]
of hermitian categories semisimplicial in $r$ (i.e. regarding the target as constant).

\begin{proposition}\label{lem:1-coskeletal}
The map $\ev_0$ just defined is $1$-coskeletal for any Poincar\'e category $(\C,\QF)$.
\end{proposition}

To establish this we in turn need the case $p=2$ of the following standard result on cubical diagrams; the case $p=1$ is spelled out in \cite[6.1.1.15]{HA}, but we do not know of a general reference. For a map $v\colon S\to [1]$, we denote $\vert v\vert:=\sum_{s\in S} v(s)$ its $\ell_1$-norm.

\begin{lemma}\label{lem:p-strongly-cocartesian}
For a cubical diagram $F\colon \lcube 1 r \to \C$, where $\C$ has finite colimits, and $p\in \mathbb N$ the following are equivalent:
\begin{enumerate}
 \item\label{item:p-strongly-cocartesian:1} $F$ is left Kan extended from $\lcube 1 r_{\leq p} = \{v \in \lcube 1 r \mid |v| \leq p\}$.
  \item\label{item:p-strongly-cocartesian:2} Each face of $F$ of dimension $\geq p+1$ that contains the initial vertex is cocartesian.
 \item\label{item:p-strongly-cocartesian:3} Each face of $F$ of dimension $p+1$ is cocartesian.
\end{enumerate}
\end{lemma}

A diagram $F\colon \lcube 1 r \to \C$ satisfying these conditions is called \emph{strongly $p$-cocartesian}, and \emph{strongly $p$-cartesian} if $F\op$ is a strongly $p$-cocartesian cubical diagram in $\C\op$.

\begin{proof}
The proof is by downward induction on $p$ for fixed $r$, where we note that for $p\geq r$ all three conditions hold trivially. Suppose now that the claim holds for some $p+1\geq 1$, and let us show the implications \eqref{item:p-strongly-cocartesian:1}$\Rightarrow$\eqref{item:p-strongly-cocartesian:2}$\Rightarrow$\eqref{item:p-strongly-cocartesian:3}$\Rightarrow$\eqref{item:p-strongly-cocartesian:1} for $p$. 

From \eqref{item:p-strongly-cocartesian:1} we conclude that $F$ is also left Kan extended from $\lcube 1 r_{\leq p+1}$ so all faces of dimension $\geq p+2$ containing the initial vertex are cocartesian by induction and the claim for the faces of dimension $p+1$ holds by the formula for Kan extensions. 

If \eqref{item:p-strongly-cocartesian:2} holds, let $v$ and $v'$ denote the initial and terminal vertex of some $p+1$-dimensional face of $\lcube 1 r$. The proof of \eqref{item:p-strongly-cocartesian:3} (for $p$) is by induction on $|v|$. For $|v|=0$ the claim holds by assumption. For the inductive step, choose a non-zero coordinate $i$ in $v$ and put $\bar v = v - \delta_i$ and $\bar v' = v'-\delta_i$. The restriction of $X$ to the $p+2$-dimensional face with initial vertex $\bar v$ and terminal vertex $v'$ is cocartesian, by assumption of the downwards induction. By assumption of the upwards induction, the further restriction to the codimension $p+1$ face spanned by $\bar v$ and $\bar v'$ is also cocartesian, and it follows that the restriction to the opposite codimension $p+1$ face, which is the original one spanned by $v$ and $v'$ is also cocartesian. 

Suppose finally that \eqref{item:p-strongly-cocartesian:3} holds for $p$. Then \eqref{item:p-strongly-cocartesian:3} also holds for $p+1$, and by inductive hypothesis \eqref{item:p-strongly-cocartesian:1} holds for $p+1$. To show \eqref{item:p-strongly-cocartesian:1} for $p$, we therefore have to show that the restriction $X \colon {\lcube 1 r_{\leq p+1}} \rightarrow \C$ is left Kan induced from its further restriction to $\lcube 1 r_{\leq p}$. But this follows from assumption \eqref{item:p-strongly-cocartesian:3}, using the pointwise formula for Kan extensions.
\end{proof}

In particular, we obtain: 

\begin{corollary}\label{cor:p-excisive-functor}
Let $\QF\colon \C\op\to \Spa$ be a hermitian structure. Then $\QF$ carries strongly cocartesian cubes in $\C$ into strongly $2$-cartesian cubes in spectra.
\end{corollary}

For the proof recall only that hermitian structures are quadratic functors, which by definition means that they carry strongly ($1$-)cocartesian $3$-cubes in $\C\op$ to cartesian cubes in spectra (and by stability a cube in $\C$ is strongly $p$-cartesian if and only if it is strongly $p$-cocartesian).

\begin{proof}[Proof of \ref{lem:1-coskeletal}]
On underlying categories, in fact the stronger statement is true that the map 
\[[r] \longmapsto \big(\mathrm{ev}_0 \colon (\Fun^\cc(\powerar, \C),\QF^{\powerar}) \longrightarrow (\C,\QF)\big)\]
is $0$-coskeletal. In this case the maps from the skeletal criterion become
\[\Fun^\cc(\powerar, \C)\longrightarrow \prod_{i=0}^r \Fun([a], \C) \times_{\prod_{i=0}^r\C} \C.\]
with the structure maps on the right being evaluation at $0 \in [a]$ and the diagonal functor. The map from left to right thus identifies with the one induced by restriction from $\powerar$ to its subposet consisting of the edges emanating from $0$, as this category is the pushout of $r+1$ copies of $[a]$ glued at their initial object. This restriction map is an equivalence by Lemma \ref{lem:strongly_cocartesian_of_sidelength_a}.

By transitivity of Kan extensions, any $0$-coskeletal map is $1$-coskeletal. Thus we find that also 
\[\xymatrix{(\Fun^\cc(\powerar, \C),\QF^{\powerar})\ar[r] \ar[d]^{\mathrm{ev}_0} & \lim_{[i] \in (\Dinj^{\leq 1}/[r])\op} (\Fun([a]^{i+1}, \C),\QF^{[a]^{i+1}}) \ar[d]^{\mathrm{ev}_0} \\
(\C,\QF) \ar[r]& \lim_{[i] \in (\Dinj^{\leq 1}/[r])\op}(\C,\QF)}\]
is cartesian at the level of underlying categories. To conclude that it is cartesian in $\Catp$, we are left to show that for any strongly cocartesian functor $F\colon \powerar\to \C$, the square
\[\xymatrix{\QF^{\powerar}(F) \ar[r] \ar[d]^{\mathrm{ev}_0} &  \lim_{[i] \in (\Dinj^{\leq 1}/[r])\op} \QF^{[a]^{i+1}}(F\vert_{[a]^{i+1}}) \ar[d]^{\mathrm{ev}_0} \\
 \QF(F(0, \dots,0)) \ar[r]& \lim_{[i] \in (\Dinj^{\leq 1}/[r])\op}\QF(F(0, \dots,0))}\] 
is cartesian as well. Now the canonical map
\[\lim_{\rcone{K}} G \longrightarrow \lim_K G \times_{\lim_K G(\ast)} G(\ast)\] 
is an equivalence for any functor $G \colon \rcone{K} \rightarrow \D$ (admitting the right hand side limits) on account of the pushout
\[\xymatrix{K \ar[r]^-{(\mathrm{id},1)} \ar[d] & K \times [1] \ar[d] \\
            \ast \ar[r] &\rcone K}\]
in $\Cat$. Asking the above square to be cartesian is thus equivalent to asking the restriction map
\[\QF^{\powerar}(F) \simeq \lim_{[i]\in \rcone{((\Dinj/[r])\op)}} \QF^{[a]^{i+1}}(F\vert_{[a]^{i+1}}) \longrightarrow \lim_{[i] \in \rcone{((\Dinj^{\leq 1}/[r])\op)}}\QF^{[a]^{i+1}}(F\vert_{[a]^{i+1}})\]
to be an equivalence, where the cone point is taken to $\QF(F(0, \dots,0))$.

Now, recalling that $(\Dinj/[r])\op$ is the poset of non-empty subsets of $[r]$ under reverse inclusion, we see that $\rcone{((\Dinj/[r])\op)}$ identifies with the cubical poset of all subsets of $[r]$ under reverse inclusion, with the cone point corresponding to the empty set. Under this identification, $\rcone{((\Dinj^{\leq 1}/[r])\op)}$ identifies with the subposet spanned by the $2$-dimensional faces emanating from $0$. The claim is thus certainly implied by the $(r+1)$-dimensional cube 
\[\rcone{((\Dinj/[r])\op)} \longrightarrow \Spa, \quad ([i]\rightarrow [r])\longmapsto \QF^{[a]^{i+1}}(F\vert_{[a]^{i+1}})\simeq \QF(F(a\cdot \chi_{[i]}))\]
being strongly $2$-cartesian; here $a\cdot \chi_{[i]}$ denotes the terminal vertex of the subcube of side-length $a$ spanned by the image of $[i]\rightarrow [r]$.

Here, the $(r+1)$-cube $a\cdot \chi_{[-]}$ in $[a]^{r+1}$ agrees with the subcube inclusion of $[1]^{r+1}\to [a]^{r+1}$ that multiplies each entry by $a$. Thus, $F$ takes this cube indeed to a strongly cocartesian cube, and so $\QF$ applied to this cube is strongly $2$-cartesian, by Corollary \ref{cor:p-excisive-functor}. 
\end{proof}

\begin{remark}
There is another beneficial perspective on the proof above, which also brings it closer in spirit to \cite{GRW}: The first part of the argument can simply be read as the translation of the coskeletality condition along the well-known equivalence
\[\Fun(\Dinj\op,\C)/\const \simeq \Fun(\eDinj\op,\C)\]
between the category of (semi-)simplicial object over constant such and the category of augmented (semi-)simplicial objects (i.e. $\eDinj$ is the category of potentially empty totally ordered finite sets and monotone injections).

Since $\rcone{((\Dinj/[r])\op)} \simeq (\eDinj/[r])\op$ the second part then verifies exactly that $[r] \mapsto (\Fun^\cc([a]^{r+1},\C),\QF^{[a]^{r+1}})$ is $1$-coskeletal as an augmented semi-simplicial Poincar\'e category.
\end{remark}

\begin{proof}[Proof of \ref{coskeletality}]
It suffices to show that the map 
\[\Pccsd a {-}(\C, \QF)\longrightarrow \const \Poinc(\C, \QF)\]
is a $1$-coskeletal map of semi-simplicial spaces for each Poincar\'e category $(\C, \QF)$; plugging in $\Q_n(\C, \QF\qshift{1})$ gives the claim from the statement.

Now, using the description of $\Pmsd{\lcube a {r+1}}(\C,\QF)$ from \ref{rem:alternativedefsofmultsurg} together with the equivalence
\[\Hom_{\Cat}(\lcube a {r+1}, \catforms(\C, \QF)) \simeq \spcforms((\C, \QF)^{\lcube a {r+1}})\]
from \cite[6.3.15]{9authI} we obtain an equivalence
\[\spcforms(\Met (\C, \QF)^{\lcube a {r+1}})\times_{\spcforms((\C, \QF)^{\lcube a {r+1}})} \Poinc(\C, \QF) \simeq \Pmsd{\lcube a {r+1}}(\C, \QF)\]
Furthermore, restricting to the collection of path-components on the right consisting of strongly cocartesian surgery data corresponds (by direct inspection) to
 \[\spcforms(\Fun^\cc(\lcube a{r+1}, \Met(\C, \QF)))\times_{\spcforms(\Fun^\cc(\lcube a{r+1}, (\C, \QF)))} \Poinc(\C, \QF)\]
 on the left. Regarding this as a semi-simplicial space in $r$, evaluation at $0$ defines a map to the constant semi-simplicial space on \[\spcforms\Met(\C,\QF) \times_{\spcforms(\C,\QF)} \Poinc(\C,\QF) = \Psd(\C,\QF).\]
By \ref{lem:1-coskeletal} this map is $1$-coskeletal, since $\spcforms \colon \Cath \rightarrow \Sps$ preserves limits by \cite[4.1.3]{9authI} (or direct inspection). Finally, the map in the statement arises from it by base change along the functor
\[\Poinc(\C,\QF) \longrightarrow \Psd(\C,\QF)\]
inserting the surgery datum $0$; formally it is induced by the identity of $\Poinc(\C,\QF)$ and the map $\Poinc(\C,\QF) \rightarrow \spcforms\Met(\C,\QF)$ arising from the the hermitian (but usually non-Poincar\'e) functor 
\[(\C,\QF\qshift{1}) \longrightarrow \Met(\C,\QF\qshift{1}), \quad X \mapsto (0 \rightarrow X).\]
Since coskeletality is clearly preserved under base change we are done.
\end{proof}

\section{Surgery on objects}\label{sec:objects}

Throughout this section we consider a Poincar\'e category $(\C, \QF)$ of dimension at least $d$. Our goal is to prove that under suitable assumptions the map
\[|\Cob^{\mc,\oc+1}(\C,\QF)| \longrightarrow |\Cob^{\mc,\oc}(\C,\QF)|\]
is an equivalence onto the path components that it hits. This is split into  two cases, surgery below the middle dimension in Sections \ref{subsec:surgery_on_component} and \ref{subsec:surgbelow}, where inter alia $\oc$ is assumed to be less than half the dimension $d$ of $(\C,\QF)$, and surgery in the middle dimension, where $\mc=\oc=d/2$ in Sections \ref{sec:middle_dimension} and \ref{subsec:surgin2}. 

We start out in \ref{subsec:surgery_on_component} by analysing when the map
\[\pi_0|\Cob^{\mc,\oc+1}(\C,\QF)| \longrightarrow \pi_0|\Cob^{\mc,\oc}(\C,\QF)|\]
is an isomorphism (assuming as we must that $2\oc < d$). Together our preparations in the previous chapter, the proof then upgrades to the space level in \ref{subsec:surgbelow} with little extra effort. That one can nevertheless see all the moving parts of our arguments come together is the reason we treat this first (despite it being the middle of three steps in the proof Theorem \ref{thm:main}). The final two sections then explain the necessary modifications in the middle dimensions.

\subsection{Surgery on components}\label{subsec:surgery_on_component}

We will begin by explaining the precise pieces of surgery data that we use in the surgery move for objects below the middle dimension, and prove their existence and uniqueness up to the disjointness relation. We use this first to derive a $\pi_0$-version of surgery below the middle dimension in Theorem~\ref{resolutionobjectspi0}, in the hope that it will serve to illustrate the overall argument. Using the identification $\pi_0|\Cob(\C,\QF)| \cong \L_0(\C,\QF\qshift{1})$ this version is essentially equivalent to the analysis of $\L$-groups in \cite[Section 1.2]{9authIII} in turn improving on results of Ranicki from \cite{rsurg2}, see Section \ref{lsection2} below for details. The reader only interested in $\L$-groups may prefer the exposition in \cite{9authIII} (which works almost verbatim in the generality of Poincar\'e categories with weight structure): Our presentation is geared for reuse in the next sections and not streamlined for the $\pi_0$-level (parametrised surgery really only appears at the level of morphisms). After unwinding notation, the actual arguments here and in \cite{9authIII} are of very similar complexity. 

To have a complete treatment of Theorem  \ref{thm:main} at the $\pi_0$-level, we also briefly discuss surgery on morphisms, which is very simple at this level, the essential difficulty only becoming visible at the level of fundamental groups. The third step in the proof of Theorem \ref{thm:main}, \emph{i.e.}\ surgery on objects in the middle dimension, is entirely vacuous at the level of components. 

 \medskip

We start by specifying the surgery data we will actually perform surgery on:

\begin{definition}\label{def:suitablebelowmid}
For $K\in \Cat$, we shall say that a $K$-fold surgery datum in $\nerv_n\Cobmsd{K}(\C, \QF)=\Pmsd{K}(\Q_n(\C, \QF\qshift 1))$, with underlying diagram $T\to X$, where $X\in \nerv_n\Cob^{\mc,\oc}(\C, \QF)$, is \emph{suitable for surgery of type $(\mc,\oc)$} if 
\begin{enumerate}
\item it is a forwards surgery datum, \emph{i.e.,} for each $k\in K$, the morphism $T(k)_{j\leq j+1}\to T(k)_{j\leq j}$ is an equivalence, 
\item for each $k\in K$, the object $T(k)_{j \leq j}$ is concentrated in degree $\oc$,
\item for each morphism $k\to k'$ in $K$, the induced morphism $T(k)_{j \leq j} \rightarrow T(k')_{j \leq j}$ is $(d-\mc)$-coconnective, 
\item also for each $k\in K$, the morphism $T(k)_{j \leq j+1} \rightarrow T(k)_{j+1 \leq j+1}$ is $d-\mc$-coconnective, and
\item for each $k\in K$, the result of surgery of $T(k)_{j \leq j} \rightarrow X_{j \leq j}$ is $\oc+1$-connective.
\end{enumerate}
\end{definition}

Note immediately that these conditions define a complete Segal space whose associated category is therefore a subcategory of $\Cobmfw K (\C, \QF)$. We first show that we can indeed use this kind of data as an input for our strategy, i.e.\ that it lies in the subcategory $\Cobmfw K^{\mc,\oc}(\C, \QF)\subset \Cobmsd K(\C, \QF)$ from Definition \ref{connectedsurgery}.

\begin{lemma}\label{lem:suitable_is_good}
Let $(\C, \QF)$ be a Poincar\'e category of dimension at least $d$. If $\oc+\mc\leq d$, then surgery data which are suitable for surgery of type $(\mc,\oc)$ are contained in (the nerve of) $\Cobmfw K^{\mc,\oc}(\C, \QF)$.
\end{lemma}

\begin{proof}
For a forwards surgery datum in $\nerv_n \Cobmfw K(\C, \QF)$ which is suitable for surgery of type $(\mc,\oc)$, denote by $F\colon [n]\times{} \lcone{K}\to \Cob(\C, \QF)$ the functor resulting from surgery. We only need to show that $F$ takes values in $\Cob^\mc(\C, \QF)$, the connectivity of objects being explicitly assumed. Thus take a morphism $\alpha$ in $[n]\times{} \lcone K$: If $\alpha=(\id_j, k'\to k)$, then the morphism $F(\alpha)$  has surgery datum  $T(k)_{j\leq j}$, if $k'=*$ is the cone point, and else $T(k)_{j\leq j}/T(k')_{j\leq j}$, both of which are $d-\mc$-coconnective, so $F(\alpha)$ is $m$-connective by \ref{lem:basic_surgery_estimate}. If $\alpha=(j\leq j+1, \id_*)$ then the connectivity condition on $F(\alpha)$ is assumed and if $\alpha=(j \leq j+1, \id_k)$ for some $k\neq *$, then $F(\alpha)$ is obtained by surgery along the datum $T(k)_{j\leq j}\leftarrow T(k)_{j\leq j+1} \to T(k)_{j+1\leq j+1}$, and so is $\mc$-connected by \ref{connsurgmor}. The lemma follows because every morphism of $[n]\times{} \lcone K$ is a composite of such morphisms.
\end{proof}

We next discuss the existence of surgery data. 

\begin{proposition}\label{existenceobjectbelowmiddle}
Let $(\C,\QF)$ be a Poincar\'e category of dimension at least $d$, and $(X,q) \in \nerv_n\Cob^{\mc,\oc}(\C,\QF) \subseteq \Poinc\Q_n(\C,\QF\qshift{1})$ for some $n \in \mathbb N$. Then $(X,q)$ extends to a surgery datum in $\Psd(\Q_n(\C, \QF\qshift 1))$, which is suitable for surgery of type $(\mc,\oc)$, provided that
\begin{enumerate}
\item $2\oc<d$, 
\item $\oc+\mc\leq d$,
\item $\oc\leq \mc$, and
\item $\Lin_\QF(Y)$ is $\oc$-connective for every $Y \in \C^\heart$.
\end{enumerate} 
\end{proposition}

As preparation we have:

\begin{observation}\label{extendingsurgdata}
Let $(\C,\QF)$ be a Poincar\'e category with a weight structure, and consider a surgery datum in $(\C, \QF)$ with underlying diagram of the form $T\to C$, with $T \in \C_{[\oc,\oc]}$. Then:
\begin{enumerate}
 \item[i)] This surgery datum extends to a surgery datum of the form $0\leftarrow 0 \to T$ on any Poincar\'e cobordism of the form $D\leftarrow W\to C$.
 \item[ii)] This surgery datum extends to a surgery datum of the form $T\xleftarrow \id T \xrightarrow \id T$ on any Poincar\'e cobordism of the form $C\leftarrow W \to D$, provided $W\to C$ is $\oc$-connective.
\end{enumerate}
\end{observation}

\begin{proof}
In both cases the underlying diagram clearly extends as described, trivially in the first and by \ref{connlift} in the second case, and since both evaluation maps
\[\lim\bigl(\QF\qshift 1(0) \to \QF\qshift 1(0) \leftarrow \QF\qshift 1(T) \bigr)\to \QF\qshift 1(T)
 \quad\mathrm{and}\quad
 \lim\bigl(\QF\qshift 1(T) \to \QF\qshift 1(T) \leftarrow \QF\qshift 1(T)\bigr)\to \QF\qshift 1(T)
\]
are equivalences, also the nullhomotopy of $q$ does.
\end{proof}

\begin{proof}[Proof of \ref{existenceobjectbelowmiddle}]
We first observe that for $(C,r) \in \Poinc(\C,\QF\qshift{1})$ and $S \in \C_{[\oc,\oc]}$ any $\oc$-connective map $f \colon S \rightarrow C$ refines to a surgery datum: Because of the fibre sequence
\[\Bil_\QF(S,S)_\hC \to \QF(S) \to \Lin_\QF(S)\]
it suffices to have $\Bil_\QF(S,S)$ and $\Lin_\QF(S)$ connective; for the bilinear term this follows from the first assumption, and for the linear term, from the last. 
The result of surgery along such a datum is then $\oc+1$-connected by Lemma \ref{lem:surgery_below_the_middle_dimension} (it is only here that we use that $2\oc$ is strictly less than $d$).

Starting with the proof proper, we can thus choose a surgery datum $T^{(0)}_{0\leq0}\to X_{0\leq 0}$ which is concentrated in degree $\oc$ and whose result of surgery is $\oc+1$-connective. By our observation above (using $\oc\leq\mc$) and the Segal condition for $\nerv \Cobsd(\C, \QF)$, it extends iteratively to a surgery datum $T^{(0)}\to X$, where $T^{(0)}$ is a constant diagram. By \ref{lem:surgery_below_the_middle_dimension} (applied to $\oc-1$) and \ref{connsurgmor}, the result of surgery along $T^{(0)}\to X$, let us call it $X^{(1)}$, again lies in $\nerv_n(\Cob^{\mc,\oc}(\C,\QF))$, and per construction $X^{(1)}_{0\leq 0}$ is even $\oc+1$-connective. 

We can thus choose another surgery datum $T^{(1)}_{1\leq 1}\to X^{(1)}_{1\leq 1}$ which is concentrated in degree $\oc$ and whose result of surgery is $\oc+1$-connective. Again by the observation above, it extends to a surgery datum $T^{(1)}\to X^{(1)}$ where $T^{(1)}$ restricts along the first face to $0\leftarrow 0\to T^{(1)}_{1\leq 1}$ and is constant on the other ($1$-)morphisms. Again, \ref{lem:surgery_below_the_middle_dimension} and \ref{connsurgmor} (using $\mc+\oc\leq d$) show that the result of surgery along $T^{(1)}\to X^{(1)}$, call it $X^{(2)}$ again lies in $\nerv_n(\Cob^{\mc,\oc}(\C,\QF))$.

Continuing this process, we obtain a chain of cobordisms in the Poincar\'e category $\Q_n(\C, \QF\qshift 1)$
\[X \rightsquigarrow X^{(1)} \rightsquigarrow \dots \rightsquigarrow X^{(n)}\]
where $X^{(j)}_{i\leq i}$ is $\oc+1$-connective for all $i \leq j$. We claim that the surgery datum $T\to X$ for the composite cobordism $X \rightsquigarrow X^{(n)}$ is suitable for surgery of type $(\mc,\oc)$.

In fact, the last condition holds by construction. For the other conditions, we note that they are met by each of the individual $T^{(k)}$'s; then it suffices to note that the surgery datum of a composite cobordism is an extension of the surgery data of its factors, and that the conditions are closed under taking extensions.
\end{proof}

We then formulate the automatic disjointness of such surgery data:

\begin{lemma}\label{lemmaautmaticdisjoint}
Let $(\C,\QF)$ be a Poincar\'e category of dimension at least $d$, and $f \colon S \rightarrow X$ and $g \colon T \rightarrow X$ two pieces of surgery data in $N_n\Cobsd(\C,\QF)$ suitable for surgery of type $(\mc,\oc)$ on the same Poincar\'e object. If $2\oc\leq d$ they can be extended to a $[1] \cup_{\{1\}} [1]$-fold surgery datum of the form
\[\xymatrix{ S \ar[r] \ar[rd]_f & S \oplus T \ar[d]^{f+g} & T \ar[l] \ar[ld]^g \\
  & X &}\]
whose restriction to the corners are the two original pieces. If $2\oc<d$ and $\oc+\mc\leq d$, then any such extension is again suitable for surgery of type $(\mc,\oc)$.
\end{lemma}

\begin{proof}
Note that since $S$ is a forwards surgery datum we have
\[\QF\qshift 1_n(S) = \lim_{\Twar[n]\op} \QF\qshift 1\circ S\op \simeq \QF\qshift 1(S_{n\leq n})
\]
as the target map $\Twar[n]\op \rightarrow [n]\op$ is a (Bousfield) localisation and thus final. The same applies to $T$ and $S\oplus T$, so by \ref{autodisj} we can find an extension of the given pieces of surgery data to $f + g \colon S \oplus T \rightarrow X$. 

It remains to check that the resulting surgery datum is again suitable for surgery of type $(\mc,\oc)$ under the stronger assumptions: The first four conditions are immediate from the construction and the assumption $\oc+\mc\leq d$, and the last one follows from \ref{lem:surgery_below_the_middle_dimension}. 
\end{proof}

Armed with this information we show:

\begin{theorem}\label{resolutionobjectspi0}
Let $(\C,\QF)$ be a Poincar\'e category of dimension at least $d$ and assume
\begin{enumerate}
\item $2\oc < d$,
\item $\oc+\mc\leq d$,
\item $\oc\leq \mc$, and
\item $\Lin_\QF(Y)$ is $\oc$-connective for every $Y \in \C^\heart$.
\end{enumerate}
Then
\[\pi_0|\Cob^{\mc,\oc+1}(\C,\QF)| \longrightarrow \pi_0|\Cob^{\mc,\oc}(\C,\QF)|\]
is an isomorphism.
\end{theorem}

\begin{proof}
As explained in Section \ref{sketch} the idea is to construct an inverse map $s$ using surgery data suitable for surgery of type $(\mc,\oc)$. So consider an object $X \in \Cob^{\mc,\oc}(\C,\QF)$. By \ref{existenceobjectbelowmiddle} (for $n=0$) we can pick such a surgery datum and can try to define $s([X])$ to be component of the result of surgery. Applying \ref{lemmaautmaticdisjoint} and \ref{lem:suitable_is_good} (with $n=0$) this component is independent of the choice of surgery datum. To see that it also only depends on the component of $X$, apply \ref{existenceobjectbelowmiddle} and \ref{lem:suitable_is_good} (for $n=1$) to any cobordism involving $X$. This establishes a map
\[s \colon \pi_0|\Cob^{\mc,\oc}(\C,\QF)| \longrightarrow \pi_0|\Cob^{\mc,\oc+1}(\C,\QF)|.\]
Using the surgery datum $0$ we immediately find that $s$ is right inverse to the map induced by inclusion, and that it is left inverse is also immediate from \ref{lem:suitable_is_good}.
\end{proof}

Let us end by briefly digressing to the case of surgery on morphisms. The aim is to prove that the map 
\[\pi_0|\Cob^{\mc+1}(\C,\QF)| \longrightarrow \pi_0|\Cob^{\mc}(\C,\QF)|\]
is an equivalence under appropriate assumptions. Note that as the two categories share objects, it is clearly surjective, so we are left to prove injectivity. Supposing that $C\leftarrow W \to D$ is a morphism in $\Cob^{\mc}(\C, \QF)$, i.e. that the left pointing arrow has a $\mc$-connective fibre, we have to show that $C$ and $D$ can be connected through a zig-zag of morphisms in $\Cob^{\mc+1}(\C,\QF)$ as well. Choose then an $\mc$-connective map 
\[T\to \fib(W\to C)\]
with $T$ concentrated in degree $\mc$. This map and the canonical nullhomotopy of the composite to $C$ uniquely extend to a commutative diagram
\[ \xymatrix{
 0 \ar[d] & T \ar[l] \ar[r]^{\mathrm{id}} \ar[d] & T \ar[d]^f\\
 C & W \ar[l] \ar[r] & D
} 
\]
The hermitian structure of $\Q_1(\C, \QF\qshift 1)$ evaluates to $0$ on the top row, so that the diagram canonically refines to a surgery datum in $\Q_1(\C, \QF\qshift 1)$. 
The result of surgery is given by 
\[  C \leftarrow W/T \to D_f
\]
so the fibre of its left pointing map, $\fib(W\rightarrow C)/T$, has become $\mc+1$-connective as desired. By \ref{corhomotopy} we also obtain a commutative diagram
\[\xymatrix{C \ar@{~>}[r]^{W/T} \ar@{~>}[d]_{\mathrm{id}}& D_f \ar@{~>}[d]^{\chi(f)} \\
            C \ar@{~>}[r]^W & D.}\]
witnessing that $D_f$ and $D$ lie in the same component of $\pi_0|\Cob^{\mc}(\C,\QF)|$. The cobordism $\chi(f)$, regarded in the direction indicated, does not lie in $\Cob^{\mc+1}(\C,\QF)$, however, (unless $T=0$) so that it may seem like we have just moved the problem one step further. But the reflection of $\chi(f)$ gives a cobordism in $\Cob^{\mc+1}(\C,\QF)$, whenever $(\C,\QF)$ has dimension at least $d$ and $d \geq 2\mc+1$: The fibre of $\chi(f) \rightarrow D$ is $\Dual_\QF(T) \in \C_{[d-\mc,\infty]}$. Thus we have found the desired zig-zag and proven:

\begin{proposition}\label{resolutionmorphismspi0}
Let $(\C,\QF)$ be a Poincar\'e category of dimension at least $d$ and assume $d \geq 2\mc+1$.Then
\[\pi_0|\Cob^{\mc+1}(\C,\QF)| \longrightarrow \pi_0|\Cob^{\mc}(\C,\QF)|\]
is an isomorphism.
\end{proposition}

We will deduce the consequences of \ref{resolutionobjectspi0} and \ref{resolutionmorphismspi0} for $\L$-groups in Section \ref{lsection2}.

\begin{remark}\label{surgerydatumformorphs}
This simple trick does not work once we consider the entirety of $|\Cob^{\mc}(\C,\QF)|$ and not just its components: We are not allowed to reflect $\chi(f)$, without breaking the commutativity of the diagram above, as the surgery datum $0 \leftarrow T \rightarrow T$ is not a forwards datum, compare Section \ref{sec:reflection_and_doubling}. We will therefore have to construct a rather more elaborate path between $C$ and $D$, see \ref{subsec:surgcompxmorph}.
\end{remark}

\subsection{Surgery below the middle dimension}\label{subsec:surgbelow}

We now upgrade Theorem \ref{resolutioninL} to the level of spaces by executing the strategy outlined in Section \ref{sketch}. The following is the relevant surgery complex:

\begin{definition}\label{def:surgerycomplexbelow}
Given a Poincar\'e category $(\C,\QF)$ of dimension at least $d$, the simplicial space of \emph{$[r]$-tuples of surgery data suitable for type $(\mc,\oc)$} on $\Cob^{\mc,\oc}(\C, \QF)$ is given in degree $n$ by the full subspace of 
\[\Pccsd 1 {r}(\Q_n(\C, \QF\qshift{1})) \simeq \nerv_n\Cobccsd 1{r}(\C,\QF)\]
 on those disjoint $[r]$-tuples of surgery data whose underlying diagram $T\to X$ satisfies:
\begin{enumerate}
\item $X$ lies in $\nerv_n\Cob^{\mc,\oc}(\C,\QF)$
\item $T$ is suitable for surgery of type $(\mc,\oc)$ at each vertex of $\lcube 1 {r+1}$ adjacent to $0$
\end{enumerate}
We denote by $\SC^{\mc,\oc}_r(\C, \QF)$ the category associated to this simplicial space. For varying $r$ this assembles into a semi-simplicial category
\[\SC^{\mc,\oc}(\C, \QF)\colon \Deltainj\op\to \Cat\]
which we call the \emph{surgery complex} for object surgery of type $(\mc, \oc)$.
\end{definition}

Recall that disjointness refers to the conditions that $T(0) \simeq 0$ and that the cube $T \colon \lcube 1 {r+1} \rightarrow \C$ is strongly cocartesian, see Definition \ref{def:cocartsurgdata}.

It is easily checked that the first part really defines a simplicial subspace of $\Pccsd 1 r\Q(\C, \QF\qshift{1})$ from Definition \ref{def:cocartsurgdata} (compare the discussion after the reflection lemma, \ref{lem:reflection_lemma}, for the definition of the semi-simplicial structure). Furthermore, it is itself a complete Segal space, thus making the second part sensible, and identifying $\SC^{\mc,\oc}_r(\C,\QF)$ as a subcategory of $\Cobccsd{1}{r}(\C,\QF)$.

Setting $\lcube 1 {r+1}_0 = \lcube 1 {r+1} \setminus \{0\}$ we see just as in \ref{lemmaautmaticdisjoint} that if $2\oc<d$ and $\oc+\mc<d$, then the whole $\lcube 1 {r+1}_0$-fold surgery datum of an element in the surgery complex is suitable for surgery of type $(\mc, \oc)$. Therefore, by the paragraph after Definition \ref{def:suitablebelowmid} (and forgetting the zero surgery datum), we have a  canonical map
\[\SC^{\mc,\oc}_r(\C,\QF) \to \Cobmfw{\lcube 1 {r+1}_0}^{\mc,\oc}(\C,\QF) .\]
Thus, we can use the surgery data from the surgery complex as an input for our surgery strategy.

\begin{proposition}\label{surgery below the middle dim}
Let $(\C,\QF)$ be a Poincar\'e category of dimension at least $d$. Then the map
\[\nerv_n\SC^{\mc,\oc}(\C,\QF) \longrightarrow \const \nerv_n \Cob^{\mc,\oc}(\C,\QF)\]
forgetting the surgery data is a trivial fibration of semi-simplicial spaces for every fixed $n \in \bbDelta$, provided
\begin{enumerate}
\item $2\oc <d$,
\item $\oc+\mc\leq d$,
\item $\oc\leq \mc$, and
\item $\Lin_\QF(Y)$ is $\oc$-connective for every $Y \in \C^\heart$.
\end{enumerate}
In particular, the forgetful map induces an equivalence
\[|\SC^{\mc,\oc}(\C,\QF)| \longrightarrow |\Cob^{\mc,\oc}(\C,\QF)|\]
in this case.
\end{proposition}

\begin{proof}
The second statement is implied by the first, see \cite[A.5.3.1]{SAG}. For the first statement, the first two of the filling conditions for a trivial fibration are handled by \ref{existenceobjectbelowmiddle} and \ref{lemmaautmaticdisjoint}. To see the higher filling conditions we show that the map in question is $1$-coskeletal: This follows from \ref{coskeletality} by considering the cartesian squares, 
\[\xymatrix{
P \ar[d] \ar[r] 
& \Cobccsd 1 {r}(\C,\QF) \ar[d] 
&& 
\SC^{\mc,\oc}_r(\C,\QF) \ar[r] \ar[d] 
& \prod_{[r]} \SC^{\mc,\oc}_0(\C,\QF) \ar[d]
\\
\Cob^{\mc,\oc}(\C,\QF) \ar[r] 
& \Cob(\C,\QF) 
&& 
P \ar[r] 
& \prod_{[r]} \Cobsd^{\mc,\oc}(\C,\QF)
}\]
which after applying $\nerv_n$ display the map in question as the composite of a $1$-coskeletal and a $0$-coskeletal map.
\end{proof}

We thus obtain:

\begin{theorem}
\label{thm:surgbelowmiddlecomplete}
Let $(\C,\QF)$ be a Poincar\'e category of dimension at least $d$ and assume
\begin{enumerate}
\item $2\oc < d$,
\item $\oc+\mc\leq d$,
\item $\oc\leq \mc$, and
\item $\Lin_\QF(Y)$ is $\oc$-connective for every $Y \in \C^\heart$.
\end{enumerate}
Then the canonical map
\[|\Cob^{\mc,\oc+1}(\C,\QF)| \longrightarrow |\Cob^{\mc,\oc}(\C,\QF)|\]
is an equivalence.
\end{theorem}

\begin{proof}
Composing the map $\SC^{\mc,\oc}_r(\C,\QF) \to \Cobmfw{\lcube 1 {r+1}_0}^{\mc,\oc}(\C,\QF)$ with the inverse of the equivalence of \ref{surgery below the middle dim} provides a section $s$ to the forgetful map
\[|[r] \longmapsto |\Cobmfw{\lcube 1 {r+1}_0}^{\mc,\oc}(\C,\QF)|| \rightarrow |\Cob^{\mc,\oc}(\C,\QF)|.\]
Furthermore, including the surgery datum $0$ provides a commutative diagram
\[\xymatrix{\const \Cob^{\mc,\oc+1}(\C,\QF) \ar[dr] \ar[rr]^{0 \Rightarrow \id} && \SC^{\mc,\oc}(\C,\QF) \ar[dl]^\fgt \\
& \const \Cob^{\mc,\oc}(\C,\QF), &}\]
which witnesses that $s$ satisfies the assumptions of \ref{generalmethodobjects} (with $\coll = |\Cob^{\mc,\oc}(\C,\QF)|$ and $I_r = \lcube 1 {r+1}_0$), the conclusion of which then proves the result.
\end{proof}

\subsection{Surgery data in the middle dimension}\label{sec:middle_dimension}

In the present and the following section, we prove a version of Theorem \ref{thm:surgbelowmiddlecomplete} for $\mc=\oc$ and a Poincar\'e category of exact dimension $d=2\oc$, which will be the final step in our proof of Theorem \ref{thm:main}. The main difference, compared to the previous situation, is that a cocartesian diagram of surgery data, whose vertices adjacent to zero are suitable for surgery of type $(\oc,\oc)$, need not be suitable for surgery of type $(\oc,\oc)$ away from zero. This is an instance of the general phenomenon that, in the middle dimension, the connectivity of the result of surgery cannot be controlled by connectivity estimates on the surgery datum. 

For the same reason, it is not always possible to achieve $\oc+1$-connectivity by surgery, even if $\Lin_\QF=0$: Note that $0$ is the only $\oc+1$-connected Poincar\'e object in $\Cob^{\oc,\oc}(\C, \QF)$ and therefore, increasing connectivity by surgery can be achieved precisely on nullcobordant Poincar\'e objects.

In the present section we show that surgery data, suitable for surgery of type $(\oc,\oc)$, exist at least in the base point component 
\[\Cob^{\oc,\oc}_0(\C, \QF)\subseteq \Cob^{\oc,\oc}(\C, \QF)\]
of $\Cob^{\oc,\oc}(\C, \QF)$, i.e.\ the wide subcategory spanned by all $(C, q)\in \Poinc(\C, \QF\qshift 1)$ with $[C]=0\in \pi_0\vert \Cob^{\oc,\oc}(\C, \QF)\vert$.

We start out by discussing the basic existence result for surgery data in the middle dimension, compare \cite[Lemma 6.21]{GRW}.

\begin{lemma}\label{lem:surgery_data_for_objects_middle_dimension_n_is_0}
Let $(\C,\QF)$ be a Poincar\'e category of exact dimension $2\oc$ and assume $\Lin_\QF(X)$ is $\oc$-connective for every $X \in \C^\heart$. Then for two objects $X, Y \in \Cob^{\oc,\oc}(\C,\QF)$ the statements
\begin{enumerate}
\item $[X] = [Y]\in \pi_0|\Cob^{\oc,\oc}(\C,\QF)|$, and
\item there is a surgery datum $f \colon T \rightarrow X$ with $T \in \C_{[\oc,\oc]}$, such that $X_f \simeq Y$
\end{enumerate}
are equivalent.
\end{lemma}

\begin{proof}
Assume the second statement. Then the opposite of the cobordism
\[X \longrightarrow \chi(f) \longrightarrow X_f\]
lies in $\Cob^{\oc,\oc}(\C,\QF)$, since the right hand map has fibre $T$.
For the more interesting converse we first show that the existence of a morphism from some $X$ to some $Y$ in $\Cob^{\oc,\oc}(\C,\QF)$ is a symmetric relation: For consider a cobordism
\[X \longleftarrow W \longrightarrow Y\]
whose left pointing map is $\oc$-connective; we claim that, after a suitable surgery on a datum of the form
\[\xymatrix{0 \ar[d]&\ar[r]\ar[l]\ar[d] T & 0\ar[d] \\
            X & W \ar[r]\ar[l]& Y}\]
(which leaves source and target unchanged) we obtain even both maps $\oc$-connective. To produce the surgery datum, choose an $\oc-1$-connective map $T \rightarrow \fib(W \rightarrow Y)$ with $T \in \C_{[-\infty,\oc-1]}$. Now from the assumptions we find $W \in \C_{[\oc,\infty]}$, so $\fib(W \rightarrow Y) \in \C_{[\oc-1,\infty]}$ and thus by \ref{connfib} $T \in \C_{[\oc-1,\oc-1]}$.

Per construction we obtain the right hand square and since $\pi_0\Hom_\C(T,X) = 0$ by connectivity, we can extend to the full diagram displayed above. 
Next we have to make sure that the form of the cobordism pulls back to
\[0 \in \pi_0\QF\qshift{1}_1(0 \leftarrow T \rightarrow 0) = \pi_0\QF(T).\]
From the fibre sequence
\[\Bil_\QF(T,T)_\hC \longrightarrow \QF(T) \longrightarrow \Lin_\QF(T)\]
we find that it suffices to show that the outer two terms are connected. For the left hand side this follows from $(\C,\QF)$ having dimension $2\oc$ (dimension at least $2\oc-1$ would even suffice), and for the right it is explicitly assumed. Finally, we analyse the result of the surgery
\[X \longleftarrow \fib(W \rightarrow \Dual_\QF(T))/T \longrightarrow Y.\]
To obtain the claim for the left hand map invoke \ref{connsurgmor} together with the dimension of $(\C,\QF)$ being (at least) $2\oc$. 
The right map, on the other hand, can be written as a composite
\[\fib(W\to \Dual_\QF(T))/T \to W/T \to Y\]
whose second map is $\oc$-connective by construction, so we are left to show that the first map is $\oc$-connective as well. But its fibre is $\Dual_\QF(T)^{[-1]}$ which is indeed $\oc$-connective because $T$ is $\oc-1$-coconnective and $(\C, \QF)$ is of dimension (at least) $2\oc$.

This concludes the proof that existence of morphisms in $\Cob^{\mc,\mc}(\C, \QF)$ is a symmetric relation on the objects. Consider now $X,Y \in \Cob^{\oc,\oc}(\C,\QF)$ in the same component. Then by the previous discussion we can find a cobordism $X \leftarrow W \rightarrow Y$ with both arrows $\oc$-connective. 
We claim, that $\fib(W \rightarrow Y) \in \C_{[\oc,\oc]}$ whence it can be used as the desired surgery datum. Per assumption $\fib(W \rightarrow Y) \in \C_{[\oc,\infty]}$ and 
\[\fib(W \rightarrow Y) \simeq \Dual_\QF(\fib(W \rightarrow X))\]
whose right hand term lies in $\C_{[-\infty,\oc]}$, since $(\C,\QF)$ has dimension exactly $2\oc$.
\end{proof}

Specialising \ref{def:suitablebelowmid} to the case where $(\C,\QF)$ has exact dimension $2\oc$, a surgery datum $T \rightarrow X$ for $X \in \nerv_n\Cob^{\oc,\oc}(\C,\QF)$ is suitable for surgery of type $(\oc,\oc)$ if
\begin{enumerate}
\item $T$ is a forwards surgery datum, \emph{i.e.}, each $T(k)_{j\leq j+1}\to T(k)_{j\leq j}$ is an equivalence,
\item each $T(k)_{j \leq j}$ is concentrated in degree $\oc$,
\item each $T(k')_{j \leq j}/T(k)_{j \leq j}$ is concentrated in degree $\oc$, 
\item also each $T(k)_{j+1 \leq j+1}/T(k)_{j \leq j+1}$ is concentrated in degree $\oc$, 
\item the result of surgery on each $T(k)_{j \leq j} \rightarrow X_{j \leq j}$ vanishes.
\end{enumerate}

We obtain:

\begin{proposition}\label{prop:surgery_data_for_objects_exist_middle_dimension}
Let $(\C,\QF)$ be a Poincar\'e category of exact dimension $2\oc$ and assume that $\Lin_\QF(Y)$ is $\oc$-connective for all $Y \in \C^\heart$. Then each element in $\nerv_n\Cob_0^{\oc,\oc}(\C,\QF)$ refines to a surgery datum suitable for surgery of type $(\oc,\oc)$.
\end{proposition}

\begin{proof}
The proof of \ref{existenceobjectbelowmiddle} applies essentially verbatim upon replacing the first paragraph with an application of \ref{lem:surgery_data_for_objects_middle_dimension_n_is_0} to generate the initial surgery datum on each object, and noting that the constructions do not leave the base point component. 
\end{proof}

See also \ref{rem:datagrw} below for an alternative argument, that is more in line with the analogous statement in \cite[Section 5.2]{GRW}.

\subsection{The surgery complex in the middle dimension}\label{subsec:surgin2}

The goal of this section is to prove the following result:

\begin{theorem}
 \label{thm:surgmid}
Let $(\C,\QF)$ be a Poincar\'e category of exact dimension $2\oc$ and assume that $\Lin_\QF(Y)$ is $\oc$-connective for every $Y \in \C^\heart$. Then the map 
\[|\Cob^{\oc,\oc+1}(\C,\QF)| \longrightarrow |\Cob^{\oc,\oc}(\C,\QF)|\]
is an equivalence onto the base point component of the target.
\end{theorem}

Before giving the proof, let us explain why the surgery argument from Section \ref{subsec:surgbelow} does not apply under in the middle dimension.
Suppose that we are given two pieces of surgery data on the same Poincar\'e object $(C,q)\in \Poinc(\C, \QF\qshift 1)$ with underlying diagrams $f\colon S\to C$ and $g \colon T\to C$, both suitable for surgery of type $(\oc,\oc)$. We note that the first part of \ref{lemmaautmaticdisjoint} still applies in the present setting 
so that there exists a direct sum surgery datum $S\oplus T\to C$. Doing surgery on the diagram
\[\xymatrix{
 0 \ar[r]\ar[d] & S\ar[d] \\ T \ar[r] & S\oplus T
}\]
yields a commutative diagram in $\Cob^{\oc,\oc}(\C, \QF)$ of the following form:
\[\xymatrix{
 C \ar@{~>}[r]^{\chi(f)} \ar@{~>}[d]_{\chi(g)} & 0 \ar@{~>}[d]^{\chi(T)}\\
 0 \ar@{~>}[r]^{\chi(S)} & C_{f+g}.
}\]
However, there is no reason any more to expect that doing surgery along the direct sum datum $S\oplus T$ renders $C$ $\oc+1$-connective (i.e.\ zero). Thus, we do not obtain a diagram of surgery data that is suitable for surgery of type $(\oc,\oc)$ and cannot directly continue with our strategy.

To solve the issue, we shall take an element of the surgery complex, given by a cocartesian diagram of surgery data, and modify it into a (no longer cocartesian) cube of surgery data, which however is suitable for surgery of type $(\oc,\oc)$ away from the initial vertex; so we can continue our surgery argument with this modified cube.

To obtain the modified cube, we observe that by \ref{lem:bicartesian_square_surgeres_to_independent}, the square above is independent so we can form its double, compare Section \ref{sec:reflection_and_doubling}, which is a diagram in $\Cob(\C,\QF)$ as follows:
\begin{equation*}\tag{$\ast$}
\xymatrix{
 C \ar@{~>}[r]^{\chi(f)}\ar@{~>}[d]^{\chi(g)} 
  & 0 \ar@{~>}[d]^{\chi(T)} \ar@{~>}[r]^{\overline{\chi(f)}}
  & C \ar@{~>}[d]^{\chi(g)}
\\
 0 \ar@{~>}[r]^{\chi(S)} \ar@{~>}[d]^{\overline{\chi(g)}} 
  & C_{f+g} \ar@{~>}[r]^{\overline{\chi(S)}} \ar@{~>}[d]^{\overline{\chi(T)}}
  & 0\ar@{~>}[d]^{\overline{\chi(g)}}
\\
 C \ar@{~>}[r]^{\chi(f)}
  & 0 \ar@{~>}[r]^{\overline{\chi(f)}}
  & C. 
}
\end{equation*}
Here we have denoted the reflection of cobordisms by overlining them. It turns out that this diagram 
lies in $\Cob^{\oc, \oc}(\C, \QF)$ and thus the subdiagram
\[\xymatrix{
 C \ar@{~>}[r]^{\chi(f)}\ar@{~>}[d]_{\chi(g)} 
  & 0 \ar@{~>}[dd]^{\overline{\chi(T)} \circ \chi(T)} 
\\
 0 \ar@{~>}[rd]_{\overline{\chi(T)} \circ \chi(S)} 
  & 
\\
  & 0 
}\]
lies in $\Cob^{\oc,\oc+1}(\C, \QF)$. Instead of using the original square of surgery data, we will thus instead use this (slanted) square to perform multiple surgery. 

When attempting the same procedure for $n>0$, there is the additional complication that the reflection of a forwards surgery datum is backwards, so that the composites in the above square end up with neither property, so we cannot use the forwards or backwards surgery equivalences. However, one can easily modify the reflected surgery datum from a backwards one to a forwards one, without changing its effect on the objects of the cobordism category. Since the reflected surgery datum is dual to the original one, this amounts to modifying the original surgery datum from a forwards one to a backwards one, without changing its effect on the objects of the cobordism category. The data that will allow us to so is the following:

\begin{definition}
A \emph{split} of a forwards surgery datum in $N_n\Cobsd(\C, \QF)=\Psd(\Q_n(\C, \QF\qshift 1))$, with underlying diagram denoted $T\to C$, consists of a map $s\colon T\to \widehat T$ in $\Q_n(\C)$ such that the following hold:
\begin{enumerate}
 \item $\widehat T$  is backwards, i.e. all maps $\widehat T_{i \leq j} \rightarrow \widehat T_{i' \leq j}$ are equivalences, and
 \item for each $i$, the map $s\colon T_{i\leq i}\to \widehat T_{i\leq i}$ is an equivalence.
\end{enumerate}
\end{definition}

The name split is justified by the fact that a split of a forwards surgery datum is the same datum as a choice of left inverse of all the maps $T_{i-1\leq i}\to T_{i\leq i}$. 

Let us outline how the extra datum of a split will be used in the argument: Consider the surgery datum underlying the diagram $(\ast)$ above in $\Cob(\C,\QF)$, i.e.
\[\xymatrix{
0 \ar[r]\ar[d] & S \ar[d]\ar[r]^{\iota_s}  & \Db(S)\ar[d] \\
T \ar[r] \ar[d]^{\iota_T} & S \oplus T \ar[r] \ar[d] & T \oplus \Db(S) \ar[d] \\
\Db(T) \ar[r] &  \Db(T) \oplus S  \ar[r]   & \Db(S \oplus T). 
}\]
By direct inspection (or Lemma \ref{lem:strongly_cocartesian_surgeres_to_independent}) it is strongly cocartesian, so determined by the top horizontal row and left vertical column. 

Now, the double being a composite, its surgery datum takes part in a fibre sequence
\[T \xrightarrow{\iota_T} \Db(T) \to \Dual_\QF(T)\]
where the last part is the surgery datum of the reflected morphism. We can therefore use the splits $t$ and $s$ of $T$ and $S$, respectively, to replace the outer two corners by the pullbacks 
\[\widehat{\Db}(T,t) = \Dual_\QF(\widehat T) \times_{\Dual_\QF T} \Db(T) \quad \text{and} \quad \widehat{\Db}(S,s) = \Dual_\QF(\widehat S) \times_{\Dual_\QF S} \Db(S),\]
respectively, resulting in a new, strongly cocartesian diagram 
\[\xymatrix{
0 \ar[r]\ar[d] & S \ar[d]\ar[r]^{(0,\iota_S)}  & \widehat{\Db}(S,s)\ar[d] \\
T \ar[r] \ar[d]^{(0,\iota_T)} & S \oplus T \ar[r] \ar[d] & T \oplus \widehat{\Db}(S,s) \ar[d] \\
\widehat{\Db}(T,t) \ar[r] &  \widehat{\Db}(T,t) \oplus S  \ar[r]   & \widehat{\Db}(S \oplus T,s\oplus t). 
}\]
which is equipped with a natural transformation to the original diagram. In particular, it refines to an element of $N_n\Cobccsd 2 2(\C,\QF)$ and by direct inspection consists entirely of forwards surgery data. Extracting the slanted subdiagram as above, then gives us a $\lcube 1 {r+1}$-fold surgery datum suitable for surgery of type $(\oc,\oc)$ away from the initial vertex. It is this square that we will use in the proof below. 

\begin{lemma}\label{lem:surgery_data_admits_split}
Let $(\C,\QF)$ be a Poincar\'e category of exact dimension $2\oc$. Then each surgery datum suitable for surgery in the middle dimension admits a split.
\end{lemma}

\begin{proof}
Let $T\to C$ be such a surgery datum. Since in the cofibre sequence 
\[T_{i-1\leq n}\to T_{i\leq n} \to T_{i\leq n}/T_{i-1\leq n}\]
all terms are concentrated in degree $\oc$, the sequence splits by \ref{heartsplit}, so we may choose a left inverse $q_i\colon T_{i\leq n}\to T_{i-1\leq n}$ of the left map. Let $\widehat T$ denote the backwards surgery datum defined by pulling back the sequence
\[T_{n\leq n} \xrightarrow{q_n} T_{n-1\leq n} \xrightarrow{q_{n-1}} \dots \xrightarrow{q_1}  T_{0\leq n}\]
along the last vertex map
\[\Twar[n]\to [n]\op, \quad (i\leq j)\mapsto j.\]
Since as a forwards surgery datum $T$ is left Kan extended from its restriction along the first vertex inclusion
\[[n]\to \Twar[n], \quad i\mapsto (i\leq n)\]
natural transformations $T\to \widehat T$ are described by natural transformations
\[(T_{0\leq n}\to \dots \to T_{n\leq n}) \to (\widehat T_{0\leq n}\to \dots \to \widehat T_{n\leq n})\]
where we note that the target chain is constant with value $T_{n\leq n}$, so that such natural transformations are again described by maps $T_{n\leq n}\to T_{n\leq n}$. Define $s\colon T\to \widehat T$ as the transformation corresponding to the identity of $T_{n\leq n}$. At the object $T_{i\leq i}$, the transformation $s$ evaluates to the map 
\[T_{i\leq i} \simeq T_{i\leq n} \to T_{n\leq n} \xrightarrow{q_{i+1}\circ\dots \circ q_n} T_{i\leq n}\]
which is an equivalence as required.
\end{proof}

\begin{remark}\label{rem:datagrw}
One can give an alternative argument for \ref{prop:surgery_data_for_objects_exist_middle_dimension}, which makes the existence and effect of splits quite transparent: Given an element $C \in \nerv_n\Cob^{\oc,\oc}(\C,\QF)$ one can start with a surgery datum on each object $f_i \colon T_i \rightarrow C_{i \leq i}$ and combine these into a surgery datum suitable for middle dimensional surgery as follows: $T_0$ can be extended to a cylindrical forwards surgery datum on the entirety of $C$, and the result extended to a $[1]$-fold surgery datum with terms $T_0 \rightarrow \Db(T_0) \rightarrow C_{i \leq j}$ (recall that $\Db(S) \simeq S \oplus \Dual_\QF(S)$ denotes the double of a surgery datum $g \colon S \rightarrow X$, i.e.\ the surgery datum of the cobordism $\overline{\chi(g)} \circ \chi(g)$). This $[1]$-fold surgery datum on $C$ restricts to one of the form
\[\xymatrix{T_0 \ar[d]^{f_0} &\ar[l] T_0\ar[d] \ar[r] & \Db(T_0)\ar[d] & \ar[l]\ar[r]\Db(T_0)\ar[d] & \Db(T_0) \ar[d]&\ar[l] \dots \\
            C_{0 \leq 0} & \ar[l] \C_{0 \leq 1} \ar[r]& C_{1 \leq 1} &\ar[l] C_{1 \leq 2} \ar[r] & C_{2 \leq 2} & \ar[l] \dots}\]
The effect of surgery is then to make $C_{0 \leq 0}$ vanish, but not affect the objects $C_{i \leq i}$ for $i > 0$. Applying a similar strategy to $T_1 \rightarrow C_{1 \leq 1}$ gives a surgery datum
\[\xymatrix{0 \ar[d] &\ar[l] \ar[d] 0\ar[r] & T_1\ar[d]_{f_1} & \ar[l]\ar[r]T_1\ar[d] & \Db(T_1) \ar[d]&\ar[l] \dots \\
            0 & \ar[l] 0 \ar[r]& C_{1 \leq 1} &\ar[l] C_{1 \leq 2} \ar[r] & C_{2 \leq 2} & \ar[l] \dots};\]
Continuing in this fashion similarly as in the proof of \ref{existenceobjectbelowmiddle}, we obtain a string of cobordisms and the surgery datum of its composite is suitable for surgery of type $(\oc,\oc)$, as we invite the reader to check. 

Each of the constituent pieces of surgery data above, then have an obvious split to
\[\xymatrix{T_0 &\ar[l] \Db(T_0) \ar[r] & \Db(T_0) & \ar[l]\ar[r]\Db(T_0) & \Db(T_0) &\ar[l] \dots}\]
\[\xymatrix{0  &\ar[l] T_1 \ar[r] & T_1 & \ar[l]\ar[r]\Db(T_1) & \Db(T_1) &\ar[l] \dots }\]
and so forth. Doing surgery on the particular kind of surgery data arising from this construction would be a more direct analogue of the strategy of Galatius and Randal-Williams, see \cite[Section 5.2]{GRW}, particularly Figure 9. However, in order to execute their strategy faithfully, one needs to show that the realisation of the surgery complex spanned by this kind of data is equivalent to its $0$-simplices, compare \cite[Theorem 5.14]{GRW}. 

We found it more convenient to isolate the splits as the pertinent datum and instead form the diagram involving the twisted doubles as explained above. For the data considered above the twisted doubles are simply given by
\[\xymatrix{\Db(T_0) &\ar[l] \Db(T_0) \ar[r] & \Db(T_0) \oplus \Db(T_0) & \ar[l]\ar[r]\Db(T_0) \oplus \Db(T_0) & \Db(T_0) \oplus \Db(T_0) &\ar[l] \dots}\]
\[\xymatrix{0 &\ar[l] 0 \ar[r] & \Db(T_1) & \ar[l]\ar[r]\Db(T_1) & \Db(T_1) \oplus \Db(T_1) &\ar[l] \dots.}\]
\end{remark}

One can easily refine the definition of the surgery complex from section \ref{subsec:surgbelow} to include the choice of splits:

\begin{definition}\label{def:surgerycomplexmiddle}
The \emph{split surgery complex} $\SSC^\oc(\C, \QF)$ is the semi-simplicial category defined in degree $r$ as the pull-back
\[\xymatrix{
 \SSC_r^{\oc}(\C, \QF) \ar[rr]\ar[d]
  && \SC_r^{\oc,\oc}(\C, \QF)\ar[d]
\\
 \prod_{[r]} \mathrm{Splt}(\C) \ar[rr]
 && \prod_{[r]}  \Span(\C)
}\]
where the right map extracts the underlying diagram of surgery data at the singletons, $\mathrm{Splt} \subseteq \Span(\Ar(\C))$ is the subcategory of spans of arrows which are forwards at the source and backwards at the target and equivalences at the objects of the span category, and the lower map evaluates arrows at the source. 
\end{definition}

\begin{proposition}
\label{prop:surgery_complex_for_objects_contractible_middle_dimension}
Let $(\C,\QF)$ be a Poincar\'e category of exact dimension $2\oc$ and assume that $\Lin_\QF(Y)$ is $\oc$-connective for each $Y \in \C^\heart$. Then for each $n$, the forgetful map 
\[\nerv_n \SSC^{\oc}(\C, \QF) \longrightarrow \const \nerv_n\Cob^{\oc, \oc}_0(\C, \QF)\]
is a trivial fibration.
In particular, the forgetful map induces an equivalence
\[|\SSC^{\oc}(\C, \QF)| \longrightarrow |\Cob^{\oc, \oc}_0(\C, \QF)|.\] 
\end{proposition}

\begin{proof}
Per construction, the forgetful map $\nerv_n\SSC^{\oc} \rightarrow \nerv_n \SC^{\oc,\oc}$ is the pullback of a $0$-coskeletal map, and thus itself $0$-coskeletal. \ref{coskeletality} combined with the same argument as in \ref{surgery below the middle dim} shows that the forgetful map $\nerv_n \SC^{\oc,\oc}\to \nerv_n \Cob^{\oc,\oc}_0(\C, \QF)$ is $1$-coskeletal. Thus, it remains to see the filling conditions in degree $0$ and $1$, It is surjective on $0$-simplices by \ref{prop:surgery_data_for_objects_exist_middle_dimension} combined with \ref{lem:surgery_data_admits_split}, and boundaries of $1$-simplices are filled by \ref{lemmaautmaticdisjoint}.
\end{proof}

We now carry out the strategy outlined in the introduction to this section and use the surgery complex above to produce a section of the map
\[\left|[r] \mapsto |\Cobmfw{\lcube 1 {r+1}_0}^{\oc,\oc}(\C,\QF)|\right| \longrightarrow |\Cob^{\oc,\oc}(\C,\QF)|\]
over the base point component, to which we can apply \ref{generalmethodobjects}.

Consider an element $(f \colon T \rightarrow C, (s_i\colon T(\delta_i)\to \widehat T(i))_{i})$ in the split surgery complex $\nerv_n\SSC^{\oc}_{n,r}(\C, \QF)$; here $\delta_i \in \lcube 1 {r+1}$ denotes the $i$-th unit vector. Doing surgery on this multiple surgery datum results in an object of
\[\Poinc(\Q_{[1]^{r+1} \times [n]}(\C,\QF\qshift{1})) \subseteq \Poinc(\Q_{\lcube 1 {r+1}}\Q_n(\C,\QF\qshift{1}))\] and we learn from \ref{lem:strongly_cocartesian_surgeres_to_independent} that it even belongs to $\Qsq_{[1], \dots, [1]}\Q_n(\C, \QF\qshift 1))$. We can thus double $C_f$ to a Poincar\'e object 
\[\Db(C_f)\in \Qsq_{[2], \dots, [2]}\Q_n(\C, \QF\qshift 1)\subseteq \Q_{[2]^{r+1}}\Q_n(\C, \QF\qshift 1).\]
Let $\Db(f) \colon \Db(T) \rightarrow C$ be the corresponding surgery datum in $\Pmsd{\lcube 2 {r+1}}\Q_n(\C,\QF\qshift{1})$ obtained from the surgery equivalence. Per construction the value of $\Db(C_f)$ at each $v \in \lcube 2 {r+1}$ with exactly one coordinate equal to $1$ coincides with the value of $C_f$ at the unit vector with the same coordinate equal to $1$. We thus find:

\begin{observation}\label{lem:condition_on_doubled_cube}
For any $v \in \lcube 2 {r+1}$ with exactly one coordinate equal to $1$, the value of $\Db(C_f)$ at $v$ lies in $\nerv_n\Cob^{\oc,\oc+1}(\C,\QF)$.\end{observation}

As discussed earlier $\Db(C_f)$ does, however, not lie in the forwards part 
\[\Poinc\Q_{[2]^{r+1} \times [n]}(\C,\QF\qshift{1}) \subseteq \Poinc\Q_{[2]^{r+1}}\Q_n(\C, \QF\qshift 1),\]
which we correct using the splits $s_i$: 
By Lemma \ref{lem:strongly_cocartesian_surgeres_to_independent}, the underlying diagram of $\Db(T)$ is a strongly cocartesian diagram $[2]^{r+1}\to \Q_n\C$. Denote by $\widehat{\Db}(T,s) \colon [2]^{r+1} \rightarrow \Q_n\C$ the left Kan extension of the diagram $[2]^{r+1}_{ax} \rightarrow \Q_n\C$ given by taking $(0,\dots,0)$ to $0$ and $\delta_i \rightarrow 2\delta_i$ to 
\[T(\delta_i) \longrightarrow \widehat{\Db}(T(\delta_i),s_i),\]
where $\delta_i$ is the $i$'th unit vector and 
\[\widehat{\Db}(T(\delta_i),s_i) = \Dual_\QF(\widehat T(i)) \times_{\Dual_\QF T(\delta_i)} \Db(T(\delta_i))\]
as in the introduction. Per construction $\widehat{\Db}(T,s)$ comes equipped with a natural transformation to $\Db(f)$, and so inherits a surgery datum by \ref{lem:cofinal_surgery_data}. As every step in this construction is functorial we obtain a map
\[\widehat{\mathrm{db}}\colon \SSC_r^\oc \longrightarrow \Cobmsd{[2]^{r+1}}(\C,\QF)\]
taking twisted doubles, that is natural in $r$, i.e. a map of semi-simplicial categories.

We next restrict to the subcube of sidelength $1$ inside $[2]^{r+1}$, indicated in the introduction for $r=1$: Consider therefore the map 
\[i_r \colon [1]^{r+1}\longrightarrow [2]^{r+1}, \quad v \longmapsto 2v - \delta_{\mathrm{min} v}\]
where
\[\delta_{\mathrm{min} v} = \begin{cases} 0 & v = 0 \\ \delta_i & i \text{ is the smallest non-zero coordinate of } v. \end{cases}\]
This map is clearly monotone and commutes with the semi-simplicial structure of source and target (which is given by inserting $0$ as a new coordinate).
Restricting along this inclusion we obtain a map
\[\msrg\colon \SSC_r^\oc \longrightarrow \Cobmsd{\lcube 1{r+1}}(\C,\QF)\]
which is natural in $r  \in \Delta_\inj$.

\begin{lemma}\label{lem:properties_of_modified_double}
After restriction along $\lcube 1{r+1}_0\subset \lcube 1{r+1}$, the  map $\msrg$ just constructed takes values in $\Cobmfw{\lcube 1{r+1}_0}^{\oc,\oc}(\C,\QF)$.
\end{lemma}

\begin{proof}
By \ref{lem:suitable_is_good}, it is enough to show that (the nerve of) $\msrg$ takes values in surgery data which are suitable for surgery of type $(\oc,\oc)$. The last condition of Definition \ref{def:suitablebelowmid} follows from \ref{lem:condition_on_doubled_cube}, since the map $i_r$ is designed so that vectors in its image have exactly one coordinate equal to $1$. We claim that the other four conditions already hold for the map $\widehat{\mathrm{db}}$ (still discarding the initial vertex). 

So consider a specific element $X = (f\colon T\to C, (s_i\colon T(\delta_i)\to \widehat T(i))_i)$ in the (nerve of the) split surgery complex. Then $\widehat{\mathrm{db}}(X)$ is a strongly cocartesian cube (of sidelength 2), and one easily checks it suffices to prove the claim after restriction to $\lcube 2{r+1}_\mathrm{ax}-\{0\}\subset \lcube 2 {r+1}_0$. This in turn we can prove by considering the values 
\[\widehat{\mathrm{db}}(X)(\delta_i)= T(\delta_i)
\quad \mathrm{and}\quad 
\widehat{\mathrm{db}}(X)(2\cdot\delta_i)/\widehat{\mathrm{db}}(X)(\delta_i)=\Dual_\QF \widehat T(i)\] 
for each $i$. The first, $T(\delta_i)$, is suitable for surgery of type $(\oc,\oc)$ by assumption and we are left to consider the second one, $\Dual_\QF\widehat T(i)$: It is forwards by construction, and since each $T(i)_{j\leq j}$ is concentrated in degree $p$, so is $(\Dual_\QF \widehat T(i))_{j\leq j}=\Dual_\QF(T(i)_{j\leq j})$, $(\C, \QF)$ being of exact dimension $2p$. It remains to see that the maps $(\Dual_\QF \widehat T(i))_{j\leq j+1}\to (\Dual_\QF \widehat T(i))_{j+1\leq j+1}$ are $p$-coconnected; this is equivalent to the maps $\widehat T(i)_{j\leq j+1}\to \widehat T(i)_{j\leq j}$ being $p$-connected. By construction these are splits of the maps $T(i)_{j\leq j+1}\to T(i)_{j+1\leq j+1}$, so the fibre of the former map is the cofibre of the latter, which is indeed concentrated in degree $p$.
\end{proof}

\begin{proof}[Proof of \ref{thm:surgmid}]
Combining the previous lemma with \ref{prop:surgery_complex_for_objects_contractible_middle_dimension} 
therefore provides us with a map 
\[s \colon |\Cob^{\oc,\oc}_0(\C,\QF)| \simeq |\SSC^{\oc}(\C,\QF)| \xrightarrow{\msrg} |[r] \mapsto |\Cobmfw{\lcube 1 {r+1}_0}^{\oc,\oc}(\C,\QF)||\]
which is easily checked to be a section of the forgetful map. Furthermore, since any twisted double of the surgery datum $0$ is necessarily $0$, it also satisfies the second assumption of \ref{generalmethodobjects}, an application of which finishes the proof.
\end{proof}

\section{Surgery on morphisms}\label{sec:morph}
Again, throughout this section we consider a Poincar\'e category $(\C, \QF)$ of dimension at least $d$, and our goal here is to address finally the case of surgery on morphisms. We thus wish to show that 
\[\Cob^{\mc+1}(\C,\QF) \longrightarrow \Cob^{\mc}(\C,\QF)\]
becomes an equivalence upon realisation under appropriate numerical assumptions, see \ref{thm:surgmor} below.

\subsection{The surgery complex for morphism surgery}\label{subsec:surgcompxmorph}

As explained at the end of Section \ref{subsec:surgery_on_component}, the surgery datum 
\[ \xymatrix{
 0 \ar[d] & T \ar[l] \ar[r]^{\mathrm{id}} \ar[d] & T \ar[d]^f\\
 C & W \ar[l] \ar[r] & D
} 
\]
can be used to improve the connectivity of $W \rightarrow C$ from $\mc$ to $\mc+1$, if the map $T\to \fib(W\to C)$ is $\mc$-connective:
The result of surgery is the diagram

\[\xymatrix{C \ar@{~>}[r]^{W/T} \ar@{~>}[d]_{\mathrm{id}}& D_f \ar@{~>}[d]^{\chi(f)} \\
            C \ar@{~>}[r]^W & D,}\]
and $W/T$ lies in $\Cob^{\mc+1}(\C,\QF)$. The cobordism $\chi(f)$, however, only lies in $\Cob^\mc(\C,\QF)$ which, once performing multiple surgery will prevent cobordisms relating different pieces of surgery data from lying in $\Cob^{\mc+1}(\C,\QF)$ (and this is required in the analogue of Proposition \ref{generalmethodobjects}). In the analysis of components in Section \ref{subsec:surgery_on_component} we solved this problem by simply reflecting $\chi(f)$, but this will break the commutativity of the above diagram, since $0 \leftarrow T \rightarrow T$ is not a forwards surgery datum. To perform parametrised surgery, following Galatius--Randal-Williams, we factor the trace of the surgery along $0\leftarrow T\rightarrow T$ into a pair of cobordisms, one which is forwards and one of which is backwards, and both of which are $\mc+1$-connected relative to the relevant boundary piece, at each object. 

More concretely, we extend the diagram of surgery data from above to a $[1]$-fold surgery datum
\[\xymatrix{
 T \ar[d] & T \ar[l]_\id \ar[r]^\id \ar[d]^\id & T \ar[d]^\id\\
 0 \ar[d] & T \ar[l] \ar[r]^{\id} \ar[d] & T \ar[d]^f\\
 C & W \ar[l] \ar[r] & D
}\]
in $\nerv_1 \Cobmbw{[1]}(\C,\QF)$. Performing backwards multiple surgery using \ref{corhomotopy} results in the left of the two commutative diagrams in $\Cob(\C,\QF)$
\[\xymatrix{C \ar@{~>}[r]^{W/T}\ar@{~>}[d] & D_f \ar@{~>}[d]^{\id} && C \ar@{~>}[r]^{W/T}\ar@{~>}[d] & D_f \ar@{~>}[d]^{\id} \\
C_T\ar@{~>}[r]^V \ar@{~>}[d]& D_f \ar@{~>}[d]^{\chi(f)} && C_T\ar@{~>}[r]^V & D_f  \\
            C \ar@{~>}[r]^W & D && C \ar@{~>}[r]^W \ar@{~>}[u]& D\ar@{~>}[u]_{\chi(f)}}\]
whose left vertical composite is the identity (being obtained by surgery along the zero surgery datum). But the lower square in this diagram is independent in the sense of Section \ref{sec:reflection_and_doubling}, since its surgery datum $T \leftarrow T \rightarrow T$ is both forwards and backwards. Thus it can be reflected resulting in the commutative diagram on the right. In this form the zig-zag from bottom to top can be regarded as parts of two natural transformations with values in $\Cob^{\mc+1}(\C,\QF)$ as required. For this reason the surgery move on morphisms will require cocartesian cubes of $[1]$-fold surgery data.

\begin{remark}
This two-step process corresponds loosely to the two-step process in the geometric surgery on morphisms from \cite{GRW}, but our actual move is simpler than theirs, resulting in an estimate in Theorem \ref{thm:surgmor} which works in one further degree in the odd dimensional case. We do not know a geometric analogue of our surgery move.
\end{remark}

We next formalise this kind of surgery data.  

\begin{definition}\label{def:of_type_kappa}
A $[1]$-fold surgery datum in $\nerv_n \Cobmsd{[1]}(\C, \QF)=\Pmsd{[1]}(\Q_n(\C, \QF\qshift 1))$ with underlying diagram  $S\to T\to C$, where $C \in \nerv_n\Cob^{\mc}(\C,\QF)$, is called \emph{suitable for surgery of type $\mc$} if 
\begin{enumerate}
\item $T$ is backwards, \emph{i.e.}, each $T_{j\leq j+1}\to T_{j+1\leq j+1}$ is an equivalence,
\item $S$ is forwards and backwards, \emph{i.e.}, a constant diagram,
\item each $S_{i \leq i}$ and $T_{i \leq i}$ are concentrated in degree $\mc$,
 \item\label{item:cofiber_concentrated} the cofibre of each map $S_{i \leq i} \rightarrow T_{i \leq i}$ is concentrated in degree $\mc+1$, and 
\item The result of surgery of $T \rightarrow C$ lies in $\nerv_n\Cob^{\mc+1}(\C,\QF)$.
\end{enumerate}
We call the surgery datum \emph{weakly of type $\mc$} provided it satisfies all conditions but possibly the last one.
\end{definition}

\begin{proposition}\label{prop:surgery_data_for_morphisms_exist}
Let $(\C,\QF)$ be a Poincar\'e category with a weight structure. Then every $C\in \nerv_n\Cob^{\mc}(\C, \QF)$ admits a surgery datum suitable for surgery of type $\mc$.
\end{proposition}

\begin{proof}
For $n=0$ the $0$ surgery datum is of type $\mc$, and for $n=1$ the $[1]$-fold surgery datum represented by the $3$-by-$3$-square considered in the discussion above is easily checked suitable for surgery of type $\mc$ (compare the discussion before \ref{resolutionmorphismspi0}). For $n>1$ one can argue as follows: For any $r\in \{0,\dots, n-1\}$, we first construct a surgery datum on $(C,q)$, with underlying diagram $S^{(r)}\to T^{(r)}\to C$, such that the result of surgery on $T^{(r)}$, call it $C^{(r)}$, has the map $C^{(r)}_{r\leq r+1}\to C^{(r)}_{r\leq r}$ $\mc$-connective: 

For this, choose a surgery datum as in the case $n=1$ on $C_{r \leq r} \leftarrow C_{r \leq r+1} \rightarrow C_{r+1 \leq r+1}$. By \ref{extendingsurgdata} we can extend its middle row to a surgery datum $T^{(r)} \rightarrow C \in \Surg_n(\C,\QF\qshift{1})$ which is constant to the right, and vanishes to the left.
Letting $S^{(r)}\colon \Twar[n]\to \C$ be the constant diagram on $T^{(r)}_{r \leq r}$ we find a canonical map $S^{(r)}\to T^{(r)}$ and the composite $S^{(r)}\to T^{(r)}\to C$ uniquely inherits the structure of a $[1]$-fold surgery datum on $C$ by means of Lemma \ref{lem:cofinal_surgery_data}. By construction, $C^{(r)}_{r\leq r+1}\to C^{(r)}_{r\leq r}$ is indeed $\mc+1$-connective.

Now, let
\[S=\bigoplus_{r=0}^{n-1} S^{(r)}, \quad T=\bigoplus_{r=0}^{n-1} T^{(r)}\]
the induced map $S\to T\to C$ canonically inherits the structure of $[1]$-fold surgery datum: From each $T^{(r)}$ being backwards, we compute
\[\lim_{\Twar[n]\op} \QF\qshift 1\circ T \simeq \QF\qshift 1 (T_{0\leq 0}) = \QF\qshift 1(0)=0,\]
so the map $T\to C$ uniquely extends to a surgery datum; and the extension to a $[1]$-fold surgery datum is again automatic by Lemma \ref{lem:cofinal_surgery_data}. We are left to show that $S \rightarrow T \rightarrow C$ is suitable for surgery of type $\mc$. The first four conditions are clear from the construction. To see that $(C_T)_{r \leq r+1} \rightarrow (C_T)_{r \leq r}$ is $\mc+1$-connective we apply \ref{extraconnectivitycalc}, which requires us to verify that $\fib(T_{r\leq r+1}\rightarrow T_{r \leq r}) \rightarrow \fib(C_{r \leq r+1} \rightarrow C_{r \leq r})$ is $\mc$-connective. But the source evaluates to $T^{(r)}_{r \leq r}$, so this is true by construction.
\end{proof}

We proceed to define the surgery complex:

\begin{definition}\label{def:surgerycomplexmorph}
Given a Poincar\'e category $(\C, \QF)$ of dimension at least $d$, the simplicial space of \emph{disjoint $[r]$-tuples of surgery data suitable for type $\mc$} on $\Cob^\mc(\C, \QF)$ is given in degree $n$ by the full subspace of 
\[\Pccsd 2 {r}(\Q_n(\C, \QF\qshift{1})) \simeq \nerv_n\Cobccsd 2{r}(\C,\QF)\]
on those $[2]^{r+1}$-fold surgery data whose underlying diagram $T\to X$ satisfies:
\begin{enumerate}
 \item $X$ lies in $\nerv_n\Cob^{\mc}(\C,\QF)$, and
 \item for each $0\leq j\leq r$, the $[1]$-fold surgery datum
\[T(\delta_j)\to T(2\delta_j)\to X\]
 obtained from by restricting along $j\colon [0]\to [r]=\{1,\dots, r+1\}$ and forgetting the entry at $0$ is suitable for surgery of type $\mc$.
\end{enumerate}
We denote by $\SC^\mc_r(\C, \QF)$ the category associated to this simplicial space. For varying $r$ these assembles into a semi-simplicial category
\[\SC^\mc(\C, \QF)\colon \Dinj\op\to \Cat\]
which we call the \emph{surgery complex} for morphism surgery of type $\mc$. It is (degreewise) a subcategory of $\Cobccsd 2 r(\C,\QF)$ (see \ref{def:cocartsurgdata}).
\end{definition}

\begin{proposition}\label{prop:surgery_complex_for_morphisms_contractible}
Suppose that $(\C, \QF)$ is a Poincar\'e category of dimension at least $d$. Then for each $n$, the forgetful map 
\[\nerv_n \SC^\mc(\C,\QF) \longrightarrow \const\nerv_n\Cob^\mc(\C,\QF)\]
is a trivial fibration provided $2\mc\leq d$, and so in particular induces an equivalence
\[|\SC^\mc(\C,\QF)| \longrightarrow |\Cob^\mc(\C,\QF)|.\]
\end{proposition}

\begin{proof}
There is an evident inclusion 
\[\SC^{\mc}(\C,\QF) \subseteq \Cobccsd{2}{\bullet}(\C,\QF)\] 
and since the conditions that isolate the source are on the edges only, it follows from \ref{coskeletality} that the forgetful map in the proposition is $1$-coskeletal just as in \ref{surgery below the middle dim}.

\ref{prop:surgery_data_for_morphisms_exist} furthermore implies that it is a surjection on $\pi_0$, so we are only left to fill boundaries of $1$-simplices. This requires us to produce an element of $\nerv_n \SC^\mc_{1}(\C,\QF)$ with boundary two given pieces of surgery data $T_0,T_1 \in \nerv_0\SC^\mc_{0}(\C,\QF)$ on the same object $(C,q)$. Per definition the diagrams underlying $T_0$ and $T_1$ assemble into 
a functor on
\[([0] \times [2]) \cup ([2] \times [0]) = [2]^2_\mathrm{ax},\]
equipped with a transformation to the constant diagram on $C$. We can left Kan extend it to a functor on $[2]\times [2]$ (which is still equipped with such a transformation) and then need to refine $T$ to a $[2]\times [2]$-fold surgery datum, compatible with the multiple surgery data on $[0]\times [2]$ and $[2]\times [0]$. By Lemma \ref{lem:cofinal_surgery_data}, this amounts to finding a nullhomotopy of the image of the Poincar\'e form $q$ under the map 
\[\Omega^\infty \QF\qshift 1_n(C) \to \Omega^\infty \QF\qshift 1_n(T(2,2)) \simeq \Omega^\infty \QF\qshift 1_n(T_0 \oplus T_1)\]
extending the given nullhomotopies of the image in $\Omega^\infty \QF\qshift 1_n(T_0)\times \Omega^\infty \QF\qshift 1_n(T_1)$. Since $T_0$ is backwards we have
\[\QF\qshift 1_n(T_0) = \lim_{\Twar[n]\op} \QF\qshift 1\circ T_0\op \simeq \QF\qshift 1((T_0)_{0\leq 0})
\]
and similarly for $T_1$ and $T_0 \oplus T_1$. Thus an extension exists by \ref{autodisj} once $2\mc \leq d$. 
\end{proof}

\subsection{Surgery on morphisms}

We seek to prove:

\begin{theorem}\label{thm:surgmor}
Let $(\C,\QF)$ be a Poincar\'e category of dimension at least $d$. Then the inclusion 
\[\Cob^{\mc+1}(\C,\QF) \longrightarrow \Cob^{\mc}(\C,\QF)\]
is an equivalence provided $2\mc+1\leq d$.
\end{theorem}

Backwards surgery defines a functor
\[\nerv_n \SC_r^\mc(\C, \QF) \to \Q_{([2]\op)^{r+1} \times [n]}(\C, \QF\qshift 1)  \subseteq \Q_{([2]\op)^{r+1}}\Q_n(\C, \QF\qshift 1).\]
In fact, by \ref{lem:strongly_cocartesian_surgeres_to_independent} it even takes values in  
\[\Qsq_{[2]\op, \dots, [2]\op} \Q_n(\C, \QF\qshift 1)\subseteq \Q_{([2]\op)^{r+1}}\Q_n(\C, \QF\qshift 1)\]
(though not in $\Qsq_{[2]\op, \dots, [2]\op,[n]}(\C, \QF\qshift 1)$, since only half the surgery datum is forwards).
Applying the Reflection Lemma \ref{lem:reflection_lemma}, and abbreviating
\[\Horn = [1]\cup_1 [1],\]
we obtain a map
\[\rsrg\colon \nerv_n \SC_{r}^\mc(\C, \QF) \longrightarrow \Poinc(\Qsq_{\Horn, \dots, \Horn} \Q_n(\C, \QF\qshift 1))\subseteq \Poinc(\Q_{\Horn^{r+1}} \Q_n(\C, \QF\qshift 1)).\]
We choose the convention that evaluation at $(0_1,\dots,0_1) \in \Horn^{r+1}$ extracts the underlying object, where $0_1 \in \Horn$ is the unique object not coming from the right summand. 

\begin{observation}
The map $\rsrg$ in fact takes values in $\Poinc\Q_{\Horn^{r+1} \times [n]}(\C,\QF\qshift{1})$ and therefore defines maps
\[\rsrg\colon \SC^\mc_r(\C,\QF) \longrightarrow \Fun(\Horn^{r+1}, \Cob(\C,\QF))\]
of semi-simplicial categories.
\end{observation}

\begin{proof}
For all $[1]\to \Horn^{r+1}$, which we further assume to be constant in all but one factor, we claim that the further restriction to $\Q_1\Q_n(\C)$ lands in $\Q_{[1]\times [n]}(\C)$; this will suffice by the Segal condition. To see the claim, we first note that restriction along the first edge $[1]\to \Horn$ corresponds to restricting the original datum along the edge $2\leq 1$ of $[2]\op$ from which the claim follows since the surgery datum is backwards. On the other hand, restriction along the second edge of $\Horn$ corresponds to  restriction along the edge $1\leq 0$ of $[2]\op$, followed by reflection in the first entry, from which the claim follows since the surgery datum is forwards. 
\end{proof}

Now let 
\[\Horn^{r+1}_0\subset\Horn^{r+1} \]
denotes the full subposet all functions that take the value $0_2\in \Horn=[1]\cup_{\{1\}} [1]$ at least once; for instance $\Horn^2_0\subset \Horn^2$ is given by 
\[\xymatrix@-1pc{  &  & (0_1,0_2) \ar[d]& &  (0_1,0_1) \ar[r] \ar[d] & (0_1,1)\ar[d] & \ar[l](0_1,0_2)\ar[d] \\
               & & (1,0_2) & \subset &  (1,0_1) \ar[r] & (1,1) & \ar[l](1,0_2) \\
              (0_2,0_1) \ar[r] &(0_2,1) & \ar[l](0_2,0_2)\ar[u]&&  (0_2,0_1) \ar[u]\ar[r]  & (0_2,1) \ar[u]& \ar[l](0_2,0_2).\ar[u]}
             \]

\begin{proposition}\label{lem:properties_of_surgery_diagram_morphisms}
Let $(\C,\QF)$ be a Poincar\'e category of dimension at least $d$ and assume that $2\mc+1\leq d$. Then functor $\rsrg$ fits into a commutative diagram
\[\xymatrix{
\SC^\mc(\C,\QF) \ar[r]^-{\rsrg}\ar[d] & \Fun(\Horn^{r+1},\Cob^\mc(\C,\QF))\ar[d]^{\mathrm{res}} \\
\Fun(\Horn^{r+1}_0,\Cob^{\mc+1}(\C,\QF)) \ar[r]& \Fun(\Horn^{r+1}_0,\Cob^\mc(\C,\QF))
}\] 
of semi-simplicial categories.
\end{proposition}

For the proof we need:

\begin{lemma}\label{lem:do_surgery_on_morphisms_estimate}
Let $X$ be a $[1]$-fold surgery datum weakly of type $\mc$ in $\Q_n(\C, \QF\qshift 1)$ with underlying diagram $S\to T\to C$, and suppose that $2\mc+1\leq d$. Then, $\rsrg(X) \colon [n]\times \Horn\to \Cob(\C, \QF)$ takes values in $\Cob^{\mc}(\C, \QF)$, and even in $\Cob^{\mc+1}(\C, \QF)$ provided $X$ lies in $\Cob^{\mc+1}(\C, \QF)$.
\end{lemma}

\begin{proof}
The statement follows by putting together the following claims (i) the trace of the surgery on each $S_{i\leq i}\to C_{i\leq i}$ is $\mc+1$-connective relative to its incoming end; (ii) the trace of the surgery on each $T_{i\leq i}/S_{i\leq i}\to (C_S)_{i\leq i}$ is $\mc+1$-connective relative to its outgoing end; (iii) $C_S \in \nerv_n\Cob^k(\C,\QF)$ if $C \in \Cob^k(\C,\QF)$; (iv) $C_T \in \Cob^\mc(\C,\QF)$ if $C_S \in \Cob^{\mc}(\C,\QF)$ and similarly for $\Cob^{\mc+1}(\C,\QF)$.

The first two are immediate from \ref{lem:basic_surgery_estimate} and last two are applications of \ref{connsurgmor}.
\end{proof}

\begin{proof}[Proof of \ref{lem:properties_of_surgery_diagram_morphisms}]
Let $T \rightarrow C$ be the underlying diagram of an element in $N_n\SC^\mc_{r}(\C, \QF)$. We shall analyse the restriction of $\rsrg(T)$ to 
\[\{a_0\} \times \dots\times  \Horn \times \dots\times\{a_r\}\subseteq \Horn^{r+1}\]
with fixed choices of $a_0, \dots a_{i-1}, a_{i+1}, \dots, a_r\in \Horn$: These evidently span the category $\Horn^{r+1}$ under composition, so showing that this span lies in $\Cob^\mc(\C,\QF)$ will give the first half of the claim. Restricting to those spans where one the $a_j$  is $0_2$, one readily checks that they span $\Horn^{r+1}_0$ under composition, so it will suffice to show that under this extra assumption the span even lies in $\Cob^{\mc+1}(\C,\QF)$. 
Under the reflection lemma, the restriction of $\rsrg(T)$ to this subposet is obtained by first applying multiple surgery along $T$, then restricting the resulting diagram $([2]\op)^{r+1} \times [n] \rightarrow \Cob(\C, \QF)$ along 
\[\{a_0\} \times \dots \times [2]\op\times\dots \times \{a_r\}\subset ([2]\op)^{r+1},\]
and reflecting the second cobordism in the resulting pair of composable cobordisms. 

Now the composable tuple of cobordisms under consideration is obtained by doing backwards $[1]$-fold surgery along 
\[\frac{T(a_0, \dots  1, \dots, a_k)}{T(a_0, \dots, 0, \dots, a_k)}
\to \frac{T(a_0, \dots  2, \dots, a_k)}{T(a_0, \dots, 0, \dots, a_k)}
\to \srg(T)(a_0, \dots, 0, \dots, a_k)\]
Since $T$ is strongly cocartesian, this is equivalent to 
\[T(0, \dots, 1, \dots, 0) \to T(0, \dots, 2, \dots, 0)\to \srg(T)(a_0, \dots, 0, \dots, a_k),\]
and the claim will follow from Lemma \ref{lem:do_surgery_on_morphisms_estimate}  if we can show that $\rsrg(T)(a_0, \dots, 0, \dots, a_k) \in \nerv_n\Cob^\mc(\C,\QF)$, and even lies in $\nerv_n\Cob^{\mc+1}(\C,\QF)$ provided that one of the $a_i$'s is $2$. 

If all $a_s$ are $0$, then  $\srg(T)(a_0, \dots, 0,\dots, a_k) = C$ indeed lies in the nerve of $\Cob^\mc(\C,\QF)$ by assumption. Then, inductively on the number of $a_s$ not being $0$, we conclude that the same is true for arbitrary values of $a_s$ by the previous argument.

Furthermore, if, say, $a_j=2$ and all other values of $a_s$ are $0$, then $\srg(T)(a_0, \dots, 0,\dots, a_k) \in \Cob^{\mc+1}(\C,\QF)$ by assumption on $T$. Again using the previous argument, we conclude inductively on the number of non-zero $a_s$, $s\neq j$, that the same is true for arbitrary values of $a_s$, $s\neq j$. 
\end{proof}

\begin{proof}[Proof of Theorem \ref{thm:surgmor}]
Consider the diagram
\[\xymatrix{
|\Cob^{\mc+1}(\C,\QF)| \ar[dd] 
& \Fun([1],|\Cob^{\mc}(\C,\QF)|) \times_{|\Cob^{\mc}(\C,\QF)|} |\Cob^{\mc+1}(\C,\QF)| \ar[d]\ar[l]_-{\ev_1} \ar@/_3pc/[ddl]_-{\ev_0} 
\\
& \Fun(\Horn^{r+1},|\Cob^{\mc}(\C,\QF)|) \times_{\Fun(\Horn^{r+1}_0,|\Cob^{\mc}(\C,\QF)|)} \Fun(\Horn^{r+1}_0,|\Cob^{\mc+1}(\C,\QF)|) \ar@/_1pc/[ld]^-{\ev_{(0_2,\dots,0_2)}} \ar@/_1pc/@{-->}[u]
\\
|\Cob^{\mc}(\C,\QF)| \ar@/_1pc/@{-->}[r]_s & |\SC^\mc(\C,\QF)|. \ar[u]^{\rsrg} \ar[l]_{\fgt}
}\]
The proof then proceeds verbatim as that of Proposition \ref{generalmethodobjects}, using contractibility of both $\Horn^{r+1}$ and $\Horn_0^{r+1}$.
\end{proof}

\section{Examples and applications}\label{sec:appl}

In this final section we assemble the various surgery moves from the previous chapters and prove the results stated in the introduction.

\subsection{Applications to Grothendieck-Witt spaces}
\label{sec:applgw}

We start out by proving Theorem \ref{thm:main}:

\begin{theorem}\label{thm:main2}
Let $(\C,\QF)$ be a Poincar\'e category of dimension at least $d$. Then the inclusion
\[\Cob^{\mc}(\C,\QF) \subseteq \Cob(\C,\QF)\]
becomes an equivalence upon realisation provided $2\mc \leq d+1$. Similarly, the inclusion
\[\Cob^{\mc,\oc}(\C,\QF) \subseteq \Cob^{\mc}(\C,\QF)\]
becomes an equivalence upon realisation provided
\begin{enumerate}
\item $\mc+\oc\leq d+1$;
\item $2\oc \leq d+1$;
\item $\oc \leq \mc+1$, and
\item $\Lin_\QF(X)$ is $(\oc-1)$-connective for every $X \in \C^\heart$. 
\end{enumerate}
In particular, if $2\oc\leq d+1$ and if $\Lin_\QF(X)$ is $(\oc-1)$-connective for every $X \in \C^\heart$, then
\[\Cob^{\oc,\oc}(\C,\QF) \rightarrow \Cob(\C,\QF)\]
is an equivalence on realisations. 

If $(\C,\QF)$ has exact dimension $2\oc$ and $\Lin_\QF(X)$ is even $\oc$-connective for all $X \in \C^\heart$, then also the inclusion
\[\Cob^{\oc,\oc+1}(\C,\QF) \subseteq \Cob(\C,\QF)_0\]
becomes an equivalence upon realisation, where the subscript denotes the component of $0$. 
\end{theorem}

\begin{proof}
We learn from Theorem \ref{thm:surgmor} that the inclusions
\[\Cob^{\mc}(\C,\QF) \longrightarrow \Cob^{\mc-1}(\C,\QF) \longrightarrow \Cob^{\mc-2}(\C,\QF) \longrightarrow \dots\]
become equivalences after realisation if $2\mc \leq d+1$. Since the extraction of morphism spaces commutes with filtered colimits of categories we find
\[\colim_{i \in \mathbb N} \Cob^{\mc-i}(\C,\QF) \simeq \Cob(\C,\QF),\]
since the weight structure on $\C$ is assumed exhaustive. This gives the first claim, as also realisation commutes with colimits. 

For the second statement we similarly apply Theorem \ref{thm:surgbelowmiddlecomplete} to the sequence
\[\Cob^{\mc,\oc}(\C,\QF) \longrightarrow \Cob^{\mc,\oc-1}(\C,\QF) \longrightarrow \Cob^{\mc,\oc-2}(\C,\QF) \longrightarrow \dots\]
to see that each of its arrows realises to an equivalence under the given assumptions. Investing additionally that $\pi_0\core \colon \Cat \rightarrow \mathrm{Set}$ commutes with filtered colimits, we find that 
\[\colim_{i \in \mathbb N} \Cob^{\mc,\oc-i}(\C,\QF) \simeq \Cob^{\mc}(\C,\QF).\]
The claim follows.

The third statement is obtained by simply plugging in $\mc=\oc$ into the first two, and the final one is a combination of this case and Theorem \ref{thm:surgmid}.
\end{proof}

With $d=0$ we find Theorem \ref{thm:resolution}, the weight theorem for Grothendieck-Witt spaces:

\begin{corollary}\label{cor:res2}
Let $(\C,\QF)$ a Poincar\'e category equipped with an exhaustive weight structure, such that $\QF(X)$ is connective for each $X \in \C^\heart$ and $\Dual_\QF$ preserves the heart of $\C$. Then the canonical map
\[\Poinc^\heart(\C,\QF)^\grp \longrightarrow \gw(\C,\QF)\]
is an equivalence.
\end{corollary}

\begin{proof}
Per assumption $(\C,\QF)$ is a Poincar\'e category of dimension exactly $0$ and from the fibre sequence 
\[\Bil_\QF(X,X)_\hC \longrightarrow \QF(X) \longrightarrow \Lin_\QF(X)\]
combined for example with Lemma \ref{weight:dualvsbil} we learn that $\Lin_\QF(X)$ takes connective values on $\C^\heart$. Thus we can apply \ref{thm:main2} to obtain that
\[|\Cob^{0,1}(\C,\QF)| \longrightarrow |\Cob(\C,\QF)|_0\]
is an equivalence. But by \ref{highlyconnectedvanishes} the source has connected groupoid core, with only object $0$, so 
\[\Cob^{0,1}(\C,\QF) \simeq \mathbb{B}\Hom_{\Cob^{0}(\C,\QF)}(0,0);\]
here we denote by $\ast/\Cat$ the category of based categories and by $\mathbb B \colon \mathrm{Mon}_{\Eone} \rightarrow */\Cat$ the left adjoint to taking endomorphisms . A priori one has to use the composition $\Eone$-structure on $\Hom_{\Cob^{0}(\C,\QF)}(0,0)$, but since $\Cob^{0}(\C,\QF)$ is closed under the symmetric monoidal structure of $\Cob(\C,\QF)$, which has unit $0$, this underlies the induced $\Einf$-structure by naturality.

Under the equivalence $\Hom_{\Cob(\C,\QF)}(0,0) \simeq \Poinc(\C,\QF)$ (per construction as $\Einf$-monoids) the subspace $\Hom_{\Cob^{0}(\C,\QF)}(0,0)$ corresponds to $\Poinc^\heart(\C,\QF)$. Thus
\[\gw(\C,\QF) \simeq \Omega|\Cob(\C,\QF)| \simeq \Omega|\Cob^{0,1}(\C,\QF)| \simeq \Omega|\mathbb{B}\Poinc^\heart(\C,\QF)| \simeq \Poinc^\heart(\C,\QF)^\grp\]
by Stasheff's recognition principle for $\Eone$-groups.
\end{proof}

By applying this to $\Hyp(\C)$ we immediately find the weight theorem in algebraic $\K$-theory as a special case:

\begin{corollary}\label{cor:fontes}
For $\C$ a stable $\infty$-category equipped with an exhaustive weight structure the natural map
\[\core(\C^\heart)^\grp \longrightarrow \k(\C)\]
is an equivalence. 
\end{corollary}

In particular (and in fact equivalently) the natural map
\[\core(\mathcal A)^\grp \longrightarrow \k(\Stab(\mathcal A))\]
is an equivalence for any additive category $\mathcal A$, since the left hand side is invariant under weak idempotent completion, e.g.\ by the group completion theorem.

\begin{proof}
From \ref{examplesofpoincweight} (6) we find that $\Hyp(\C)$ inherits an weight structure of dimension $0$ with $\Hyp(\C)^\heart = \C^\heart \times (\C^\heart)\op$. The equivalence $\Poinc(\Hyp(\C)) \simeq \core(\C)$ from \cite[2.2.5 Proposition]{9authI} therefore restricts to an equivalence
\[\Poinc^\heart(\Hyp(\C)) \simeq \core(\C^\heart)\]
on the one hand and provides an equivalence
\[\gw(\Hyp(\C)) \simeq \k(\C)\]
on the other, compare \cite[4.1.5 Corollary]{9authIII}. Now apply \ref{cor:res2}.
\end{proof}

\begin{remark}\label{rem:fontes} Corollary \ref{cor:fontes} is probably well-known to experts but has a somewhat chequered history in the literature:
\begin{enumerate}
\item In the case where $\mathcal A$ is an ordinary category, the result can be obtained from the literature by combining Quillen's `$+=\Q$'-theorem (identifying the group completion of $\core(\mathcal A)$ with his $\Q$-construction of $\mathcal A$ regarded as a split exact category, see \cite[p. 228]{QuillenHigherK2} or \cite[Theorem 7.1]{Kbook}), the Gillet-Waldhausen theorem (which allows one to pass to the $\Q$- or $\mathcal S$-construction applied to the Waldhausen category bounded chain complexes in $\mathcal A$, see \cite[Theorem 1.11.7]{tt} or \cite{raptisgw}) and finally comparison results to the $\K$-theory of higher categories (which equates this with the $\Q$- or $\mathcal S$- construction of the $\infty$-category of finite chain complexes in $\mathcal A$, which models $\Stab(\mathcal A)$ in this case, see \cite[Section 7.2]{BGT} or \cite[10.10 Proposition]{barwickhigher}).
\item In a very similar vein, the case where $\C = (\Spa/B)^\omega$, where $\K((\Spa/B)^\omega) \simeq \mathrm{A}(B)$, was essentially established by Waldhausen in \cite[Theorem 2.2.1]{waldhausen} as a consequence of his sphere theorem. 
\item In \cite[Theorem VI.7.1]{EKMM} Elmendorf, Kriz, Mandell and May proved a `$+\text{-equals-}\mathcal S$'-theorem for an arbitrary connective ring spectrum $R$ (in particular recovering the previous point using $R = \mathbb S[\Omega B]$ for connected pointed $B$). 
\item For $\C=\Mod^{\omega}_R$ (or $\mathcal A=\Proj(R)$) with $R$ a connective $\Eone$-ring, Lurie provided a direct $\infty$-categorical proof as part of a lecture course on algebraic K-theory in \cite[Lecture 19]{Lurie-K}. This proof was then massaged into the general statement of \ref{cor:fontes} by Clausen in an appendix to \cite{clausen-artin}, which was unfortunately removed before publication. We first learned about the general statement from a version of this appendix that Clausen kindly shared with us. 
\item One can in fact deduce the general statement of \ref{cor:fontes} from the special case of perfect modules over an $\Eone$-ring: Both functors in the statement preserve filtered colimits, and every idempotent complete stable category is the filtered colimit of its monogenic idempotent complete stable subcategories. By the Schwede-Shipley theorem these are equivalent categories of perfect modules over an $\mathbb E_1$-ring, whence this special case yields the result for idempotent complete $\C$; the general case then follows from an application of the cofinality theorem to the right and an explicit analysis on the left (see e.g.\ \cite[Lemma 7.5]{HLS} and the surrounding discussion).
\item Raptis recently gave a general version of Waldhausen's sphere theorem, that simultaneously applies to the case of an ordinary additive category, see \cite[Theorem 6.11]{Raptisdevissage}. Since the $\K$-theory of every stable $\infty$-category can also be written as the $\K$-theory of a Waldhausen category (by the strictification results mentioned above) one can use this result to approach \ref{cor:fontes} as briefly indicated in \cite[Example 6.15]{Raptisdevissage}.
\item To the best of our knowledge the full statement of \ref{cor:fontes} first appeared in print as the main result of \cite{Fontes2018weight}. Fontes' proof strategy closely follows that in the case where $\mathcal A$ is an ordinary additive category, but using Barwick's Waldhausen $\infty$-categories to avoid the need for strictification, but many of the arguments in \cite{Fontes2018weight} lack details. More importantly, we also feel obliged to point out that the applications Fontes gives in \cite[Section 6.2 \& 6.3]{Fontes2018weight} are flawed: In \cite[Corollary 6.5]{Fontes2018weight} he claims to deduce the full Gillet-Waldhausen theorem, which concerns the relation between $\K$-spectra of general (i.e. non-split) exact categories $\mathcal E$ and their Waldhausen-categories of finite chain complexes $\CCbdd(\mathcal E)$. However, his main result (i.e.\ \ref{cor:fontes} above) only compares the $\K$-theory of $\CCbdd_\infty(\mathcal E)$ with that of $\mathcal E$, where the latter is regarded with its split-exact structure and not the given one. Since $\K$-spectra are generally sensitive to changes in the exact structure this invalidates \cite[Corollary 6.3]{Fontes2018weight}. The deduction of Quillen's resolution theorem suffers the same defect, see \cite[Corollary 6.7]{Fontes2018weight}.\footnote{During the revision process for the present paper, complete proofs of these statements (appropriately corrected), in particular of Corollary \ref{cor:fontes}, appeared in work of Saunier \cite{Saunier}, who introduced heart structures on stable categories, that generalise weight structures to allow for non-split exact structures on hearts.}
\item Another place where the full statement of \ref{cor:fontes} recently appeared is \cite[Theorem 1.35]{heleodoro}. However, in \cite[Definition 1.29]{heleodoro} Heleodoro crucially (mis)defines a weight structure on $\C$, so that for stable $\C$ the spectra $\hom_\C(x,y)$ are required to be coconnected if $x \in \C_{[-\infty,0]}$ and $y \in \C_{[1,\infty]}$ (whereas \ref{def:weightstr} requires them to be connected). His argument makes explicit use of this faulty requirement, which in fact forces $\C \simeq 0$. He has informed us that an erratum is forthcoming.
\item Finally, let us warn the reader that \cite[Corollary 4.1]{Sosniloweight}, the statement of which reads $\K(\C^\heart) \simeq \K(\C)$, while looking deceptively similar is in fact not directly related to \ref{cor:fontes}, because Sosnilo simply defines the $\K$-spectrum of an additive category $\mathcal A$ as $\K(\Stab(\mathcal A))$. This statement is therefore an immediate consequence of the equivalence $\Stab(\C^\heart) \simeq \C$, which he provides in \cite[Section 3]{Sosniloweight} (and which forms the basis for \ref{prop:poinconstab} above).
\end{enumerate}
\end{remark}

We now turn to hermitian examples. 

\begin{corollary}\label{moduleresolution}
If $A$ is a connective $\Eone$-ring, $(M,N,\alpha)$ an invertible module with genuine involution over $A$ and $c \subseteq \K_0(A)$ a subgroup, such that 
\begin{enumerate}
\item[i)] both $M$ and $\hom_A(M,A)$ are connective, 
\item[ii)] $N$ is connective,
\item[iii)] $c$ is closed under the involution on $\K_0(A)$ induced by $M$, 
and 
\item[iv)] $[A] \in c$.
\end{enumerate}
Then the natural map
\[\Poinc^\heart(\Mod^c_A,\QF_M^\alpha)^\grp \longrightarrow \gw(\Mod^c_A,\QF_M^\alpha)\]
is an equivalence.
\end{corollary}

Following the discussion in \ref{examplesofweight} (2) the fourth assumption can be relaxed; it is only used to guarantee that the weight structure on $\Modp_A$ restricts to $\Mod^c_A$.

\begin{proof}
The third condition ensures that $(\Mod^c_A,\QF_M^\alpha)$ is indeed a Poincar\'e category and the fourth that connective $A$-modules take part in a weight structure on $\Mod^c_A$, see \ref{examplesofweight} (2). By the characterisation of $(\Mod^c_A)^\heart$, condition i) implies that the weight structure has exact dimension $0$, and via the fibre sequence
\[\hom_A(X \otimes X, M)_\hC \longrightarrow \QF_M^\alpha(X) \longrightarrow \hom_A(X,N)\]
from the construction of $\QF_M^\alpha$ the second condition ensures that it takes connective values on the heart.
\end{proof}

The component of some $(P,q)$ in $\Poinc(\Mod^c_A,\QF_M^\alpha)$ is the classifying space of the $\QF_M^\alpha$-orthogonal $\Eone$-group $\mathrm{O}^\alpha(P,q)$, i.e. it is the component of the quotient $\Omega^\infty\QF_M^\alpha(P) \sslash \mathrm{Gl}(P)$
containing $q$. The $\Eone$-group $\mathrm{O}^\alpha(P,q)$ itself then lies in the fibre sequence 
\[\mathrm{O}^\alpha(P,q) \longrightarrow \mathrm{Gl}(P) \xrightarrow{\mathrm{act}_q} \Omega^\infty\QF_M^\alpha(P)\]
exhibiting it as the stabiliser of $q$. Abbreviating $\mathrm{O}^\alpha(\hyp(A^g))$ by $\mathrm{O}_{g,g}^\alpha(A)$ we deduce:

\begin{proposition}\label{homologygw}
If $A$ is a connective $\Eone$-ring, $(M,N,\alpha)$ an invertible module with genuine involution over $A$, such that 
\begin{enumerate}
\item[i)] both $M$ and $\hom_A(M,A)$ are connective, and 
\item[ii)] $N$ is connected.
\end{enumerate}
Then the natural map
\[\colim_{g \in \mathbb N} \mathrm H_*(\mathrm{BO}_{g,g}^\alpha(A);\mathbb Z) \longrightarrow \mathrm H_*(\gw(\Modp_A,\QF_M^\alpha)_0;\mathbb Z)\]
is an isomorphism.
\end{proposition}

For the proof we need:

\begin{lemma}\label{met=hypinheart}
If $(\C,\QF)$ is a Poincar\'e category of dimension at most
$0$, and $\Lin_\QF$ takes connected values on the heart of $\C$, then any $(P,q) \in \Poinc(\C, \QF)$ admitting a connective Lagrangian $L$  
is equivalent to $\hyp(L)$.
\end{lemma}

For $\C = \Dperf(R)$ this is \cite[1.2.11]{9authIII}, and the argument given there works in general. We repeat it for the sake of completeness.

\begin{proof}
Ranicki's algebraic Thom construction, see \cite[2.3.20]{9authI}, provides us with an equivalence
\[\Poinc(\Met(\C, \Q))\simeq \grpcr\catforms(\C, \QF\qshift{-1})\]
which sends a Lagrangian $L\to P$ to a form on $\Dual_\QF(L)^{[-1]}$, in a way so the canonical Lagrangian $L\to \hyp(L)$  corresponds to the $0$-form on $\Dual_\QF(L)^{[-1]}$. 

Thus it suffices to show that, under our hypotheses, $\pi_0  \QF\qshift{-1}(\Dual_\QF(L)^{[-1]})=0$, but this is immediate from the fibre sequence 
\[\Bil_\QF(\Dual_\QF L^{[-1]},\Dual_\QF L^{[-1]})_\hC \longrightarrow \QF((\Dual_\QF L)^{[-1]}) \longrightarrow \Lin_\QF(\Dual_\QF L^{[-1]})\]
since $\Dual_\QF(L)$ is coconnective.
\end{proof}

\begin{proof}[Proof of Proposition \ref{homologygw}]
The claim will follow from Corollary \ref{moduleresolution} and the group completion theorem \cite{mcduff-segal,nik-grp} once we establish that inverting the class of $\hyp(A)$ in the commutative monoid $\pi_0\Poinc^\heart(\Modp_A,\QF_M^\alpha)$ yields a group.

But for every Poincar\'e form $(P,q)$ the form $(P,q) \oplus (P,-q)$ admits the diagonal of $P$ as a Lagrangian, see \cite[Construction 2.3.8]{9authII}. If $P \in (\Modp_A)^\heart$ this forces $(P,q) \oplus (P,-q) \simeq \hyp(P)$ by Lemma \ref{met=hypinheart}. But by the characterisation of the heart in \ref{examplesofweight} (2), $\hyp(P)$ is a retract of $\hyp(A^g)$ for some $g$, in total making every object of $\Poinc^\heart(\Modp_A,\QF^\alpha_M)$ a direct summand in $\hyp(A^g)$ for some $g \in \mathbb N$, which clearly suffices.
\end{proof}

We can now specialise this result to discrete rings, in particular proving Theorem \ref{thm:main_special}.

\begin{corollary}\label{main_special2}
Let $R$ be a ring and 
\[\fpm = (M_\Ct \xrightarrow{\tau} Q \xrightarrow{\rho} M^\Ct)\]
be an invertible form parameter with associated Poincar\'e structure $\QF_M^{\g\fpm}$ on $\Dperf(R)$. Then the natural map
\[\gw^\fpm_{\mathrm{cl}}(R;M) = \Unimod^\fpm(R;M)^\grp \longrightarrow \gw(\Dperf(R),\QF_M^{\g\fpm})\]
is an equivalence. If $\tau$ is in addition surjective we find that the natural map
\[\colim_{g \in \mathbb N} \mathrm H_*(\mathrm{O}^\fpm_{g,g}(R);\mathbb Z) \longrightarrow \mathrm H_*(\gw(\Dperf(R),\QF^{\g\fpm}_M)_0;\mathbb Z)\]
an isomorphism.
\end{corollary}

Here the left hand side denotes the group homology of the discrete group $\mathrm{O}^\fpm_{g,g}(R)$. Unwinding definitions, it is the group of $R$-linear automorphisms of $\hyp(R)^g=R^g\oplus M^g$, where $M$ is considered as a $R$-module via the inclusion $R \rightarrow R \otimes_\mathbb Z R$ into the second factor, that preserve the $\lambda$-form
\begin{align*}
 R^g\oplus M^g \times R^g\oplus M^g \to M, \quad & ((r,m), (s,n))\mapsto r\cdot n+\sigma(s\cdot m),\\
 R^g\oplus M^g\to Q, \quad & (r,m) \mapsto \tau(r\cdot m)
\end{align*}
with $r\cdot m= \sum_i (r_i\otimes 1) m_i$ and $\sigma$ denoting the involution on $M$; recall also that $M$ is assumed projective over $R$. With this result established we shall follow the conventions of \cite{9authIII} and denote the common value of the two spaces under consideration by $\gw^{\g\fpm}(R;M)$.

\begin{proof} Per construction we have 
\[\Poinc^\heart(\Dperf(R),\QF_M^{\g\fpm}) \simeq \Unimod^\fpm(R;M).\]
To apply \ref{moduleresolution} and \ref{homologygw} recall only that $\Lin_{\QF^{\g\fpm}_M} \simeq \hom_A(-,N)$, where the underlying spectrum of $N$ is $\cof(\tau \colon M_\hC \rightarrow Q)$. This is clearly always connective and even connected if $\tau$ is surjective.
\end{proof}

Form parameters satisfying the surjectivity condition of this corollary are precisely those describing quadratic modules in Bak's language \cite{Bak}. Examples are $\fpm \in \{\pm \q,\pm \ev\}$, but usually not $\fpm = \pm \s$. 
The second part of Corollary \ref{main_special2} has an analogue whenever there is a unimodular $\lambda$-form $(Q,q)$ such that every unimodular form is a direct summand of some $(Q,q)^g$. Note that $\hyp(\mathbb Z)\in \Unimod^s(\mathbb Z)$ does indeed not have this property since no odd form can appear as a direct summand of an even form. Forms as required can, however, often be found even through partial classification results. 

To illustrate this point recall that indefinite unimodular symmetric forms over $\mathbb Z$ are classified by their parity, rank and signature, see \cite[Section 7]{Serre}, so every element of $\Unimod^\s(\mathbb Z)$ is a direct summand of some
\[\left(\mathbb Z^2, \begin{pmatrix} 1 & 0 \\ 0 & - 1 \end{pmatrix}\right)^g \cong \left(\mathbb Z^{2g}, \begin{pmatrix} \mathrm{I} & 0 \\ 0 & - \mathrm{I} \end{pmatrix}\right)\]
whose isometry group we denote $\mathrm{O}_{\langle g,g\rangle}(\mathbb Z)$.
The group completion theorem combined with \ref{main_special2} thus yields
\[\H_*(\gw^{\g\s}(\mathbb Z)_0; \mathbb Z) = \colim_{g \in \mathbb N} \H_*(\mathrm{O}_{\langle g,g\rangle}(\mathbb Z);\mathbb Z).\]
To finish the illustration let us use this identification to derive the stable group homology modulo $2$ of the groups $\mathrm{O}_{\langle g,g\rangle}(\mathbb Z)$; more elaborate calculations, based on the same method, will appear in joint work with Land and Nikolaus \cite{homologyo}. To state the outcome, note that there are maps
\[\gw^{\g\s}(\mathbb Z) \longrightarrow \gw^{\g\s}(\mathbb F_3) \quad \text{and} \quad \gw^{\g\s}(\mathbb Z) \longrightarrow \gw^\mathrm{top}(\mathbb R)\]
induced by the evident ring homomorphisms; here the topological version $\gw^\mathrm{top}(\mathbb R)$ can for example be constructed as the group completion of $\Unimod^{\mathrm{top}}(\mathbb R)$, the $\Einf$-space underlying the topologically enriched groupoid of unimodular forms (via the usual topology on matrix groups over $\mathbb R$), or by realising the simplicial space
\[n \longmapsto \gw^{\g\s}(\mathrm{C}(\Delta^n,\mathbb R)).\]
In either case the result is equivalent to $\mathrm{ko} \times \mathrm{ko}$ by extracting the positive and negative definite parts; using the former definition this follows from the group completion theorem, and using the latter it follows for example from \cite[Proposition 10.2]{schlichting-derived} and the discussion thereafter. On the other hand, Fiedorowicz, Friedlander and Priddy used Quillen's techniques to produce an equivalence
\[\gw^{\g\s}(\mathbb F_3)_0 \simeq \fib(\psi^3 - \id \colon \mathrm{BO} \rightarrow \mathrm{BSO}),\]
where $\psi^3$ denotes the third Adams operation, see \cite[Corollary 1.6]{friedlander-computations} or \cite[Theorem III.3.1 (d)]{FP-classical}. 

\begin{theorem}\label{hlgyo}
Base changing to $\mathbb F_3$ and extracting the positive part over $\mathbb R$ yields a $2$-adic equivalence
\[\gw^{\g\s}(\mathbb Z)_0 \longrightarrow \fib(\psi^3 - \id \colon \mathrm{BO} \rightarrow \mathrm{BSO}) \times \mathrm{BO}.\]
In particular, we obtain isomorphisms
\[\colim_{g \in \mathbb N} \mathrm{H}_*(\mathrm{O}_{\langle g,g\rangle}(\mathbb Z);\mathbb F_2)\longrightarrow \mathbb F_2[e_i, f_i | i \geq 1] \otimes \Lambda[b_i| i \geq 1]\]
and
\[\mathbb F_2[v_i,w_i,a_i | i \geq 1] \longrightarrow \lim_{g \in \mathbb N\op}\mathrm{H} ^*(\mathrm{O}_{\langle g,g\rangle}(\mathbb Z);\mathbb F_2)\]
where the $e_i,f_i,v_i, w_i$ and $b_i$ have degree $i$, whereas the $a_i$ have degree $2i-1$. Both are compatible with products (i.e. the block-sum product in the former and the cup product in the latter case).
\end{theorem}

The (co)homology isomorphisms in the second part can also directly be described using the homomorphisms\[\mathrm{O}_{\langle g,g\rangle}(\mathbb Z) \longrightarrow \mathrm{O}_{\langle g,g\rangle}(\mathbb F_3) \quad \text{and} \quad \mathrm{O}_{\langle g,g\rangle}(\mathbb Z) \longrightarrow \mathrm{O}^{\mathrm{top}}_{\langle g,g\rangle}(\mathbb R) \simeq \mathrm{O}(g) \times \mathrm{O}(g)\]
induced by the evident ring maps (and finally extracting the positive part in the second case). In this language Fiedorowicz and Priddy construct classes in the (co)homology of $\mathrm{O}_{\langle g,g\rangle}(\mathbb F_3)$ that can serve as the $f_i,v_i, b_i$ and $a_i$ purely group theoretically in \cite[Section IV.2]{FP-classical}; the $w_i$ are simply the Stiefel-Whitney classes of the positive definite part, and the $e_i$ given as the images of 
\[\mathbb F_2 = \mathrm{H}_i(\mathrm{BO}(1);\mathbb F_2) \longrightarrow \mathrm{H}_i(\mathrm{BO}(g);\mathbb F_2).\]

\begin{proof}
The proof most readily available from the literature has two inputs. On the one hand the analogous statement for $\mathbb Z[1/2]$ in place of $\mathbb Z$ is essentially due to Berrick and Karoubi: Their main result in \cite{berrick-karoubi} provides a $2$-adically cartesian square 
\[\xymatrix{ \gw^\gs(\mathbb Z[1/2])_0 \ar[r]\ar[d] & \gw^{\mathrm{top}}(\mathbb R)_0 \ar[d] \\
             \gw^\gs(\mathbb F_3)_0 \ar[r]^{\mathrm{Br}} & \gw^{\mathrm{top}}(\mathbb C)_0,}\]
             whose lower horizontal map is produced by Brauer lifting. Now, recalling $\mathrm{O}_{g,g}(\mathbb C) \simeq \mathrm{O}(2g)$ the group completion theorem implies $\gw^{\mathrm{top}}(\mathbb C) \simeq \mathrm{ko}$, and the right vertical map identifies with the addition $\mathrm{ko} \times \mathrm{ko} \rightarrow \mathrm{ko}$ restricted to the unit component. In particular, it is split by the inclusion into either factor. Combining this with the identification of the lower left term mentioned above gives that
\[\gw^{\g\s}(\mathbb Z[1/2])_0 \longrightarrow \fib(\psi^3 - \id \colon \mathrm{BO} \rightarrow \mathrm{BSO}) \times \mathrm{BO}\]
is a $2$-adic equivalence. The second ingredient is the fact that the map
\[\gw^{\g\s}(\mathbb Z)_0 \longrightarrow \gw^{\g\s}(\mathbb Z[1/2])_0\]
is a $2$-adic equivalence, which was conjectured by Berrick and Karoubi and is proven in \cite[3.1.10]{9authIII}.
 
A somewhat more direct argument, which simultaneously handles the level of path-components, was given by Nikolaus in \cite[Section 15]{nik-lec}, based on the equivalence $\gw^{\g\s}(\mathbb Z)_2^\wedge \simeq (\k(\mathbb Z)_2^\wedge)^\hC$ established in \cite[3.1.6]{9authIII}. \\

The second statement has one further input, namely the computation of the (co)homology of the fibre term by Fiedorowicz and Priddy in \cite[Section I.3]{FP-classical} (though of course this already enters implicitly into their identification of $\gw^{\g\s}(\mathbb F_3)$). They observe that the Serre spectral sequences of the fibre sequence
\[\mathrm{SO} \longrightarrow \fib(\psi^3 - \id \colon \mathrm{BO} \rightarrow \mathrm{BSO}) \longrightarrow \mathrm{BO}\]
with coefficients in $\mathbb F_2$ collapse and the (co)homology is free enough to give
\[\mathrm{H}_*(\fib(\psi^3 - \id \colon \mathrm{BO} \rightarrow \mathrm{BSO});\mathbb F_2) \cong \mathrm{H}_*(\mathrm{BO};\mathbb F_2) \otimes \mathrm{H}_*(\mathrm{SO};\mathbb F_2) \cong \mathbb F_2[f_i |i \in \mathbb N] \otimes \Lambda[b_i | i \in \mathbb N]\]
and 
\[\mathrm{H}^*(\fib(\psi^3 - \id \colon \mathrm{BO} \rightarrow \mathrm{BSO});\mathbb F_2) \cong \mathrm{H}^*(\mathrm{BO};\mathbb F_2) \otimes \mathrm{H}^*(\mathrm{SO};\mathbb F_2) \cong \mathbb F_2[v_i |i \in \mathbb N] \otimes \mathbb F_2[a_i | i \in \mathbb N].\]
\end{proof}

\begin{remark}\label{homstabsymmrmk}

Contrary to a claim in an earlier version of this paper, homological stability results for the groups $\mathrm{O}_{\langle g,g\rangle}(\mathbb Z)$ do not seem known, despite identifications
\[\left(\mathbb Z^{2},\begin{pmatrix} 1 & 0 \\ 0&  -1\end{pmatrix}\right)^g \cong \left(\mathbb Z^{2},\begin{pmatrix} 1 & 0 \\ 0&  -1\end{pmatrix}\right) \oplus \left(\mathbb Z^2, \begin{pmatrix} 0 & 1 \\ 1 & 0\end{pmatrix}\right)^{g-1}\]
which for example follow from the classification of indefinite forms over the integers in terms of rank, signature and parity, or directly from 
\[\begin{pmatrix} 1 & 1& 1 \\ -1 & 0 & -1 \\ 0& 1 & 1\end{pmatrix} \in \mathrm{Gl}_3(\mathbb Z) \quad \text{conjugating} \quad \begin{pmatrix} 1 & 0 & 0 \\ 0 & 0 & 1 \\ 0 & 1 & 0\end{pmatrix} \quad \text{into} \quad \begin{pmatrix} 1& 0 & 0 \\ 0 & 1 & 0 \\ 0 & 0 & -1 \end{pmatrix}.\]
The most generally applicable stabilisation result seems to be \cite[Theorem B]{nina}, where Friedrich shows that stabilisation by $\mathrm{hyp}(\mathbb Z)$ induces isomorphisms 
\[\mathrm{O}^\mathrm{q}(M,q) \rightarrow \mathrm{O}^\mathrm{q}((M,q) \oplus \mathrm{hyp}(\mathbb Z))\]
for an arbitrary quadratic form $q$ on an abelian group $M$ in a certain range of degrees depending of $(M,q)$. And while generally $\mathrm{O}^\mathrm{q}(P,q) = \mathrm{O}^\mathrm{s}(P,q)$ for quadratic forms over the integers, Friedrich's arguments are not applicable for two reasons:
\begin{enumerate}
\item The forms under consideration are odd, so do not admit quadratic refinements at all but their existence is needed in the proof of the crucial \cite[Lemma 3.5]{nina}.
\item It is unclear whether stabilisation by the standard odd form and by $\mathrm{hyp}(\mathbb Z)$ induce the same map on homology under any identification of forms as above, but it is the former which occur in the colimit of Theorem \ref{hlgyo}.
\end{enumerate}
Given the ubiquity of homological stability in the cohomology of arithmetic groups we find it hard to imagine that the required stability statement for the groups $\mathrm{O}_{\langle g,g\rangle}(\mathbb Z)$ does not hold; but as the standard techniques for addressing such questions are quite different from the methods of this paper, we shall not pursue this point here.
\end{remark}

\subsection{Applications to the classical $\Q$-construction}\label{subsec:applQ}

As already mentioned in the introduction, \ref{cor:res2} contains the strongest assertion of \ref{thm:main2} also for a Poincar\'e category $(\C,\QF)$ of arbitrary even dimension $d = 2\oc$ by considering $(\C^{\langle \oc \rangle},\QF)$. In this section we will explore the applications of \ref{thm:main2} in odd dimensions, where again it suffices to treat a fixed case, say the case $(\C,\QF\qshift{-1})$ where $(\C,\QF)$ has dimension $0$. In this case \ref{thm:main2} implies that the natural map
\[|\Cob^{0,0}(\C,\QF\qshift{-1})| \longrightarrow |\Cob(\C,\QF\qshift{-1})|\]
is an equivalence, provided $\Lin_\QF$, or equivalently $\QF$, takes connective values on $\C^\heart$. We start by analysing the left-hand side.

\begin{lemma}\label{lem:cob00}
Let $(\C, \QF)$ is a Poincar\'e category of dimension $0$. Then, a morphism $c\leftarrow w\to d$ of $\Cob(\C, \QF\qshift{-1})$ belongs to $\Cob^{0,0}(\C, \QF\qshift{-1})$ if and only if $c,d,w \in \C^\heart$ and additionally   $w \rightarrow c$ is the projection to a direct summand or equivalently $w \rightarrow d$ is the inclusion of a direct summand.
\end{lemma}

\begin{proof}
By \ref{leftarrightar} $w \rightarrow d$ is $0$-coconnective if and only if $w \rightarrow c$ is $0$-connective. If furthermore source and target lie in $\C^\heart$ then connectivity is equivalent to being the projection to a direct summand by \ref{heartsplit} and dually for coconnective maps. Together these facts imply the statement by inspecting definitions. 
\end{proof}

\begin{lemma}\label{lem:Q_is_1_category}
Let $(\C, \QF)$ is a Poincar\'e category of dimension $0$ and suppose further that 
$\C^\heart$ is an ordinary category, and that 
$\QF$ takes $1$-truncated values on $\C^\heart$. Then, $\Cob^{0,0}(\C, \QF\qshift{-1})$ is also an ordinary category.
\end{lemma}

\begin{proof}
Recall that the mapping spaces in the category associated to a Segal space $X$ are described by the fibres of the map $(d_1,d_0)\colon X_1\to X_0\times X_0$. Denote then by 
\[\Poinc(\Q_1(\C, \QF))^{0,0}\subseteq \Poinc(\Q_1(\C, \QF))\]
the components of cobordisms that lie in $\Cob^{0,0}(\C, \QF\qshift{-1})$. Using the identification of these in the previous lemma it suffices to show that the top horizontal map in the diagram
\[\xymatrix{\Poinc(\Q_1(\C, \QF))^{0,0}\ar[rr]^-{(d_1,d_0)} \ar[d]^{\fgt} && \Poinc^\heart(\C, \QF)^2 \ar[d]^\fgt \\
\core(\Q_1(\C))^{0,0}\ar[rr]^-{(d_1,d_0)} && \core(\C^\heart)^2}\]
has discrete fibres. To see this, we compare them to the fibres of the lower horizontal map, whose source consists of those path components in $\core\Q_1(\C)$ spanned by objects $c \leftarrow w \rightarrow d$ with all three of $c,w,d$ lying in $\C^\heart$, $w \rightarrow c$ the projection to a direct summand and $w \rightarrow d$ the inclusion of one. The assumption that $\core(\C^\heart)$ is discrete makes the lower horizontal map a faithful functor between ordinary groupoids, whence it has discrete fibres. 

The fibres of the comparison map (\emph{i.e.,} the total homotopy fibre of the square) can also be computed by first taking vertical fibres and so are given by the full subspace of Poincar\'e forms in the limit of
\[\{q\} \rightarrow \Omega^\infty \QF(w) \leftarrow \{q'\},\]
where $q$ and $q'$ are Poincar\'e forms on $c$ and $d$, respectively. This limit is either empty or equivalent to $\Omega^{\infty+1}\QF(w)$ and therefore discrete as well. 
\end{proof}

Let us now restrict to the case where $(\C,\QF)$ has dimension $0$ and $\QF$ takes discrete values on $\C^\heart$. In this situation, the same is true for the bilinear part of $\QF$, and since by the dimension assumption the duality $\Dual_\QF$ on $\C$ restricts to a duality on $\C^\heart$, we find that $\C^\heart$ automatically has discrete mapping spaces, \emph{i.e.}\ it is an ordinary category. The most important example is the one from Section~\ref{sec:genuine_GW}, namely $\C=\Modp_R$ for a discrete commutative ring $R$ equipped with its canonical weight structure and the Poincar\'e structure $\QF^{\g\fpm}$ for $\fpm\in \{\s, \sks, \q, \skq\}$; in this case $\C^\heart$ is the category of finitely generated projective $R$-modules, equipped with its usual duality. We find:

\begin{corollary}
Let $(\C, \QF)$ is a Poincar\'e category of dimension $0$ and suppose further that 
$\QF$ takes discrete values on $\C^\heart$. Then, $\Cob^{0,0}(\C, \QF\qshift{-1})$ is also an ordinary category.
\end{corollary}

With these assumptions in place $\QF$ restricts to a functor $\QF\colon \C^\heart\to \Ab$, which is additively quadratic (see the discussion before \ref{prop:poinconstab}) and whose bilinear part is corepresented (in $\C^\heart$) by a duality. Thus the pair $(\C^\heart, \QF\vert_{\C^\heart})$ canonically gives rise to an additive form category $(\C^\heart,(\Dual_\QF)_{|\C^\heart},\mathrm{can},\QF_{|\C^\heart})$ as considered by Schlichting in \cite{Schlichting2019}; in fact, by Yoneda's lemma, the additional data required for a form category in \cite[Definition 2.1]{Schlichting2019} is always uniquely determined (up to canonical equivalence) by its quadratic functor, compare \cite[Remark 2.4]{Schlichting2019}.

We will now set out to identify $\Cob^{0,0}(\C, \QF\qshift{-1})$ with a version of the classical hermitian $\Q$-construction of $\C^\heart$. 
To facilitate the discussion note that the assumptions on a fixed $(\C,\QF)$ are strict enough to give a model of $\Poinc^\heart(\C,\QF)$ as an ordinary groupoid, which we shall denote by $\Poinc(\C^\heart,\QF)$, that is well-defined up to unique isomorphism (and not just equivalence): Its objects are pairs $(c,q)$ with $c \in \mathrm{Ho}(\C^\heart)$ and $q \in \pi_0\Omega^\infty\QF(c)$, and its morphisms $(c,q) \rightarrow (c',q')$ are isomorphisms $f \colon c \rightarrow c'$ in $\mathrm{Ho}(\C^\heart)$ with $f^*q' = q$.

\begin{definition}\label{def:classhermq}
For $(\C,\QF)$ a Poincar\'e category of dimension $0$, such that $\QF_{|\C^\heart}$ factors through $\mathrm{H} \colon \Ab \rightarrow \Spa$, let $\Qhclass(\C^\heart,\QF)$ be the category whose objects are those of $\Poinc(\C^\heart,\QF)$ and whose morphisms $\xymatrix{(c,q) \ar@{~>}[r]&(c', q')}$ are isomorphism classes of spans
\[c\xleftarrow{p} w \xrightarrow{i} c'\]
in $\mathrm{Ho}(\C^\heart)$, with $p$ the projection to a direct summand and $i$ the inclusion of one, and that additionally satisfy the following two conditions:
\begin{enumerate}
 \item\label{item:same_form}the restrictions of $q$ and $q'$ to $w$ agree, and
 \item\label{item:Poincare_condition} $(q')_\sharp \circ i\colon w\to \Dual_{\QF}(c')$ induces an isomorphism $\ker(p)\to \ker(\Dual_{\QF}(i))$.
\end{enumerate}
Composition is by forming pullbacks of spans.
\end{definition}

\begin{remark}
\begin{enumerate}
\item In \cite{Schlichting2019} Schlichting in fact works with the subdivided $\mathcal S$- rather than the $\Q$-construction. The translation is that the former is equivalent to the Rezk-nerve of the latter, as explained in \cite[Appendix B.1]{9authII} in the setting of stable $\infty$-categories with duality; using Lemma \ref{lem:cob_vs_clasical_Q} below, the proof given there works almost verbatim to compare the definition above with \cite[Definition 6.3]{Schlichting2019}. For the reader tempted to fill in the remaining details, let us mention that the only additional information required is that Schlichting's quadratic functor on $\mathcal S_{2k+1}(\C^\heart)$ agrees with the restriction of our $\QF^{\Twar[k]} \colon \Q_k(\C)\op \rightarrow \Spa$ to diagrams in $\C^\heart$. In the notation of \cite[Definition 6.1]{Schlichting2019} this boils down to the transformation $\varphi \colon A \rightarrow A^\#$ being uniquely determined by the family $\xi_i \in \mathcal Q(A_i)$ in case $A \colon \Ar[n] \rightarrow \mathcal A$ is an element of $\mathcal{S}_n(\mathcal A)$. Indeed, the equation $\rho(\xi_i) = \varphi_{i'} \circ A_{i \leq i'}$ determines $\varphi_i$ for all $i \in \mathcal P$ with $i=i'$ and for $\mathcal P = \Ar[n]$ this is the case for $i = (0 \leq n)$. But all components of $\varphi$ are subquotients of $\varphi_{(0 \leq n)}$ by naturality, so uniquely determined as well.
\item In case the natural map $\pi_0\QF(X) \rightarrow (\pi_0\Bil_\QF(X,X))^\Ct$ is an isomorphism for each $X \in \C^\heart$, the exact form category $(\C^\heart,(\Dual_\QF)_{|\C^\heart},\mathrm{can},\QF_{|\C^\heart})$ is associated to its underlying exact category with duality as in \cite[Example 2.24]{Schlichting2019} and the above definition reduces to Schlichting's implementation of the hermitian $\Q$-construction for exact categories with duality in \cite[Definition 4.1]{Schlichting2010}. This is the case for the Poincar\'e structures $\QF_M^{\g\s} \colon \Dperf(R)\op \rightarrow \Spa$ and their generalisations appearing in Section \ref{sec:krcomp} below.
\item We shall refrain from giving a full discussion of the functoriality of the $\Q$-construction described above, as we do not need it. Let us, however, warn the reader that our differentiation between $\C^\heart$ and $\mathrm{Ho}(\C^\heart)$ is not meant to suggest that $\Qhcl$ can be made into a functor with target the ordinary category of ordinary categories, but rather only to indicate the place, where we access $\C^\heart$ up to isomorphism, rather than equivalence. 

It is, however, easy to check that $\Qhcl$ gives a functor into the homotopy category of ordinary categories, and, indeed Definition \ref{def:classhermq} can be canonically upgraded to a functor from the homotopy bicategory of Poincar\'e categories with a weight structure satisfying its assumptions to the bicategory of ordinary categories. With that extension in place the maps occuring in the next results are each part of a natural transformation of bicategories; of course can also simply use \ref{lem:cob_vs_clasical_Q} below to \emph{define} such functoriality for $\Qhcl$.
\end{enumerate}
\end{remark}

\begin{lemma}\label{lem:cob_vs_clasical_Q}
For $(\C,\QF)$ a Poincar\'e category of dimension $0$, such that  $\QF_{|\C^\heart}$ factors through $\mathrm{H} \colon \Ab \rightarrow \Spa$, there is a canonical equivalence 
\[\Cob^{0,0}(\C, \QF\qshift{-1}) \simeq \Qhclass(\C^\heart, \QF).\]
\end{lemma}

\begin{proof}
Since $\core\Cob^{0,0}(\C,\QF\qshift{-1}) \simeq \Poinc^\heart(\C,\QF)$ by \ref{lem:cob00}, we may assume that the objects of $\Cob^{0,0}(\C,\QF)$ are given by those of $\Poinc(\C^\heart,\QF)$. In this case, we claim that there is an isomorphism
\[r\colon \Ho(\Cob^{0,0}(\C, \QF\qshift{-1})) \to \Qhclass(\C^\heart, \QF)\]
of ordinary categories induced by the identity on objects and by sending an equivalence class of Poincar\'e cobordisms to the equivalence class of its underlying span in $\mathrm{Ho}(\C^\heart)$.

Indeed, the underlying span of any Poincar\'e cobordism
\[c \xleftarrow{p} w \xrightarrow{i} c'\]
satisfies condition (\ref{item:same_form}) by definition, and the Poincar\'e condition says that 
\[(q')_\# \colon \fib(p) \longrightarrow \Dual_{\QF\qshift{-1}}\fib(i)\]
is an equivalence (where the fibres are formed in the surrounding category $\C$). But if $p$ is the projection to a direct summand in $\C^\heart$ we have $\fib(p) \cong \ker(p)$ and similarly $\Dual_{\QF\qshift{-1}}(\fib(i)) \cong \ker(\Dual_\QF(i))$ if $i$ is the inclusion of a direct summand (where the kernels are formed in $\mathrm{Ho}(\C^\heart)$).

Thus a Poincar\'e cobordism also satisfies (\ref{item:Poincare_condition}). We finally check that $r$ is fully faithful: Since a span satisfying the two conditions from \ref{def:classhermq} can be refined to a Poincar\'e cobordism by the computations above it is certainly surjective on morphism sets, and since $\QF$ evaluates to discrete spectra on the heart of $\C$ an equivalence of Poincar\'e cobordisms from $(X,q)$ to $(X',q')$ is nothing but an equivalence of its underlying spans, which shows that it is injective as well.
\end{proof}

Thus, summarising our discussion, we obtain:

\begin{theorem}\label{thm:Q=Cob}
Let $(\C, \QF)$ be Poincar\'e category of dimension $0$ such that $\QF$ takes values in discrete spectra. Then, there is a canonical inclusion
\[\Qhclass(\C^\heart, \QF)\to \Cob(\C, \QF\qshift{-1})\]
as a subcategory, which induces an equivalence on realisations. 
\end{theorem}

Our first application again pertains to plain algebraic $\K$-theory. Suppose $\C$ is a stable $\infty$-category equipped with an exhaustive weight structure, such that $\C^\heart$ is an ordinary category. Then $\Hyp(\C)$ satisfies the assumptions of Theorem \ref{thm:Q=Cob} and by inspection $\Qhclass(\Hyp(\C)) \simeq \Qclass(\C^\heart)$ is simply Quillen's classical $\Q$-construction of $\C^\heart$ regarded as a split-exact category (i.e. it is the wide subcategory of $\Span(\C^\heart)$, where the left pointing leg in a morphism is required to be the projection to a direct summand and the right pointing leg is the inclusion of one). 

But $\Hyp(\C)\qshift{-1} \simeq \Hyp(\C)$ via the functor 
\[\C \times \C\op \longrightarrow \C \times \C\op, \quad (X,Y) \longmapsto (X^{[1]},Y),\]
so from $\Cob(\Hyp(\C)) \simeq \Span(\C)$ and \ref{cor:res2} we find:

\begin{corollary}\label{q=+stable}
If $\C$ is a stable $\infty$-category equipped with an exhaustive weight structure such that $\C^\heart$ is an ordinary category, then both inclusions
\[\Qclass(\C^\heart) \longrightarrow \Span(\C) \longleftarrow \mathbb B\core(\C^\heart)\]
induce equivalences upon realisation. 
\end{corollary}

As explained in Section \ref{subsec:weightstable} an additive category occurs as the heart of a weight structure on a stable $\infty$-category if and only if it is $\flat$-additive, \emph{i.e.} if every map in $\mathcal A$ that admits a retraction is in fact the inclusion of a direct summand. We thus recover:

\begin{corollary}[Quillen]\label{cor:+=qwic}
If $\mathcal A$ is a $\flat$-additive, ordinary category, there is a canonical equivalence
\[\Omega|\Qclass(\mathcal A)| \simeq \core(\mathcal A)^\grp,\]
 i.e. the direct sum and split-exact $\K$-spaces of $\mathcal A$ agree.
\end{corollary}

\begin{remark}
\begin{enumerate}
\item Let us emphasise, that the restriction to ordinary categories in \ref{cor:+=qwic} is solely for ease of exposition. For any $\flat$-additive $\infty$-category $\mathcal A$ it is easily checked that $\Cob^{0,0}(\Hyp(\Stab(\mathcal A))\qshift{-1})$ agrees with Barwick's extension of Quillen's $\Q$-construction from \cite{BarwickQ} applied to $\mathcal A$ (again regarded with its split exact structure). In Barwick's notation the argument above then applies verbatim to give
\[\Omega|\Q(\mathcal A)| \simeq \core(\mathcal A)^\grp;\]
we warn the reader explicitly of the clash with our use of $\Q(\C)$ for stable $\C$.
\item For the reader interested in unwinding our proof for this special case, let us remark that the filtration $\Cob^{\mc,\oc}(\Hyp(\C)\qshift{d})$
which we employed for $d=0,-1$ to obtain \ref{q=+stable}, translates to the filtration on $\Span(\C)$, where an object $X$ is assumed both $\oc$-connective and $d+1-\oc$-coconnective, and where a morphism $c \leftarrow w \rightarrow d$ has left pointing arrow $\mc$-connective and right pointing arrow $d+1-\mc$-coconnective.

Furthermore, a surgery datum on an object $X$ of $\Poinc(\Hyp(C)) \simeq \core(\C)$ is simply a null-composite sequence $T \rightarrow X \rightarrow S$ and the trace of a surgery translates to
\[X \longleftarrow \fib(X \rightarrow S) \longrightarrow \fib(X \rightarrow S)/T.\]
Our arguments therefore systematically alter objects and morphisms in $\Span(\C)$ by embedding them into appropriately chosen null-composite sequences and then replacing them with the `homology' thereof.
\item The statement of \ref{cor:+=qwic} remains valid upon dropping the requirement that $\mathcal A$ be weakly idempotent complete, see e.g.\ \cite[Theorem IV.7.1]{Kbook}. In fact both sides are invariant under forming weak idempotent closures; this follows from the group completion theorem for the right hand side and is contained in the cofinality theorem, see e.g.\ \cite[Theorem IV.8.9]{Kbook}, for the left hand side. Thus the restriction to $\flat$-additive categories loses no material content. 

We also suspect that our methods can be applied in a slightly more finessed manner, to obtain the full statement, but since we largely included \ref{cor:+=qwic} for illustrative purposes only, let us refrain from attempting this here. 

A similar comment applies to the next corollary.
\end{enumerate}
\end{remark}

We finally combine the statements proven so far to obtain what is colloquially known as Giffen's theorem:

\begin{corollary}[Schlichting]\label{cor:Qfib}
Let $(\C, \QF)$ be a Poincar\'e category of dimension $0$, such that $\QF$ takes discrete values on $\C^\heart$. Then there is a canonical fibre sequence of spaces
\[\Poinc(\C^\heart,\QF)^\grp \to |\Qhclass(\C^\heart, \QF)| \xrightarrow{\fgt} |\Qclass(\C^\heart)|.\]
\end{corollary}

Using \ref{prop:poinconstab} this result can indeed be applied to Schlichting's additive form categories $\mathcal A$ (at least the weakly idempotent complete ones, but one can again check that all terms are invariant under weak idempotent completion) by extending a suitable quadratic functor $\mathcal A\op \rightarrow \Ab$ to a Poincar\'e structure $\Stab(\mathcal A)\op \rightarrow \Spa$. 

After an erroneous proof in \cite{CL} and a partial fix in \cite{schlichting-giffen}, even the case of ordinary categories with duality was only solved in full by Hesselholt and Madsen in \cite{HM}, and the general case is the main result of \cite{Schlichting2019}. Classically, the fibre of the right hand map would serve as the foundation for developing hermitian $\K$-theory in a framework of ordinary categories as implemented by Schlichting in \cite{Schlichting2010} if $2$ is assumed invertible, and initiated in \cite{HM, Schlichting2019} in general. However, the greater flexibility of Poincar\'e categories (or for that matter any derived set-up such as those in \cite{schlichting-derived}, if $2$ is assumed invertible, or \cite{Spitzweckreal}) allows for a direct definition of Grothendieck-Witt spaces and spectra without recourse to the non-hermitian case (since these make it possible to shift the duality). The sequence above then appears a posteriori in the guise of what we have termed the Bott-Genauer sequence
\[|\Cob(\C,\QF\qshift{-1})| \xrightarrow{\fgt} |\Span(\C)| \xrightarrow{\hyp} |\Cob(\C,\QF)|\]
due to its resemblance to Genauer's fibre sequence
\[|\Cob_{d+1}| \longrightarrow |\Cob_{d+1}^\partial| \longrightarrow |\Cob_d|,\]
from geometric topology, see \cite[Corollary 2.5.2]{9authII}. 

\begin{proof}
Rotating the Bott-Genauer sequence once to the left gives a fibre sequence 
\[\Omega|\Cob(\C,\QF)| \longrightarrow |\Cob(\C,\QF\qshift{-1})| \xrightarrow{\fgt} |\Span(\C)|\]
which Corollary \ref{cor:res2}, Theorem \ref{thm:Q=Cob} and Corollary \ref{q=+stable} identify with the sequence from the statement.
\end{proof}

\begin{remark}
\begin{enumerate}
\item Just as with \ref{cor:+=qwic} the assumption in \ref{cor:Qfib} that $\QF$ take discrete values on $\C^\heart$ (which in particular forces $\C^\heart$ to be an ordinary category) can be relaxed to $\QF$ taking connective values on $\C^\heart$ at the cost of treating $\Cob^{0,0}(\C,\QF\qshift{-1})$ as a $\infty$-categorical version of the hermitian $\Q$-construction of $(\C^\heart,\QF)$.
\item Interestingly, the faulty argument for \ref{cor:Qfib} in \cite{CL} can sometimes be salvaged in the non-additive setting of proto-exact categories: There it can happen that every split mono- or epimorphism splits \emph{uniquely} (a split exact category with this property automatically vanishes) and this allows the argument from \cite{CL} to go through, see \cite{f1stuff}. This result was used to analyse $\K$- and $\GW$-spectra of monoid schemes in \cite{f2stuff}, in particular, showing that this non-additive version of the theory behaves quite differently from anything that appears in the additive/stable set-up.
\item We hope to combine, in future work, our parametrised algebraic surgery with its geometric counterpart to better understand the relationship between (suitably stabilised) diffeomorphism groups and cobordism categories in odd dimensions. As mentioned in the introduction the program of Galatius and Randal-Williams meets an obstruction in this regime, as does the algebraic version, and the work in \cite{HP} suggests a close connection between the two. 
\end{enumerate}
\end{remark}

\subsection{Applications to L-groups}\label{lsection2}

In this section we explain the implication of Theorem \ref{thm:main2} at the level of path-components, where they give results on $\L$-groups. To this end recall that the main result of \cite{9authII} provides in particular a fibre sequence
\[\k(\C,\QF)_\hC \xrightarrow{\hyp} \gw(\C,\QF) \xrightarrow{\mathrm{bord}} \l(\C,\QF),\]
where $\l(\C,\QF)$ is defined as the realisation of the simplicial object formed by the Poincar\'e $n$-ads, see \cite[Section 4.4]{9authII}. To mine this statement for information on $\gw^{\g\s}(R;M)$, say, it remains to understand the right hand term. As part of the main results of \cite{9authIII}, we achieved this by identifying the homotopy groups of $\l(\Dperf(R),\QF^{\g\s}_M)$ with the bordism groups of short symmetric chain complexes introduced by Ranicki, and similar results hold for general form parameters. The purpose of this final section is to explain, that these identifications are also contained in Theorem \ref{thm:main2}, in fact, already in the discussion in Section \ref{subsec:surgery_on_component} (which in fact do not rely on our parametrised version of surgery, as explained at the beginning of that section). 

The starting point for our approach is that, by design, the homotopy groups of $\l(\C,\QF)$ are always the cobordism groups of the shifts of $(\C,\QF)$, more precisely
\[\L_n(\C,\QF) \cong \pi_0|\Cob(\C,\QF\qshift{-n-1})|,\]
see \cite[Proposition 2.3.7 \& Corollary 4.4.6]{9authII}. We then put
\[\L^{\mc,\oc}_n(\C,\QF) = \pi_0|\Cob^{\mc,\oc}(\C,\QF\qshift{-n-1})|;\]
let us warn the reader immediately that this notation clashes with the corresponding one from \cite{9authIII}, where similar superscripts are used to indicate different bounds; unfortunately we did not find a good solution for this clash of notation and hope no confusion will arise.

\begin{theorem}\label{resolutioninL}
Let $(\C,\QF)$ be a Poincar\'e category of dimension at least $d$ and assume 
\begin{enumerate}
\item $2\oc \leq d-n$, 
\item $\mc+\oc\leq d-n$  and
\item $\Lin_\QF(Y)$ is $\oc+n$-connective for every $Y \in \C^\heart$.
\end{enumerate}
Then the map
\[\L^{\mc,\oc}_n(\C,\QF) \longrightarrow \L_n(\C,\QF)\]
is an isomorphism. 
\end{theorem}

Such connectivity estimates on cycles for $\L$-groups are of course by no means new (save for the generality afforded by our framework). In fact, the primordial case is Ranicki's description of quadratic $\L$-groups, originally defined by Wall using only unimodular quadratic forms, as cobordism groups of quadratic Poincar\'e chain complexes in \cite{rsurg2, rsurg1}. Since $\Lin_\QF \simeq 0$ for a quadratic Poincar\'e structure $\QF$ his results are special cases of Theorem \ref{resolutioninL}: Using the canonical weight structure on $\Dperf(R)$ of dimension $0$, the case of the even quadratic $\L$-groups $\L^\mathrm q_{2k}(R,M)$ follows via $\oc = -k = \mc$ (compare also Corollary \ref{Lheart=L} below) and in the case of the odd $\L$-groups $\L_{2k+1}^\mathrm{q}(R,M)$ the reduction to Ranicki's unimodular formations is covered by taking $\oc = -k-1$ and $\mc = -k$; we shall refrain from spelling this out further, since we have nothing to add to Ranicki's treatment.

\begin{proof}
From Theorem \ref{thm:main2} we immediately learn that the claim is true for $\mc+\oc=d-n$ and $\oc \leq \mc +1$ (and this is in fact all we will need below); alternatively, the claim is also immediate from the far simpler \ref{resolutionobjectspi0} and \ref{resolutionmorphismspi0}. For smaller $\mc$ consider the factorisation 
\[\pi_0|\Cob^{d-n-\oc,\oc}(\C,\QF\qshift{-n-1})| \longrightarrow \pi_0|\Cob^{\mc,\oc}(\C,\QF\qshift{-n-1})| \longrightarrow \pi_0|\Cob(\C,\QF\qshift{-n-1})|\]
whose composite is an isomorphism. But the former map is clearly surjective so all maps are bijective.
\end{proof}

\begin{corollary}\label{prop:short}
If $(\C,\QF)$ is a Poincar\'e category of exact dimension $0$, and $\Lin_\QF(Y)$ is connective for every $Y$ in $\C^\heart$, then the map
\[\L_n^\short(\C,\QF) \rightarrow \L_n(\C,\QF)\]
is an isomorphism for all $n \geq 0$.
\end{corollary}

Here, we denote by $\L_n^\short(\C,\QF)$ the Poincar\'e bordism group spanned by short Poincar\'e objects of $(\C,\QF\qshift{-n})$, i.e.\ those that lie in $\C_{[-n,0]}$, modulo short cobordisms $C \leftarrow W \rightarrow D$, i.e.\ with $W \in \C_{[-n-1,0]}$. We leave the simple check that it really is an equivalence relation to the reader. Note that this is only a sensible definition if $\QF$ has exact dimension $0$ and $n \geq 0$. 

For $(\C,\QF) = (\Dperf(R),\QF^{\g\s}_M)$ this statement is precisely \cite[1.2.7]{9authIII} (beware of the notation difference between the papers here), and the proof given there also applies essentially verbatim in the present situation. We can now simply appeal to \ref{resolutioninL}, however.

\begin{proof}
Apply the theorem with $\mc=-n-1$ and $\oc=-n$ to obtain that 
\[\L_n^{-n-1,-n}(\C,\QF) \longrightarrow \L_n(\C,\QF)\]
is an isomorphism. The left hand side has the correct cycles, so it remains to analyse the induced relation. But any morphism $C \leftarrow W \rightarrow D$ in $\Cob(\C,\QF\qshift{-n})$ between objects concentrated in $\C_{[-n,0]}$ satisfies $W \in \C_{[-n-1,0]}$ if and only if one (and then both) of its constituent maps are $(-n-1)$-connective.
\end{proof}

Short $\L$-groups can be compared for different Poincar\'e structures much more easily. We first used this in \cite{9authIII} to analyse the relation between $\QF^{\g\s}$ and $\QF^{\s}$. Essentially the same method as employed there proves:

\begin{proposition}\label{prop:shortcomparison}
Let $(\C,\QF)$ be a Poincar\'e category of dimension $0$ and $\eta\colon \QF \Rightarrow \Phi$ a transformation to another Poincar\'e structure on $\C$. If for some $n \geq 0$
\begin{enumerate}
\item the induced transformation $\Bil_\QF \Rightarrow \Bil_\Phi$ is an equivalence, and
\item the induced map $\Omega^\infty \Lin_\QF(X) \rightarrow \Omega^\infty\Lin_\Phi(X)$ is $(n+1)$-connective for all $X \in \C^\heart$,
\end{enumerate}
then the map
\[\eta_* \colon \L^\short_k(\C,\QF) \longrightarrow \L^\short_k(\C,\Phi)\]
is an isomorphism for all $0 \leq k \leq n$ and surjective for $k=n+1$.
\end{proposition}

In particular, we recover \cite[1.2.18]{9authIII}, \emph{i.e.} that there is an isomorphism
\[\L^{\g\s}_k(R;M) \longrightarrow \L^\short_k(R;M)\]
for all $k \geq 0$ by combining the statement with \ref{prop:short}.

Note also that if we a priori assume $(\C,\Phi)$ to also have dimension $0$, the two assumptions can be combined into the statement that 
\[\Omega^\infty\QF(X) \rightarrow \Omega^\infty\Phi(X)\]
 is $(n+1)$-connective for all $X \in \C^\heart$: This implies that also $\Omega^\infty\Bil_\QF(X,Y) \rightarrow \Omega^\infty \Bil_\Phi(X,Y)$ is $(n+1)$-connective, so in particular
\[\pi_0\Hom_{\C^\heart}(X,\Dual_\QF Y) \longrightarrow \pi_0\Hom_{\C^\heart}(X,\Dual_\Phi Y)\]
is an isomorphism. But by Yoneda's lemma this implies that $\eta_*\colon \Dual_\QF Y \rightarrow \Dual_\Phi Y$ is an equivalence for all $Y \in \C^\heart$. But the subcategory of all $\C$ spanned by all objects $Y$ for which this holds is stable, so has to be all of $\C$, as the weight structure is exhaustive. This proves condition (1) above, and the consequent equivalence to condition (2) follows from the decomposition of a quadratic functor into its linear and bilinear parts.
\begin{proof}
Using this translation it is easy to check that the map $\L_0^\short(\C,\QF) \longrightarrow \L_0^\short(\C,\Phi)$ is surjective even under the weaker assumption of (1) together with $\pi_0\Lin_\QF(X) \rightarrow \pi_0\Lin_\Phi(X)$ being surjective for all $X \in \C^\heart$. Using this fact it clearly suffices to establish surjectivity in degree $n+1$ and injectivity in degree $n$.

Given assumption (1) it suffices for surjectivity to show that the map
\[\pi_0\Omega^{\infty+n+1}\QF(C) \longrightarrow \pi_0\Omega^{\infty+n+1}\Phi(C)\]
is surjective for all $C \in \C_{[-n-1,0]}$, or equivalently that the boundary homomorphism
\[\pi_{n+1}\Phi(C) \longrightarrow \pi_n\fib(\QF(C) \rightarrow \Phi(C))\]
vanishes. But by another application of (1) the fibre of $\QF(C) \rightarrow \Phi(C)$ agrees with that of $\Lin_\QF(C) \rightarrow \Lin_\Phi(C)$, so we will be done if we can show that the analogous statement for the linear parts. Now if $C \in \C^\heart$ the infinite loop space of this fibre is assumed $(n+1)$-connective. This clearly implies that the target of the boundary map vanishes not only for such $C$ but also for all $C\qshift{-k}$ with $0 \leq k \leq n$. Since this vanishing statement is closed under extensions it follows for all $C \in \C_{[-n,0]}$. For arbitrary $C \in \C_{[-n-1,0]}$ pick a weight decomposition $X \rightarrow C \rightarrow Y$ with $X \in \C_{[-n-1,-n-1]}$ and $Y \in \C_{[-n,0]}$ using \ref{connfib}. It induces a diagram
\[\xymatrix{\pi_{n+1} \Lin_\Phi(Y) \ar[d] \ar[r] & \pi_{n+1} \Lin_\Phi(C) \ar[d] \ar[r] & \pi_{n+1} \Lin_\Phi(X) \ar[d] \\
            \pi_n\fib(\Lin_\QF(Y) \rightarrow \Lin_\Phi(Y)) \ar[r] & \pi_n\fib(\Lin_\QF(C) \rightarrow \Lin_\Phi(C))\ar[r] &\pi_n\fib(\Lin_\QF(X) \rightarrow \Lin_\Phi(X))}\]
            with exact rows. But its lower left hand corner vanishes by the previous argument, so we are reduced to showing that the right vertical map vanishes. But its kernel is the image of 
            \[\pi_{n+1}\Lin_\QF(X) =\pi_0\Lin_\QF(X\qshift{n+1}) \longrightarrow \pi_0\Lin_\Phi(X\qshift{n+1}) =\pi_{n+1}\Lin_\Phi(X)\]
            which is surjective by assumption.

Injectivity in degree $n$ is similarly implied by the map
\[\pi_n(\QF(X) \times_{\QF(W)} \QF(Y)) \longrightarrow \pi_n(\QF(X) \times_{\Phi(W)} \QF(Y))\]
being surjective, whenever $X \leftarrow W \rightarrow Y$ is a diagram in $\C$ with $X,Y \in \C_{[-n,0]}$ and $W \in \C_{[-n-1,0]}$. From the associated long exact sequence one finds this equivalent to the map 
\[\pi_n(\QF(X) \times_{\Phi(W)} \QF(Y)) \longrightarrow \pi_n(\fib(\QF(W) \rightarrow \Phi(W)))\]
vanishing. This follows by a very similar argument: The target agrees with the fibre on linear parts, so it suffices to treat those. Then choose a decomposition $U \rightarrow W \rightarrow V$ with $U \in \C_{[-n-1,-n-1]}$ and $V \in \C_{[-n,0]}$. The maps $W \rightarrow X$ and $W \rightarrow Y$ both factor through $V$ (since $\pi_0\Hom_\C(U,X) = 0$ straight from the axioms of a weight structure). Choosing such factorisations we find a diagram
\[\xymatrix{\pi_{n} \Lin_\Phi\QF(X) \times_{\Lin_\Phi(V)} \Lin_\QF(Y) \ar[d] \ar[r] & \pi_{n} \Lin_\QF(X) \times_{\Lin_\Phi(W)} \Lin_\QF(Y) \ar[d] \ar[r] & \pi_{n}(\Lin_\Phi(U)\qshift{-1}) \ar[d] \\
            \pi_n\fib(\Lin_\QF(V) \rightarrow \Lin_\Phi(V)) \ar[r] & \pi_n\fib(\Lin_\QF(W) \rightarrow \Lin_\Phi(W))\ar[r] &\pi_n\fib(\Lin_\QF(U) \rightarrow \Lin_\Phi(U))}\]
            with exact rows. By the same arguments as before, the lower left corner vanishes, as does the right vertical map, and together this gives the claim.
\end{proof}

We can also give a mild simplification and generalisation of \cite[Proposition 1.2.15]{9authIII}.

\begin{corollary}\label{Lheart=L}
If $(\C,\QF)$ is a Poincar\'e category of exact dimension $0$, and $\Lin_\QF(Y)$ is connective for every $Y$ in $\C^\heart$, then the evident map
\[\L_0^\heart(\C,\QF) \longrightarrow \L_0(\C,\QF)\]
is an isomorphism.
\end{corollary}

Here, we denote by $\L_0^\heart(\C,\QF)$ the abelian group obtained by taking the monoid of Poincar\'e objects in the heart of $\C$ as cycles, and dividing them by the congruence relation generated by declaring all strictly metabolic forms equivalent to $0$, where a strictly metabolic form $M$ is one that admits a Lagrangian which also lies in the heart of $\C$; such a Lagrangian is nothing but a morphism $\xymatrix{0 \ar@{~>}[r]& M}$ in $\Cob^{0,0}(\C,\QF)$. Note that for $\QF$ the animation of a classical form parameter $\fpm$ over some ring $R$ this is nothing but the ordinary Witt group $\mathrm{W}^\fpm(R,M)$ of $\fpm$, see \cite[Definition 2.27]{Schlichting2019}, so we indeed obtain:

\begin{corollary}
Let $R$ be a ring and
\[\fpm = (M_\Ct \xrightarrow{\tau} Q \xrightarrow{\rho} M^\Ct)\]
be an invertible form parameter with associated Poincar\'e structure $\QF_M^{\g\fpm}$ on $\Dperf(R)$. Then the tautological map
\[\mathrm{W}^\fpm(R,M) \longrightarrow \L_0(\Dperf(R),\QF_M^{\g\fpm})\]
is an isomorphism.
\end{corollary}

\begin{proof}[Proof of Corollary \ref{Lheart=L}]
Since a strict Lagrangian $D$ in some Poincar\'e object $C$ in $\C^\heart$ defines a morphism $C \leftarrow D \rightarrow 0$ in $\Cob^{0,0}(\C,\QF\qshift{-1})$, there is an obvious factorisation 
\[\L_0^\heart(\C,\QF) \longrightarrow \L_0^{0,0}(\C,\QF) \longrightarrow \L_0(\C,\QF).\]

We will argue that both maps are bijective. 
For the second map use \ref{resolutioninL} with $\oc=\mc=d=n=0$. The first map is evidently surjective, so it remains to show that it is also injective. To see this we need to show that the existence of a morphism $(C,q)\cobmor (C', q)$ in $\Cob^{0,0}(\C, \QF)$ implies $[(C,q)]=[(C', q')]\in \L_0^\heart(\C, \QF)$. But any such morphism defines a strict Langrangian in $(C, q)\oplus (C', -q')$ and hence $[(C, q)]=-[(C', -q')]=[(C', q')]$ (applying the same reasoning to the identity morphism of $(C', q')$ for the second equality).
\end{proof}

Finally, to obtain the first corollary of Theorem \ref{thm:main_special} in the introduction, we quote from \cite[1.3.7]{9authIII} that the comparison map
\[\L_n^\short(\Dperf(R),\QF^{\s}_M) \longrightarrow \L_n(\Dperf(R),\QF^\s_M)\]
is an isomorphism if $n \geq d-1$ for every noetherian ring $R$ of global dimension $d$. This is achieved by a different kind of surgery to which the set-up of this paper does not apply. The statement in particular implies that
\[\l^{\g\s}(R;M) \longrightarrow \l^{\s}(R;M)\]
is an equivalence if $R$ is a Dedekind domain. Combining this with Corollary \ref{main_special2} and the fibre sequence in \cite[Corollary 4.4.13]{9authII}, we finally find:

\begin{corollary}\label{cor:2corintro}
For any Dedekind ring $R$ and an invertible module with involution $M$ over $R$ there is a fibre sequence
\[\k(R;M)_\hC \xrightarrow{\hyp} \gw^{\s}_\mathrm{cl}(R;M) \xrightarrow{\mathrm{bord}} \l^{\s}(R;M),\]
that canonically splits after inverting $2$. 
\end{corollary}

The final statement advertised in the introduction, namely that the forgetful map
\[\fgt \colon \gw^{\s}_\mathrm{cl}(R;M) \longrightarrow \k(R;M)^\hC\]
is a $2$-adic equivalence, whenever the fraction field of $R$ is in addition a number field, follows by combining \ref{main_special2} with \cite[Theorem 3.1.6]{9authIII}, which proves the statement for $\gw(\Dperf(R),\QF_M^{\g\s})$. Let us also mention that for such $R$ and $M=R$ with trivial involution \cite[Corollary 2.2.4]{9authIII} expresses the homotopy groups of $\l^{\s}(R)$ entirely in terms of the Picard group, the symmetric Witt-group, and the number of dyadic primes in $R$.

\subsection{Applications to real $\K$-spectra}\label{sec:krcomp}

Corollary \ref{cor:res2} also allows us to compare the real $\K$-theory spectrum constructed by Heine, Lopez-Avila and Spitzweck in \cite{Spitzweckgroupcompletion} with that constructed in \cite{9authII}, and thus in particular, compute the geometric fixed points of the former.

To set the stage, recall that \cite[Section 7.4]{9authI} constructs for every Poincar\'e category $(\C,\QF)$ a Mackey-object $\mathrm{gHyp}(\C,\QF)$ in $\Catp$, \emph{i.e.} a product preserving functor $\mathrm{Span}({\operatorname{\Ct-\mathrm{Fin}}}) \rightarrow \Catp$ from the span category of finite $\Ct$-sets. Its values on the $\Ct$-sets $\Ct$ and $\ast$ are given by $\Hyp(\C)$ and $(\C,\QF)$, respectively. Taking Grothendieck-Witt spectra then results in a genuine $\Ct$-spectrum, \emph{i.e.} a product preserving functor $\Span({\operatorname{\Ct-\mathrm{Fin}}}) \rightarrow \Spa$, the real algebraic $\K$-spectrum $\KR(\C,\QF)$, per construction its underlying spectrum (the value at $\Ct$) is $\GW(\Hyp(\C)) \simeq \K(\C)$ and its genuine fixed points (the value at $\ast$) are $\GW(\C,\QF)$.

Similarly, \cite[Section 4]{Spitzweckgroupcompletion} constructs for every additive category with duality $(\mathcal A,\Dual) \in \Cat^\hC$ (where the $\Ct$-action on $\Cat$ is by opponing) a product preserving functor $\Span({\operatorname{\Ct-\mathrm{Fin}}}) \rightarrow \mathrm{Mon}_{\E_\infty}(\Sps)$ whose values at $\Ct$ and $\ast$ are $\core(\mathcal A)$ and $\core(\mathcal A)^\hC$, respectively, where the $\Ct$-action on $\core(\mathcal A)$ is via the duality $\Dual$. Following this functor by group completion and delooping
\[\mathrm{Mon}_{\E_\infty}(\Sps) \longrightarrow \mathrm{Grp}_{\E_\infty}(\Sps) \longrightarrow \Spa\]
defines the genuine $\Ct$-spectrum $\kr(\mathcal A, \Dual)$.

By the results of \cite[Section 7.1]{9authI} or \cite[Section 5]{Spitzweckgroupcompletion} the duality $\Dual$ is classified by a symmetric functor $\Bil_\Dual \colon \mathcal A\op \times \mathcal A\op \rightarrow \mathrm{Grp}_{\E_\infty}(\Sps)$ with
\[\Bil_\Dual(x,y) \simeq \Hom_\mathcal A(x,\Dual y).\] 
After weak idempotent completion, the pair $(\mathcal A, (\Bil_\Dual \circ \Delta)^\hC)$, where the homotopy fixed points are formed in $\mathrm{Grp}_{\E_\infty}(\Sps)$ and then delooped (instead of first delooping and then forming them in $\Spa$), is easily checked to lie in $\Catpad$, which receives an equivalence $(-)^\heart \colon \Catpw \rightarrow \Catpad$ by Proposition \ref{prop:poinconstab}. We shall denote a preimage by $(\Stab(\mathcal A), \QF_\Dual^{\gs})$, as this fits with the example $\mathcal A = \PProj(R)$. In particular, $\QF_\Dual^{\gs}$ is usually different from the symmetric Poincar\'e structure $\QF_\Dual^\s$ associated with its duality.

The goal of the present section is to show:
\begin{proposition}\label{prop:krcomparison}
There is a canonical comparison map $\kr(\mathcal A,\Dual) \rightarrow \KR(\Stab(\mathcal A),\QF_\Dual^{\gs})$
which induces equivalences 
\[\kr(\mathcal A,\Dual)^{\g\Ct} \rightarrow \tau_{\geq 0}\GW(\Stab(\mathcal A),\QF_\Dual^{\gs}) \quad \text{and} \quad \kr(\mathcal A,\Dual)^{\varphi\Ct} \rightarrow \tau_{\geq 0}\L(\Stab(\mathcal A),\QF_\Dual^{\gs})\]
and on underlying spectra. In particular, 
\[\pi_*(\kr(\mathcal A,\Dual)^{\varphi\Ct}) \cong \L^\short_*(\Stab(\mathcal A),\QF_\Dual^{\s}).\]
\end{proposition}
We start out by constructing the comparison map. To this end we shall need to recall the definition of the two genuine $\Ct$-spectra in question in greater detail, and this in turn requires the language of $\Ct$-categories. Recall that a $\Ct$-category is a functor $\D \colon \mathcal O(\Ct)\op \rightarrow \Cat$ with $\mathcal O(\Ct)$ the orbit category of $\Ct$. This unwinds to a category $\D(\Ct)$ with $\Ct$-action, together with a functor $\D(\ast) \rightarrow \D(\Ct)^\hC$. The construction of $\KR(\C,\QF)$ begins with the observation that the functor $\Catp \rightarrow (\Catx)^\hC$ extracting the duality, where $\Catx$ carries the opponing action, defines such a $\Ct$-category, which is semi-additive; for a discussion of that notion we refer the reader to \cite[Section 5]{nardin} or \cite[Section 7.4]{9authI}.

Now by \cite[Proposition 5.11 \& Theorem 6.5]{nardin}, there is a semi-additive $\Ct$-category $\mathcal Span(\Ct)$ (there called $\underline{\bf{\mathrm A}}^{\mskip-2mu\textit{eff}}(\Ct)$) with
\[(\mathcal Span(\Ct))(\ast) \simeq \Span({\operatorname{\Ct-\mathrm{Fin}}}) \quad \text{and} \quad (\mathcal Span(\Ct))(\Ct) \simeq \Span(\mathrm{Fin}),\]
such that evaluation at $\ast \in \Span(\Ct)$ defines an equivalence
\[\Hom^{{\operatorname{\Ct-\times}}}(\mathcal Span(\Ct),\D) \simeq \core(\D(\ast))\]
whenever $\D$ is semi-additive; here the left hand side denotes the space of $\Ct$-functors (i.e. natural transformations between the defining diagrams) commuting with finite $\Ct$-products (in fact in \cite[Theorem 6.5]{nardin} the existence of all finite $\Ct$-limits in $\D$ is assumed, but the proof only makes use of finite $\Ct$-products). The action of $\Ct$ on $\Span(\mathrm{Fin})$ that is part of the definition of this $\Ct$-category  is trivial and the requisite functor
\[\Span(\Fun(\mathrm B\Ct,\mathrm{Fin})) = \Span({\operatorname{\Ct-\mathrm{Fin}}}) \longrightarrow \Span(\mathrm{Fin})^\hC = \Fun(\mathrm B\Ct,\Span(\mathrm{Fin}))\]
is easily unwound to be the fully faithful inclusion used in \cite[Section 4]{Spitzweckgroupcompletion} for the construction of $\kr(\mathcal A,\Dual)$.

The Mackey object $\mathrm{gHyp}(\C,\QF) \colon \Span({\operatorname{\Ct-\mathrm{Fin}}}) \rightarrow \Catp$, and consequently $\KR(\C,\QF)$ is now simply defined by taking the preimage of $(\C,\QF) \in \Catp$ under Nardin's equivalence and then evaluating the resulting $\Ct$-functor at $\ast \in \mathcal O(\Ct)$.

On the other hand, every functor $\D \colon \mathrm{B}\Ct \rightarrow \Cat$ gives rise to a $\Ct$-category by right Kan extension along the fully faithful inclusion $\mathrm B\Ct \rightarrow \mathcal O(\Ct)\op$. The value at $\ast \in \mathcal O(\Ct)$ is then simply $\D^\hC$ with structure map the identity. Starting with a semi-additive category, this again results in a semi-additive $\Ct$-category. The construction of the Mackey object group completing to $\kr(\mathcal A,\Dual)$ can be interpreted in terms of this construction as follows: Starting with the category $\Cata$, which is semi-additive, and the action of $\Ct$ via opponing, we obtain from Nardin's equivalence above for every $(\mathcal A,\Dual) \in (\Cata)^\hC$ a functor $\mathrm{ghyp}(\mathcal A,\Dual) \colon \Span({\operatorname{\Ct-\mathrm{Fin}}}) \rightarrow (\Cata)^\hC$. Applying $\core(-)^\hC \colon (\Cata)^\hC \rightarrow \mathrm{Mon}_{\Einf}(\Sps)$ then recovers the Mackey object from \cite{Spitzweckgroupcompletion}. 

To obtain the desired comparison map, we now observe that also $\Catpw \rightarrow (\Catw)^\hC$ gives rise to a semi-additive $\Ct$-category which comes equipped with maps to both the $\Ct$-category arising from $\Catp \rightarrow (\Catx)^\hC$ by forgetting weight structures, and to that arising from the identity of $(\Cata)^\hC$ by applying the forgetful functor $\Catp \rightarrow (\Catx)^\hC$ and then taking hearts. Applying Nardin's result a third time, we obtain from $(\Stab(\mathcal A),\QF^\gs_\Dual) \in \Catpw$ a Mackey-object $\Span({\operatorname{\Ct-\mathrm{Fin}}}) \rightarrow \Catpw$. Per construction it is underlain by $\mathrm{gHyp}(\Stab(\mathcal A),\QF^\gs_\Dual)$ and on the heart we find $\mathrm{ghyp}(\mathcal A,\Dual)$. We thus obtain maps
\[\Poinc(\mathrm{gHyp}(\Stab(\mathcal A),\QF^\gs_\Dual)) \Longleftarrow
 \Poinc^\heart(\mathrm{gHyp}(\Stab(\mathcal A),\QF^\gs_\Dual)) \Longrightarrow \core(\mathrm{ghyp}(\mathcal A,\Dual))^\hC.\]
But the right hand transformation is an equivalence when evaluated on both $\ast$ and $\Ct \in \Span({\operatorname{\Ct-\mathrm{Fin}}})$, and since they generate $\Span({\operatorname{\Ct-\mathrm{Fin}}})$ under products, it is an equivalence. Since Grothendieck-Witt spaces are grouplike there is then a unique way to complete the arising square
\[\xymatrix{\core(\mathrm{ghyp}(\mathcal A,\Dual))^\hC \ar[r] \ar[d] & \Poinc(\mathrm{gHyp}(\Stab(\mathcal A),\QF^\gs_\Dual)) \ar[d] \\
(\core(\mathrm{ghyp}(\mathcal A,\Dual))^\hC)^\grp \ar@{-->}[r] & \gw(\mathrm{gHyp}(\Stab(\mathcal A),\QF^\gs_\Dual)), }\]
and using $\gw \simeq \Omega^\infty \GW$ the lower transformation finally adjoins to the desired map
\[\kr(\mathcal A,\Dual) \longrightarrow \KR(\Stab(\mathcal A),\QF^\gs_\Dual),\]
whenever $\mathcal A$ is weakly idempotent complete. But the left hand side is invariant under weak idempotent completion by the group completion theorem, so the general case follows.

\begin{proof}[Proof of Proposition \ref{prop:krcomparison}]
Unwinding definitions, we find the map on underlying spectra to induce the canonical map
\[\core(\mathcal A)^\grp \rightarrow \k(\Stab(\mathcal A))\]
on the associated infinite loopspaces. It is an equivalences by \ref{cor:fontes}. Since the underlying spectra of both $\kr(\mathcal A,\Dual)$ and $\KR(\Stab(\mathcal A),\QF_\Dual^\gs)$ are connective, this proves the third claim. Similarly, the induced map on the infinite loop spaces of the genuine fixed points is
\[(\core(\mathcal A)^\hC)^\grp \rightarrow \gw(\Stab(\mathcal A),\QF_\Dual^\gs)\]
which is an equivalence by \ref{cor:res2}. Per construction the genuine fixed points of $\kr(\mathcal A,\Dual)$ are connective, so we obtain the first claim. The second one, follows by comparing the long exact sequences of the fibre sequences
\[\xymatrix{\kr(\mathcal A,\Dual)_\hC \ar[r] \ar[d] & \kr(\mathcal A,\Dual)^{\mathrm g\Ct} \ar[r] \ar[d] &\kr(\mathcal A,\Dual)^{\varphi\Ct} \ar[d] \\
\KR(\Stab(\mathcal A),\QF_\Dual^\gs)_\hC \ar[r] & \KR(\Stab(\mathcal A),\QF_\Dual^\gs)^{\mathrm g\Ct} \ar[r] & \KR(\Stab(\mathcal A),\QF_\Dual^\gs)^{\varphi\Ct}}\]
together with the equivalence of the lower right corner to $\L(\Stab(\mathcal A),\QF_\Dual^\gs)$ from \cite[Corollary 4.6.2]{9authII}.

The final claim about the homotopy groups of the geometric fixed points now follows from \ref{prop:short} and \ref{prop:shortcomparison}.
\end{proof}

\appendix
\section{The weight theorem for $\L$-spaces (by Y.\ Harpaz)}

The goal of the present appendix is twofold: In the first part we provide an $\L$-theoretic counterpart to Theorem \ref{thm:resolution}, by comparing the $\L$-space of a $0$-dimensional Poincar\'e category $(\C,\QF)$ with that of its heart (suitably defined), see Theorem \ref{prop:comp-L}. The method is the parametrised algebraic surgery developed in the body of the paper, though it can be significantly simplified in the case at hand.

In the second part we then explain how this result can be used to approach Theorem \ref{thm:resolution} itself. More precisely, we show how Theorem \ref{prop:comp-L}, together with a certain additivity property of the direct sum Grothendieck-Witt-functor $(\C,\QF) \mapsto \Poinc^\heart(\C,\QF)^\grp$, reduces Theorem \ref{thm:resolution} to the case of ordinary $\K$-theory, i.e., to the statement that $\core(\C^\heart)^\grp \simeq \k(\C)$
is an equivalence. Restricting to the case where $\C^\heart$ is an ordinary additive category and $\QF$ takes discrete values on the heart, the requisite additivity statement was established by Schlichting in \cite{Schlichting2019} and the $\K$-theoretic statement is essentially the classical theorem of Gillet and Waldhausen. Thus our discussion in particular provides an alternative proof of Theorem \ref{thm:main_special}, arguably the most important special case, that is logically independent from the body of the paper, see Corollary \ref{cor:bordinary}.

Given the additivity results from \cite{9authII} (which pertain to the functor $\gw$) Theorem \ref{thm:resolution} in fact implies the requisite additivity statement for direct sum $\GW$-spaces in full generality. An independent verification of this fact could thus be used to provide a proof of Theorem \ref{thm:resolution} in full, though we do not pursue this here. Note also, that the proof provided in the body of the text simply derives the $\K$-theoretic statement as a special case of the general machinery and requires no a priori development of direct sum $\K$- or $\GW$-spaces, both features not shared by the strategy we present below.

\subsection{Parametrised algebraic surgery on $\L$-spaces}

Recall from Proposition~\ref{prop:poinconstab} above that extracting the heart induces an equivalence from the category \(\Catpw\) of $0$-dimensional Poincaré categories (and Poincaré functors of degree $0$ between them), and a certain category \(\Catpad\) of weakly idempotent complete additive categories, equipped with a suitable additive analogue of the notion of a Poincaré structure. Recall from \S\ref{subsec:weightstable} the convention that we abbreviate weakly idempotent complete additive categories to \(\wic\)-additive categories. We shall similarly call the objects of \(\Catpad\) \(\wic\)-additive Poincaré categories.

We now proceed to define the \(\L\)-theory space of such a \(\wic\)-additive Poincaré category \((\A,\QF)\). 
We will say that a map \(f\colon \x \to \y\) in a \(\wic\)-additive category \(\A\) is a \emph{split projection} of it admits a section \(\y \to \x\), and a \emph{split inclusion} if it admits a retraction \(\y \to \x\). Being weakly idempotent complete, any split projection in \(\A\) has a fibre \(\z\) and is (non-canonically) isomorphic to a projection \(\z \oplus \y \to \y\), and similarly any split inclusion has a cofibre \(\z\) and is isomorphic to the inclusion \(\x \to \z \oplus \x\). We then also note that pullbacks along split projections and pushouts along split inclusions exist in any \(\wic\)-additive category \(\A\).

The hermitian \(\Q\)-construction can be defined in this setting as follows: for \([n] \in \Del\) we consider the Poincaré \(\wic\)-additive category \(\Q_n(\A,\QF)\) whose underlying \(\wic\)-additive category \(\Q_n(\A)\) is the full subcategory of \(\Fun(\Twar[n],\A)\) spanned by those functors \(\vphi\colon \Twar[n] \to \A\) such that for every \(i \leq j \leq k\) the arrow \(\vphi(i \leq k) \to \vphi(i \leq j)\) is a split projection, the arrow \(\vphi(i \leq k) \to \vphi(j \leq k)\) is a split inclusion, and the square 
\begin{equation}\label{e:exactness} 
\xymatrix{
\vphi(i \leq k) \ar[r]\ar[d] & \vphi(j \leq k) \ar[d] \\
\vphi(i \leq j) \ar[r] & \vphi(j,j) \\
}
\end{equation}
is bicartesian. The quadratic functor is defined by \(\QF_n(\vphi) = \lim_{(i \leq j)\in\Twar[n]\op} \QF(\vphi(i \leq j))\). These assemble to form a simplicial object \(\Q(\A,\QF)\) of Poincaré \(\wic\)-additive categories and Poincaré additive functors between them. 
Given a functor \(\F\colon \Catpad \to \Grp_{\Einf}\) we then define
\[ \F\Q^{(n)}(\A,\QF) \coloneqq \displaystyle\mathop{\colim}_{([k_1],\dots,[k_n]) \in \Del^n}\F(\Q_{k_1}\cdots \Q_{k_n}(\A,\QF)) ,\] 
where the transition maps are induced by the identification \(\Q_0=\id\) in the exterior \(\Q\)-term. We refer to \(\F\Q^{(n)}(\A,\QF)\) as the \(n\)-fold \(\Q\)-construction of \(\F\) (evaluated at \((\A,\QF)\)). 
We then define the \(\L\)-theory space of \((\A,\QF)\) to be
\[ \Lspace^{\pl}(\A,\QF) \coloneqq \colim_n \GWspace^{\pl}\Q^{(n)}(\A,\QF).\]

\begin{aremark}\label{rem:L-with-poinc}
It can be verified that the \(\Einf\)-monoid
\[\Poinc\Q^{(n)}(\A,\QF) \coloneqq \displaystyle\mathop{\colim}_{([k_1],\dots,[k_n]) \in \Del^n}\Poinc(\Q_{k_1}\cdots \Q_{k_n}(\A,\QF)) \]
is group-like. Since group completion commutes with colimits 
it then follows that
\[ \GWspace^{\pl}\Q^{(n)}(\A,\QF) \simeq \Poinc\Q^{(n)}(\A,\QF),\]
and so we could have also defined \(\Lspace^{\pl}\) using \(\Poinc\) instead of \(\GWspace^{\pl}\). 
\end{aremark}

In the setting of ordinary additive categories (and more generally, ordinary exact categories) with duality, Schlichting considered this construction in~\cite[Remark 8]{schlichting-mv}, where it was called the \emph{Witt space}. On the other hand, the above construction can also be performed verbatim in the setting of Poincaré \emph{stable} categories, using the relevant notions of \(\Q\)-construction and \(\GW\)-space defined in~\cite[Sections 3 \& 4]{9authII}. By~\cite[Remark 4.4.12]{9authII} for a Poincaré (stable) category \((\C,\QF)\) we indeed have a natural equivalence
\[\colim_n \GWspace\Q^{(n)}(\C,\QF) = \colim_n\Om^{\infty-n}\GW(\C,\QF\qshift{-n}) \simeq  \Om^{\infty}\L(\C,\QF) = \Lspace(\C,\QF)\]
between this construction and the associated \(\L\)-space, where \(\GW\) is the \(\GW\)-spectrum and the first equality is a definitional one (\(\GWspace\Q^{(n)}(\C,\QF\qshift{n})\) being the \(n\)'th term in the \(\Omega\)-spectrum model for \(\GW\) used in~\cite[Definitions 4.2.1]{9authII}, see also~\cite[Proposition 3.4.5(ii)]{9authII}).

We now come to our main result:

\begin{atheorem}
\label{prop:comp-L}
Let \((\D,\QF)\) be a $0$-dimensional Poincaré (stable) category. Suppose that \(\QF^{\heart} \coloneqq \QF|_{\D^{\heart}}\) takes values in connective spectra, so that we can consider it as taking values in \(\Einf\)-groups. Then the natural map
\[ \Lspace^{\pl}(\D^{\heart},\QF^{\heart}) \lto \Lspace(\D,\QF) \]
issued by the preceding discussion, is an equivalence of spaces.
\end{atheorem}

The remainder of this section is dedicated to the proof of Proposition~\ref{prop:comp-L} using surgery arguments as in the body of the paper. The main point is that the surgery move becomes simpler in the \(\L\)-theory setting, essentially due to the fact that one does not need to distinguish between forward and backward surgery, but just perform surgery on morphisms (and higher simplices) simply as Poincaré objects in the \(\Q\)-construction.

We begin by observing that if \((\D,\QF)\) is a weighted Poincaré category then for every \(n\) the Poincaré category \(\Q_n(\D,\QF)\) inherits a weight structure:

\begin{alemma}\label{lem:weight-on-Q}
Let \(\C\) be a stable category equipped with a weight structure. Then there exists a weight structure on \(\Q_n(\C)\) such that an object \(\vphi\colon \Twar[n] \to \C\) of \(\Q_n(\C)\) lies in \(\Q_n(\C)_{[a,b]}\) (for \(-\infty \leq a,b \leq \infty\)) if and only if
\begin{enumerate}
\item
for every \(i \leq j \in [n]\) the object \(\vphi(i \leq j)\) lies in \(\C_{[a,b]}\);
\item
for every \(i \leq j \leq k \in [n]\) the fibre of the map \(\vphi(i \leq k) \to \vphi(i \leq j)\) lies in \(\C_{[a,b]}\);
\item
for every \(i \leq j \leq k \in [n]\) the cofibre of the map \(\vphi(i \leq k) \to \vphi(j \leq k)\) lies in \(\C_{[a,b]}\).
\end{enumerate}
\end{alemma}

\begin{proof}
We show that the full subcategories obtained by taking \((a,b) = (0,\infty)\) and \((a,b) = (-\infty,0)\) satisfy the axioms of a weight structure.
For this, note that \(\Q_n(\C)\) is naturally equivalent to the \(\bS\)-construction \(\bS_{2n+1}(\C)\) and the latter can in turn be identified with \(\Fun([2n],\C)\), compare e.g.\ \cite[Section B.1]{9authII}. Tracing these equivalences it will suffice to show that for every \(m \geq 0\) there is a weight structure on \(\Fun([m],\C)\) whose connective objects are those \(x_0 \to \dots \to x_m\) such that each \(x_i\) is connective and whose coconnective objects are those for which each \(x_i\) is coconnective and each \(x_i \to x_{i+1}\) has a coconnective cofibre. We now argue by induction on \(m\). For \(m = 0\) we just have \(\C\) itself with its underlying weight structure, so the claim is clear. Now suppose the claim is known for \(m-1 \geq 0\) and let us prove for \(m\). 
Clearly the classes of connective and coconnective objects are closed under retracts. 
Let \(x_0 \to \dots \to x_m\) be a chain of maps between coconnective objects with coconnective cofibres and \(y_0 \to \dots \to y_m\) a chain of maps between connective objects. We may then consider the pushout square
\[\xymatrix{
[0 \to \dots \to 0 \to x_{m-1}] \ar[r]\ar[d] & [x_0 \to \dots \to x_{m-1} = x_{m-1}] \ar[d] \\
[0 \to \dots \to 0 \to x_{m}] \ar[r] & [x_0 \to \dots \to x_{m-1} \to x_m]
}\]
in \(\Fun([m],\C)\), where we observe that the top right corner is left Kan extended from \([m-1] \subseteq [m]\) and the objects in the left column are left Kan extended from \(\{m\} \subseteq [m]\). We then have that 
\[ \hom_{\Fun([m],\C)}([x_0 \to \dots \to x_m],[y_0 \to \dots \to y_m]) =\] \[\hom_{\Fun([m-1],\C)}([x_0 \to \dots \to x_{m-1}],[y_0 \to \dots \to y_{m-1}]) \times_{\hom_{\C}(x_{m-1},y_m)} \hom_{\C}(x_m,y_m), \]
and the last term is connective by the induction hypothesis since 
\[\fib[\hom_{\C}(x_{m},y_m) \to \hom_{\C}(x_{m-1},y_m)] = \hom_{\C}(\cof[x_{m-1}\to x_m],y_m)]\] 
is connective. 

It is left to establish the factorisation axiom. Given a general object \([x_0 \to \dots \to x_m]\), we have by the induction hypothesis an exact sequence 
\[[y_0 \to \dots \to y_{m-1}] \lto [x_0 \to \dots \to x_{m-1}] \lto [z_0 \to \dots \to z_{m-1}]\] 
such that each \(y_i\) is coconnective, each \(y_i \to y_{i+1}\) has a coconnective cofibre, and each \(z_i\) is 1-connective cofibre. Write \(x'\) for the cofibre of the composite map \(y_{m-1} \to x_{m-1} \to x_m\). Invoking the factorisation axiom for \(\C\) we may find an exact sequence \(y'\to x'\to z_m\) such that \(y'\) is coconnective and \(z_m\) is 1-connective. We then obtain an exact sequence
\[ [y_0 \to \dots \to y_{m-1} \lto y_{m+1} \oplus y'] \lto [x_0 \to \dots \to x_m] \to [z_0 \to \dots \to z_m] \]
whose left most object is coconnective in \(\Fun([m],\C)\) and whose right most object is 1-connective there.
\end{proof}

We now observe that if \((\D,\QF)\) is a weighted Poincaré category of dimension 0 then \(\Q_n(\D,\QF)\) is again 0-dimensional with respect to the weight structure of Lemma~\ref{lem:weight-on-Q}, and \(\QF^{\heart} \coloneqq \Om^{\infty}\QF_{\D^{\heart}}\) is an additive Poincaré structure on \(\D^{\heart}\). In addition, in this case the heart of the above weight structure on \(\Q_n(\D,\QF)\) identifies with the image of \(\Q_n(\D^{\heart},\QF^{\heart})\), where the latter is defined as above for the Poincaré \(\wic\)-additive category \((\D^{\heart},\QF^{\heart})\). Here we use the fact that a map in \(\D\) whose source and target lie in \(\D^{\heart}\) has its fibre lying in \(\D^{\heart}\) if and only if it admits a section, and has is cofibre lying in \(\D^{\heart}\) if and only if it admits a retraction. We also note that, as is visible by the explicit description of Lemma~\ref{lem:weight-on-Q}, all the simplicial transition functors in the simplicial object \(\Q(\D,\QF)\) are 0-dimensional, that is, they preserve hearts.

\begin{adefinition}
For a weighted Poincaré category \((\D,\QF)\) of dimension 0,
we will say that an 
object in \(\D\) is \emph{\(m\)-bounded} if it lies in \(\D_{[-m,m]}\). Similarly, we will say that a Poincaré object is \(m\)-bounded if its underlying object is so. We then write \(\Poinc_m(\D,\QF) \subseteq \Poinc(\D,\QF)\) for the full subspace spanned by the \(m\)-bounded Poincaré objects.
\end{adefinition}

\begin{aremark}\label{rem:bounded-is-connective}
For a 0-dimensional weighted Poincaré category \((\D,\QF)\), the duality switches between the classes of \((-m)\)-connective and \(m\)-coconnective objects. It then follows that a Poincaré object is \(m\)-bounded if and only if its underlying object is \((-m)\)-connective. Applying this for the 0-dimensional weighted Poincaré category \(\Q_1(\D,\QF)\) we see that a cobordism is \(m\)-bounded if and only if its source and target are \(m\)-connective and its left leg is \(m\)-connective.
\end{aremark}

\begin{aremark}\label{rem:weight-cob}
A Poincaré object of \(\Q_n(\D,\QF)\) is \(m\)-bounded exactly when it corresponds to a composable sequence of cobordisms, each of which is \(m\)-bounded in \(\Q_1(\D,\QF)\). In particular, the simplicial object \(\Poinc_m\Q(\D,\QF)\) is again a (complete) Segal space, which by Remark~\ref{rem:bounded-is-connective} corresponds to the subcategory \(\Cob^{m,m}(\D,\QF) \subseteq \Cob(\D,\QF)\) spanned by the \(m\)-connective Poincaré objects and cobordisms with left leg \(m\)-connective between them, see \S\ref{sketch}.
\end{aremark}

Given \(n \geq 0\) we now denote by
\[ \Poinc_m\Q^{(n)}(\D,\QF) \coloneqq \displaystyle\mathop{\colim}_{([k_1],\dots,[k_n]) \in \Del^n}\Poinc_m(\Q_{k_1}\cdots \Q_{k_n}(\D,\QF)) \]
the space obtained by performing the iterated \(\Q\)-construction and taking \(m\)-bounded Poincaré objects.

\begin{alemma}\label{l:comp}
Under the assumptions of Proposition~\ref{prop:comp-L}, for every \(m,n \geq 0\) the square of spaces
\begin{equation}\label{e:comp-4} 
\xymatrix{
\Poinc_m\Q^{(n)}(\D,\QF) \ar[r]\ar[d] & \Poinc_m\Q^{(n+1)}(\D,\QF) \ar[d] \\
\Poinc_{m+1}\Q^{(n)}(\D,\QF) \ar[r]\ar@{-->}[ur] & \Poinc_{m+1}\Q^{(n+1)}(\D,\QF) \\
}
\end{equation}
admits a dotted lift as indicated.
\end{alemma}

The proof of our main result Proposition~\ref{prop:comp-L} can be reduced to the statement of Lemma~\ref{l:comp}:
\begin{proof}[Proof of Proposition~\ref{prop:comp-L} given Lemma~\ref{l:comp}]
Using Remark~\ref{rem:L-with-poinc} we may factor the map in Proposition~\ref{prop:comp-L} as a transfinite composition
\[\Lspace^{\pl}(\D^{\heart},\QF^{\heart}) = 
\colim_n \Poinc_0\Q^{(n)}(\D,\QF) \to \colim_n \Poinc_1\Q^{(n)}(\D,\QF) \to \dots \to \colim_n\Poinc\Q^{(n)}(\D,\QF) = \Lspace(\D,\QF).\]
It will then suffice to show that the map
\begin{equation}\label{e:comp-3} 
\colim_n \Poinc_m\Q^{(n)}(\D,\QF) \lto \colim_n \Poinc_{m+1}\Q^{(n)}(\D,\QF)
\end{equation}
is an equivalence of spaces for every \(m \geq 0\).

Using the dotted lifts in~\eqref{e:comp-4} we may weave the sequences \(\{\Poinc\Q^{(n)}_m(\D,\QF)\}_{n \in \NN}\) and \(\{\Poinc\Q^{(n)}_{m+1}(\D,\QF)\}_{n \in \NN}\) into a single sequence
\begin{equation}\label{e:big-sequence}
\Poinc_m(\D,\QF) \lto \Poinc_{m+1}(\D,\QF) \lto \Poinc_m\Q^{(1)}(\D,\QF) \to \Poinc_{m+1}\Q^{(1)}(\D,\QF) \to \dots ,
\end{equation}
where the inclusions of both \(\{\Poinc_m\Q^{(n)}(\D,\QF)\}_{n \in \NN}\) and \(\{\Poinc_{m+1}\Q^{(n)}(\D,\QF)\}_{n \in \NN}\) induce an equivalence on colimits by cofinality. Let \(P\) denote the colimit of~\eqref{e:big-sequence}. Then under this equivalence the map~\eqref{e:comp-3} identifies with the map \(P \to P\) induced on the colimit of~\eqref{e:big-sequence} by a shift, 
and is hence an equivalence. 
\end{proof}

We now turn to the proof of Lemma~\ref{l:comp}. We first note that the square~\eqref{e:comp-4} (without the dotted lift) is induced by a square of the form
\begin{equation}\label{e:comp-5} 
\xymatrix{
\Poinc_m(-) \ar[r]\ar[d] & \Poinc_m\Q^{(1)}(-) \ar[d] \\
\Poinc_{m+1}(-) \ar[r] & \Poinc_{m+1}\Q^{(1)}(-) \\
}
\end{equation}
by evaluating on a certain diagram of weighted Poincaré categories and taking colimits. Here we understand~\eqref{e:comp-5} as a square of functors defined on the category \(\Catpw\) of 0-dimensional weighted Poincaré categories and 0-dimensional Poincaré functors between them. 
To construct a lift in~\eqref{e:comp-4} it will hence suffice to construct a lift in~\eqref{e:comp-5} as a square of functors.

\begin{aremark}
Suppose that \(\D = \Dperf(R)\) for some ring \(R\), equipped with its standard weight structure, and that \(\QF\) is a Poincaré structure on \(\Dperf(R)\) with duality associated to some invertible module with involution \(M\) over \(R\). Then \((\D,\QF)\) is a \(0\)-dimensional Poincaré category, and the condition that \(\QF\) takes connective values on the heart is equivalent to \(\QF\) being 0-quadratic in the sense of~\cite[Definition 1.1.2]{9authIII}. For \(m \geq 0\) the group 
\[ \pi_0\Poinc_{m}\Q^{(1)}(\Dperf(R),\QF) = |\pi_0\Poinc_{m}\Q(\Dperf(R),\QF)|\] 
can then be identified with the bounded L-group \(\L^{2m,2m}(R,\QF)\) of~\cite[Definition 1.2.1]{9authIII}, and~\cite[Propositions 1.2.6 \& 1.2.14]{9authIII} shows that the maps $\L^{2m,2m}_0(R,\QF) \to \L_0(R,\QF)$ are isomorphisms under these assumptions for every $m \geq 0$. 
It then follows that the right vertical map in~\eqref{e:comp-5} is an isomorphism on \(\pi_0\) when evaluated on such \((\Dperf(R),\QF)\).
In fact, the arguments of~\cite[Section 1.2]{9authIII} can be extended in a relatively straightforward manner to an arbitrary 0-dimensional Poincaré category \((\D,\QF)\) with \(\QF\) taking connective values on the heart, and the right vertical map in~\eqref{e:comp-5} is an isomorphism on \(\pi_0\) when evaluated on any such \((\D,\QF)\). In particular, the fact that a section on the level of \(\pi_0\) exists and is unique is essentially contained in~\cite{9authIII}. The crux of what follows is really to show that this section can be chosen coherently on the level of spaces. 
\end{aremark}

To proceed, we now replace~\eqref{e:comp-5} by an equivalent model. 
Let \(\Poset\) denote the category of posets and let \(\Beta\colon \Del \to \Poset\) (resp. \(\Beta^{\circ}\colon\Del \to \Poset$) be the functor which sends \([n]\) to the poset of subsets (resp.\ non-empty subsets) of \(\{0,\dots,n\}\). We have a natural transformation \(\Beta^{\circ} \to \Beta\) which is levelwise a full inclusion.

\begin{aconstruction}\label{c:cubes}
Given a Poincaré category \(\D\) and \([n] \in \Del\) let us denote by 
\[ \R_n(\D,\QF) \coloneqq \Q_{\Beta[n]}(\D,\QF) \subseteq (\D^{\Twar(\Beta[n])},\lim \circ \QF) .\]
We consider \(\R(\D,\QF)\) as a simplicial object in \(\Catp\). Similarly, we denote by \(\R^{\circ}_n(\D,\QF) \coloneqq \Q_{\Beta^{\circ}(n)}(\D,\QF)\), so that we have a natural map of simplicial objects
\[
\R(\D,\QF) \to \R^{\circ}(\D,\QF) .
\]
\end{aconstruction}

\begin{aremark}\label{r:cubes}
Explicitly, \(\R_n(\D,\QF)\) is the full subcategory of \(\D^{\Twar(\Beta[n])}\) spanned by those diagrams \(\vphi\colon \Twar(\Beta[n]) \to \D\) such that the square
\[ \xymatrix{
\vphi(S \subseteq U) \ar[r]\ar[d] & \vphi(T \subseteq U) \ar[d] \\
\vphi(S \subseteq T) \ar[r] & \vphi(T \subseteq T) \\
}\]
is exact for every \(S \subseteq T \subseteq U \in \Beta[n]\). It is endowed with the hermitian structure restricted from the diagram structure on \(\D^{\Twar(\Beta[n])}\), and is Poincaré because it is a limit of Poincaré categories. Poincaré objects in \(\R_n(\D,\QF)\) then correspond to \(\Beta[n]\)-diagrams in \(\Cob(\D,\QF)\) by \cite[Proposition 2.3.6]{9authII}. Similarly, \(\R^{\circ}_n(\D,\QF)\) is a Poincaré category and Poincaré objects in \(\R^{\circ}_n(\D,\QF)\) correspond to \(\Beta^{\circ}([n])\)-diagrams in \(\Cob(\D,\QF)\).
\end{aremark}

We note that by \cite[Lemma 2.2.5]{9authII}, the simplicial object \(\R^{\circ}(\D,\QF)\) is the simplicial subdivision of the simplicial object \(\Q(\D,\QF)\), i.e.\ the evaluation of the limit preserving extension of $\Q(\D,\QF) \colon \bbDelta\op \rightarrow \Catp$ to a functor $\sSps\op \rightarrow \Catp$ at the barycentric subdivision of $\Delta^n$, which is by definition the nerve of $\Beta^{\circ}([n])$. There is always a canonical map from a simplicial object to its subdivision which is induced by the natural transformation \(\tau_n\colon \Beta^{\circ}([n]) \to [n]\) given by \(\tau_n(S) = \max(S)\). 
The commutative square
\begin{equation}\label{e:rho}
\xymatrix{
\Beta^{\circ}([n]) \ar[r]\ar[d] & [n] \ar[d] \\
\Beta[n] \ar[r] & \ast \\
}
\end{equation}
of functors \(\Del \to \Poset\) then induces a commutative square of simplcial objects 
\begin{equation}\label{e:rho-2}
\xymatrix{
(\D,\QF) \ar[r]\ar[d] & \Q(\D,\QF) \ar[d] \\
\R(\D,\QF) \ar[r] & \R^{\circ}(\D,\QF) \\
}
\end{equation}
in \(\Catp\).

\begin{adefinition}\label{d:bounded-4}
We will say that an object \(\vphi \in \R_n(\D,\QF)\) (resp. \(\vphi \in \R^{\circ}_n(\D,\QF)\)) is \emph{locally \(m\)-bounded} if for every simplex \(\sig\colon[k] \to \Beta[n]\) (resp. \(\sig\colon[k] \to \Beta^{\circ}([n])$), the image of \(\vphi\) in 
\(\Q_k(\D,\QF)\) is \(m\)-bounded. We will say that a Poincaré object in \(\R_n(\D,\QF)\) or \(\R_n^{\circ}(\D,\QF)\) is locally \(m\)-bounded if its underlying object is so. We will denote by \(\Poinc_m(\R_n(\D,\QF)) \subseteq \Poinc(\R_n(\D,\QF))\) and \(\Poinc_m(\R^{\circ}_n(\D,\QF)) \subseteq \Poinc(\R^{\circ}_n(\D,\QF))\) the respective subspaces of locally \(m\)-bounded Poincaré objects.
\end{adefinition}

The square~\eqref{e:rho-2} determines a commutative cube of spaces
\begin{equation}\label{e:rho-3}
\xymatrix{
\Poinc_m(\D,\QF) \ar[rr] \ar[dd] \ar[dr] && |\Poinc_m(\Q(\D,\QF))| \ar[dr]^-{\simeq} \ar[dd] |!{[dl];[dr]}\hole \\
& |\Poinc_m(\R(\D,\QF))| \ar[rr] \ar[dd] && |\Poinc_m(\R^{\circ}(\D,\QF))| \ar[dd] \\
\Poinc_{m+1}(\D,\QF) \ar[rr] |!{[ur];[dr]}\hole \ar[dr] && |\Poinc_{m+1}(\Q(\D,\QF))| \ar[rd]^{\simeq} \\
& |\Poinc_{m+1}(\R(\D,\QF))| \ar[rr] && |\Poinc_{m+1}(\R^{\circ}(\D,\QF))| \\
}
\end{equation}
which is natural in \(\D\) as a weighted Poincaré category, and where the indicated equivalences are a consequence of the general fact that simplicial subdivision $\mathrm{Ex} \colon \sSps \rightarrow \sSps$ preserves geometric realisations; for simplicial sets this is a classical result of Kan, and the case of general simplicial spaces easily reduces to this by picking a trivial fibration from a simplicial set \cite[Lecture 7, Proposition 7]{Lurie-L}. We may then consider~\eqref{e:rho-3} as a natural transformation between the square of functors~\eqref{e:comp-5} and the square of functors  
\begin{equation}\label{e:comp-6} 
\xymatrix{
|\Poinc_m(\R(-))| \ar[r]\ar[d] & |\Poinc_m(\R^{\circ}(-))| \ar[d] \\
|\Poinc_{m+1}(\R(-))| \ar[r] & |\Poinc_{m+1}(\R^{\circ}(-))|  \\
}
\end{equation}

\begin{alemma}\label{l:translation}
In the cube~\eqref{e:rho-3} the maps \(\Poinc_m(\D,\QF) \to |\Poinc_m(\R(\D,\QF))|\) and \(\Poinc_{m+1}(\D,\QF) \to |\Poinc_{m+1}(\R(\D,\QF))|\) are equivalences. In other words, when considered as a natural transformations between~\eqref{e:comp-5} and~\eqref{e:comp-6}, the cube~\eqref{e:rho-3} is an equivalence.
\end{alemma}
\begin{proof}
Both equivalences follow from the fact that \([n] \mapsto (\Q(\D,\QF))^{\Beta([n])}\) extends to a split simplicial object with augmentation \((\D,\QF)\). To see this, simply note that \(\Beta(-)\) extends to a functor 
\[ \Del_{-\infty} \lto \Poset \quad\quad S \longmapsto \Beta(S \setminus \{-\infty\}) \]
where we recall that \(\Del_{-\infty}\) is the category consisting of the linearly ordered sets 
\[ \{-\infty\}, \{-\infty,0\},\{-\infty,1\},\dots,\{-\infty,1,\dots,n\},\dots \] and order preserving maps which preserve \(-\infty\) between then, to which \(\Del\) maps via the association \([n] \mapsto \{-\infty,\dots,n\}\).
\end{proof}

\begin{proof}[Proof of Lemma~\ref{l:comp}]
In light of Lemma~\ref{l:translation}, it will suffice to show that the square~\eqref{e:comp-6} admits a lift. We now begin the surgery move argument. Recall the surgery equivalence
\begin{equation}\label{e:surgery-2}
\Cob(\D,\QF)_{(\x,q)/} \xrightarrow{\simeq} \Sur(\x,q), 
\end{equation}
where \(\Sur(\x,q)\) is the total space of the right fibration \(\pi\colon \Sur(\x,q) \to \D_{/\x}\) classified by the functor \(\chi(f\colon \z \to \x) = \{0\} \times_{\Om^{\infty}\Q(\z)}\{f^*q\}\).  This equivalence records the fact that cobordisms out of \(\x\) can be described in terms of surgery data. 

As mentioned in Remark~\ref{rem:weight-cob}, the simplicial object \(\Poinc_m\Q(\D,\QF)\) is the complete Segal space corresponding to the subcategory \(\Cob^{m,m}(\D,\QF) \subseteq \Cob(\D,\QF)\) spanned by the \(m\)-bounded Poincaré objects and the \(m\)-bounded cobordisms between them. 
Given an \(m\)-bounded Poincaré object \((\x,q)\), the equivalence~\eqref{e:surgery-2} then identifies the subcategory \(\Cob^{m,m}(\D)_{(\x,q)/} \subseteq \Cob(\D,\QF)_{(\x,q)/}\) with the subcategory of \(\Sur(\x,q)\) whose objects are the pairs \((f\colon \z \to \x,\eta\colon f^*q \sim 0)\) for which \(\Sig \z\) is \(m\)-bounded and whose maps are those maps \(g\colon\z \to \z'\) over \(\x\) for which \(\Sig\cof(g)\) is \(m\)-bounded. We note that the condition that \(\x\) and \(\Sig \z\) are \(m\)-bounded implies that the trace of the surgery on \(\x\) along \(\z\) is \(m\)-bounded as well.

Now by Remark~\ref{r:cubes} Poincaré objects in \(\R_n(\D,\QF)\) can be identified with diagrams \(\Beta[n] \to \Cob(\D,\QF)\). In particular, if \(\vphi \colon \Twar\Beta([n]) \to \D\) is an object of \(\R_n(\D,\QF)\) equipped with a Poincaré structure \(q_{\vphi}\), and we set \((\x,q) = (\vphi(\emptyset \subseteq \emptyset),q_{\vphi}(\emptyset\subseteq\emptyset))\) to be the Poincaré object obtained by evaluating 
at \((\emptyset \subseteq \emptyset) \in \Twar\Beta([n])\), then \((\vphi,q_{\vphi})\) determines a diagram \(\Beta[n] \to \Cob(\D)_{(\x,q)/}\) (sending \(\emptyset\) to the identity on \((\x,q)$), and hence a diagram
\[ g_\vphi\colon \Beta[n] \lto \Sur(\x,q) \]
(sending \(\emptyset\) to the trivial surgery datum). In particular, for every \(S \subseteq [n]\) we have a surgery datum 
\[ g_{\vphi}(S) = (f\colon\z_{S} \to \x,\eta_S\colon f_S^*q\sim 0) .\]
For \(m, n \geq 0\) we now define \(\SC^m_n(\D,\QF) \subseteq \Poinc_{m+1}(\R_n(\D,\QF))\) to be the full subspace spanned by those locally \((m+1)\)-bounded Poincaré objects \((\vphi,q) \in \Poinc_{m+1}(\R_n(\D,\QF)\) such that
\begin{enumerate}
\item
the restriction \(\vphi|_{\Beta^{\circ}(n)} \in \R^{\circ}_n\D\) is locally \(m\)-bounded; and
\item
the composite cube diagram
\begin{equation}\label{e:composed} 
\Beta[n] \xrightarrow{g_\vphi} \Sur(\x,q) \xrightarrow{\pi} \D_{/\x} \lto \D 
\end{equation}
is strongly cocartesian and takes values in \(\Sig^{-m-1}\D^{\heart} \subseteq \D\). 
\end{enumerate}
We note that in the above situation, the cube diagram~\eqref{e:composed} sends \(\emptyset \in \Beta[n]\) to \(\z_{\emptyset} = 0\). The condition that this cube is strongly cocartesian is then equivalent to the condition that for every non-empty \(S \subseteq [n]\), the maps \(\z_{\{i\}} \to \z_S\) exhibit \(\z_S\) as the direct sum of \(\oplus_{i \in S}\z_{\{i\}}\) in \(\D\). For any \(S \subseteq S'\), the associated cobordism \(\x_S \leftarrow \cob_{S,S'} \to \x_{S'}\) 
corresponds to a surgery data of the form \(\z_{S,S'} \to \x_S\), 
with 
\[ \z_{S,S'} = \cof(\z_S \to \z_{S'}) \simeq \oplus_{i \in S' \setminus S}\z_{\{i\}} .\]
Since all the \(\z_{\{i\}}\) are assumed to lie in \(\Sig^{-m-1}\D\), the same holds for each of these direct sums, so the suspension of each of these surgery data is in particular \(m\)-bounded. In light of this, Condition (i) above is equivalent to \(\x_{\{i\}} = \vphi(\{i\} \subseteq \{i\})\) being \(m\)-bounded for every \(i \in [n]\). 

We note that Conditions (i) and (ii) above are preserved by the face operators in \(\Poinc_m\R(\D,\QF)\), so that the spaces \(\SC^m_n(\D,\QF)\) assemble into a \emph{semi-simplicial space} \(\SC^m(\D,\QF)\), the analogue of the \emph{surgery complex} in the present setting. It fits into a commutative diagram of the form
\begin{equation}\label{e:r-1} 
\xymatrix{
\Poinc_m(\D,\QF) \ar^-{\simeq}[r]\ar[d] & |\Poinc_m(\R(\D,\QF))| \ar[r]\ar[d] & |\Poinc_m(\R^{\circ}(\D,\QF))| \ar[d] \\
|\SC^m(\D,\QF)| \ar[r]\ar[urr] & |\Poinc_{m+1}(\R(\D,\QF))| \ar[r] & |\Poinc_{m+1}(\R^{\circ}(\D,\QF))|  \ .
}
\end{equation}
Here the top left horizontal map is the equivalence of Lemma~\ref{l:translation},
and the left vertical map is the inclusion of constant cubes of cobordisms on \(m\)-bounded objects. 
The diagonal lift in~\eqref{e:r-1} is simply given by restriction along \(\Beta^{\circ}([n]) \subseteq \Beta[n]\).  
The diagram~\eqref{e:r-1} is visibly natural in \((\D,\QF)\) with respect to 0-dimensional Poincaré functors. In particular, we may consider the external rectangle as a square of functors \(\Catpw \to \Sps\), equipped with a natural transformation to the right square, which is just~\eqref{e:comp-6}. In addition, the external rectangle comes with a lift by construction. We would hence be done if we knew that this natural transformation is an equivalence. By Lemma~\ref{l:translation} it will suffice to show that the forget-the-surgery map 
\begin{equation}\label{e:aug} 
\SC^{m}_n(\D,\QF) \lto \Poinc_{m+1}(\D,\QF) ,
\end{equation}
sending \(\vphi\colon\Twar\Beta([n]) \to \D\) to \(\vphi(\emptyset\subseteq \emptyset)\), induces an equivalence
\[ |\SC^{m}(\D,\QF)| \xrightarrow{\simeq} \Poinc_{m+1}(\D,\QF) .\]
This is proven in Proposition~\ref{prop:core} just below. 
\end{proof}

The following proposition encompasses the essence of the surgery argument, showing that the choice of surgery data becomes essentially unique upon realisation, as described in \S\ref{coskeletalsection}             in the body of the paper. We carry out the argument in the present context for completeness.
\begin{aproposition}\label{prop:core}
Let \((\D,\QF)\) be a weighted 0-dimensional Poincaré category such that \(\QF\) sends \(\D^{\heart}\) to connective spectra. Then the map 
\begin{equation}\label{e:key-map}
\SC^{m}(\D,\QF) \lto \Poinc_{m+1}(\D,\QF),
\end{equation}
where the target is considered as a constant semi-simplicial space,
induces an equivalence after geometric realisation. 
\end{aproposition}
\begin{proof}
Since colimits in spaces are universal we can check this on fibres. 
In other words, it will suffice to show that for every \((m+1)\)-bounded Poincaré object \((\x,q) \in \Poinc_{m+1}(\D,\QF)\), the geometric realisation of the simplicial space 
\[ X \coloneqq \SC^{m}(\D,\QF)_{(\x,q)} \] 
is contractible, where \(\SC^{m}(\D,\QF)_{(\x,q)}\) denotes the (homotopy) fibre of~\eqref{e:aug} over \((\x,q)\). 
To have an explicit picture in mind we recall that the space of \(n\)-simplices of \(X\) corresponds, under the surgery equivalence, to the groupoid of those functors \(\Beta[n] \to \Sur(\x,q)\) whose underlying cube \(\Beta[n] \to \D\) is strongly cocartesian of the form \(S \mapsto \oplus_{i \in S} \z_i\) with \(\z_i \in \Sig^{-m-1}\D^{\heart}\), and such that the trace of doing surgery on \(\x\) along \(\z_i\) is \((-m)\)-connective for every \(i\) (see Remark~\ref{rem:bounded-is-connective}). This trace lies in an exact sequence between \(\cof[\z_i \to \x]\) and \(\Om \Dual z_i\), and since the latter is \(m\)-connective (and in particular \((-m)\)-connective, since \(m \geq 0\)), we have that 
the trace of the surgery is \((-m)\)-connective if and only if the cofibre of the map \(\z_i \to \x\) is \((-m)\)-connective. 
Such a cube in \(\D\) is uniquely (and freely) determined by the collection of maps \(f_i\colon \z_i \to\x\). In addition, since \(\Beta[n]\) has a final object \([n] \in \Beta[n]\), and \(\Sur(\x,q) \to \D_{/\x}\) is a right fibration, we have that lifting a diagram \(\Beta[n] \to \D_{/\x}\) to \(\Sur(\x,q)\) is equivalent to lifting the image of \([n] \in \Beta[n]\). It then follows that the data of a point in \(X_n\) can be equivalently described by the following:
\begin{enumerate}
\item 
A collection of objects \(\z_0,\dots,\z_n \in \Sig^{-m-1}\D^{\heart}\). 
\item
A collection of maps \(f_i\colon \z_i \to \x\) whose cofibres are \((-m)\)-connective.
\item
A null-homotopy of the image of \(q\) in \(\Om^{\infty}(\QF(\z_0 \oplus \dots \oplus \z_n))\).
\end{enumerate}
We now claim that \(X\) is Kan contractible, that is, the maps \(X_n = X(\Delta^n) \to X(\partial \Delta^n)\) 
are surjective on components, and so its geometric realisation is contractible as desired; see the beginning of \S\ref{coskeletalsection} for the notation. By the discussion there it will suffice to prove the following:
\begin{enumerate}
\item
The space \(X_0\) is non-empty.
\item
The fibres of the map \(X_1 \to X_0 \times X_0\) are all non-empty.
\item
\(X\) is \(1\)-coskeletal, that is, the map \(X \to \cosk_1\sk_1(X)\) is an equivalence.
\end{enumerate}
To verify property (i) we note that by the definition of a weight structure there exists a cofibre sequence of the form \(\z \st{f}{\to} \x \to \y\) such that \(\z\) is \((-m-1)\)-coconnective and \(\y\) is \((-m)\)-connective. Since \(\x\) is \((-m-1)\)-connective it then follows that \(\z\) is \((-m-1)\)-connective as well, and so \(\z \in \Sig^{-m-1}\D^{\heart}\). To make \(\z\) into a surgery data we need to produce a null-homotopy of \(f^*q\) in \(\Om^{\infty}\QF(\z)\). This is possible since our assumption that \(\QF|_{\D^{\heart}}\) is connective implies the same for the linear part \(\Lam_{\QF}\) of \(\QF\), and so \(\Lam_{\QF}(\z)\) is \((m+1)\)-connective. In addition, the standing 0-dimensionality assumption means that the bilinear part takes connective values on \(\D^{\heart}\), and so \(\Bil_{\QF}(\z,\z)\) is \((2m+2)\)-connective. Together this means that \(\QF(\z)\) is \((m+1)\)-connective, and in particular \(1\)-connective, so that the desired null-homotopy exists. 
This shows (i).

To verify property (ii) we need to check that if \((f_0\colon \z_0 \to \x,\eta_0\colon f_0^*q \sim 0)\) and \((f_1\colon \z_1 \to \x,\eta_1\colon f_1^*q \sim 0)\) are two surgery data with \(\z_0,\z_1 \in \Sig^{-m-1}\D^{\heart}\), then there exists a null-homotopy of \(f_{01}^*q\) in \(\Om^{\infty}\QF(\z_0\oplus \z_1)\) which extends both \(\eta_0\) and \(\eta_1\), where \(f_{01}\colon \z_0 \oplus \z_1 \to \x\) is the map induced by \(f_0\) and \(f_1\). Indeed, \(\QF(\z_0\oplus \z_1)\) canonically decomposes as
\[ \QF(\z_0\oplus \z_1) \simeq \QF(\z_0) \oplus \QF(\z_1) \oplus \Bil_{\QF}(\z_0,\z_1), \]
and  
\(\Bil_{\QF}(\z_0,\z_1)\) is \((2m+2)\)-connective as mentioned above, since \(\z_0,\z_1 \in \Sig^{-m-1}\D^{\heart}\). 

Let us now show (iii). For brevity let us write 
\(Z_n \coloneqq X(\sk_2\Delta^n)\) (the target is nothing but the morphism space from $\sk_2\Delta^n$ to $X$), so that the claim we need to prove is that the natural map \(X_n \to Z_n\) is an equivalence. Let \(\I_{\leq 2} \subseteq \Beta[n]\) be the subposet spanned by the sets of cardinality \(\leq 2\). Let \(\Fun_{\scc}(\Beta[n],\D_{/\x}) \subseteq \Fun(\Beta[n],\D_{/\x})\) be the full subcategory spanned by those cube diagrams which are strongly cocartesian and let \(\Fun'(\Beta[n],\Sur(\x,q)) \subseteq \Fun(\Beta[n],\Sur(\x,q))\) be the inverse image of \(\Fun_{\scc}(\Beta[n],\D_{/\x})\) in \(\Fun(\Beta[n],\Sur(\x,q))\). Consider the commutative diagram
\[ \xymatrix{
X_n \ar@{^{(}->}[r]\ar[d] & \core(\Fun'(\Beta[n],\Sur(\x,q))) \ar[r]\ar[d] & \core(\Fun_{\scc}(\Beta[n],\D_{/\x})) \ar[d]^{\simeq} \\
Z_n \ar@{^{(}->}[r] & \core(\Fun(\I_{\leq 2},\Sur(\x,q))) \ar[r] & \core(\Fun(\I_{\leq 2},\D_{/\x})) \\
}\]
Then the right vertical map is an equivalence by uniqueness of left Kan extensions, and the left square is cartesian since the conditions determining \(X_n\) inside \(\core(\Fun'(\Beta[n],\Sur(\x,q)))\) are given in terms of the values on subsets of cardinality \(\leq 1\): indeed, the objects \(\oplus_{i \in S}\z_i\) belong to \(\Sig^{-m-1}\D^{\heart}\) iff each \(\z_i \in \Sig^{-m-1}\D^{\heart}\), and the maps \(\oplus_{i \in S}\z_i \to \x\) all have \((-m)\)-connective cofibre if and only if each \(\z_i \to \x\) has. To finish the proof it will hence suffice to show that the right square is cartesian. Unwinding the definitions, what we need to show is that if \(\vphi\colon\Beta[n] \to \D_{/\x}\) is a strongly cocartesian cube then the natural map
\[ \lim_{\Beta[n]\op}\chi\vphi \lto \lim_{\I_{\leq 2}\op}\chi\vphi \]
is an equivalence, where \(\chi\colon \D_{/\x} \to \Sps\) is the functor \(\chi(f\colon \z \to \x) = \{0\} \times_{\Om^{\infty}\QF(\z)}\{f^*q\}\). But this is a direct consequence of the fact that \(\chi\) is \(2\)-excisive as a functor to spaces (since \(\Om^{\infty}\QF\) is such) and hence \(\chi\vphi\) is a right Kan extension of \(\chi\vphi|_{\I_{\leq 2}}\) whenever \(\vphi\colon \Beta[n] \to \D_{/\x}\) is strongly cocartesian.
\end{proof}

\begin{aremark}
Finally, observe that Theorem \ref{prop:comp-L} is also a direct consequence of Theorem \ref{thm:resolution} using the weight structures from Lemma \ref{lem:weight-on-Q} and the fact that the colimit over $\GW\Q^{(n)}$ is indeed a model for the $\L$-spectra of Poincar\'e categories as remarked at the start of this section, see \cite[Corollary 3.6.16 \& Theorem 4.4.11]{9authII}.
\end{aremark}

\subsection{An alternative approach to Theorem~\ref{thm:main_special}}

To explain this approach 
let us begin by reviewing algebraic \(\K\)-theory and \(\GW\)-theory in the setting of \(\wic\)-additive categories.  
In this context it is natural to work with the group completion version of algebraic \(\K\)-theory, and so we will simply define the algebraic \(\K\)-space
\[ \Kspace^{\pl}(\A) \coloneqq \core(\A)^{\grp} \]
to be the group completion of the \(\Einf\)-monoid \(\core\A\), that is, the core groupoid of \(\A\) with its direct sum monoidal structure. 
Similarly, given a \(\wic\)-additive Poincaré structure \(\QF \colon \A\op \to \Grp_{\Einf}\), we will define its Grothendieck--Witt space 
\[ \GWspace^{\pl}(\A) \coloneqq \Poinc(\A,\QF)^{\grp}\]
to be the group completion of the \(\Einf\)-monoid \(\Poinc(\A,\QF)\) of Poincaré objects in \(\A\). 

The main property we will need to know about the functors \(\Kspace^{\pl}\) and \(\GWspace^{\pl}\) is a suitable form of \emph{additivity}. 
To express this, we first set up some terminology. 
We will use the term \emph{split adjunction} to mean an adjunction between \(\wic\)-additive categories whose unit is a split inclusion and whose counit is a split projection. By a left-split (resp. right-split) subcategory of a \(\wic\)-additive category \(\B\) we will mean a full subcategory inclusion \(\A \hrar \B\) which is the left (resp. right) functor in a split adjunction. We note that in that case \(\A\) inherits from \(\B\) the property of being \(\wic\)-additive. Given a left-split subcategory \(i\colon \A \hrar \B\) and a right-split subcategory \(j\colon \C \hrar \B\), we will say that \(\A\) and \(\C\) are \emph{orthogonal} if \(\Hom(x,y) = \ast\) for every \(x \in \A\) and \(y \in \C\). We note that this is equivalent to the condition the right adjoint \(p\) of \(i\) vanishes on \(\C\) and also to the condition that the left adjoint \(q\) of \(j\) vanishes on \(\A\).

\begin{adefinition}\label{def:additive}
We will say that a functor \(\F\colon \Catwic \to \Grp_{\Einf}\) to \(\Einf\)-groups is \emph{additive} if the following holds: whenever we have a \(\wic\)-additive category \(\B\) and a pair of orthogonal split subcategories \(\A,\C \subseteq \B\) with adjoints \(p\colon \B \to \A\) and \(q\colon \B \to \C\) as above, the induced map
\[ \F(\A) \oplus \F(\C) \oplus \F(\ker(p) \cap \ker(q)) \lto \F(\B) \]
is an equivalence of \(\Einf\)-groups.
\end{adefinition}

\begin{remark}
In the situation of Definition~\ref{def:additive}, if \(\A \subseteq \B\) is a left-split subcategory with right adjoint \(p\colon \B \to \A\) then the subcategory \(A^{\perp} \coloneqq \ker(p) \subseteq \B\) is a right-split subcategory with left adjoint \(q'(x) = \cof[ip(x) \to x]\), and this is essentially by definition the maximal right-split subcategory of \(\B\) orthogonal to \(\A\), so we may call it the \emph{right-orthogonal complement} of \(\A\). In this case one has \(\ker(p) \cap \ker(q') = 0\), and so for any additive functor \(\F\colon \Catwic \to \Grp_{\Einf}\) one has \(\F(\B) = \F(\A) \oplus \F(\A^{\perp})\). On the other hand, \(\C\) can now be considered as a right-split subcategory of \(\A^{\perp}\), and has a similarly defined left-orthogonal complement in \(\A^{\perp}\), given by \(\ker(q|_{\A^{\perp}}) = \ker(p) \cap \ker(q)\). Additivity for this pair of orthogonal subcategories then yields \(\F(\A^{\perp}) = \F(\C) \oplus \F(\ker(p) \cap \ker(q))\). The upshot of this is that the additivity property as formulated in Definition~\ref{def:additive} can be deduced from the a priori weaker property obtained by restricting only to pairs \(\A,\C \subseteq \B\) which are \emph{maximal orthogonal}, that is, \(\C = \A^{\perp}\) and \(\A = \C^{\perp}\) (these two conditions are in fact equivalent to each other, and both are equivalent to the condition that \(\ker(p) \cap \ker(q) = 0\)).
\end{remark}

\begin{aproposition}\label{prop:additivity}
The functor \(\Kspace^{\pl}\colon \Catwic \to \Grp_{\Einf}\) is additive.
\end{aproposition}

In the setting of \emph{ordinary} \(\wic\)-additive categories, Proposition~\ref{prop:additivity} can be deduced from Waldhausen additivity and Quillen's comparison of algebraic \(\K\)-theory and the group completion of \(\core\A\) in the split-exact case (see~\cite[p. 228]{QuillenHigherK2}). In general, Waldhausen additivity is also known in the higher categorical setting by the work of Barwick~\cite{barwickhigher}, and the comparison with direct sum \(\K\)-theory is discussed in detail in Remark \ref{rem:fontes}. Alternatively, a direct proof of Proposition~\ref{prop:additivity} in the higher categorical setting is also possible by suitably adapting Quillen's original arguments.

In the setting of Poincaré \(\wic\)-additive categories, we suggest the following as the analogue of the above additivity property. Given a Poincaré \(\wic\)-additive category \((\A,\QF)\), we will say that a left-split additive subcategory \(\Lag \subseteq \A\) is \emph{isotropic} if \(\QF\) vanishes on \(\Lag\). 
Given an isotropic subcategory \(\Lag \subseteq \A\) with split right adjoint \(p\colon \A \to \Lag\), the full subcategory of \(\A\) spanned by the objects \(\x \in \A\) such that \(p(\x) = p(\Dual \x) = 0\) is closed under the duality, and so inherits from \(\A\) an additive Poincaré structure. We will refer to this full subcategory as the \emph{homology} of \(\Lag\), and denote it by \(\Hlgy(\Lag)\). On the other hand, if we let \(\Hyp(\Lag) \coloneqq (\Lag \times \Lag\op,\Hom)\) then we have a Poincaré additive functor \(\Hyp(\Lag) \to (\A,\QF)\) sending \((\x,\y)\) to \(p(\x) \oplus \Dual p(\y)\).

\begin{adefinition}
We will say that a functor \(\F\colon \Catpad \to \Grp_{\Einf}\) is \emph{additive} if for every Poincaré \(\wic\)-additive category \((\A,\QF)\) with isotropic subcategory \(\Lag \subseteq \A\), the induced functor
\[ \F(\Hyp(\mathcal E)) \oplus \F(\Hlgy(\mathcal E)) \lto \F(\A,\QF) \]
is an equivalence 
of \(\Einf\)-groups.
\end{adefinition}

\begin{aproposition}\label{prop:poincare-additivity}
The functor \(\GWspace^{\pl}\colon \Catpad \to \Grp_{\Einf}\) is additive.
\end{aproposition}

\begin{aremark}
Proposition~\ref{prop:poincare-additivity} implies Proposition~\ref{prop:additivity}. This follows from the identification \(\GWspace^{\pl}(\Hyp(-)) \simeq \Kspace^{\pl}(-)\) and the fact that if \(\A,\C \subseteq \B\) are a pair of orthogonal split subcategories with adjoints \(p,q\) as in Definition~\ref{def:additive}, then the Poincaré \(\wic\)-additive category 
\(\Hyp(\B)\) admits \(\A \times \C\op\) as an isotropic subcategory whose homology is \(\Hyp(\ker(p) \cap \ker(q))\). 
\end{aremark}

A posteriori, Proposition~\ref{prop:poincare-additivity} follows from Theorem \ref{thm:resolution} and the additivity of the \(\GW\)-space for stable Poincaré categories, specifically \cite[Theorem 3.2.10]{9authII}. On the other hand, if one restricts attention to ordinary \(\wic\)-additive categories equipped with \(\Ab\)-valued Poincaré structures, then Proposition~\ref{prop:poincare-additivity} follows from the work of Schlichting~\cite{Schlichting2019}. More precisely, the statements of both Proposition~\ref{prop:additivity} and Proposition~\ref{prop:poincare-additivity} can be reduced to certain universal cases by purely formal arguments. In the case of algebraic \(\K\)-theory, it suffices to check the claim for the orthogonal subcategory pair
\[ \A \hrar \Seq(\A) \hlar \A ,\]
where \(\Seq(\A)\) is the category of split exact sequences \(\x \hrar \y \twoheadrightarrow \z\) in \(\A\), the left-split copy of \(\A\) consists of the sequences \([\x = \x \to 0]\) and the right-split copy of \(\A\) consists of the sequences \([0 \to \z = \z]\), see e.g.\ \cite[Proposition 1.3.2]{waldhausen} and \cite[Remark 2.7.6]{9authII} for versions of the argument in the language of Waldhausen categories and stable $\infty$-categories, respectively.
Similarly, in the hermitian setting it suffices to prove additivity in the universal case of \(\Seq(\A,\QF)\), which has \(\Seq(\A)\) as an underlying \(\wic\)-additive category and its Poincaré structure sends 
a split exact sequence \(\x \hrar \y \twoheadrightarrow \z\) to 
\(\fib[\QF(\y) \to \QF(\x)]\); in this case the first reduction step is \cite[Proposition 6.7 (ii)]{Schlichting2019}, which shows that maps induced by cobordant Poincar\'e functors differ by a hyperbolic and the second consists of observing that the proof of the isotropic decomposition theorem \cite[Theorem 3.2.10]{9authII} requires only this input and works mutatis mutandis for $\flat$-additive in place of stable Poincar\'e categories. In the case when \(\A\) is an ordinary category and \(\QF\) takes values in abelian groups, this universal case of additivity was proven by Schlichting in~\cite[Theorem 5.1]{Schlichting2019}.

\begin{aproposition}\label{prop:l-theory}
Assuming Proposition~\ref{prop:poincare-additivity}, there is a natural bifibre sequence of \(\Einf\)-groups
\[ \Kspace^{\pl}(\A,\QF)_{\hC} \lto \GWspace^{\pl}(\A,\QF) \lto \Lspace^{\pl}(\A,\QF)\]
in which the first arrow is the one induced by the hyperbolic map \(\Kspace^{\pl}(\A) \to \GWspace^{\pl}(\A,\QF)\).
\end{aproposition}

\begin{aremark}\label{rem:l-theory-ordinary}
When \(\A\) is ordinary and \(\QF\) takes values in \(\Ab\subseteq \Spa\), Proposition~\ref{prop:poincare-additivity} is known from~\cite{Schlichting2019}. Running the proof below in that case works verbatim since the \(\Q\)-construction preserves ordinary categories. In particular, if \(\A\) is ordinary then Proposition~\ref{prop:l-theory} holds without extra assumptions. 
\end{aremark}

\begin{proof}[Proof of Proposition~\ref{prop:l-theory}]
For every \(n\) the Poincaré \(\wic\)-additive category \(\Q_n(\A,\QF)\) has an isotropic subcategory spanned by the objects corresponding to sequences of ``purely backward'' spans starting from \(0\), 
whose homology is the full subcategory of constant diagrams, that is, the image of the degeneracy operator \((\A,\QF) = \Q_0(\A,\QF) \to \Q_n(\A,\QF)\). For every additive functor \(\F\colon \Catpad \to \Grp_{\Einf}\) one consequently has
\[ \F(\Q_n(\A,\QF)) \simeq \F(\A,\QF) \oplus \F(\Hyp(\A))^n .\]
This means that if \(\F \Rightarrow \G\) is a natural transformation of additive functors \(\F,\G\colon \Catpad \to \Grp_{\Einf}\) such that \(\F(\Hyp(\A)) \to \G(\Hyp(\A))\) is an equivalence for every \(\A\), then the square
\[ \xymatrix{
\F(\A,\QF) \ar[r]\ar[d] & \F\Q^{(1)}(\A,\QF) \ar[d] \\
\G(\A,\QF) \ar[r] & \G\Q^{(1)}(\A,\QF)
}\]
is bicartesian in \(\Grp_{\Einf}\). Using the identification \(\Q_n(\Hyp(\A)) = \Hyp(\Q_n(\A))\) one obtains more generally that the square
\[ \xymatrix{
\F(\A,\QF) \ar[r]\ar[d] & \F\Q^{(n)}(\A,\QF) \ar[d] \\
\G(\A,\QF) \ar[r] & \G\Q^{(n)}(\A,\QF)
}\]
is bicartesian for every \(n\), and passing to transfinite compositions we obtain that 
\[ \xymatrix{
\F(\A,\QF) \ar[r]\ar[d] & \colim_n\F\Q^{(n)}(\A,\QF) \ar[d] \\
\G(\A,\QF) \ar[r] & \colim_n\G\Q^{(n)}(\A,\QF)
}\]
is bicartesian. Applying this to the natural transformation \(\Kspace^{\pl}_{\hC}(-) \Rightarrow \GWspace^{\pl}(-)\) it will hence suffice to show that \(\colim_n\Kspace(\Q^{(n)}(\A))\) vanishes. Indeed, it is well known (see~\cite{BarwickQ}) that the \(\Q\)-construction \(\Q\A\) is equivalent to the edgewise subdivision of the \emph{\(\bS\)-construction} \(\bS\A\). For any additive functor \(\F\) one has that \(\F(\bS_n\A) \simeq \F(A)^n\), and more precisely that \(\F(\bS\A)\) is the bar construction of the \(\Einf\)-group \(\F(\A)\). It follows that \(|\F(\Q\A)| \simeq |\F(\bS\A)| \simeq \mathrm B\F(\A)\) and similarly \(\F\Q^{(n)}(\A) \simeq \mathrm{B}^n\F(\A)\). We consequently obtain that \(\colim_n\Kspace^{\pl}\Q^{(n)}(\A) = \colim_n\mathrm{B}^n\Kspace^{\pl}(\A)  \simeq \ast\), as desired.
\end{proof}

\begin{acorollary}\label{cor:bordinary}
Let \((\D,\QFD)\) be a $0$-dimensional Poincaré (stable) category whose heart \(\A = \D^{\heart}\) is an ordinary category and such that \(\QF^{\heart} \coloneqq \QF|_{\A}\) takes values in \(\Ab \subset \Spa\). Then the map
\[ \GWspace^{\pl}(\A,\QF^{\heart}) \lto \GWspace(\D,\QF) \]
is an equivalence of spaces.
\end{acorollary}
\begin{proof}
By Remark~\ref{rem:l-theory-ordinary} we can apply Proposition~\ref{prop:l-theory} unconditionally when \(\A\) is ordinary. In the commutative diagram of \(\Einf\)-groups
\[ \xymatrix{
\Kspace^{\pl}(\A)_{\hC} \ar[r]\ar[d] & \GWspace^{\pl}(\A,\QF^{\heart})\ar[r]\ar[d] & \Lspace^{\pl}(\A,\QF^{\heart}) \ar[d] \\
\Kspace(\D)_{\hC} \ar[r] & \GWspace(\D,\QF) \ar[r] & \Lspace(\D,\QF)
}\]
whose rows are bifibre sequences the right most vertical map is then an equivalence by Proposition~\ref{prop:comp-L} and the left vertical map is an equivalence by the weight theorem for \(\K\)-theory (see Remark \ref{rem:fontes}). It then follows that the middle vertical map is an equivalence, as desired.
\end{proof}

\bibliographystyle{amsalpha} 
 
\newcommand{\etalchar}[1]{$^{#1}$}
\providecommand{\bysame}{\leavevmode\hbox to3em{\hrulefill}\thinspace}
\providecommand{\MR}{\relax\ifhmode\unskip\space\fi MR }
\providecommand{\MRhref}[2]{%
  \href{http://www.ams.org/mathscinet-getitem?mr=#1}{#2}
}
\providecommand{\href}[2]{#2}

\end{document}